\documentclass[twoside,letterpaper,draft,10pt]{amsart}

\usepackage{amssymb,amsmath,amsthm,latexsym}
\usepackage{amsfonts}
\usepackage[twlrm]{rawfonts}
\usepackage[english]{babel}
\usepackage[utf8]{inputenc}
\usepackage{indentfirst}
\usepackage{graphicx}
\usepackage[usenames]{color}
\usepackage[matrix,arrow,curve]{xy}
\usepackage{lscape}
\usepackage{pdflscape}
\usepackage{multicol}

\usepackage{mathtools}
\usepackage{multirow}
\usepackage[fixlanguage]{babelbib}
\usepackage[font=small,format=plain,labelfont=bf,up,textfont=it,up]{caption}
\usepackage[usenames,svgnames,dvipsnames]{xcolor}

\usepackage[a4paper,top=4.0cm,bottom=2.0cm,left=2.5cm,right=2.5cm]{geometry}

\newtheorem{teo}{Theorem}[section]
\newtheorem{cor}[teo]{Corollary}
\newtheorem{lem}[teo]{Lemma}
\newtheorem{pro}[teo]{Proposition}
\newtheorem{obs}[teo]{Remark}
\newtheorem{afr}[teo]{Claim}
\newtheorem{defi}[teo]{Definition}
\newtheorem{exe}[teo]{Example}
\newenvironment{dem}[1][Proof]{\noindent\textbf{#1.} }{\hfill \rule{0.5em}{0.5em}}

\newcommand{\dsum}{\displaystyle\sum}

\newcommand{\Rr}{\mathfrak{r}}
\newcommand{\Ll}{\mathfrak{l}}
\newcommand{\D}{\Delta}
\newcommand{\N}{\mathbb{N}}
\newcommand{\F}{\mathcal{F}}
\newcommand{\T}{\Theta}
\newcommand{\C}{\mathcal{C}}
\newcommand{\Om}{\Omega}
\newcommand{\om}{\omega}
\newcommand{\Gm}{\Gamma}

\newcommand{\bp}{\bar{p}}

\newcommand{\p}{\mathcal{P}}
\newcommand{\Z}{\mathcal{Z}}
\newcommand{\R}{\Rightarrow}
\newcommand{\al}{\alpha}
\newcommand{\U}{\mathcal{U}}

\newcommand{\m}{{-1}}

\newcommand{\Aut}{\mathbf{Aut}}
\newcommand{\End}{\mbox{End}}
\newcommand{\Pic}{\mathbf{Pic}}
\newcommand{\Pics}{\mathbf{PicS}}
\newcommand{\Hom}{\mbox{Hom}}

\newcommand{\s}{\sigma}
\newcommand{\et}{\theta}
\newcommand{\e}{\varepsilon}
\newcommand{\ot}{\otimes}
\usepackage{float}

\begin{document}

\title[Partial generalized crossed products and a seven term exact sequence (expanded version)]{Partial generalized crossed products and a seven term exact sequence (expanded version)}

\author[M. Dokuchaev]{Mikhailo Dokuchaev}
\address{Instituto de
Matem\'atica e Estat\'\i stica\\
Universidade de S\~ao Paulo \\
Rua do Mat\~ao, 1010\\
05508-090, S\~ao Paulo, SP, Brasil}
\email{dokucha@ime.usp.br}

\author[I.\ Rocha]{Itailma Rocha}
\address{Unidade Acad\^emica de Matem\'atica, Universidade Federal de Campina Grande - UAMat/UFCG, 
Avenida Apr\'{\i}gio Veloso, 882, 58429-900, Campina Grande, PB, Brasil.}
\email{itailma@mat.ufcg.edu.br}

\thanks{ The first named author was partially supported by FAPESP of Brazil (Process 2020/16594-0) and by CNPq of Brazil (Process 312683/2021-9). The second named author was partially supported by Grant 3223/2021 of Para\'{\i}ba State Research Foundation (FAPESQ) and CNPq.}

\keywords{Partial action, partial representation, crossed product, partial group cohomology.}


\subjclass[2000]{Primary 16W22,     16D20,  16S35; Secondary  20M18.}

\date{}

\begin{abstract} Given a non-necessarily commutative unital ring $R$ and a unital partial representation $\Theta $ of a group $G$  into the Picard  semigroup 
$\Pics(R)$ of the isomorphism classes of partially invertible $R$-bimodules, we construct an abelian group $\C(\T/R) $ formed by  the isomorphism classes of partial
generalized  crossed products related to $\Theta  $ and identify an 
appropriate second partial cohomology group   of $G$ with a naturally defined subgroup $\C _0(\T/R) $ of 
$\C(\T/R).$ Then we use the obtained results to give an analogue   of the  Chase-Harrison-Rosenberg exact sequence associated with an extension of non-necessarily commutative rings $R\subseteq S$ with the same unity and a unital partial representation $  G  \to \mathcal{S}_R(S)$ of an arbitrary  group  $G$ into the monoid    $\mathcal{S}_R(S)$  of the   $R$-subbimodules of $S.$ This generalizes
the works by Kanzaki and Miyashita.
\end{abstract}

\maketitle





\section{Introduction}\label{sec:introduction}

The purpose of this paper is to give a partial action version of the Y. Miyashita's non-commutative analogue 
\cite{miyashita1973exact} of the 
 Chase-Harrison-Rosenberg sequence \cite{chase1965galois}. The latter  is related to a Galois extension 
$R^G\subseteq R$ of commutative 
unital rings with a finite Galois group $G$ and was obtained applying  a seven term Amitsur cohomology exact sequence 
of S. U. Chase and A. Rosenberg  \cite{chase1965amitsur}, proved using  spectral sequences. A constructive proof of the 
Chase-Harrison-Rosenberg sequence was given by T. Kanzaki \cite{kanzaki1968generalized}, employing the novel notion of a generalized crossed product. The sequence  is of the form
\begin{align*}
&0\to H^1(G, \U(R)) {\to} {\bf Pic}(R^G) {\to}  {\bf Pic}(R)^G {\to }
H^2(G, \U(R)  ) {\to} B(R/R^G) {\to} H^1(G,{\bf Pic}(R)){\to} H^3(G, \U(R)),
\end{align*}
where $ {\bf Pic}(R)$ is the Picard group of the isomorphism classes of finitely generated projective $R$-modules of rank $1,$ $H^n(G,_{-})$ is the cohomology group of $G$ and $B(R/R^G)$ is the relative Brauer group, whose elements are the
 equivalence classes of the Azumaya  $R^G$-algebras split by $R.$ There is an induced action { of $G$} on $ {\bf Pic}(R),$ and
 $ {\bf Pic}(R)^G$ denotes the  subgroup of the elements fixed by the action.
 If $R$ is a field, then the first homomorphism of the sequence immediately gives Emmy Noether's theorem $H^1(G, \U(R))=1,$ which extends  the E. Kummer's result of 1855 on cyclic field extensions, known as Hilbert's Theorem 90 due to the fact of being listed under number 90 in D. Hilbert's
Zahlbericht. Another important fact, being generalized by the sequence, 
is the crossed product theorem, stating that if $R$ is a field, then the map $H^2(G, \U(R)  ) \ni [f] \mapsto [R\ast _f G] \in  B(R/R^G)$ is an isomorphism of groups, where   $R\ast _f G$ stands for the crossed product associated with the 
$2$-cocycle $f$ and the Galois action of $G$ on $R.$

Besides  Galois theory, crossed products are relevant  in ring theory, in the theory of von Neumann algebras and in that of   $C^*$-algebras.   In the latter area more general types of crossed products were introduced, in particular, crossed products by partial actions, which are useful to deal with important classes of $C^*$-algebras generated by partial isometries. Prominent early examples of partial crossed product descriptions  include those for the Cuntz-Krieger $C^*$-algebras  \cite{E2} and for the more general Exel-Laca  $C^*$-algebras   \cite{EL}. Among the more recent relevant examples we mention    the tame $C^*$-algebras related to finite bipartite separated graphs   \cite{AraE1}, the full and reduced $C^*$-algebras  of $E$-unitary or strongly $E^*$-unitary inverse semigroups  \cite{MiSt},  
  the Carlsen-Matsumoto $C^*$-algebras associated to arbitrary subshifts  on finite alphabets \cite{DE2} and 
 certain ultragraph $C^*$-algebras \cite{GR3}. Highly interesting  particular cases include $C^*$-algebras associated to dynamical systems of type $(m,n)$ (see \cite{AraEKa})  and  graph $C^*$-algebras (see \cite{E6}).

In algebra, partial crossed products turned out be applicable to group graded algebras \cite{dokuchaev2005associativity}, \cite{dokuchaev2008crossed}, to   Hecke algebras \cite{Exel2008}, to Leavitt path algebras \cite{GR} and to Steinberg algebras \cite{BeuCor}, \cite{BeuGon2}, \cite{HazratLi}.  In particular, it was shown in  \cite{dokuchaev2008crossed} that any group graded algebra $A$, satisfying a natural mild condition, after passing  to finite matrices over $A$ of an appropriate infinite size,  becomes
a crossed product by a twisted partial action.  This is an algebraic analogue of R. Exel's result  \cite{exel1997twisted} on $C^*$-algebraic bundles. In addition,  partial crossed product analogues are useful for semigroups, as it can be seen in \cite{KL}, \cite{CornGould}, \cite{Khry1}, \cite{Kud},   \cite{Kud2}.

Among the diverse algebraic developments on partial actions we mention two steps related to our theme: the partial 
action treatment of Galois theory of  commutative rings  in \cite{dokuchaev2007partial} and the introduction of a group cohomology  based on partial actions in  \cite{dokuchaev2015partial}. While the former  inspired a Hopf theoretic approach to partial actions by  S. Caenepeel and K. Janssen in   \cite{CaenJan}, which in its turn became a starting point for interesting and fruitful Hopf theoretic developments, the latter is useful for partial projetive group representations \cite{DoSa},  for extensions of semi-lattices of groups by groups \cite{DKh3} and to the study of the  ideal structure of reduced   (global) $C^*$-crossed products  \cite{KennedySchafhauser}. It also influenced the introduction of partial cohomology of Hopf algebras in  \cite{BaMoTe} and that of  groupoids in \cite{NysOinPin3}.
For more information on advances around partial actions and applications the reader is referred to R. Exel's book \cite{E6} and to the survey article \cite{D3}.

In view of the above mentioned progresses,  it was natural to generalize the   Chase-Harrison-Rosenberg sequence to partial Galois extensions of commutative rings.  This was done in a sequence of two papers \cite{DoPaPi} and \cite{DoPaPi2}.
New ingredients came into the picture:  the first one is the inverse semigroup ${\bf PicS}(R)$ of the finitely generated projective $R$-modules of rank $\leq 1,$ on which the Galois partial action $\alpha $ of the finite group $G$  on the commutative 
unital ring $R$  induces a partial action $\alpha ^*;$ the second one is certain partial representations needed to introduce a partial version of generalized crossed products.  All (global) cohomology groups in the Chase-Harrison-Rosenberg sequence  are substituted  by their partial theoretic analogues, and the sequence takes the form:
\begin{align*}
&0\to H^1(G,\alpha , R){\to} {\bf Pic}(R^\alpha ){\to }{\bf PicS}(R)^{\alpha^*}\cap {\bf Pic}(R) {\to }H^2(G,\alpha, R) {\to} B(R/R^\alpha){\to}\\& H^1(G,\alpha^*,{\bf PicS}(R)){\to}   
H^3 (G,\alpha , R),
\end{align*}  where  $R^\alpha $ and 
${\bf PicS}(R)^{\alpha ^*}$ are the subring and the submonoid of  partial invariants, respectively (see  
Section~\ref{sec:ParAcParRepParHom}
for definitions).  Two immediate consequences of the sequence are a partial version of the Hilbert's Theorem 90 \cite[Corollary 6.8]{DoPaPi2} and the partial  crossed product theorem \cite[Corollary 6.9]{DoPaPi2}.

T. Kanzaki's  proof of the     Chase-Harrison-Rosenberg sequence inspired  Y. Miyashita to use generalized crossed products to produce  a non-commutative analogue of the sequence \cite{miyashita1973exact}. Instead of a Galois extension, Miyashita's sequence is  related to a rather general setting of a ring extension $R\subseteq S$ with the same unity and a fixed homomorphism $\Theta :  G \to  {\bf Inv} _R(S),$ where $G$ is  an arbitrary non-necessarily finite  group and 
$ {\bf Inv} _R(S)$ is the group of invertible $R$-subbimodules of $S.$ 
The sequence  has a  somewhat different form, but  the case of a Galois extension of commutative rings $R^G \subseteq R$ can  be recovered by taking $S=  R\star G $ and a rather natural homomorphism 
$G \to   {\bf Inv} _R( R\star G) ,$  where  $R\ast G$ is the skew group ring  by the Galois action of $G$ on $R.$
One of the crucial ingredients of the non-commutative version of the sequence is an appropriate analogue  
${\mathcal  B} (\Theta /R)$ of the Brauer group. It is defined as a  quotient of the abelian group ${\mathcal C} (\Theta /R)$ of the isomorphism classes of certain generalized crossed products related to $\Theta .$ 
A subgroup  ${\mathcal C}_0 (\Theta /R)$ of  ${\mathcal C} (\Theta /R)$ is shown to be isomorphic to a second cohomology group of $G,$ and the restriction of the natural map  
${\mathcal C} (\Theta /R)\to {\mathcal  B} (\Theta /R)$ to     ${\mathcal C}_0 (\Theta /R)$ provides one of the homomorphisms of the sequence.
Miyashita's sequence was extended for a technically more difficult case of rings with local units by   L. El Kaoutit and   J. G\'omez-Torrecillas in \cite{el2012invertible} using their results from  \cite{el2010invertible}.  
 
 Several other articles related to the Chase-Harrison-Rosenberg sequence were published. In particular, in \cite{CroRaeWil}  D. Crocker, I. Raeburn and D. Williams obtained  a $C^*$-analogue of the   sequence, using the equivariant version of the Picard group 
\cite[Proposition 2]{CroRaeWil},  that of the Brauer  group \cite{CroKumRaeWil} and  the cohomology theory developed by C. C. Moore \cite{Moore} for actions of locally compact groups on Polish modules.

We begin by giving some preliminary facts on projective modules over a
non-necessarily commutative   ring $R$ with $1$ in Section~\ref{sec:Proj}, followed by  the notion of the generalized Picard semigroup  $ \Pics(R)$ of the isomorphism classes of the partially invertible $R$-bimodules   in Section~\ref{sec: Pics} and by some technical facts on the bimodule relation  $M|N$ and isomorphisms of the form  $M\ot_RN\simeq N\ot_RM$ in Section~\ref{section: modulos similares}. Furthermore,  given an  extension of rings $R\subseteq S$ with the same unity,  { we recall    in Section~\ref{sec: pSR} the definition of the group $\p(S/R)   $ of the isomorphism classes of certain} $R$-bimodule maps $P\to X,$ 
where  $P$ is an $R$-bimodule, whose  isomorphism class $[P]$ lies in the Picard group  $\Pic(R)$ and $[X]\in \Pic(S).$ The group $ \p(S/R)  $ participates an exact sequence 
{ given  in  \cite[Theorem 1.5]{miyashita1973exact} (see Theorem~\ref{sequenciaparaodiagramadoprimeiromorfismo})},  
which is used to construct the first two homomorphisms of our main  sequence, whose second term is  an appropriate subgroup 
 of $\p(S/R)$.

Section~\ref{sec:ParGenCrosProd} deals with  partial generalized  crossed products. After recalling in Section~\ref{sec:ParAcParRepParHom} some background on partial actions, partial representations and partial cohomology of a group $G,$  we give  in Section~\ref{sec:ParRep G to Pics} some auxiliary facts related to a unital partial representation of the form $G\to  \Pics(R)$ and, given such a representation, construct  the associated partial actions $\alpha$ and $\alpha ^*$ on $\Z$ and   $\Pics(R),$ respectively, 
where  $\Z$ stands for the centre of $R.$ In particular, this allows us to consider partial cohomology groups of a group $G$ with values in $\Z, $  which we denote by $H_\T^{n}(G,\al,\Z).$ Notice that along with the usual  concept of a partial action of a group on a semigroup, in which  the domains of the partial isomorphisms are assumed to be two-sided ideals, it may be useful to keep in mind a relaxed version allowing them to be only subsemigroups (see Definition~\ref{def:ParAcOnSemigr} and  Proposition~\ref{proparaacaoparcialsobrePics}).

The main technical part of Section~\ref{sec:ParGenCrosProd} is concentrated in  Section~\ref{sec:FactSet GenParCrossProd}, where we give the concept   of a factor set and that of a generalized crossed product related to a fixed unital partial   
representation $\T : G \to \Pics(R)$,  and prove several results on them. Factor sets are needed to produce generalized { partial} crossed products and they may arise from partial representations $\T $  as  families of $R$-bimodule isomorphisms  of the form $f^\T=\{f_{x,y}^\T: \T_x\ot_R\T_y\longrightarrow 1_x\T_{xy}, \ x,y\in G \},$
where $1_x,$ $(x\in G),$ is a central idempotent in $R.$   
 An obstruction for $f^\T$ to be a factor set  is a partial $3$-cocycle (see Proposition~\ref{3cocicloquecorrigeaassociatividade} and Corollary~\ref{obs3cocicloparasimilares}), which  is automatically  trivial if the partial representation comes from that of the  form $\Theta:  G  \to \mathcal{S}_R(S).$  Here   $\mathcal{S}_R(S)$ denotes the set of the   $R$-subbimodules of $S,$ equipped with the natural  multiplication
$MN=\left\lbrace  \dsum  m_in_i; \ m_i \in M, n_i \in N\right\rbrace ,$ $M,N\in \mathcal{S}_R(S),$ turning $\mathcal{S}_R(S)$ into  a monoid with neutral element $R.$   Section~\ref{subsection: SRS} gives some preliminaries on unital partial representations $G  \to \mathcal{S}_R(S).$ 

As it can be seen already in \cite{DoPaPi} and \cite{DoPaPi2},  it is more laborious to deal with the partial action setting. The non-commutative case brings new technical challenges, and a big part of the work is related to the partial action version of the group ${\mathcal C} (\Theta /R),$ to which 
 Section~\ref{sec:groupC} is dedicated. More specifically, 
 given a unital partial representation 
  $\T:G\longrightarrow\Pics(R)$   we   construct in  Theorem~\ref{grupoC} the abelian group $\C(\T/R)$ of the isomorphism
classes of certain generalized crossed products  related to $\T .$  Its subgroup   ${\mathcal C}_0 (\Theta /R)$ is given in  Proposition~\ref{pro:subgrC_0}  and  an isomorphism
between       ${\mathcal C}_0 (\Theta /R)$ and a second partial cohomology group of $G$ is established in Theorem~\ref{C0isomorfoH2}.

 Our seven term exact sequence is obtained in a series of subsections of  Section~\ref{cap: sequencia}, and it  is produced starting with
 an extension of rings $R\subseteq S$ with he same unity and  
   a fixed    unital partial representation
   $\T : G  \to \mathcal{S}_R(S).$   An important step is the construction  in Section~\ref{sec:secondExSeq} of a  group homomorphism 
${\mathcal L}: \Pic_\Z(R)^{(G)} \to \C(\T/R)$ (see Theorem~\ref{teo:homL}), where  $\Pic_\Z(R)^{(G)}$ is a subgroup of 
$\Pic_\Z(R)\subseteq \Pic (R)$
(see Section~\ref{sec: Pics} and Lemma~\ref{lemma:PicZ^G}). The homomorphism ${\mathcal L} $ is used to define our   analogue  of the Brauer group  in Section~\ref{sec:B}  as the quotient $\mathcal{B}(\T/R)=\dfrac{\C(\T/R)}{\mbox{Im}(\mathcal{L})}.$ 

The sixth term  $\overline{H^1}(G,\al^*,\Pics_0(R))$ of our sequence is introduced in  Section~\ref{sec:overline{H}} and it  involves $1$-cocycles with values in the subsemigroup $\Pics_0(R)\subseteq \Pics(R)$ provided by  Lemma~\ref{propriedadesdePics0}. The same lemma implies that  the partial action $\al^*$ of $G$ on  $ \Pics(R)$ can be restricted to  $\Pics_0(R).$ Then the  group $\overline{H^1}(G,\al^*,\Pics_0(R))$ is  defined as the quotient 
$\dfrac{Z^1(G,\al^*,\Pics_0(R))}{(\zeta\circ\mathcal{L})(\Pic_\Z(R)^{(G)})},$ where $	\zeta:  \C(\T/R)  \longrightarrow  Z^1(G,\al^*,\Pics_0(R))$ is the group homomorphism given in 
Lemma~\ref{definicaodezeta}. It follows by Remark~\ref{rem:B^1} that $\overline{H^1}(G,\al^*,\Pics_0(R))$ is a quotient of the partial cohomology group ${H^1}(G,\al^*,\Pics_0(R)).$

Besides being used in the definition of our version  $\mathcal{B}(\T/R)$ of a Brauer group,   the map  ${\mathcal L} $ is involved in the construction of three consecutive homomorphisms of our sequence:
\begin{equation*}
	 \Pic_\Z(R)\cap\Pics_\Z(R)^{\al^*} \to H_\T^2(G,\al,\Z) 
		\to \mathcal{B}(\T/R) \to 	 \overline{H}^1(G,\al^*,\Pics_0(R))  
	\end{equation*}
(see  Proposition~\ref{terceiraseqexata},  Remark~\ref{phi5} and the definition of $\varphi _3$ after the proof of  Theorem~\ref{teo:homL}).
The entire sequence is given in Theorem~\ref{theo:main}, and it is of the form:

\begin{align*}
	1 & \to   H_\T^{1}(G,\al,\Z) \to  \p_\Z(\D(\T)/R)^{(G)}  \to \Pic_\Z(R)\cap\Pics_\Z(R)^{\al^*} \to H_\T^2(G,\al,\Z) 
		\to \mathcal{B}(\T/R)\\ &\to 	 \overline{H}^1(G,\al^*,\Pics_0(R)) \to H^3_\T(G,\al,\Z),  
	\end{align*}
where the subgroup   $\p_\Z(\D(\T)/R)^{(G)} $ of $\p (\D(\T)/R) $ is defined in Section~\ref{seq1} (see Lemma~\ref{lemma:SubgrP^G}).

\section{Background}\label{sec:backkground}
%

In all what follows  $R$ will stand for  an associative ring  with unity element $1.$ We denote by $\Z$ the center of  $R$ and by $\U(R)$ the unit group of $R$. A right  $R$-module $M$ is called unital if $m\cdot 1=m$, for all $m \in M$. In this case, we have that $MR=R$ and $M\ot_RR\simeq M$ via the right $R$-action on $M$, i.e. $m\ot r\mapsto m\cdot r$. Analogously we define a left unital $R$-module.   We shall  sometimes write $M_R$ to say that $M$ is a 
right module over $R$. The meaning of $_RM $ will be similar.  Given an $R$-bimodule $M$ we denote by 
$\Aut_{R-R}(M)$ the group of all  $(R,R)$-bimodule automorphisms of $M.$ Moreover, the bimodule $M$ is called central if $r\cdot m=m\cdot r$, for all   $m \in M$ and $r \in R$. 

Let $e$ be a central idempotent of $R$ and $M$ an $R$-bimodule. Write $eM=\{m\in M \, :\,  em=m\}.$ Obviously,  $eM$ is an  $R$-subbimodule of $M,$ since $e$ is a central element. Moreover, there is an $R$-bimodule isomorphism $Re\ot_RM\simeq eM$ defined by $r\ot m\mapsto rm$, whose inverse is given by $m\mapsto e\ot m$.
Analogously, $Me=\{m \in M\, :\,  \ me=m\}$ is a $R$-subbimodule of $M$ and $M\ot_RRe\simeq Me$, as $R$-bimodules.
 If $F:M\rightarrow N$ is an isomorphism of $R$-bimodules, then the restriction of $F$ to $eM$ is an $R$-bimodule isomorphism onto $eN$. In some cases, for simplicity of notation, we will denote this restriction by the same symbol $F$.
  

\subsection{Projective modules over non-commutative rings}\label{sec:Proj}

 In this section  we list some  properties of finitely generated  projective $R$-modules, whose proofs do not seem to 
be easy to find in general references like books.


Let $M$ and $N$ be $R$-bimodules. Then the sets $\Hom(M_R,N_R)$ and $\Hom(_RM,_RN)$ are $R$-bimodules via:
\begin{equation}
(r\cdot f)(m)=rf(m) \ \ \mbox{and} \ \ (f\cdot r)(m)=f(rm), \ \ f\in \Hom(M_R,N_R), m \in M, r \in R, \label{estrelaMebimodulo}
\end{equation} 
\begin{equation}
(r\cdot g)(m)=g(mr) \ \ \mbox{and} \ \ (g\cdot r)(m)=g(m)r, \ \ g \in \Hom(_RM,_RN), m\in M, r\in R .\label{Mestrelaebimodulo} 
\end{equation}

In particular, taking $N=R,$ we have that the sets  ${^*M}=\Hom(M_R,R_R)$ and 
${M^*}=\Hom(_RM,_RR)$ are $R$-bimodules via  (\ref{estrelaMebimodulo}) and (\ref{Mestrelaebimodulo}).

\begin{lem}\label{lemisomorfismoPPestrela} Let $P$ be an $R$-bimodule. 
	\begin{itemize}
		\item[(i)] If $P$ is  {  finitely generated projective as a right} $R$-module, then the map $\varphi:P\ot_R {^*P}  \longrightarrow  \End(P_R),$ defined by $\varphi(p\ot f)(p')=pf(p'),$
		is an $R$-bimodule isomorphism.\\
		
		\item[(ii)]  If $P$ is  {  finitely generated  projective 
		as a left} $R$-module, then the map $\varphi': P^{*}\ot_RP \longrightarrow  \End(_RP),$ defined by $\varphi(f\ot p)(p')=f(p')p,$
		is an $R$-bimodule isomorphism. 
	\end{itemize}
	
\end{lem}
\begin{dem} 
	(i) Clearly $\varphi$ is a well-defined  left $R$-linear map.  Given $r \in R,$ we have
	$$\varphi(p\ot f\cdot r)(p')=p(f\cdot r)(p')=pf(rp')=\varphi(p\ot f)(rp')=(\varphi(p\ot f)\cdot r)(p'), $$
	for all $p' \in P$. Thus $\varphi$ is right $R$-linear and therefore $R$-bilinear. 
	Let $\{p_i,f_i\}_{i=1,2,...,n}$ be a dual base of $P$. Given  $f\in \End(P_R)$ we have  $\dsum_{i=1}^nf(p_iw)\ot f_i \in P\ot_R{^*P}$ and 
	$$\varphi\left( \dsum_{i=1}^nf(p_i)\ot f_i\right)(p')= \dsum_{i=1}^nf(p_i)f_i(p')=\dsum_{i=1}^nf(p_if_i(p'))=f(p'),$$
	for all $p' \in P$. Thus $\varphi\left( \dsum_{i=1}^nf(p_i)\ot f_i\right)=f$ and therefore $\varphi$ is onto.
	
	Now, let $p_l' \in P$ and $f_l' \in {^*P}$, with $l=1,2,...,k$, such that $\dsum_{l=1}^kp_l'f_l'(p)=0$, for all $p \in P.$ Then
	$$\dsum_{l=1}^kp_l'\ot f_l'=\dsum_{l,i}p_if_i(p_l')\ot f_l'=\dsum_{l,i}p_i\ot f_i(p_l')\cdot f_l'=0,$$
	since 
	$$\left( \dsum_{l=1}^kf_i(p_l')\cdot f_l'\right)(p)=\dsum_{l=1}^kf_i(p_l')f_l'(p)=\dsum_{l=1}^kf_i(p_l'f_l'(p))=0, $$
	for all $p \in P$. This implies that  $\varphi$ is injective. Therefore, $\varphi$ is an $R$-bimodule isomorphism.
	
	(ii) is obtained analogously to   (i).
	%
	
\end{dem}

\begin{lem}\label{tensordemodulospfg} Let  $P$ and $Q$ be $R$-bimodules. If $P$ and $Q$ are {  finitely generated  projective as right (left)} $R$-modules, then $P\ot_RQ$ is { a  finitely generated  projective right (left) } $R$-module.
\end{lem}
\begin{dem} Let $P$ and $Q$ be { $R$-bimodules, which are   finitely generated  projective as right  $R$-modules.} Let $\{f_i,p_i\}_{i=1,2,...,n}$ and $\{g_j,q_j\}_{j=1,2,...,m}$ be dual bases of { $P_R$ and $Q_R$,} respectively, and  let
	$$\begin{array}{c c c l}
	F_{i,j}: & P\ot_RQ & \longrightarrow & R\\ 
	& p\ot q & \longmapsto & g_j(f_i(p)q).
	\end{array}$$
	Clearly $F_{i,j}$ is a well-defined right $R$-linear map, that is $F_{i,j}\in {^*(P\ot_RQ)}$. Moreover, 
	\begin{eqnarray*}
		\dsum_{i,j=1}^{n,m}(p_i\ot q_j)F_{i,j}(p\ot q) & = & \dsum_{i,j=1}^{n,m}(p_i\ot q_j)g_j(f_i(p)q)=\dsum_{i,j=1}^{n,m}p_i\ot q_jg_j(f_i(p)q)\\
		& = & \dsum_{i=1}^np_i\ot f_i(p)q= \dsum_{i=1}^np_if_i(p)\ot q= p\ot q,
	\end{eqnarray*}
	for all $p \in P$ e $q\in Q$. Thus, $\{F_{i,j}, p_i\ot q_j\}_{i,j}$ is dual base  { for the right $R$-module} $P\ot_RQ,$ and consequently   $P\ot_RQ$ is a {  finitely generated  projective right $R$-module. 
	The left-hand side version is symmetric.}

	
\end{dem}

\begin{lem}\label{PestrelaQestrelaisoaQPestrela} Let $P$ and $Q$ be $R$-bimodules. 
	\begin{itemize}
		\item[(i)]  If $P$ and $Q$ are {  $R$-bimodules, which are finitely generated  projective as
		right $R$-modules,} then the map $\eta:{^*P}\ot_R{^*Q}  \rightarrow  {^*(Q\ot_RP)},$ defined by $\eta(f\ot g)(p\otimes q)=f(g(q)p),$
		is an $R$-bimodule isomorphism.
		
		\item[(ii)] If $P$ and $Q$ are { $R$-bimodules, which are  finitely generated  projective as left $R$-modules,} then the map $\eta':  P^*\ot_RQ^*  \rightarrow (Q\ot_RP)^*,$ defined by $\eta'(f\ot g)(q\ot p)=g(qf(p)),$
		is an $R$-bimodule isomorphism.
		
	\end{itemize}
\end{lem}
\begin{dem} (i) Let us check first that  $\eta$ is a well-defined homomorphism of  $R$-bimodules. Indeed, given  $r \in R,$ we have
	\begin{equation*}
	\eta(f\ot g)(q\ot pr)=f(g(q)pr)=f(g(q)p)r=\eta(f\ot g)(q\ot p)r.
	\end{equation*}
Thus, $\eta(f\ot g) \in {^*(Q\ot_RP)}$. Now, for $r \in R, p \in P$ and $q \in Q$
	\begin{equation*}
	\eta(f\cdot r\ot g)(q\ot p)=(f\cdot r)(g(q)p)=f(rg(q)p)=\eta(f\ot r\cdot g)(q\ot p).
	\end{equation*}
 Hence $\eta$ is  $R$-balanced and therefore it is a well-defined map. Clearly,  $\eta$ is left $R$-linear. On the other hand,
	\begin{equation*}
	\eta(f\ot g\cdot r)(q\ot p)=f((g\cdot r)(q)p)=f(g(rq)p)=(\eta(f\ot g)\cdot r)(q\ot p),
	\end{equation*}
	for $r\in R$. Thus, $\eta$ is right $R$-linear. Therefore, $\eta$ is an $R$-bimodule morphism.
	
	Let $\{p_i, f_i \}_{i=1,2,...,n}$ and $\{q_j,g_j \}_{j=1,2,...,m}$ be dual bases of $P$ and  $Q$, respectively. We check that
	\begin{equation*}
	\begin{array}{c c l}
	{^*(Q\ot_RP)} & \longrightarrow & {^*P}\ot_R{^*Q}\\
	F & \longmapsto & \dsum_{i,j}F(q_j\ot p_i)\cdot f_i\ot g_j
	\end{array}
	\end{equation*}
	is the inverse for $\eta$. Indeed, given $F\in {^*(Q\ot_RP)}$ we have
	\begin{eqnarray*}
		\eta\left( \dsum_{i,j}F(q_j\ot p_i)\cdot f_i\ot g_j\right)(q\ot p) & = & \dsum_{i,j}F(q_j\ot p_i)f_i(g_j(q)p)=  \dsum_{i,j}F(q_j\ot p_if_i(g_j(q)p))\\
		& = & \dsum_{j}F(q_j\ot g_j(q)p)= \dsum_{j}F(q_jg_j(q)\ot p)= F(q\ot p),
	\end{eqnarray*}
	for all $p \in P$ and $q \in Q$. 
	
	On the other hand, given $f \in {^*P}$ and $g\in {^*Q}$ we have
	\begin{eqnarray*}
		f\ot g  & \longmapsto &   (q\ot p \mapsto f(g(q)p)) \longmapsto \dsum_{i,j}f(g(q_j)p_i)f_i\ot g_j\\
		& = & \dsum_{i,j}(f\cdot g(q_j))(p_i)f_i\ot g_j =\dsum_{j}f\cdot g(q_j)\ot g_j\\
		& = &  \dsum_{j}f\ot g(q_j)\cdot g_j=f\ot g.
	\end{eqnarray*} 
	
	Therefore,  $\eta$ is an $R$-bimodule isomorphism.  
	
	The proof of (ii) is analogous to that of (i).  \end{dem}

\subsection{The generalized Picard semigroup}
\label{sec: Pics}
{We recall that an  $R$-bimodule $P$ is called {\it invertible} if  there  exists an $R$-bimodule $Q$ such that $$P\ot_RQ\simeq R\simeq Q\ot_RP,$$ { as $R$-bimodules (see, for example, \cite[Chapter II]{Bass}).} We shall need the following known fact.
	\begin{pro}\label{proModInv}{ \cite[Chapter II]{Bass}} Let $P$ be an  $R$-bimodule. Then the following are equivalent:
		\begin{itemize}
			\item[$(a)$] $P$ is invertible.
			\item[$(b)$] $P$ is a finitely generated  projective  right  $R$-module, which is a generator and  $R\simeq End(P_R)$. 
			\item[$(c)$] $P$ is a finitely generated  projective  left  $R$-module, which is a generator and $R\simeq End(_RP).$  
		\end{itemize}	
	\end{pro}
	We  also recall that the  \textit{Picard group} $\Pic(R)$ of a ring $R$ is the set of the isomorphism classes  of  invertible $R$-bimodules with multiplication given by $\otimes_R$. 
	\begin{obs}\label{isomorfismoemPic} { It can be readily verified that if} $[P]\in \Pic(R),$ then there exist $R$-bimodule isomorphisms:
		\begin{center}	
			\begin{tabular}{c c c }
				$\begin{array}{c c l}
					{^*P}\ot_RP & \stackrel{\Rr}{\longrightarrow}& R\\
					f\ot p & \longmapsto  & f(p)
				\end{array},$   & &  $\begin{array}{c c l}
					P\ot_R{^*P} & \stackrel{\Ll}{\longrightarrow} & \End(P_R)\simeq R\\
					p\ot f & \longmapsto  & (p'\mapsto pf(p'))
				\end{array},$ \\
				&  & \\
				& & \\
				$\begin{array}{c c l}
					P\ot_RP^* & \stackrel{\Ll'}{\longrightarrow} & R\\
					p\ot f & \longmapsto  & f(p)
				\end{array}$  & and & $\begin{array}{c c l}
					{P^*}\ot_RP & \stackrel{\Rr'}{\longrightarrow} & \End(_RP)\simeq R\\
					f\ot p & \longmapsto  & (p'\mapsto f(p')p)
				\end{array}.$ \\
				&  & \\
			\end{tabular}		
		\end{center}	
	\end{obs}
	By  \cite[Theorem 1.1]{morita1967endomorphism},  ${^*P}\simeq P^*$ as $R$-bimodules, and hence   $[^*P]=[P^*]=[P]^\m$ in $\Pic(R)$. We shall write $[P]^\m=[P^\m]$.   Note that  we can choose the isomorphisms   
	$\Rr$ and $\Ll$ so that  $(R,R,P,P^*,\Rr,\Ll)$ is a Morita  context, and consequently:
	\begin{equation}
		\Ll(p\ot q)p'=p\Rr(q\ot p') \ \ \mbox{and} \ \ \Rr(q\ot p)q'=q\Ll(p\ot q'), \label{morita}
	\end{equation}
	for all $p,p' \in P$ and $q,q'\in P^*$  { (see \cite[Chapter II]{Bass}).} 
	{ For a $K$-algebra $R$ over a commutative ring $K$ we denote} by  $\Pic_K(R)$ the set of the isomorphism classes of  invertible  $R$-bimodules which are central  over $K,$  i.e.,  $kp=pk$ for all  $k \in K$ and $p \in P$. In the case when  $R$  is commutative, we have that  $[P]\in \Pic_R(R)$ if and only if $P$ is a finitely generated projective $R$-module of rank $1$  {(see, for example, \cite[II, \S 5]{demeyer1971separable})}.\footnote{If $R$ is commutative, then  $\Pic_R(R)$ is usually denoted simply by $\Pic(R)$.}}

  We say that an  
$R$-bimodule $P$ is \textit{partially invertible} if 

\begin{itemize}
	\item[(i)] $P$ is finitely generated projective left and right  $R$-module;
	\item[(ii)] The maps
	\begin{center}
		\begin{tabular}{c c c}
			$\begin{array}{c c l}
			R& \longrightarrow & \End(P_R)\\
			r & \longmapsto & (p\mapsto rp)
			\end{array}$ & and & $\begin{array}{c c l}
			R& \longrightarrow & \End(_RP)\\
			r & \longmapsto & (p\mapsto pr)
			\end{array}$
		\end{tabular}
	\end{center}
	are epimorphisms. 
\end{itemize}

{\begin{exe}\label{exe:ParInvMod} Let $e\in R$ be a central idempotent and $P$  an invertible $eR$-bimodule. Define a structure of an $R$-bimodule on $P$ by setting $r\cdot p = (er)p$ and $p\cdot r = p(er), r\in R, p \in P.$ Then $\End(P_R)=\End(P_{eR}), \End(_RP) = \End(_{eR}P)$ and it is readily seen that $P$ is a partially invertible $R$-bimodule.
\end{exe}
}

We denote by $\mathbf{PicS}(R)$ the set of the isomorphism classes $[P]$ of partially invertible $R$-bimodules, 
that is, 
\begin{equation*}
\mathbf{PicS}(R)=\{[P], P \ \mbox{is a partially invertible $R$-bimodule}\}.
\end{equation*}

\begin{pro} $\Pics(R)$ is a monoid with multiplication defined by $[P][Q]=[P\ot_RQ]$.
\end{pro}
\begin{dem}
	By Lemma \ref{tensordemodulospfg}, we have that $P\ot_RQ$ is a left and right finitely generated projective   
$R$-module. Let  $\psi$ be the map defined by the following  chain of $R$-bimodule isomorphisms:
	$$\xymatrix{   R\ar[r]\ar@/_1.5cm/@{-->}[rrrrddd]_{\psi} & \End(P_R)\ar[r]^{\varphi_P^\m} & P\ot_R{^*P}\ar[r] & P\ot_RR\ot_R{^*P}\ar[r] &  P\ot_R\End(Q_R)\ot_R{^*P}\ar[d]^{P\ot \varphi_Q^\m\ot {^*P}} \\
		& & & & P\ot_RQ\ot_R{^*Q}\ot_R{^*P}\ar[d]^{P\ot Q\ot \eta}  \\
		& & & & P\ot_RQ\ot_R{^*(P\ot Q)}\ar[d]^{\varphi_{P\ot Q}}\\
		& & & & \End(P\ot_RQ)   }$$
	where $\varphi_P,\varphi_Q$ and $\varphi_{P\ot Q}$ are as in  Lemma \ref{lemisomorfismoPPestrela} and $\eta$ is
 the isomorphism given by Lemma \ref{PestrelaQestrelaisoaQPestrela}. 
	Then $\psi$ is onto and we have that 
	$$\psi(r)(p\ot q)=rp\ot q.$$
Indeed, let $\{p_i,f_i\}_{i=1,...,n}$ and $\{q_j,g_j\}_{j=1,...,m}$ be dual bases for $P$ and $Q,$ respectively,  as right 
$R$-modules. Given $r\in R$,  the above sequence of morphism  is:
	\begin{eqnarray*}
		r& \longmapsto & (p\mapsto rp)\longmapsto \dsum_{i=1}^nrp_i\ot f_i \longmapsto \dsum_{i=1}^nrp_i\ot 1 \ot f_i\\
		& \longmapsto & \dsum_{i=1}^nrp_i\ot Id_{Q} \ot f_i\longmapsto \dsum_{i,j} rp_i\ot q_j\ot g_j\ot f_i\\
		& \longmapsto & \dsum_{i,j} rp_i\ot q_j \ot ((p\ot q)\mapsto g_j(f_i(p)q))\\
		& \longmapsto &\left[  p\ot q\mapsto \dsum_{i,j}rp_i\ot q_jg_j(f_i(p)q)\right] = \left[ p\ot q \mapsto \dsum_{i=1}^nrp_i\ot f_i(p)q\right] \\
		& = & \left[ p\ot q\mapsto \dsum_{i=1}^nrp_if_i(p)\ot q\right] = [p\ot q \mapsto rp\ot q].
	\end{eqnarray*}
	
	Analogously, $\psi':R\longrightarrow \End(_RP\ot Q),$ determined by
	\begin{equation*}
	 \psi'(r)(p\ot q)=p\ot qr,
	\end{equation*}
	is onto. 
	Thus, $[P\ot_RQ] \in \Pics(R)$. Clearly, $[R]$ is the unity element  of  $\Pics(R)$. Therefore, $\Pics(R)$ is a monoid. 
	
	\end{dem}

%

\begin{obs}
{It is immediately  seen that} 	$	\mathcal{U}(\Pics(R))=\Pic(R). $
\end{obs}

{With respect to Example~\ref{exe:ParInvMod} note that if the $eR$-bimodule $Q$ is an inverse for the $eR$-bimodule $P,$ then considering $P$ and $Q$ as $R$-bimodules it is easily seen that
$[P][Q][P]=[Q]$ and $[Q][P][Q]=[P]$ in $\Pics(R),$ so that the $R$-bimodule $Q$ is a partial (or weak) inverse of the $R$-bimodule $P.$
}

Analogously to the  case of the Picard groups, if $R$ is a $K$-algebra, then  $\Pics_K(R)$ stands for  the 
 subsemigroup of those  isomorphism classes 
$[P]\in \Pics(R),$ in which the bimodules $P$ are central over $K.$  If $R$ is a commutative ring, then 
$\Pics(R)=\Pics_R(R)$ consists of the isomorphism classes of the finitely generated projective  central 
$R$-bimodules of rank less than or equal to one (see \cite[Proposition 3.6]{DoPaPi}). Moreover, by 
\cite[Proposition  3.8]{DoPaPi}, we have that  $\Pics(R)$ is an inverse  semigroup, with  $[P]^*=[P^*]$ 
for all  $[P]\in \Pics(R)$. 


\subsection{{ The bimodule relation $M|N$ and a tensor product commuting  isomorphism}}
\label{section: modulos similares}

{The central  idempotents of the ring  $R$ will be important for us.  We shall use the following  easy   isomorphism. Let   $e$ be a central  idempotent in $R$ and $M$ an $R$-bimodule. Write $eM=\{m\in M; em=m\}$. Since $e$ is central in $R$, then $eM$ is an  $R$-subbimodule of $M$. Moreover, we have the following  isomorphism of $R$-bimodules:
\begin{equation}
	\begin{array}{c c l}
		Re\ot_RM& \longrightarrow & eM,\\
		r\ot m & \longmapsto & rm, 
	\end{array} \label{isomorfismocomidempotentes}
\end{equation}
whose   inverse is given by   
\begin{equation*}
	\begin{array}{c c l}
		eM & \longrightarrow & Re\ot_R M,\\
		m & \longmapsto & e\ot m.
	\end{array}
\end{equation*}
Analogously,  $Me=\{m \in M; \ me=m\}$ is an  $R$-subbimodule of $M$ and we have the    $R$-bimodule isomorphism:
\begin{equation*}
	\begin{array}{c c l}
		M\ot_RRe& \longrightarrow & Me\\
		m\ot re & \longmapsto & mre
	\end{array}.
\end{equation*}
Let $M$ and $N$ be $R$-bimodules, $F:M\longrightarrow N$ an  $R$-bimodule isomorphism  and $e$ a central idempotent in $R$. The restriction of $F$ to the  $R$-subbimodule $eM$ is also an $R$-bimodule isomorphism between  $eM$ and $eN,$ which, with a slight abuse of notation, will be denoted by the same symbol $F.$  }

{Let $M$ and $N$ be $R$-bimodules. We shall write  $M|N$ if $M$ is isomorphic, as an $R$-bimodule, to a direct summand 
of some direct power  of $N$, that is, if there exists  an $R$-bimodule $M'$ such that  $N^{(n)}\simeq M\oplus M'$, for some   $n \in \N$. This is a  reflexive and transitive relation which is compatible with the tensor product, in the sense that,
if    $M|N$ and $Q$ is an  $R$-bimodule, then 
\begin{equation}
	(M\ot_RQ)|(N\ot_RQ) \ \ \mbox{and} \ \ (Q\ot_RM)|(Q\ot_RN). \label{combatibilidadecomopt}
\end{equation}}

 \begin{obs}\label{MdivideReRdivideM} We immediately have:
	\begin{itemize}
		\item[(i)] If $M_R|R_R$, then $M$ is a finitely generated projective right  $R$-module.
		
		\item[(ii)] If $R_R|M_R$, then $M_R$ is a generator of the category of the right $R$-modules.
	\end{itemize}
\end{obs}
%
%
%
%
%

\vspace{0.3cm}
\begin{lem}\label{fgdedivide} Let $M$ and $N$  be $R$-bimodules.  Then $M|N$ if and only if there 
exist  $R$-bimodule homomorphisms $f_i:M\rightarrow N$ and $g_i:N\rightarrow M$,  $i=1,2,...,n$, such that 
$\sum_{i=1}^ng_i\circ f_i=Id_M$.
\end{lem}   
\begin{dem}
	Indeed, if $M|N$, then there exist $n \in \N$ and an $R$-bimodule  $T$ such that $N^{(n)}=M\oplus T$. Consider 
	$$f_i:M\stackrel{\iota}{\hookrightarrow} N^{(n)}\stackrel{\pi_i}{\twoheadrightarrow} N \ \ \ \mbox{and} \ \ \ g_i:N\stackrel{\iota_i}{\hookrightarrow}N^{(n)}\stackrel{\pi}{\twoheadrightarrow}M,$$
	for $i=1,2,...,n$,  where $\pi_i$ and $\iota_i$ are canonical projections and injections. 
Then, $f_i$ and $g_i$ are $R$-bilinear maps and  
	\begin{equation*}
		\dsum_{i=1}^ng_i\circ f_i =  \dsum_{i=1}^n(\pi\circ \iota_i\circ \pi_i\circ \iota)=\pi\circ\dsum_{i=1}^n(\iota_i\circ\pi_i) \circ \iota
		 =  \pi\circ id_{N^{(n)}}\circ \iota= Id_M. 
	\end{equation*}
	
	Conversely, suppose that there exist $f_i:M\rightarrow N$ and $g_i:N\rightarrow M$, for $i=1,2,...,n$, such
 that $\dsum_{i=1}^ng_i\circ f_i=Id_M$. Consider the exact sequence
	\begin{equation}
	\xymatrix{1 \ar[r] & M\ar[r]^{f} & N^{(n)}\ar[r] & \dfrac{N^{(n)}}{\mbox{Im}(f)}\ar[r] & 1}, \label{seqexataMdivideN}
	\end{equation}
	where $f:M\longrightarrow N^{(n)}$ is defined by $f(m)=(f_1(m),...,f_n(m))$.  Let $g:N^{(n)}\rightarrow M$ be  defined by
	$$g(x_1,...,x_n)=\dsum_{i=1}^ng_i(x_i).$$
	Then,
	\begin{equation*}
	(g\circ f)(m)  =  g(f_1(m),...,f_n(m))=\dsum_{i=1}^ng_i(f_i(m))=m,
	\end{equation*}
	and the sequence (\ref{seqexataMdivideN}) splits. Therefore, $M$ is a direct summand of  $N^{(n)}$.
	\end{dem}

{\begin{obs}\label{MdivideRimplicaemZbimodulocentral}  If $M$ is an  $R$-bimodule such that   $M|R$, then  $M$ is a central $\Z$-bimodule. 
\end{obs}
\noindent Indeed, there exist $R$-bilinear maps $f_i:M\longrightarrow R$  $g_i:R\longrightarrow M$,  $i=1,2,...,n$, with $\dsum_{i=1}^ng_if_i=Id_M$. Take $m \in M$ and $r \in \Z .$ Then
\begin{equation*}
	mr=\dsum_{i=1}^ng_i(f_i(m))r=\dsum_{i=1}^ng_i(\underbrace{f_i(m)}_{\in R}r)=\dsum_{i=1}^ng_i(rf_i(m))=\dsum_{i=1}^nrg_i(f_i(m))=rm.
\end{equation*}
Hence, $M$ is  central over $\Z .$}

\begin{pro}\cite[Corollary 3]{miyashita1973exact}\label{isomorfismoT} Let $M$ and $N$ be $R$-bimodules. If $M|R$ and $N|R$, then there is an $R$-bimodule isomorphism $M\ot_RN\simeq N\ot_RM$ given by
$$\begin{array}{c c c l}
T_{M,N}: & M\ot_RN & \longrightarrow & N\ot_RM\\
& x\ot y & \longmapsto & \dsum_{i=1}^nf_i(x)y\ot g_i(1),
\end{array}$$
where $f_i:M\rightarrow R$, $g_i:R\rightarrow M$, for $i=1,2,...,n$, are $R$-bilinear maps,     such that  $\sum_{i=1}^ng_if_i=Id_M,$ as given by Lemma~\ref{fgdedivide}.
\end{pro}

\begin{obs}\label{remarkT}
The isomorphism $T_{M,N}$ can also be defined by 
$$\begin{array}{c c c l}
T_{M,N}: & M\ot_RN & \longrightarrow & N\ot_RM\\
& x\ot y & \longmapsto & \dsum_{j=1}^mg'_j(1)\ot xf_j'(y),
\end{array}$$
where  $f_j':N\rightarrow R$, $g_j':R\rightarrow N$, for $j=1,2,...,m$, are $R$-bilinear maps such that  $\dsum_{j=1}^mg_j'f_j'=Id_N$.
\end{obs}

 Given an $R$-bimodule $M ,$ we write  $C_M(R)=\{m \in M, \ rm=mr \ \mbox{for all }\ r \in R\}$. The set
 $C_M(R)$ is a central $\Z$-bimodule, and we may endow 
 { $R\ot_{\Z} C_M(R)$} with a structure of  an $R$-bimodule 
 via  the $R$-bimodule structure of $R$, that is,
{ $$r_1\cdot (r\ot _{\Z} m)\cdot r_2=r_1rr_2\ot _{\Z }m,$$}
 for all $r,r_1,r_2 \in R$ and $m \in C_M(R)$.

 \begin{lem}\cite[Lemma 2.4]{miyashita1973exact}\label{MdivideBisomorfismocomcentralizador} Let $M$ be an  $R$-bimodule such that  $M|R$ (as bimodules), then there  exists an $R$-bimodule  isomorphism  $M\simeq R\ot_{\Z}C_M(R).$
 	
 \end{lem}
 %

 \begin{cor}\label{TcomZbomodcentral} Let $M$ and $N$ be  $R$-bimodules such that  $M|R$  and $N$ is central  over 
 	$\Z .$ Then there  exists an $R$-bimodule  isomorphism  $M\ot_RN\simeq N\ot_RM$.
 \end{cor}
 \begin{dem} Since $M|R$, by  Lemma~\ref{MdivideBisomorfismocomcentralizador} we have the sequence of $R$-bimodule  isomorphisms: 
 	\begin{eqnarray*}
 		N\ot_RM & \simeq&  N\ot_RR\ot_\Z C_M(R)\simeq N\ot_\Z C_M(R)\\
 		& \simeq & C_M(R)\ot_\Z N\simeq C_M(R)\ot_\Z R\ot_RN\\
 		& \simeq&  M\ot_RN.
 	\end{eqnarray*}
 \end{dem}
 
 \begin{obs}
 	Analogously to  Proposition~\ref{isomorfismoT}, the isomorphism in Corollary~\ref{TcomZbomodcentral} is given by  
 	$$\begin{array}{c c c l}
 		T_{N,M}: & N\ot_RM & \longrightarrow & M\ot_RN,\\
 		& n\ot m & \longmapsto & \dsum_{i=1}^ng_i(1)\ot nf_i(m),
 	\end{array}$$
 	where $f_i:M\longrightarrow R$ and  $g_i:R\longrightarrow M$ are as in  Lemma~\ref{fgdedivide}.
 \end{obs}
 
 In the following examples we construct  the isomorphism  $T_{-,-}$ in  some particular cases, which will be used later.
 
 \begin{exe}\label{execomRe1} Let $e$ be a central  idempotent in $R .$ Then $Re$ is a direct summand  of $R$, in particular, $Re|R$. If $M$ is a central  $\Z$-bimodule, then there exists  an $R$-bimodule isomorphism   $Re\ot_RM\simeq M\ot_RRe$. Observe that the canonical  inclusion  and the canonical  projection  $i:Re\longrightarrow R$ and $\pi:R\longrightarrow Re$ satisfy  the conditions of  Lemma~\ref{fgdedivide} and we have the isomorphisms  
 	\begin{center}
 		\begin{tabular}{c c c }
 			$\begin{array}{c c l}
 				Re\ot_RM & \longrightarrow & M\ot_RRe,\\
 				re\ot m & \longrightarrow & rm\ot e
 			\end{array}$ & and  & $\begin{array}{c c l}
 				M\ot_RRe & \longrightarrow & Re\ot_RM,\\
 				m\ot re & \longrightarrow & e\ot mr .
 			\end{array}$ 
 		\end{tabular}
 	\end{center} 
 \end{exe}

 \begin{exe}\label{TcomMeNisomorfosaR} Let $M$ and $N$ be $R$-bimodules such that  $N|R$ and there exists  an $R$-bimodule  isomorphism  $f:M\longrightarrow R$. In particular, $M|R$. Observe that the isomorphism  $f:M\longrightarrow R$ and its inverse $f^\m:R\longrightarrow M$ satisfy the condition of Lemma~\ref{fgdedivide} and, consequently, the isomorphism  $T_{M,N}$ can be written as 
 	$$\begin{array}{c c c l}
 		T_{M,N}: & M\ot_RN & \longrightarrow & N\ot_RM, \\
 		& m\ot n & \longmapsto & f(m)n\ot f^{-1}(1). 
 	\end{array}$$
 \end{exe}

\begin{exe}\label{execomRe2} Let $e$ be a central idempotent of $R$ and $M$ an $R$-bimodule such that $M\simeq Re$. Since $Re$ is a direct summand of $R$, then $Re|R$ and, by transitivity, $M|R$. If $N$ is an $R$-bimodule such that  $N|R$, then there exists an  $R$-bimodule isomorphism  $M\ot_RN\simeq N\ot_RM$. Let $h:M\rightarrow Re$ be an $R$-bimodule isomorphism and consider  $i:Re\longrightarrow R$ and $\pi:R\longrightarrow Re,$ the  canonical projection and  injection, respectively. Let $f:M\rightarrow R$ and $g:R\rightarrow M$ be defined by $f=i\circ h$ and $g=h^\m\circ \pi$. It is easy 
to see that $gf=Id_M$ and 
	$$\begin{array}{c c l}
	M\ot_RN & \longrightarrow & N\ot_RM,\\
	m\ot n & \longmapsto & f(m)n\ot g(1).
	\end{array} $$

\end{exe}

%

\subsection{The group $\p(S/R)$}\label{sec:Group P}
\label{sec: pSR}

Let $R\subseteq S$ be an extension of rings with the same unity element.  Denote by $\mathcal{M}(S/R)$ the set  of the elements    $\xymatrix@C=1.2cm{ P\ar@{=>}[r]|-{[\phi]}& X}$, where $P$ is an  $R$-bimodule, $X$ is an  $S$-bimodule and $\phi:P\longrightarrow X$ is an $R$-bilinear map such that the maps  
\begin{equation}
	\begin{tabular}{c c c}
		${\begin{array}{c c c l}
				\bar{\phi_r}: & P\otimes_RS & \longrightarrow & X ,\\
				& p\otimes_Rs & \longrightarrow & \phi(p)s \end{array}}$ & \ \ \mbox{and} \ \ & \ \ \ 	${\begin{array}{c c c l}
				\bar{\phi_l}:&	S\otimes_RP & \longrightarrow & X ,\\
				&	s\otimes_Rp & \longrightarrow & s\phi(p) \end{array}} $
	\end{tabular} \label{phirephil}
\end{equation}
are isomorphisms of $R$-$S$-bimodules and $S$-$R$-bimodules, respectively.  Given   $\xymatrix@C=1.2cm{ P\ar@{=>}[r]|{[\phi]}& X}$ and $\xymatrix@C=1.2cm{ Q\ar@{=>}[r]|{[\psi]}& Y}$ in $\mathcal{M}(S/R)$, define the element  $\xymatrix@C=1.9cm{ P\ot_RQ\ar@{=>}[r]|-{[\phi\ot \psi]} & X\ot_S Y,}$ where 
$$\begin{array}{c c c l}
	\phi\ot\psi: & P\ot_RQ& \longrightarrow & X\ot_S Y,\\
	& p\ot q& \longmapsto & \phi(p)\ot\psi(q).
\end{array}$$

A morphism from $\xymatrix@C=1.2cm{ P\ar@{=>}[r]|{[\phi]}& X}$ to $\xymatrix@C=1.2cm{ Q\ar@{=>}[r]|{[\psi]}& Y}$  is a pair $(\al, \beta)$, where $\al:P\longrightarrow Q$ is $R$-bilinear, $\beta:X\longrightarrow Y$ is $S$-bilinear and the  diagram 
$$\xymatrix{P\ar[rrr]^{\phi}\ar[dd]_{\al} & & &  X\ar[dd]^{\beta}\\
	& & & \\
	Q\ar[rrr]_{\psi} & &  & Y}$$
is commutative. If $\al$ and $\beta$ are isomorphisms of $R$-bimodules and $S$-bimodules, respectively, then  $(\al, \beta)$ is an  isomorphism between  $\xymatrix@C=1.2cm{ P\ar@{=>}[r]|{[\phi]}& X}$ and $\xymatrix@C=1.2cm{ Q\ar@{=>}[r]|{[\psi]}& Y.}$

We denote by $\p(S/R)$ the set of the isomorphism classes  of elements $\xymatrix@C=1.2cm{ [P]\ar@{=>}[r]|{[\phi]}& [X]}$,  where $[P]\in \Pic(R)$ and $[X]\in \Pic(S)$. { By a slight abuse of notation  $\xymatrix@C=1.2cm{ P\ar@{=>}[r]|{[\phi]}& X}$ will also denote the isomorphim class of an element  $\xymatrix@C=1.2cm{ P\ar@{=>}[r]|{[\phi]}& X}$ of
	$\mathcal{M}(S/R).$}  We recall the following fact. 

\begin{pro}\label{grupopSR} \cite[Lemma 3.1]{yoichi1971galois}. Let $R\subseteq S$ be an extension of rings with the same unity.  Let $[P]\in \Pic(R)$, $[X]\in \Pic(S)$ and $\phi:P\longrightarrow X$ be an $R$-bilinear map.   Suppose that  $\bar{\phi}_l$ is an isomorphism. Then 
	\begin{itemize}
		\item[(i)]  $\bar{\phi}_r$ is an  isomorphism. 
		\item[(ii)]  $\phi$ is injective.
		\item[(iii)]  $\xymatrix@C=1.2cm{ [P^*]\ar@{=>}[r]|{[\phi^*]}& [X^*]} \in \p(S/R)$, where $\phi^*(f)(s\phi(p))=sf(p),$ with $s\in S,$  $p \in P$. 
		\item[(iv)]  $\xymatrix@C=1.2cm{ [^*P]\ar@{=>}[r]|{[\phi^+]}& [X^*]} \in \p(S/R)$, where $\phi^+(g)(\phi(p)s)=g(p)s,$ with $s\in S ,$  $p \in P$. 
		\item[(v)]  $(\xymatrix@C=1.2cm{ [P^*]\ar@{=>}[r]|{[\phi^*]}& [X^*]})=(\xymatrix@C=1.2cm{ [^*P]\ar@{=>}[r]|{[\phi^+]}& [^*X]})$ in $\p(S/R)$.
	\end{itemize}
\end{pro}

\begin{obs}\label{obsphirouphil}
	\begin{itemize}
		\item[(i)] Proposition~\ref{grupopSR} still holds if we interchange  $\bar{\phi}_r$ with  $\bar{\phi}_l.$
		
		\item[$(ii)$] It follows from  Proposition~\ref{grupopSR} that if  $[P]\in \Pic(R)$ and $[X]\in \Pic(S)$, then $\xymatrix@C=1.2cm{ [P]\ar@{=>}[r]|{[\phi]}& [X]}\in \p(S/R),$ provided that  $\bar{\phi}_l$ (or $\bar{\phi}_r$) is an isomorphism.
	\end{itemize}  
\end{obs}

{By \cite[Theorem 1.3]{miyashita1973exact}, $\p(S/R)$ is a group,} { in which the product of the isomorphim class of  $\xymatrix@C=1.2cm{ P\ar@{=>}[r]|{[\phi]}& X} \in \mathcal{M}(S/R)$  by that of  $\xymatrix@C=1.2cm{ Q\ar@{=>}[r]|{[\psi]}& Y} \in \mathcal{M}(S/R)$ is given by the isomorphism class of 
	$$\xymatrix@C=1.9cm{ P\ot_RQ\ar@{=>}[r]|-{[\phi\ot \psi]} & X\ot_SY}\in \mathcal{M}(S/R),$$ and the inverse of the class of  $\xymatrix@C=1.2cm{ P\ar@{=>}[r]|{[\phi]}& X}$ is  that of  
	$\xymatrix@C=1.2cm{ [P^*]\ar@{=>}[r]|{[\phi^*]}& [X^*].}$}

%
%

If $R$ is a  $K$-algebra, we define   $\p_K(S/R)=\{\xymatrix@C=1.2cm{ [P]\ar@{=>}[r]|{[\phi]} & [X]} \in \p(S/R); [P]\in \Pic_K(R)\}$. It is easy to see that  $\p_K(S/R)$ is a subgroup of  $\p(S/R)$. Indeed, it is enough to check that if  $P$ is a central  $K$-bimodule, then so is $P^*.$ Let $k \in K$ and $f\in P^*$, then
\begin{equation*}
(f\cdot k)(p)=f(p)k=k f(p)=f(k p)=f(pk)=(k\cdot f)(p),
\end{equation*}
for all $p \in P$. Therefore, $[P^*]\in \Pic_K(R),$ as desired.


We shall next present  an exact sequence which was given in  \cite[Theorem 1.5]{miyashita1973exact}. We shall use it  in Section~\ref{seq1} { and we include a proof for the reader's convenience.}
Denote by $\mathbf{Aut}_{R-\mbox{rings}}(S)$ the  group of those automorphisms of the ring $S$  which fix each element of $R.$ Then  
$$\begin{array}{c c c l}
\mathcal{F}: & \U(\Z) & \longrightarrow & \mathbf{Aut}_{R-\mbox{rings}}(S),\\
& r & \longmapsto & (s\mapsto rsr^\m)
\end{array}$$
is a group homomorphism.

Given $f \in \mathbf{Aut}_{R-\mbox{rings}}(S)$, we define an $S$-bimodule $S_f$ by putting $S_f=S$ as groups, with the bimodule  actions determined  by:
\begin{equation*}
s'*t=s't \ \mbox{and}  \ \ t*s=tf(s), \ \mbox{for}\ \  s,s' \in S,\ t \in S_f.
\end{equation*}
Then
$$\begin{array}{c c c l}
\mathcal{E}: & \mathbf{Aut}_{R-\mbox{rings}}(S) & \longrightarrow & \p(S/R),\\
& f & \longmapsto & (\xymatrix@C=1.2cm{ [R]\ar@{=>}[r]|{[\iota_f]} & [S_f]}),
\end{array}$$
is a group homomorphism, where 
$\iota_f:R\longrightarrow S_f$ is the  inclusion.

\begin{teo}\label{sequenciaparaodiagramadoprimeiromorfismo} The sequence
$$\U(\Z)\stackrel{\mathcal{F}}{\longrightarrow}\mathbf{Aut}_{R-\mbox{rings}}(S)\stackrel{\mathcal{E}}{\longrightarrow}\p(S/R)\stackrel{\vartheta}{\longrightarrow}\Pic(R),$$
where $\vartheta(\xymatrix@C=1.2cm{ [P]\ar@{=>}[r]|{[\phi]} & [X]})=[P]$, is exact.
\end{teo}
\begin{dem} Let us check the   exactness at $\mathbf{Aut}_{R\mbox{-rings}}(S)$. Let $f \in \mathbf{Aut}_{R\mbox{-rings}}(S)$ be such that  $(\xymatrix@C=1.2cm{ [R]\ar@{=>}[r]|{[\iota_f]} & [S_f]})=(\xymatrix@C=1.2cm{ [R]\ar@{=>}[r]|{[\iota]} & [S]})$ in $\mathcal{P}(S/R)$. Then we have the following commutative diagram: 
$$\xymatrix{R\ar[rr]^{\iota}\ar[dd]_{\beta} & & S\ar[dd]^{\alpha}\\
	& & \\
	R\ar[rr]_{\iota_f} & & S_f}$$
where $\alpha$ is an $S$-bimodule  isomorphism. By the  commutativity of the diagram, we have that $R=\al(R)$. Let $d=\alpha^{-1}(1)\in R$. Since $\alpha$ is an $S$-bimodule  isomorphism and $f(r)=r$ for  $r\in R$, it follows that $d \in \U(\Z)$. Moreover,
$$df(s)=d*s=\alpha^\m(1)*s=\alpha^\m(s)=s\alpha^\m(1)=sd, $$
for all $s \in S$. Hence, $f(s)=d^\m sd=\mathcal{F}(d^\m)(s)$ for all   $s \in S.$ This means that  $\mathcal{F}(d^\m)=f$ and,  consequently, $f \in \mbox{Im}(\mathcal{F})$.

Let $f \in \mbox{Im}(\mathcal{F}).$ Then there exists  $r \in \U(\Z)$ such that $f(s)=rsr^\m$ for all  $s\in S$. Let $\om:S\longrightarrow S_f$ be defined by 
$\om(s)=sr^\m$. Clearly,  $\om$ is a left {$S$-linear  isomorphism}. If $s,s' \in S$, then
\begin{equation*}
	\om(ss')  =  ss'r^\m=sr^\m rs'r^\m=sr^\m f(s')=\om(s)*s'.
\end{equation*} 
Therefore, $\om$ is an $S$-bimodule  isomorphism. Consider the $R$-bimodule isomorphism  $\lambda:R\longrightarrow R$, defined by  $\lambda(x)=xr^\m$. Then the   diagram
\begin{equation*}
	\xymatrix{ R \ar[rrr]^{i}\ar[dd]_{\lambda} & & &  S\ar[dd]^{\om}\\
		& & & \\
		R\ar[rrr]_{i_f} & & & S_f}
\end{equation*}
is commutative. Hence, $\mathcal{E}(f)=(\xymatrix@C=1.2cm{ [R]\ar@{=>}[r]|{[\iota_f]} & [S_f]})=(\xymatrix@C=1.2cm{ [R]\ar@{=>}[r]|{[\iota]} & [S]})$ in $\p(S/R)$. Thus $f\in \ker(\mathcal{E})$.

For the  exactness  at $\p(S/R)$, take  $(\xymatrix@C=1.2cm{ [P]\ar@{=>}[r]|{[\phi]} & [X]}) \in \mathcal{P}(S/R)$ with $[P]=[R]$ in $\mathbf{Pic}(R)$. Then there exists an $R$-bilinear   isomorphism  $\gamma:R\longrightarrow P$. We have the following isomorphisms  of  $(R,S)$-bimodules and $(S,R)$-bimodules
$$\al: S\stackrel{\simeq}{\longrightarrow}R\ot_RS\stackrel{\gamma\ot_RS}{\longrightarrow} P\ot_RS\stackrel{\bar{\phi}_r}{\longrightarrow} X,$$
$$\beta: S\stackrel{\simeq}{\longrightarrow}S\ot_RR\stackrel{S\ot_R\gamma}{\longrightarrow} S\ot_RP\stackrel{\bar{\phi}_l}{\longrightarrow} X.$$
Observe that $\al(s)=\phi(\gamma(1))s$ and $\beta(s)=s\phi(\gamma(1))$, for $s\in S$. Define $f=\beta^{-1}\circ \al:S\longrightarrow S$. Then
\begin{equation}
	f(s)\phi(\gamma(1))=\beta(f(s))=\al(s)=\phi(\gamma(1))s, \label{eq9}
\end{equation}
for all   $s \in S$. Hence,  for  $s,t \in S$ we obtain 
\begin{eqnarray*}
	\al(st) & = & \phi(\gamma(1))st \stackrel{(\ref{eq9})}{=}  f(s)\phi(\gamma(1))t \\
	& \stackrel{(\ref{eq9})}{=} &f(s)f(t)\phi(\gamma(1))\\
	& = & \beta(f(s)f(t)).
\end{eqnarray*}
Applying  $\beta^{-1}$ we get $f(st)=f(s)f(t).$ Moreover, since $\phi$ and $\gamma$ are $R$-bilinear, we see that  
$$\alpha(r)=\phi(\gamma(1))r=\phi(\gamma(r))=r\phi(\gamma(1))=\beta(r),$$
for all $r \in R$. Applying $\beta^{-1},$ we obtain that  $f(r)=r$, for all  $r \in R$. Hence, $f \in \mathbf{Aut}_{R\mbox{-rings}}(S)$.  { Since $S$ and $S_f$ coincide as sets, we may consider $\beta:S_f\longrightarrow X .$} Obviously,  $\beta$ is left  $S$-linear. For $t \in S$ and $s \in S_f$ we have that
\begin{eqnarray*}
	\beta(s* t) & = & \beta(sf(t))=sf(t)\phi(\gamma(1))\\
	& \stackrel{(\ref{eq9})}{=} & s\phi(\gamma(1))t=\beta(s)t. 
\end{eqnarray*}
Consequently, $\beta:S_f \longrightarrow X$ is an $S$-bilinear isomorphism. Moreover, since $\beta(r)=(\phi\circ \gamma)(r)$, for all $r \in R$, the following diagram is commutative
$$\xymatrix{R\ar[rrr]^{i_f}\ar[dd]_{\gamma}& & & S_f\ar[dd]^{\beta}\\
	& & \\
	P\ar[rrr]_{\phi} & & & X.}$$
Thus, $\mathcal{E}(f)=(\xymatrix@C=1.2cm{ [R]\ar@{=>}[r]|{[\iota_f]} & [S_f]})=(\xymatrix@C=1.2cm{ [P]\ar@{=>}[r]|{[\phi]} & [X]})$ em $\mathcal{P}(S/R)$. Therefore, $\ker(\vartheta)\subseteq \mbox{Im}(\mathcal{E})$. The opposite inclusion is immediate. 
\end{dem}

\section{Partial generalized crossed products}\label{sec:ParGenCrosProd}

\subsection{Partial {actions, partial representations} and partial cohomology}\label{sec:ParAcParRepParHom}

{ We shall need a concept of a partial group action on a semigroup which we give below and which  is more general than that  considered so far (see, in particular, \cite{dokuchaev2015partial}).}   

\begin{defi}\label{def:ParAcOnSemigr}
	Let $G$ be a group and $S$ a semigroup. A \underline{partial action} $\alpha $ of $G$ on $S$ is a family of subsemigroups $S_x$,  $(x\in G)$, and  semigroup isomorphisms  $\al_x:S_{x^\m}\longrightarrow S_x$ which satisfy the following  conditions:
	\begin{itemize}
		\item[(i)] $S_1=S$ and $\al_1=Id_S$,
		\item[(ii)] $\al_y^\m(S_y\cap S_{x^\m})\subseteq S_{(xy)^\m}$,
		\item[(iii)] $\al_x\circ \al_y(s)=\al_{xy}(s)$, for each $s \in \al_y^\m(S_y\cap S_{x^\m}).$
	\end{itemize} 
\end{defi}
We shall write for simplicity  $\al=(S_x,\al_x).$  
As it can be seen in   \cite{dokuchaev2005associativity}, conditions   (ii) and (iii) of the above definition imply   that $\al_x^\m=\al_{x^\m}$ and  
\begin{equation}
	\al_x(S_{x^\m}\cap S_{y})=S_x\cap S_{xy}, \ \ \mbox{for all} \ x,y \in G. \label{alxemSxinversoSy}
\end{equation}

{ We say that the partial action   $\al=(S_x,\al_x)$ is {\it unital} if each   $S_x$ is an ideal in $S$ generated by an   idempotent which is central in $S$, that is, $S_x=S1_x$,  $1_x\in \Z (S),$ for all $x \in G.$ In this case,}  $S_x\cap S_y=S1_x1_y$ and (\ref{alxemSxinversoSy}) implies that  $\al_x(1_y1_{x^\m})=1_x1_{xy},$ for all   $x,y \in G.$ As a consequence, 
\begin{equation*}
	\al_{xy}(s1_{y^\m x^\m})1_x=\al_x(\al_y(s1_{y^\m})1_{x^\m}), \ \ \mbox{for all } \ x,y \in G \ \mbox{and} \ s \in S.
\end{equation*}
{ The subsemigroup of the  invariants of $S$ with respect to the unital partial action   $\al$ is defined by}
\begin{equation*}
	S^\al=\{s\in S; \ \al_x(s1_{x^\m})=s1_x, \ \mbox{for all } x\in G \}.
\end{equation*}

%

We recall the next:
\begin{defi} \label{defipartialrepresentation} A partial  representation of $G$ in a monoid $S$ is a map  
	$$\begin{array}{c  c c l}
		\et:& G& \longrightarrow & S\\
		& x   & \longmapsto & \et_x
	\end{array}$$
	which satisfies the following   properties:
	\begin{itemize}
		\item[(i)] $\et_{1_G}=1_S$,
		\item[(ii)] $\et_x\et_y\et_{y^{-1}}=\et_{xy}\et_{y^{-1}}$, for all  $x,y \in G$. 
		\item[(iii)] $\et_{x^{-1}}\et_x\et_y=\et_{x^{-1}}\et_{xy}$,  for all  $x,y \in G$.
	\end{itemize}
\end{defi}
It follows from Definition \ref{defipartialrepresentation} that
\begin{equation}
\et_x\et_{x^\m}\et_x=\et_x,   \ \  \mbox{for all } \  x \in G. \label{etx}
\end{equation} 
 By \cite{dokuchaev2000partial} we know that the {  $\e_x$ are   idempotents such that  
	\begin{equation}
		\e_x \e_y = \e_y \e_x \ \  \mbox{and} \ \ \et_x\e_y=\e_{xy}\et_x  \ \  \mbox{for all } \ x,y\in G. \label{Txeeycomutam}
\end{equation}	}
In particular, 
\begin{equation}
	\e_x\et_x=\et_x \ \ \mbox{ e} \ \ \  \et_x\e_{x^\m}=\et_x, \ \mbox{for all} \  x\in G. \label{exetx}
\end{equation}
Moreover, by (\ref{etx})
\begin{equation}
\et_x\et_y=\et_x\et_{x^\m}\et_x\et_y=\et_x\et_{x^\m}\et_{xy}=\e_x\et_{xy}, \ \ \mbox{for all} \ \ x,y\in G. \label{extxty=extxy}
\end{equation}
Applying (\ref{Txeeycomutam})  to (\ref{extxty=extxy}), we obtain  
\begin{equation}
	\et_x\et_y=\e_x\et_{xy}=\et_{xy}\e_{y^\m},  \ \mbox{{ for all} } \ x,y \in G. \label{txtyeym=txyeym}
\end{equation}

{ The need of a more general concept of a partial group action on a semigroup is justified by the following fact, which will be useful for us.} 
\begin{pro}\label{proparaacaoparcialsobrePics} Let $S$ be a   monoid and let $\et:G\longrightarrow S$ be a partial representation. Write  $S_x=\e_x S\e_x$,  $(x\in G)$. Then,  $\al^*=(S_x,\al^*_x)$, where 
	$$\begin{array}{c c c l}
		\al_x^*: & S_{x^{-1}} & \longrightarrow & S_x\\
		& s & \longmapsto & \et_xs\et_{x{^{-1}}},
	\end{array}$$ 
	is a partial action of $G$ on $S$.	
\end{pro}
\begin{dem} Clearly,  $S_1=S$ and $\al_1^*=Id_S$.  Observe that $s \in S_x $ if and only if  $\e_xs\e_x=s$. Indeed, obviously, if $s=\e_xs\e_x$, then $s\in S_x$. On the other hand, if $s \in S_x$, then  $s=\e_xs'\e_x$, for some $s'\in S$, and $\e_xs\e_x=\e_x\e_xs'\e_x\e_x=\e_xs'\e_x=s,$ as desired. 
	
	By (\ref{exetx}), we have \begin{equation*}
		\al^*_x(s)=\et_xs\et_{x^\m}=\e_x\et_xs\et_{x^\m}\e_x=\e_x\al_x^*(s)\e_x,
	\end{equation*}
	and thus,  $\al_x^*(s)\in S_x$. For $s,s' \in S_{x^\m}$, we see that
	\begin{equation*}
		\al_x^*(ss')=  \et_xss'\et_{x^\m}=\et_xs\e_{x^\m}s'\et_{x^\m}=\et_xs\et_{x^\m}\et_xs'\et_{x^\m}=\al_x^*(s)\al_x^*(s').
	\end{equation*} 
	Hence, $\al_x^*$ is a semigroup { homomorphism.} 
	Given $s \in S_{x^\m},$ we obtain
	\begin{eqnarray*}
		(\al_{x^\m}^*\circ\al_x^*)(s) & = & \al_{x^\m}^*(\et_x s\et_{x^\m}) = \et_{x^\m}\et_x s\et_{x^\m}\et_{x}\\
		& = & \e_{x^\m}s\e_{x^\m}=s.
	\end{eqnarray*}
	Analogously,   $(\al_x^*\circ\al_{x^\m}^*)(s)=s$, for all  $s \in S_x$. Therefore, $\al_x^*$ is an isomorphism, whose inverse is  $\al_{x^\m}^*$. 
	
	If $s \in S_y\cap S_{x^\m}$, then
	\begin{eqnarray*}
		\al_{y^\m}^*(s) & = & \et_{y^\m}s\et_y \stackrel{s\in S_{x^\m}}{=}  \et_{y^\m}\e_{x^\m}s\e_{x^\m}\et_y\\
		& \stackrel{(\ref{Txeeycomutam})}{=} & \e_{(xy)^\m}\et_{y^\m}s\et_{y}\e_{(xy)^\m}\in S_{(xy)^\m}.
	\end{eqnarray*}
	Hence,  $\al_{y^\m}^*(S_y\cap S_{x^\m}) \subseteq S_{(xy)^{-1}}.$ If $s \in \al_{y^\m}^*(S_y\cap S_{x^\m}) \subseteq S_{(xy)^{-1}}\cap S_{y^\m}$, then using  (\ref{extxty=extxy}) and (\ref{txtyeym=txyeym}) we compute 
	\begin{eqnarray*}
		(\al_x^*\circ \al_y^*)(s) & = & \al_x^*(\et_ys\et_{y^\m})=\et_x\et_ys\et_{y^\m}\et_{x^\m}\\
		& = & \et_{xy}\e_{y^\m}s\e_{y^\m}\et_{y^\m x^\m}\\
		& = & \et_{xy}s\et_{y^\m x^\m}\\
		& = & \al_{xy}^*(s).
	\end{eqnarray*}
		Consequently, $\al^*$ is a partial action. 
\end{dem}

Let $K$ be a commutative ring (or a commutative monoid) and $\al=(D_x,\al_x)$ a unital partial action of $G$ on $K$, where $D_x$ is generated by the central idempotent  $1_x$. An $n$-\textit{cochain}, with $n\in \N$, of $G$ with values in $K$ is a function $f:G^n\longrightarrow K$ such that $f(x_1,...,x_n)\in \U(K1_{x}1_{x_1x_2}...1_{x_1x_2...x_n})$. For $n=0$, we define a $0$-cochain as  an element in $\U(K)$. 

Let $C^n(G,\al,K)$ be the set of all $n$-cochains of $G$ with values in  $K$. Then, $C^n(G,\al,K)$ is an abelian group with the multiplication defined point-wise, that is, 
\begin{equation*}
fg(x_1,...,x_n)=f(x_1,...,x_n)g(x_1,...,x_n), \ \ \mbox{with} \ x_1,...,x_n \in G.
\end{equation*}
Clearly,  $I(x_1,...,x_n)=1_{x_1}1_{x_1x_2}...1_{x_1x_2...x_n}$ is the unit for this multiplication and \begin{center}
{ $f^\m(x_1,...,x_n)=f(x_1,...,x_n)^\m \in 
	\U(K1_{x}1_{x_1x_2}...1_{x_1x_2...x_n})$,} for $x_1,...,x_n \in G$.
\end{center} 
%

\begin{pro}\label{prodeltan} \cite[Proposition 1.5]{dokuchaev2015partial} The map  
$\delta^n:C^n(G,\al,K)\longrightarrow C^{n+1}(G,\al,K)$ defined by 
	\begin{eqnarray}
	(\delta^nf)(x_1,...,x_{n+1}) & = & \al_{x_1}(f(x_2,...,x_{n+1})1_{x_1^\m})\prod_{i=1}^nf(x_1,...,x_ix_{i+1},...,x_{n+1})^{(-1)^i}\nonumber\\
	& & f(x_1,...,x_n)^{(-1)^{n+1}}. \label{definicaodedeltapequeno}
	\end{eqnarray}
	where the inverse elements are taken in the corresponding ideals, is a group morphism such that
	\begin{equation*}
	(\delta^{n+1}\circ \delta^n)(f)(x_1,...,x_{n+2})=1_{x_1}1_{x_1x_2}...1_{x_1x_2...x_{n+2}},
	\end{equation*}
	for all $f \in C^n(G,\al,K)$ and $x_1,...,x_{n+1}\in G$.
\end{pro}

\begin{obs}
	If $n=0$, then $\delta^0: \U(K)\longmapsto C^1(G,\al,K)$ is defined by \begin{equation*}
	(\delta^0k)(x)=\al_{x}(k1_{x^\m})k^\m, \ \ \mbox{for each} \ k \in \U(K).
	\end{equation*}
\end{obs} 
We define the groups
\begin{equation*}
Z^n(G,\al,K)=\ker(\delta^n) \ \ \mbox{and} \ \ B^n(G,\al,K)=\mbox{Im}(\delta^{n-1})
\end{equation*}
that are  \textit{the group of the partial $n$-cocycles and that of the partial $n$-coboundaries}, respectively.  
By Proposition~\ref{prodeltan}, we have $B^n(G,\al,K)\subseteq Z^n(G,\al,K)$. Thus, we define the group of the partial $n$-cohomologies by 
\begin{equation*}
H^n(G,\al,K)=\dfrac{Z^n(G,\al,K)}{B^n(G,\al,K)}.
\end{equation*}
For $n=0$ we set $H^0(G,\al,K)=Z^0(G,\al,K)=ker(\delta^0)$.


\begin{exe}
	For $n=0$ we have $$H^0(G,\al,K)=Z^0(G,\al,K)=\{k \in \U(K), \al_x(k1_{x^\m})=k1_x, \ \forall \ x\in G\}$$
	$$B^1(G,\al,K)=\mbox{Im}(\delta^0)=\{f \in C^1(G,\al,K); \exists \ k \in \U(K) \ \ f(x)=\al_x(k1_{x^\m})k^\m, \ \forall \ x\in G\}.$$
	
	For $n=1$,  $(\delta^1f)(x,y)=\al_x(f(y)1_{x^\m})f(xy)^\m f(x)$, where $f \in C^1(G,\al,K)$. Then
	$$Z^1(G,\al,K)=\{ f\in C^1(G,\al,K); \al_x(f(y)1_{x^\m})f(x)=f(xy)1_x, \ \forall \ x,y \in G \}.$$
	$$ {\small B^2(G,\al,K)=\{ f\in C^2(G,\al,K); \exists \ \s\in C^1(G,\al,K),  \   f(x,y)=\al_x(\s(y)1_{x^\m})\s(xy)^\m\s(x), \ \forall \ x,y \in G \}.}$$
	For $n=2$, $$(\delta^2f)(x,y,z)=\al_x(f(y,z)1_{x^\m})f(xy,z)^\m f(x,yz) f(x,y)^\m,$$ where $f \in C^2(G,\al,K)$ and $x,y \in G.$ Then,
	$$Z^2(G,\al,K)=\{f \in C^2(G,\al,K); \al_x(f(y,z)1_{x^\m})f(x,yz)=f(xy,z)f(x,y), \ \forall \ x,y,z \in G \}.$$
	\end{exe}
Let $f,f'\in Z^n(G,\al,K).$ We say that $f$ and $f'$ are \textit{cohomologous} if there is $g \in C^{n-1}(G,\al,K)$ such that $f=f'(\delta^ng)$. In this case,  $[f]=[f']$ in $H^n(G,\al,K)$.  

An $n$-cocycle $f$ is called \textit{normalized} if $$f(1,x_1,...,x_{n-1})=f(x_1,1,...,x_{n-1})=\cdots =f(x_1,...,x_{n-1},1)=1_{x_1}1_{x_1x_2}...1_{x_1...x_{n-1}},$$ 
 for all $x_1, x_2, \ldots, x_{n-1} \in G.$

\begin{obs}
Any partial $1$-cocycle $\s$ is normalized. Indeed, we have  $\al_x(\s_y1_{x^\m})\s_{xy}^\m\s_x=1_x1_{xy},$
 for all $x,y \in G.$ 
  In particular, taking $x=y=1,$ we obtain $\s_1^2=\s_1,$ which  implies 
that $\s_1=1$. Therefore, $\s$ is normalized. 
By \cite[Remark 2.6]{dokuchaev2015partial} if  $\s\in Z^2(G,\al, K)$, then there is  
a normalized $\widetilde{\s}\in Z^2(G,\al,K)$  such that $\s=\widetilde{\s}(\delta^1\beta)$, 
for some $\beta \in C^1(G,\al, K)$. Therefore, $\s$ is cohomologous to a normalized partial  $2$-cocycle.	
	\end{obs}

Observe also that by \cite[Remark 2.6]{dokuchaev2015partial},  if $\s\in Z^2(G,\al,\Z)$, then there exists a normalized  $\widetilde{\s}\in Z^2(G,\al,\Z)$ such that  $\s=\widetilde{\s}(\delta^1\beta)$, for some  $\beta \in C^1(G,\al,\Z)$. Hence, $\s$ is cohomologous to a normalized partial  $2$-cocycle.

\subsection{Unital partial representations  $G \to \Pics(R)$}\label{sec:ParRep G to Pics}

{Consider a partial representation $$\T:  G \longrightarrow  \Pics(R)$$ of a group $G$  into the monoid  
$\Pics(R).$ For each $x\in G$ fix a representative $\T _x$ of the isomorphism class $\T (x) \in  \Pics(R),$ i.e. 
$\T (x) = [ \T _x ].$}
 We say that $\T$ is \textit{unital} if  $[\e_x]=[\T_x][\T_{x^{-1}}]=[R1_x]$, where $1_x$ is a central idempotent of $R$, for all $x \in G$. Thus, there is  an $R$-bimodule isomorphism $\T_x\ot_R\T_{x^\m}\simeq R1_x$, for all $x \in G$.
In the following auxiliary fact we gather some isomorphisms involving $\T _x$ for further use.

\begin{lem}\label{isomorfismoscomT} Let $\T:G\longrightarrow\Pics(R)$ be a unital partial representation   with $\T_x\ot_R\T_{x^\m}\simeq R1_x,$   $(x \in G)$. Then, for all  $x,y\in G$, we have { the following $R$-bimodule isomorphisms:}
	\begin{itemize}
	\item[(i)] $\T_x\ot_RR1_y\simeq R1_{xy}\ot_R\T_x$. In particular,
	\begin{equation}
		R1_x\ot_R \T_x\simeq \T_x \ \mbox{and} \ \ \T_x\ot_R R1_{x^\m}\simeq \T_x,  \ \mbox{for all} \ x\in G. \label{TxeR1x}
	\end{equation}
	\item[(ii)] $\T_x\ot_R\T_y\simeq R1_x\ot_R\T_{xy}$  and $\T_x\ot_R \T_y\simeq \T_{xy}\ot_R R1_{y^\m}$,
	\item[(iii)] $\T_x\ot_R\T_y\ot_R\T_{(xy)^{\m}}\simeq R1_x1_{xy}$, 
	\item[(iv)] $\T_x\ot_R\T_y\simeq \T_x\ot_R\T_y\ot_RR1_{y^\m x^\m}$.
\end{itemize}
\end{lem}
\begin{dem} Item (i) is a  direct consequence of  (\ref{Txeeycomutam}), whereas    (ii) follows from  (\ref{extxty=extxy}) and (\ref{txtyeym=txyeym}). 
	For (iii) we have
	\begin{eqnarray*}
		\T_x\ot_R\T_y\ot_R\T_{(xy)^{\m}} & \stackrel{(ii)}{\simeq} & R1_x\ot_R \T_{xy}\ot_R\T_{(xy)^{-1}}\\
		& \simeq & R1_x\ot_RR1_{xy}\simeq R1_x1_{xy}.
	\end{eqnarray*}
	As to (iv):
	\begin{eqnarray*}
		\T_x\ot_R\T_y\ot_RR1_{y^\m x^\m} & \simeq & \T_x\ot_R\T_y\ot_R\T_{(xy)^{\m}}\ot_R\T_{xy}\\
		& \simeq & \T_x\ot_R\T_y\ot_R\T_{y^{\m}}\ot_R\T_{x^{-1}}\ot_R\T_{xy}\\
		& \simeq & \T_x\ot_R\T_y\ot_R\T_{y^{\m}}\ot_R\T_{x^{-1}}\ot_R\T_{x}\ot_R\T_y\\
		& \simeq & \T_x\ot_RR1_{y}\ot_RR1_{x^{-1}}\ot_R\T_y\\
		& \simeq & \T_x\ot_RR1_{x^{-1}}\ot_RR1_{y}\ot_R\T_y\\
		& \simeq & \T_x\ot_R\T_y.\\
	\end{eqnarray*}
\end{dem}

Using the isomorphism in (\ref{isomorfismocomidempotentes}), and     (i) and (ii) of Lemma~\ref{isomorfismoscomT}, we have:
\begin{equation}
	1_{xy}\T_x\simeq \T_x1_y \ \ \mbox{and}  \ \ \T_x\ot_R\T_y\simeq 1_x\T_{xy}, \ \ \mbox{for all} \ x,y \in G, \label{1xyTxisomorfismo}
\end{equation} as $R$-bimodules.
{In particular,}
 \begin{equation}
1_x\T_x\simeq \T_x \quad \mbox{and} \quad \T_x1_{x^\m}\simeq \T_x, \label{unidadescomTx}
 \end{equation}

\begin{lem}\label{lemTxcomutacom1y} Let $\T:G\longrightarrow \Pics(R)$ be a unital partial representation with  $\e_x=\T_x\ot_R\T_{x^{-1}}\simeq R1_x$, where $\T(x)=[\T_x]$ for all $x\in G$. Then,
	\begin{equation}
		u_x1_y=1_{xy}u_x1_y=1_{xy}u_x, \ \ \mbox{for all } x,y \in G \ \mbox{and}\ u_x\in \T_x. \label{uxe1y}
	\end{equation}
	In particular, $1_{xy}\T_x=\T_x1_y$, for all $x,y \in G$.
\end{lem}
\begin{dem} By (\ref{1xyTxisomorfismo}) there exists an $R$-bimodule isomorphism $\kappa_{x,y}:\T_x1_y\longrightarrow 1_{xy}\T_x$. Given $u_x\in \T_x$, there is $v_x \in \T_x$ such that $\kappa_{x,y}(u_x1_y)=1_{xy}v_x$. Then
	$$\kappa_{x,y}(1_{xy}u_x1_y)=1_{xy}v_x=\kappa_{x,y}(u_x1_y).$$
	Since $\kappa_{x,y}$ is an isomorphism, then $1_{xy}u_x1_y=u_x1_y.$ On the other hand, there exists $v_x'\in \T_x$ such that $\kappa^{-1}_{x,y}(1_{xy}u_x)=v_x'1_y$ and 
	$$\kappa^{-1}_{x,y}(1_{xy}u_x1_y)=v_x'1_y=\kappa^{-1}_{x,y}(1_{xy}u_x).$$
	Thus,  $1_{xy}u_x1_y=1_{xy}u_x$. Therefore,
	$$u_x1_y=1_{xy}u_x1_y=1_{xy}u_x, \ \mbox{for all } \ u_x \in \T_x.$$ \end{dem}


It follows from Lemma \ref{lemTxcomutacom1y} that 
\begin{equation}
1_xu_x=u_x \ \ \mbox{and} \ \ u_x1_{x^\m}=u_x, \ \mbox{for all} \ x\in G, u_x \in \T_x. \label{uxe1x}
\end{equation}
The maps $R1_x\ot_R\T_x\simeq \T_x$ and  $\T_x\ot_RR1_{x^\m}\simeq \T_x$ defined via the right $R$-action and 
the left $R$-action on $\T_x$, respectively, are isomorphisms.  Thus, $\T_x$ is a unital 
$(R1_x,  R1_{x^\m})$-bimodule 
 and  there exist  $R$-bimodule isomorphisms
\begin{equation}
\tau_x:\T_x\ot_R\T_{x^{-1}}\longrightarrow R1_x \ \ \mbox{and} \ \ \tau_{x^\m}:\T_{x^\m}\ot_R\T_x\longrightarrow R1_{x^\m}. \label{isoTxTxm}
\end{equation}

Note that $\T_x\ot_R\T_{x^\m}=\T_x\ot_{R1_{x^\m}}\T_{x^\m}$  as  $(R1_x, R1_{x^\m})$-bimodules. In fact, by (\ref{uxe1x})
$$u_xr\ot_{R1_{x^\m}}u_{x^\m}=u_x1_{x^\m}r\ot_{R1_{x^\m}}u_{x^\m}=u_x\ot_{R1_{x^\m}}r1_{x^\m}u_{x^\m}=u_x\ot_{R1_{x^\m}}ru_{x^\m},$$
for all $u_x \in \T_x, u_{x^\m}\in \T_{x^\m}$ and  $r\in R$.
Thus, we can choose the isomorphisms $\tau_x$ and $\tau_{x^\m}$ such that $(R1_x,R1_{x^\m},\Theta_x,\Theta_{x^\m},\tau_x,\tau_{x^\m})$ is a Morita context, that is, $\tau_x$ and $\tau_{x^\m}$ satisfy: 
\begin{equation}
\tau_x(u_x\otimes u_{x^\m})u_x=u_x\tau_{x^\m}(u_{x^\m}\otimes u_x'), \label{contextomonita1}
\end{equation}
\begin{equation}
\tau_{x^\m}(u_{x^\m}\otimes u_x)u_{x^\m}'=u_{x^\m}\tau_x(u_x\otimes u_{x^\m}'), \label{contextomonita2}
\end{equation}
for all $u_x,u_x' \in \T_x$. For simplicity, we write $\tau_{x}(u_x\ot u_{x^\m})=u_xu_{x^\m}$, where $u_x\in \T_x$ and $u_{x^\m}\in \T_{x^{-1}}$. Thus,  by  (\ref{contextomonita1}) and (\ref{contextomonita2}) we have the associativity:
\begin{equation}
u_xu_{x^\m}u_{x}'=(u_xu_{x^\m})u_x'=u_x(u_{x^\m}u_x')  \label{moritasemtau}
\end{equation}
for $u_x,u_{x}'\in \T_x$ and $u_{x^\m}\in\T_{x^{-1}}$.

\begin{lem} For all $u_x \in \T_x$, with $x \in G$, and $r \in \Z$, the following equality holds:
	\begin{equation}
	(u_x(u_yu_{y^\m})u_{x^\m})u_{xy}ru_{(xy)^\m}=(u_x(u_yru_{y^\m})u_{x^\m})u_{xy}u_{(xy)^\m}. \label{ruxuyuxym1}
	\end{equation}
\end{lem}
\begin{dem}
	Firstly, by (\ref{1xyTxisomorfismo}) we have
	$$\T_x\ot\T_y\ot\T_{(xy)^\m}\simeq R1_{x}\ot \T_{xy}\ot \T_{(xy)^\m}\simeq R1_{x}\ot R1_{xy}\simeq R1_x1_{xy}.$$
  	Let $\varphi$ be the composition of the above  $R$-bimodule isomorphisms.  For each $r \in \Z$ we have $r\varphi(u_x\ot u_y\ot u_{(xy)^\m})=\varphi(u_x\ot u_y\ot u_{(xy)^\m})r$. Since $\varphi$ is $R$-bilinear, it follows that $$\varphi(ru_x\ot u_y\ot u_{(xy)^\m})=\varphi(u_x\ot u_y\ot u_{(xy)^\m}r).$$
	Moreover, since $\varphi$ is an isomorphism, we have
	\begin{equation}
	ru_x\ot u_y\ot u_{(xy)^\m}=u_x\ot u_y\ot u_{(xy)^\m}r \label{ruxuyuxym}
	\end{equation}
	for all $r\in \Z$.
	Now, note that
	$$(u_x(u_yu_{y^\m})u_{x^\m})u_{xy}ru_{(xy)^\m}=\tau_{xy}(\tau_x(u_x\tau_y(u_y\ot u_{y^\m})\ot u_{x^\m})u_{xy}\ot u_{(xy)^\m}).$$
	Rewriting the  right-hand side of the above equality and using (\ref{ruxuyuxym}), we have
	\begin{eqnarray*}
	(u_x(u_yu_{y^\m})u_{x^\m})u_{xy}ru_{(xy)^\m} & = & \tau_{xy}(\tau_x(u_x\tau_y(u_y\ot u_{y^\m})\ot u_{x^\m})u_{xy}r\ot u_{(xy)^\m})\\
	& = & [\tau_{xy}\circ (\tau_x\ot \T_{xy}\ot \T_{(xy)^\m})\circ (\T_x\ot  \tau_y\ot \T_{x^\m}\ot \T_{xy}\ot \T_{(xy)^\m})]\\
	& & (u_x\ot u_y\ot u_{y^\m}\ot u_{x^\m}\ot u_{xy}r\ot u_{(xy)^\m})\\
	& = & [\tau_{xy}\circ (\tau_x\ot \T_{xy}\ot \T_{(xy)^\m})\circ (\T_x\ot  \tau_y\ot \T_{x^\m}\ot \T_{xy}\ot \T_{(xy)^\m})]\\
	& & (u_x\ot u_yr\ot u_{y^\m}\ot u_{x^\m}\ot u_{xy}\ot u_{(xy)^\m})\\
	& = & \tau_{xy}(\tau_x(u_x\tau_y(u_yr\ot u_{y^\m})\ot u_{x^\m})u_{xy}\ot u_{(xy)^\m})\\
	& = & (u_x(u_yru_{y^\m})u_{x^\m})u_{xy}u_{(xy)^\m}.
	\end{eqnarray*}

\end{dem}

For the rest of this section, we will fix a unital partial representation $\T:G\longrightarrow \Pics(R)$ with
 $\T(x)=[\T_x]$ and $\e_x=\T_x\ot_R\T_{x^{-1}}\simeq R1_x$, for all $x\in G$. Denote by $\dsum_{(x)}\om_x\ot \om_{x^\m} \in \T_x\ot_R\T_{x^\m}$ the inverse image of $1_x$ under the isomorphism $\tau_x$, that is, $\dsum_{(x)}\om_x\om_{x^\m}=1_x$. 
\begin{lem}\label{gammatil} There is a group homomorphism
	$$	\begin{array}{c c l}
	\Aut_{R-R}(1_y\T_x) & \stackrel{\widetilde{(-)}}{\longrightarrow} & \U(\Z1_x1_y)\\
	\s & \longmapsto & \widetilde{\s},
	\end{array}$$
	where
	\begin{equation}
	\widetilde{\s}=\dsum_{(x)}\s(1_y\omega_x)\omega_{x^{-1}}. \label{defideftil}
	\end{equation}
 Moreover,
	\begin{equation}
	\s(1_yu_x)=\widetilde{\s}u_x, \ \ \mbox{for all} \ u_x \in \T_x. \label{relacaodefcomftil}
	\end{equation}	
\end{lem}
\begin{dem}	
Firstly,  since $\s$ is $R$-bilinear, we have
	\begin{equation}
		\s(1_yu_x)u_{x^\m} v_xv_{x^\m}=\s(1_yu_xu_{x^\m}v_x)v_{x^\m}, \label{asssigma}
	\end{equation}
	for all $u_x,v_x \in \T_x$. Let us  check first  that $\widetilde{\s}$ does not depend on the choice of the decomposition of $1_x$. Let $\dsum_{(\widetilde{x})}\widetilde{\omega}_x\widetilde{\om}_{x^\m}=1_x,$ with $\widetilde{\om}_x \in \T_x$ and $\widetilde{\om}_{x^\m}\in \T_{x^\m},$ be another decomposition of $1_x.$ Then by (\ref{asssigma})
		\begin{eqnarray*}
			\widetilde{\s} & = & \dsum_{(x)}\s(1_y\om_x)\om_{x^\m}=\dsum_{(x)}\s(1_y\om_x)\om_{x^\m}1_x= \dsum_{(x),(\widetilde{x})}\s(1_y\om_x)\om_{x^\m}{\om}_x\widetilde{\om}_{x^\m}\\
			& = &  \dsum_{(x),(\widetilde{x})}\s(1_y\om_x\om_{x^\m}\widetilde{\om}_x)\widetilde{\om}_{x^\m}= \dsum_{(\widetilde{x})}\s(1_y1_x\widetilde{\om}_x)\widetilde{\om}_{x^\m}= \dsum_{(\widetilde{x})}\s(1_y\widetilde{\om}_x)\widetilde{\om}_{x^\m}.	\end{eqnarray*}
		Clearly  $\widetilde{\s} \in R1_x1_y$. Now, for $r\in R$:
		\begin{eqnarray*}
			\widetilde{\s}r & = & \dsum_{(x)} \s(1_y\om_x)\om_{x^\m}r=\dsum_{(x)} \s(1_y\om_x)\om_{x^\m}r1_x= \dsum_{(x),(\widetilde{x})} \s(1_y\om_x)\om_{x^\m}r\widetilde{\om}_x\widetilde{\om}_{x^\m}\\
			& = & \dsum_{(x),(\widetilde{x})} \s(1_y\om_x\om_{x^\m}r\widetilde{\om}_x)\widetilde{\om}_{x^\m}=\dsum_{(\widetilde{x})} \s(1_y1_xr\widetilde{\om}_x)\widetilde{\om}_{x^\m}\\
			& = &  \dsum_{(\widetilde{x})} r\s(1_y\widetilde{\om}_x)\widetilde{\om}_{x^\m}=r\widetilde{\s}.
		\end{eqnarray*}
		Thus, $\widetilde{\s}\in \Z1_x1_y$. It remains to verify that $\widetilde{\s}$ is an invertible element of $\Z1_x1_y$. Given $\s^\m \in \Aut_{R-R}(1_y\T_x),$ take $\widetilde{\s^\m} \in \Z1_x1_y.$ Then
		\begin{eqnarray*}
			\widetilde{\s}\widetilde{\s^\m}& = & \dsum_{(x),(\widetilde{x})}\s(1_y\om_x)\om_{x^\m}\s^{-1}(1_y\widetilde{\om}_x)\widetilde{\om}_{x^\m}= \dsum_{(x),(\widetilde{x})}\s(1_y\om_x\om_{x^\m}\s^{-1}(1_y\widetilde{\om}_x))\widetilde{\om}_{x^\m}\\
			& = &  \dsum_{(\widetilde{x})}\s(1_y1_x\s^{-1}(1_y\widetilde{\om}_x))\widetilde{\om}_{x^\m}= \dsum_{(\widetilde{x})}\s(\s^{-1}(1_y\widetilde{\om}_x))\widetilde{\om}_{x^\m}= \dsum_{(\widetilde{x})}1_y\widetilde{\om}_x\widetilde{\om}_{x^\m}=1_y1_x.
		\end{eqnarray*}
		Analogously, $\widetilde{\s^\m}\widetilde{\s}=1_x1_y$. Therefore, $\widetilde{\s}\in \U(\Z1_x1_y)$.

	Moreover, if  $u_x\in \T_x$, then:
	\begin{equation*}
		\widetilde{\s}u_x =  \dsum_{(x)}\s(1_y\om_x)\om_{x^\m}u_x=\dsum_{(x)}\s(1_y\om_x\om_{x^\m}u_x)= \s(1_y1_xu_x)=\s(1_yu_x).
	\end{equation*}
	and (\ref{relacaodefcomftil}) follows.
	
	Given $\s,\gamma\in \Aut_{R-R}(1_y\T_x)$, we have
\begin{eqnarray*}
			\widetilde{\s}\widetilde{\gamma} & = & \dsum_{(x),(\widetilde{x})} \s(1_y\om_x)\om_{x^\m}\gamma(1_y\widetilde{\om}_x)\widetilde{\omega}_{x^\m}=\dsum_{(x),(\widetilde{x})} \s(1_y\om_x\om_{x^\m}\gamma(1_y\widetilde{\om}_x))\widetilde{\omega}_{x^\m} \\
			& = &  \dsum_{(\widetilde{x})} \s(1_y1_x\gamma(1_y\widetilde{\om}_x))\widetilde{\omega}_{x^\m}= \dsum_{(\widetilde{x})} \s(\gamma(1_y\widetilde{\om}_x))\widetilde{\omega}_{x^\m} = \widetilde{\s\circ \gamma}.
	\end{eqnarray*}
	Therefore, $\widetilde{(-)}$ is a group morphism.

	Now, given $r\in \U(\Z1_x1_y)$, consider the map 
	$$\begin{array}{c c c l}
		\s_r: & 1_y\T_x & \longrightarrow & 1_y\T_x ,\\
		& v_x & \longmapsto & rv_x.
	\end{array}$$
	Clearly  $\s_r$ is right $R$-linear. Since $r \in \Z$ then $\s_r$ is left $R$-linear. Indeed, if $r' \in R,$ then 
	\begin{equation*}
		\s_r(r'v_x)=rr'v_x=r'rv_x=r'\s_r(v_x).
	\end{equation*}
	Let $r^\m$ be the inverse element of $r$ in $\Z1_x1_y$ and let $\s_{r^\m}:1_y\T_x\rightarrow 1_y\T_x$ be defined by $\sigma_{r^\m}(u_x)=r^\m u_x$.
	Given $v_x\in 1_y\T_x$,  we have that 
	\begin{equation*}
		\s_{r^\m}(\s_r(v_x))=\s_{r^\m}(rv_x)=r^\m (rv_x)=(r^\m r)v_x=1_x1_yv_x=v_x.
	\end{equation*}
	Analogously, $\s_r(\s_{r^\m}(v_x))=v_x$, for all $v_x \in 1_y\T_x$. Thus, $\sigma_r$ is a bijection with $\s_r^\m=\sigma_{r^\m}$. It follows that $\s_r \in \Aut_{R-R}(1_y\T_x)$.
	Notice now that
	\begin{eqnarray*}
		\widetilde{\s_r}=\dsum_{(x)}\s_r(1_y\om_x)\om_{x^\m}=\dsum_{(x)}r1_y\om_x\om_{x^\m}=r1_y1_x=r.
	\end{eqnarray*}
	Therefore, $\widetilde{(-)}$ is onto.
	Let $\gamma,\s \in \Aut_{R-R}(1_y\T_x)$ such that $\widetilde{\gamma}=\widetilde{\s}$. For any $v_x \in 1_y\T_x$, by (\ref{relacaodefcomftil}) it follows that
	$$\gamma(v_x)=\widetilde{\gamma}v_x=\widetilde{\s}v_x=\s(v_x).$$
	Thus, $\gamma=\s$ and, consequently  $\widetilde{(-)}$ is injective and, therefore, it is a group isomorphism.

\end{dem}

%


\begin{pro}\label{acaoparcialsobreZ} Let  $D_x=\Z1_x$ be the ideal of $\Z$ generated by the central idempotent $1_x$ and 
let   $\al_x:D_{x^\m}\mapsto D_x$ be defined by $$\al_x(r)=\sum_{(x)}\om_xr\om_{x^\m}.$$ Then, $\al=(D_x, \al_x)_{x\in G}$ is a partial action of $G$ on $\Z$.
\end{pro}

\begin{dem}
	First, let us  check that  $\al_x$ does not depend on the choice of the decomposition of $1_x$. Let  $\dsum_{(\widetilde{x})}\widetilde{\omega}_x\widetilde{\omega}_{x^\m}=1_x$, with  $\widetilde{\omega}_x\in \T_x$ and $\widetilde{\omega}_{x^\m}\in \T_{x^\m}$, be another decomposition of $1_x$. Given $r\in \Z , $  using (\ref{uxe1x}) and (\ref{moritasemtau}), we have
	\begin{eqnarray*}
			\al_x(r1_{x^\m}) & = & \dsum_{(x)}\omega_xr\omega_{x^\m}\stackrel{(\ref{uxe1x})}{=}\dsum_{(x)}\omega_xr\omega_{x^\m}1_{x}= \dsum_{(x),(\widetilde{x})}\omega_xr\omega_{x^\m}\widetilde{\omega}_x\widetilde{\omega}_{x^\m}\\
			& \stackrel{r \in \Z}{=}  &   \dsum_{(x),(\widetilde{x})}\omega_x\omega_{x^\m}\widetilde{\omega}_xr\widetilde{\omega}_{x^\m}= \dsum_{(\widetilde{x})}1_x\widetilde{\omega}_xr\widetilde{\omega}_{x^\m} \stackrel{(\ref{uxe1x})}{=}  \dsum_{(\widetilde{x})}\widetilde{\omega}_xr\widetilde{\omega}_{x^\m}.
		\end{eqnarray*}
		Now, for $r \in \Z$ and for any $r' \in R$, we  see that  
		\begin{eqnarray*}
			\al_x(r1_{x^{-1}})r'& = & \dsum_{(x)}\omega_xr \omega_{x^{-1}}r'=\dsum_{(x)}\omega_xr\omega_{x^{-1}}1_xr'=  \dsum_{(x),(\widetilde{x})}\omega_xr \omega_{x^{-1}}r'\widetilde{\omega}_x \widetilde{\omega}_{x^{-1}}\\
			& \stackrel{r \in \Z}{=}  & \dsum_{(x),(\widetilde{x})}\omega_x \omega_{x^{-1}}r'\widetilde{\omega}_xr \widetilde{\omega}_{x^{-1}}
			= \dsum_{(\widetilde{x})}1_xr'\widetilde{\omega}_xr \widetilde{\omega}_{x^{-1}}=\dsum_{(\widetilde{x})}r'\widetilde{\omega}_xr \widetilde{\omega}_{x^{-1}}=r'\al_x(r1_{x^\m}).
	\end{eqnarray*}
	Thus, $\al_x(r1_{x^{-1}}) \in D_x=\Z1_x$ and $\al_x$ is well-defined. Let $\sum_{(\overline{x})}\overline{\omega}_{x^\m}\ot \overline{\om}_{x}\in \T_{x^\m}\ot_R\T_{x}$ be the inverse image of $1_{x^\m}$ under the isomorphism $\tau_{x^\m}$, that is, $\dsum_{(\overline{x})}\overline{\omega}_{x^\m}\overline{\omega}_x=1_{x^\m}$. Given  $r\in D_x , $ we have
		\begin{eqnarray*}
			(\al_x\circ \al_{x^\m})(r) & = & \al_x\left(\dsum_{(\overline{x})}\overline{\omega}_{x^{-1}}r \overline{\omega}_x \right)=\dsum_{(x),(\overline{x})}\omega_x\overline{\omega}_{x^{-1}}r \overline{\omega}_x\omega_{x^\m}
			\stackrel{r \in \Z}{=}  \dsum_{(x),(\overline{x})}r\omega_x\overline{\omega}_{x^{-1}} \overline{\omega}_x\omega_{x^\m}\\
			& = &\dsum_{(x)}r\omega_x1_{x^\m}\omega_{x^\m}= \dsum_{(x)}r\omega_x\omega_{x^\m}=r1_x=r.
	\end{eqnarray*}
	
	Analogously, given $r \in D_{x^\m}, $ we obtain $(\al_{x^\m}\circ \al_x)(r)=r$. Therefore, $\al_x^{\m}=\al_{x^\m}$. 
	
	For $r,s \in D_x$ we see that 
	\begin{eqnarray*}
			\al_x(r)\al_x(s) & = &\dsum_{(x),(\widetilde{x})}\omega_xr\omega_{x^{-1}}\widetilde{\omega}_xs\widetilde{\omega}_{x^\m} \stackrel{s\in \Z}{=} \dsum_{(x),(\widetilde{x})}\omega_xrs\omega_{x^{-1}}\widetilde{\omega}_x\widetilde{\omega}_{x^\m} \\
			& = & \dsum_{(x)}\omega_xrs\omega_{x^{-1}}1_x  = \dsum_{(x)}\omega_xrs\omega_{x^{-1}} = \al_x(rs).
	\end{eqnarray*}
	Thus, $\al_x$ is a ring isomorphism.  
	
	Write  $1_{y^\m}=\dsum_{(\overline{y})}\overline{\om}_{y^\m}\overline{\om}_y$. Given  $r \in D_{y}\cap D_{x^\m}$, using (\ref{uxe1y}),  we obtain
	\begin{equation*}
		\al_{y^\m}(r)=\dsum_{(\overline{y})}\overline{\om}_{y^\m}r\overline{\om}_{y} = \dsum_{(\overline{y})}\overline{\om}_{y^\m}r1_{x^\m}\overline{\om}_{y}=\dsum_{(\overline{y})}\overline{\om}_{y^\m}r\overline{\om}_{y}1_{(xy)^\m} \in D_{(xy)^\m}.
	\end{equation*}
	Hence, $\al_{y^\m}(D_y\cap D_{x^\m})\subseteq D_{(xy)^\m}$. 
	Using again (\ref{uxe1y}),
	\begin{equation}
		1_x1_{xy}=\dsum_{(x)}\om_x\om_{x^\m}1_{xy}=\dsum_{(x)}\om_x1_y\om_{x^\m}=\dsum_{(x),(y)}\om_x(\om_y\om_{y^\m})\om_{x^\m}.\label{1x1xy2}
	\end{equation}
	Write  $\dsum_{(xy)}\omega_{xy}\omega_{(xy)^\m}=1_{xy}$ and take  $r \in D_{y^{-1}}\cap D_{(xy)^{-1}}$. By (\ref{ruxuyuxym1}) and (\ref{1x1xy2}) we compute
	\begin{eqnarray*}
		\al_{xy}(r) & = & \dsum_{(xy)}\omega_{xy}r\omega_{(xy)^{-1}}
		\stackrel{r\in D_{y^\m}}{=}   \dsum_{(xy)}1_{xy}\omega_{xy}1_{y^\m}r\omega_{(xy)^{-1}}\stackrel{(\ref{uxe1y})}{=}  \dsum_{(xy)}1_{xy}1_x\omega_{xy}r \omega_{(xy)^{-1}}\\
		& \stackrel{(\ref{1x1xy2})}{=} & \dsum_{\stackrel{(x),(y),}{(xy)}} (\omega_x(\omega_y \omega_{y^{-1}}) \omega_{x^{-1}})\omega_{xy}r\omega_{(xy)^{-1}}=\dsum_{\stackrel{(x),(y),}{(xy)}} (\om_x(\om_yr\om_{y^\m})\om_{x^\m})\om_{xy}\om_{(xy)^\m}\\
		& = & \dsum_{(x),(y)} (\om_x(\om_yr\om_{y^\m})\om_{x^\m})1_{xy} \stackrel{(\ref{uxe1y})}{=}  \dsum_{(x),(y)} (\om_x(\om_yr\om_{y^\m})1_y\om_{x^\m})\\
		& = & \dsum_{(x),(y)} (\om_x(\om_yr\om_{y^\m})\om_{x^\m})=\al_x\circ \al_y(r).
	\end{eqnarray*}
	Therefore, $\al=(D_x,\al_x)_{x \in G}$ is a partial action of $G$ on $\Z$.  \end{dem}

\begin{lem} Let $\al=(D_x,\al_x)$ be the partial action of $G$ on $\Z$ constructed in Proposition~\ref{acaoparcialsobreZ} {and $M$ an $R$-bimodule}. Then:
	\begin{itemize}
	\item[(i)] For all $u_x \in \T_x$ and $r \in \Z$ we have
		\begin{equation}
		\al_x(r1_{x^\m})u_x=u_xr. \label{uxercomutam}
		\end{equation}
	\item[(ii)] If $M| \T_x$, then 
		\begin{equation}
		\al_x(r1_{x^\m})m=m r \label{mercomutam}
		\end{equation}
for all $m \in M$ and $r \in \Z.$
	\end{itemize}
\end{lem}
\begin{dem}
	(i) Given $r \in \Z$ and $u_x \in \T_x$, by (\ref{uxe1x}) and (\ref{moritasemtau}) we have
	\begin{eqnarray*}
		\al_x(r1_{x^\m})u_x & = & \dsum_{(x)}(\omega_xr\omega_{x^{-1}})u_x =  \dsum_{(x)}\omega_xr( \omega_{x^{-1}}u_x)\\
		& \stackrel{r \in \Z}{=} & \dsum_{(x)}\omega_x( \omega_{x^{-1}}u_x)r = \dsum_{(x)}(\omega_x  \omega_{x^{-1}})u_xr\\
		& = & 1_xu_xr\stackrel{(\ref{uxe1x})}{=}u_xr.
	\end{eqnarray*}

	(ii) Let $f_i:M\rightarrow \T_x$ and $g_i:\T_x\rightarrow M$, $i=1,2,...,n$, be $R$-bimodule morphisms such that $\sum_{i=1}^ng_if_i=Id_M$. Given $m \in M$ and $r\in \Z , $ using (\ref{uxercomutam}) we obtain
	\begin{eqnarray*}
		mr & = & \dsum_{i=1}^ng_i(f_i(m))r=\dsum_{i=1}^ng_i(\underbrace{f_i(m)}_{\in \T_x}r)
		=  \dsum_{i=1}^ng_i(\al_x(r1_{x^\m})f_i(m))\\
		& = & \dsum_{i=1}^n\al_x(r1_{x^\m})g_i(f_i(m))
		=  \al_x(r1_{x^\m})m.
	\end{eqnarray*}
	
\end{dem}

We denote by $H^n_\T(G,\al,\Z)$ the group of partial $n$-cohomologies, where $\al$ is the partial action induced by the  untial partial represetation $\T$.


We shall define now a partial action of $G$ on $\Pics(R) .$ Let $\mathcal{X}_x=[R1_x]\Pics(R)[R1_x]$ and
\begin{equation*}
	\begin{array}{c c c l}
		\al_x^*: & \mathcal{X}_{x^\m} & \longrightarrow & \mathcal{X}_x ,\\
		& [P] & \longmapsto & [\T_x\ot_R P\ot_R \T_{x^\m}].
	\end{array}
\end{equation*}
It follows by Proposition~\ref{proparaacaoparcialsobrePics}  that  $\al^*=(\mathcal{X}_x,\al_x^*)$ is a partial action of $G$ on  $\Pics(R)$.  
Observe that if   $[P]\in \Pics_\Z(R)$, then by  Corollary~\ref{TcomZbomodcentral}, we obtain that  $P\ot_RR1_{x}\simeq R1_x\ot_RP$, as $R$-bimodules,  for all $x \in G$. Hence, $[R1_x]$ is a { central} idempotent  in $\Pics_\Z(R)$. Moreover, if $[P]\in \Pics_\Z(R)$, then $[\T_x\ot_R P\ot_R \T_{x^\m}]\in \Pics_\Z(R)$. Indeed, it is enough to verify that   $\T_x\ot_R P\ot_R \T_{x^\m}$ is also a central $\Z$-bimodule. Let $p\in P$, $u_x \in \T_x, u_{x^\m}\in \T_{x^\m}$ and $r \in \Z .$ Since $P$ is central over  $\Z ,$  using  (\ref{uxercomutam}) we obtain:
\begin{eqnarray*}
	u_x \ot p \ot u_{x^\m}r & = & u_x\ot p\ot \al_{x^\m}(r1_x)u_{x^\m} =u_x\ot p\al_{x^\m}(r1_x)\ot u_{x^\m}\\
	& = & u_x\ot \al_{x^\m}(r1_x)p\ot u_{x^\m}=u_x\al_{x^\m}(r1_x)\ot p\ot u_{x^\m}\\
	& = &  \al_x(\al_{x^\m}(r1_x))u_x\ot p\ot u_{x^\m}=r1_xu_x\ot p \ot u_{x^\m}\\
	& = & ru_x\ot p \ot u_{x^\m}.
\end{eqnarray*}		
Thus, restricting the above partial action to  $\Pics_\Z(R)$, we have that $\mathcal{X}_x=\Pics_\Z(R)[R1_x],$ i.e.  
$\mathcal{X}_x$ is the  ideal in $\Pics_\Z(R)$ generated by the central  idempotent  $[R1_x].$ Furthermore,  
\begin{equation*}
	\begin{array}{c c c l}
		\al_x^*: & \Pics_\Z(R)[R1_{x^\m}] & \longrightarrow & \Pics_\Z(R)[R1_x],\\
		& [P][R1_{x^\m}] & \longmapsto & [\T_x\ot_R P\ot_R \T_{x^\m}]
	\end{array}
\end{equation*}
is well-defined and, consequently,  $\al^*=(\mathcal{X}_x,\al^*_x)$ is a partial action of $G$ on $\Pics_\Z(R)$. 

The subsemigroup of invariants of  $\Pics_\Z(R)$ is given by 
\begin{equation}
	\Pics_\Z(R)^{\al^*}=\{[P]\in \Pics_\Z(R); \ \T_x\ot_RP\simeq P\ot_R\T_{x}, \ \mbox{for all} \ x \in G\}. \label{invariantesdePicsZ}
\end{equation}
Indeed, $[P]\in \Pics_\Z(R)^{\al^*},$ if and only if  $\al_x^*([P][R1_{x^\m}])=[P][R1_x]$, for all  $x \in G$. It follows that 
\begin{eqnarray*}
	[P]\in \Pics_\Z(R)^{\al^*} & \Leftrightarrow & \T_x\ot_RP\ot_R\T_{x^{-1}}\simeq P\ot_RR1_x,  \ \ \forall  x\in G;\\
	& \Leftrightarrow & \T_x\ot_RP\ot_R\T_{x^{-1}}\ot_R\T_x\simeq P\ot_RR1_x\ot_R\T_x, \  \ \forall x\in G;\\
	& \Leftrightarrow & \T_x\ot_RP\ot_RR1_{x^{-1}}\simeq P\ot_R\T_x, \  \ \forall x\in G;\\
	& \Leftrightarrow & \T_x\ot_RR1_{x^{-1}}\ot_RP\simeq P\ot_R\T_x, \ \ \forall x\in G;\\
	& \Leftrightarrow & \T_x\ot_RP\simeq P\ot_R\T_x, \  \ x\in G.\\
\end{eqnarray*}

\subsection{Partial Generalized Crossed Products}
\label{sec:FactSet GenParCrossProd}
\label{sec: pcgp}
Let 
$$\begin{array}{c c c c}
\T: & G& \longrightarrow & \Pics(R),\\
& x & \longmapsto & [\T_x],
\end{array}$$ be a unital partial representation with $\T_x\ot\T_{x^{-1}}\simeq R1_x$,  for all $x\in G$. By (\ref{1xyTxisomorfismo}) there is a family of $R$-bimodule isomorphisms
\begin{equation}
f^\T=\{f_{x,y}^\T:\T_x\ot_R \T_y\longrightarrow 1_x\T_{xy}, \ x,y\in G\}. 
\end{equation}

Following \cite{DoPaPi2}, we say that $f^\T$ is a \textit{factor set} for $\T$ if $f^\T$ satisfies the  following commutative diagram: 
\begin{equation}
\xymatrix{ \T_x\ot_R \T_y\ot_R\T_z\ar[r]^{\T_x\ot f^\T_{y,z}}\ar[dd]_{f^\T_{x,y}\ot \T_z} & \T_x\ot_R 1_y\T_{yz}\ar@{=}[r] & 1_{xy}\T_x\ot_R \T_{yz}\ar[dd]^{f^\T_{x,yz}}\\
	& & \\
	1_x\T_{xy}\ot_R \T_z\ar[rr]_{f^\T_{xy,z}} & & 1_x1_{xy}\T_{xyz}  } \label{diagramaassociatividadeprodutocruzado}
\end{equation}

%

Let $f^\T=\{f_{x,y}^{\T}:\T_x\ot_R\T_y\longrightarrow 1_x\T_{xy},\ x,y \in G\}$ be a factor set for $\T .$ The set $\D(\T)=\bigoplus_{x\in G}\T_x$ with multiplication defined by
\begin{equation}
u_x\stackrel{\T}{\circ} u_y = f_{x,y}^\T(u_x\ot u_y)\in 1_x\T_{xy}, \ \ u_x\in \T_x,u_y\in \T_y,
\end{equation}
is called a \textit{partial generalized crossed product}. 

\begin{pro}\label{pcgpeanelcomutativocom1} Let $\T:G\longrightarrow \Pics(R)$ be a unital partial representation with $\varepsilon_x=\T_x\ot_R\T_{x^\m}\simeq R1_x$, for all $x\in G$, and let $f^\T=\{f_{x,y}^\T:\T_x\ot_R\T_y\longrightarrow 1_x\T_{xy}, \ x,y \in G\}$ be   a factor set for  $\T$. Then, the partial generalized crossed product  $\D(\T)$ is an associative ring with unity and $R\simeq \T_1$ is a subring of $\D(\T)$.
\end{pro}
\begin{dem} The commutativity of   (\ref{diagramaassociatividadeprodutocruzado}) implies the associativity of the multiplication in  $\D(\T) .$  
 It is easy to see that $\T_1$ is a subring of ${\D(\T)}$ since $f_{1,1}^\T:\T_1\ot\T_1\rightarrow \T_1$ is an $R$-bimodule isomorphism.  
	Let  $j:R\longrightarrow\T_1$ be a $R$-bimodule isomorphism;  it exists because $\T$ is a partial representation and $[\T _1]=[R].$ Write $v=j(1) .$ Then $\T_1=j(R)=Rv=vR$ and
	$$rv=rj(1)=j(r)=j(1)r=vr, \ \ \mbox{for all} \ r \in R.$$
	Let $c \in R$ be such that  $f_{1,1}^\T(v\ot v)=cv$. 
	\begin{afr} $c\in \U(\Z)$.
	\end{afr} Indeed, let $r \in R$ be such that $rv=0 .$  Then $0=rv=rj(1)=j(r)$. Since $j$ is an isomorphism,  we have  $r=0.$  For any $r \in R$ we see that
	\begin{equation*}
	rcv=rf_{1,1}^\T(v\ot v)=f_{1,1}^\T(rv\ot v)=f_{1,1}^\T(v\ot vr)=f_{1,1}^\T(v\ot v)r=vcr=crv.
	\end{equation*}
	Thus, $(rc-cr)v=0,$ which implies   that  $rc=cr$, for any $r\in R.$ 
	
	Since $f_{1,1}^\T$ is an isomorphism, there is $d \in R$ such that $(f_{1,1}^\T)^\m(v)=dv\ot v.$ Then
	\begin{equation*}
	dcv=df_{1,1}^\T(v\ot v)=f_{1,1}^\T(dv\ot v)=f_{1,1}^\T(f_{1,1}^\T)^\m(v)=v.
	\end{equation*}
	Thus, $dc=cd=1$. Therefore $c\in \U(\Z)$.
	
	Denote $u=c^{\m}v\in \T_1 .$ Then $ru=rc^\m v=c^\m vr=ur$, for all $r \in R$ and
	\begin{equation}
	f_{1,1}^\T(u\ot u)=f_{1,1}^\T(c^{-1}v\ot c^{-1}v)=c^{\m}f^{\T}_{1,1}(v\ot v)c^\m=c^\m cvc^\m=vc^\m=u. \label{uidempotente}
	\end{equation}
	Let $\iota:R\longrightarrow\T_1$ be the $R$-bimodule isomorphism defined by  $\iota(r)=ru$. Given $r,r' \in R$ we have
	\begin{equation*}
	\iota(r)\stackrel{\T}{\circ}\iota(r')=f_{1,1}^\T(ru\ot r'u)=rr'f^{\T}_{1,1}(u\ot u)\stackrel{(\ref{uidempotente})}{=}rr'u=\iota(rr').
	\end{equation*}
	so, $\iota$ is a ring isomorphism.
	
	Finally, we check that  $\iota(1)=u \in \T_1$ is the unity of $\D(\T)$. Let $u_x\in \T_x$, $x\in G$. Since $f_{1,x}^\T:\T_1\ot_R\T_x\longrightarrow \T_x$ is an $R$-bimodule isomorphism,  there is $\sum_{i=1}^nur_i\ot u_x^i=u\ot \widetilde{u}_x\in \T_1\ot_R\T_x$, where $\widetilde{u}_x=\sum_{i=1}^nr_iu_x^i$, such that $f_{1,x}^\T(u\ot \widetilde{u}_x)=u_x$. Then,
	\begin{eqnarray*}
		(u\stackrel{\T}{\circ} u_x)& = & f_{1,x}^\T(u\ot u_x)=f_{1,x}^\T(u\ot f_{1,x}^\T(u\ot \widetilde{u}_x))\\
		& = & f_{1,x}^\T(f_{1,1}^\T(u\ot u)\ot \widetilde{u}_x)\\
		& \stackrel{(\ref{uidempotente})}{=} & f_{1,x}^\T(u\ot \widetilde{u}_x)=u_x. 
	\end{eqnarray*} 
	Analogously, since $f_{x,1}^\T:\T_x\ot_R\T_1\longrightarrow \T_x$ is an isomorphism, there is 
$\overline{u_x}\in \T_x$ such that $f_{x,1}^\T(\overline{u_x}\ot u)=u_x .$  Then
	\begin{eqnarray*}
		(u_x\stackrel{\T}{\circ} u) & = & f_{x,1}^\T(u_x\ot u) = f_{x,1}^\T(f_{x,1}^\T(\overline{u_x}\ot u)\ot u)\\
		& = & f_{x,1}^\T(\overline{u_x}\ot f_{1,1}^\T(u\ot u))\\
		& = & f_{x,1}^\T(\overline{u_x}\ot u)=u_x.
	\end{eqnarray*}
	Therefore, $u$ is the unity of $\D(\T)$.
	\end{dem}

We identify $R$ with $\T_1$ and $1$ with $\iota(1)=u \in \T_1$, via the ring isomorphism $\iota$. Then,  $R\subseteq \D(\T)$ is an extension of rings with the same unity.

\begin{obs}\label{unidadeediagramasemDT} Let $\D(\T)$ be a partial generalized crossed product, 
  $\iota:R\longrightarrow \T_1$ an isomorphism of rings and $R$-bimodules and   $\iota(1)=u \in \T_1$ be the unity of $\D(\T)$. Then, the following diagrams are commutative:
	\begin{equation}
	\begin{tabular}{ccc}
	$\xymatrix{ R\ot_R\T_x \ar[rd]_{\iota\ot \T_x}\ar[rr]^{\simeq}& & \T_x\\
		& \T_1\ot_R\T_x\ar[ur]_{f_{1,x}^{\T}} & }$ &  \mbox{and} & $\xymatrix{ \T_x\ot_RR \ar[rd]_{\T_x\ot \iota}\ar[rr]^{\simeq}& & \T_x\\
		& \T_x\ot_R\T_1\ar[ur]_{f_{x,1}^\T} & }$
	\end{tabular} \label{diagramasdaunidade}
	\end{equation}
	 Indeed, given  $r \in R$ and $u_x \in \T_x ,$ we have
	\begin{equation*}
	(f_{1,x}^\T\circ(\iota\ot \T_x))(r\ot u_x)  =  f_{1,x}^\T(\iota(r)\ot u_x)= rf_{1,x}^\T(\iota(1)\ot u_x)=ru_x.
	\end{equation*} 
	Analogously, 
	\begin{equation*}
	(f_{x,1}^\T\circ (\T_x\ot \iota))(u_x\ot r) = f^{\T}_{x,1}(u_x\ot \iota(r))=f^{\T}_{x,1}(u_x\ot \iota(1))r=u_xr.
	\end{equation*}

\end{obs}

By  (\ref{diagramasdaunidade}), we can write 
\begin{equation}
u_x\stackrel{\T}{\circ} r=u_xr \quad \mbox{e} \quad r\stackrel{\T}{\circ} u_x=ru_x, \label{produtoporr}
\end{equation}
for all $u_x \in \T_x$ and  $r \in R$.

\begin{pro}\label{3cocicloquecorrigeaassociatividade} Let  $f^\T=\{f_{x,y}^\T: \T_x\ot_R\T_y\longrightarrow 1_x\T_{xy}, \ x,y\in G \}$ be a family of $R$-bimodule isomorphisms which satisfies the commutative diagrams   (\ref{diagramasdaunidade}). Given $x,y,z \in G ,$ define an  $R$-bimodule   isomorphism $\beta_{x,y,z}:1_x1_{xy}\T_{xyz}\longrightarrow 1_x1_{xy}\T_{xyz}$ via the commutative  diagram
	\begin{equation*}
		\xymatrix{ \T_x\ot_R\T_y\ot_R\T_z\ar[rr]^{f_{x,y}^\T\ot \T_z} \ar[dd]_{\T_x\ot f^\T_{y,z}}& & 1_x\T_{xy}\ot_R\T_z\ar[dd]^{f_{xy,z}^\T} \\ 
			& & \\
			1_{xy}\T_x\ot_R\T_{yz}\ar[rd]_{f_{x,yz}^\T} & & 1_x1_{xy}\T_{xyz}\\
			& 1_{x}1_{xy}\T_{xyz}\ar@{-->}[ru]_{\beta_{x,y,z}} &                 }
	\end{equation*}
	that is, 
	\begin{equation}
		\beta_{x,y,z}\circ f^\T_{x,yz}\circ (\T_x\ot f^\T_{y,z})=f_{xy,z}^\T\circ (f_{x,y}^\T\ot \T_z), \ \forall \ x,y,z\in G. \label{definicaodebetaem3cociclo}
	\end{equation}
	Then,
	\begin{equation*}
		\begin{array}{c c c l}
			\widetilde{\beta_{-,-,-}}: & G\times G\times G & \longrightarrow & \U(\Z)\\
			& (x,y,z) & \longmapsto & \widetilde{\beta_{x,y,z}}
		\end{array},
	\end{equation*}
	is a normalized element of $Z^3_{\T}(G,\al,\Z),$ where $\widetilde{\beta_{x,y,z}}$ is defined as in  Lemma~\ref{gammatil}. Furthermore, if $\widetilde{\beta_{-,-,-}}\in B_\T^3(G,\al,\Z)$, that is, if there exists  $\sigma_{-,-}:G\times G \longrightarrow \U(\Z)$, with $\sigma_{x,y}\in \U(\Z1_x1_{xy})$ and
	\begin{equation}
		\widetilde{\beta_{x,y,z}}=\al_x(\sigma_{y,z}1_{x^\m})\sigma_{xy,z}^\m\sigma_{x,yz}\sigma_{x,y}^\m, \ \ \mbox{for all } x,y,z \in G,\label{betaedelta2desigma}
	\end{equation}
	then the family of $R$-bimodule isomorphisms  
	\begin{equation*}
		\begin{array}{c c c l}
			\bar{f}_{x,y}^\T: & \T_x\ot_R\T_y & \longrightarrow & 1_x\T_{x,y},\\
			& u_x\ot u_y & \longmapsto & \sigma_{x,y}f_{x,y}^\T(u_x\ot u_y)
		\end{array}
	\end{equation*}
	is a factor set for  $\T$.
\end{pro}
\begin{dem}
	Since $\beta_{x,y,z} \in \Aut_{R-R}(1_x1_{xy}\T_{xyz})$,    Lemma~\ref{gammatil} implies that  $\widetilde{\beta_{x,y,z}}\in \U(\Z1_x1_{xy}1_{xyz})$, for all  $x,y,z \in G$. Write $f_{x,y}^\T(u_x\ot u_y)=(u_x\stackrel{\T}{\circ}u_y)$. By (\ref{definicaodebetaem3cociclo}) we have that
	\begin{equation*}
		\beta_{x,y,z}(u_x\stackrel{\T}{\circ}(u_y\stackrel{\T}{\circ}u_z))=((u_x\stackrel{\T}{\circ}u_y)\stackrel{\T}{\circ} u_z),
	\end{equation*}
	for all  $u_x \in \T_x,u_y \in \T_y$ and $u_z \in \T_z$. Since $(u_x\stackrel{\T}{\circ}(u_y\stackrel{\T}{\circ}u_z)) \in 1_x1_{xy}\T_{xyz}$, then by  (\ref{relacaodefcomftil}),
	\begin{equation}
		\widetilde{\beta_{x,y,z}}(u_x\stackrel{\T}{\circ}(u_y\stackrel{\T}{\circ}u_z))=((u_x\stackrel{\T}{\circ}u_y)\stackrel{\T}{\circ} u_z). \label{betatilefxyz}
	\end{equation}
	Let $x,y,z,t \in G$ e $u_x \in \T_x,u_y \in \T_y, u_z \in \T_z$ and $u_t \in \T_t$. Then,
	\begin{eqnarray*}
		(((u_x\stackrel{\T}{\circ}u_y)\stackrel{\T}{\circ} u_z)\stackrel{\T}{\circ}u_t) & \stackrel{(\ref{betatilefxyz})}{=} & \widetilde{\beta_{x,y,z}}((u_x\stackrel{\T}{\circ}(u_y\stackrel{\T}{\circ}u_z))\stackrel{\T}{\circ}u_t)\\
		& \stackrel{(\ref{betatilefxyz})}{=}& \widetilde{\beta_{x,y,z}}\widetilde{\beta_{x,yz,t}}(u_x\stackrel{\T}{\circ}((u_y\stackrel{\T}{\circ}u_z)\stackrel{\T}{\circ}u_t)) \\
		& \stackrel{(\ref{betatilefxyz})}{=} & \widetilde{\beta_{x,y,z}}\widetilde{\beta_{x,yz,t}}(u_x\stackrel{\T}{\circ}\widetilde{\beta_{y,z,t}}(u_y\stackrel{\T}{\circ}(u_z\stackrel{\T}{\circ}u_t)))\\
		& = & \widetilde{\beta_{x,y,z}}\widetilde{\beta_{x,yz,t}}(u_x\widetilde{\beta_{y,z,t}}\stackrel{\T}{\circ}(u_y\stackrel{\T}{\circ}(u_z\stackrel{\T}{\circ}u_t)))\\
		& \stackrel{(\ref{uxercomutam})}{=} & \widetilde{\beta_{x,y,z}}\widetilde{\beta_{x,yz,t}}\al_x(\widetilde{\beta_{y,z,t}}1_{x^\m})(u_x\stackrel{\T}{\circ}(u_y\stackrel{\T}{\circ}(u_z\stackrel{\T}{\circ}u_t)))\\
		& \stackrel{(\ref{betatilefxyz})}{=} & \widetilde{\beta_{x,y,z}}\widetilde{\beta_{x,yz,t}}\al_x(\widetilde{\beta_{y,z,t}}1_{x^\m})\widetilde{\beta_{x,y,zt}}^{\m}((u_x\stackrel{\T}{\circ}u_y)\stackrel{\T}{\circ}(u_z\stackrel{\T}{\circ}u_t))\\
		& \stackrel{(\ref{betatilefxyz})}{=} & \widetilde{\beta_{x,y,z}}\widetilde{\beta_{x,yz,t}}\al_x(\widetilde{\beta_{y,z,t}}1_{x^\m})\widetilde{\beta_{x,y,zt}}^{\m}\widetilde{\beta_{xy,z,t}}^{\m}(((u_x\stackrel{\T}{\circ}u_y)\stackrel{\T}{\circ}u_z)\stackrel{\T}{\circ}u_t).
	\end{eqnarray*}
	Thus, 
	\begin{equation*}
		(((u_x\stackrel{\T}{\circ}u_y)\stackrel{\T}{\circ} u_z)\stackrel{\T}{\circ}u_t) =\widetilde{\beta_{x,y,z}}\widetilde{\beta_{x,yz,t}}\al_x(\widetilde{\beta_{y,z,t}}1_{x^\m})\widetilde{\beta_{x,y,zt}}^{\m}\widetilde{\beta_{xy,z,t}}^{\m}(((u_x\stackrel{\T}{\circ}u_y)\stackrel{\T}{\circ}u_z)\stackrel{\T}{\circ}u_t),
	\end{equation*}
	for all $u_x \in \T_x,u_y \in \T_y, u_z \in \T_z$ and $u_t \in \T_z$. 
{Write
		$$1_x=\dsum_{(x)}(\omega_x\stackrel{\T}{\circ}\omega_{x^{-1}}), 1_y=\dsum_{(y)}(\omega_y\stackrel{\T}{\circ}\omega_{y^{-1}}), 1_z=\dsum_{(z)}(\omega_z\stackrel{\T}{\circ}\omega_{z^{-1}}) \ \mbox{and} \ 1_t=\dsum_{(t)}(\omega_t\stackrel{\T}{\circ}\omega_{t^{-1}}),$$
	where $\om_l \in \T_l$, for all  $l \in G.$ Then,
		\begin{equation*}
			(((\om_x\stackrel{\T}{\circ}\om_y)\stackrel{\T}{\circ} \om_z)\stackrel{\T}{\circ}\om_t) =\widetilde{\beta_{x,y,z}}\widetilde{\beta_{x,yz,t}}\al_x(\widetilde{\beta_{y,z,t}}1_{x^\m})\widetilde{\beta_{x,y,zt}}^{\m}\widetilde{\beta_{xy,z,t}}^{\m}(((\om_x\stackrel{\T}{\circ}\om_y)\stackrel{\T}{\circ}\om_z)\stackrel{\T}{\circ}\om_t),
			\end{equation*}
		and, applying  $f_{xyzt,t^\m}^\T$,  we obtain:
		\begin{equation*}
			(((\om_x\stackrel{\T}{\circ}\om_y)\stackrel{\T}{\circ} \om_z)\stackrel{\T}{\circ}\om_t)\stackrel{\T}{\circ}\om_{t^\m} =\widetilde{\beta_{x,y,z}}\widetilde{\beta_{x,yz,t}}\al_x(\widetilde{\beta_{y,z,t}}1_{x^\m})\widetilde{\beta_{x,y,zt}}^{\m}\widetilde{\beta_{xy,z,t}}^{\m}((((\om_x\stackrel{\T}{\circ}\om_y)\stackrel{\T}{\circ}\om_z)\stackrel{\T}{\circ}\om_t)\stackrel{\T}{\circ}\om_{t^\m}).
			\end{equation*}
		By (\ref{betatilefxyz}),
		\begin{eqnarray*}
				& & \widetilde{\beta_{xyz,t,t^\m}}(((\om_x\stackrel{\T}{\circ}\om_y)\stackrel{\T}{\circ} \om_z)\stackrel{\T}{\circ}(\om_t\stackrel{\T}{\circ}\om_{t^\m})) \\
				& = & \widetilde{\beta_{x,y,z}}\widetilde{\beta_{x,yz,t}}\al_x(\widetilde{\beta_{y,z,t}}1_{x^\m})\widetilde{\beta_{x,y,zt}}^{\m}\widetilde{\beta_{xy,z,t}}^{\m} \widetilde{\beta_{xyz,t,t^\m}}(((\om_x\stackrel{\T}{\circ}\om_y)\stackrel{\T}{\circ} \om_z)\stackrel{\T}{\circ}(\om_t\stackrel{\T}{\circ}\om_{t^\m})).
			\end{eqnarray*}
		Multiplying by  $\widetilde{\beta_{xyz,t,t^\m}}^\m$:
		\begin{eqnarray*}
				& & 1_{xyz}1_{xyzt}(((\om_x\stackrel{\T}{\circ}\om_y)\stackrel{\T}{\circ} \om_z)\stackrel{\T}{\circ}(\om_t\stackrel{\T}{\circ}\om_{t^\m})) \\
				& = & \widetilde{\beta_{x,y,z}}\widetilde{\beta_{x,yz,t}}\al_x(\widetilde{\beta_{y,z,t}}1_{x^\m})\widetilde{\beta_{x,y,zt}}^{\m}\widetilde{\beta_{xy,z,t}}^{\m} 1_{xyz}1_{xyzt}(((\om_x\stackrel{\T}{\circ}\om_y)\stackrel{\T}{\circ} \om_z)\stackrel{\T}{\circ}(\om_t\stackrel{\T}{\circ}\om_{t^\m}))\\
				& = & \widetilde{\beta_{x,y,z}}\widetilde{\beta_{x,yz,t}}\al_x(\widetilde{\beta_{y,z,t}}1_{x^\m})\widetilde{\beta_{x,y,zt}}^{\m}\widetilde{\beta_{xy,z,t}}^{\m} (((\om_x\stackrel{\T}{\circ}\om_y)\stackrel{\T}{\circ} \om_z)\stackrel{\T}{\circ}(\om_t\stackrel{\T}{\circ}\om_{t^\m}))
			\end{eqnarray*}
		Summing over  $(t)$ and using the commutativity of  (\ref{diagramasdaunidade}) and (\ref{uxe1y}), we compute: 
		\begin{equation*}
			1_{xyz}1_{xyzt}(((\om_x\stackrel{\T}{\circ}\om_y)\stackrel{\T}{\circ} \om_z))1_t \\
			=   \widetilde{\beta_{x,y,z}}\widetilde{\beta_{x,yz,t}}\al_x(\widetilde{\beta_{y,z,t}}1_{x^\m})\widetilde{\beta_{x,y,zt}}^{\m}\widetilde{\beta_{xy,z,t}}^{\m} (((\om_x\stackrel{\T}{\circ}\om_y)\stackrel{\T}{\circ} \om_z))1_t
			\end{equation*}
		\begin{equation*}
			1_{xyz}1_{xyzt}1_{xyzt}(((\om_x\stackrel{\T}{\circ}\om_y)\stackrel{\T}{\circ} \om_z)) 
			=   \widetilde{\beta_{x,y,z}}\widetilde{\beta_{x,yz,t}}\al_x(\widetilde{\beta_{y,z,t}}1_{x^\m})\widetilde{\beta_{x,y,zt}}^{\m}\widetilde{\beta_{xy,z,t}}^{\m} 1_{xyzt}(((\om_x\stackrel{\T}{\circ}\om_y)\stackrel{\T}{\circ} \om_z))
			\end{equation*}
		\begin{equation*}
			1_{xyz}1_{xyzt}(((\om_x\stackrel{\T}{\circ}\om_y)\stackrel{\T}{\circ} \om_z))
			=   \widetilde{\beta_{x,y,z}}\widetilde{\beta_{x,yz,t}}\al_x(\widetilde{\beta_{y,z,t}}1_{x^\m})\widetilde{\beta_{x,y,zt}}^{\m}\widetilde{\beta_{xy,z,t}}^{\m}(((\om_x\stackrel{\T}{\circ}\om_y)\stackrel{\T}{\circ} \om_z)).
			\end{equation*}
		Repeating the same   argument for   $z,y,x$, we get:
	$$\widetilde{\beta_{x,y,z}}\widetilde{\beta_{x,yz,t}}\al_x(\widetilde{\beta_{y,z,t}}1_{x^\m})\widetilde{\beta_{x,y,zt}}^{\m}\widetilde{\beta_{xy,z,t}}^{\m}=1_x1_{xy}1_{xyz}1_{xyzt},$$
for all $x,y,z,t \in G$.}
%
%
%
	Hence, $\widetilde{\beta_{-,-,-}}\in Z^3_\T(G,\al, \Z)$. 
	
	Let us check that  $\widetilde{\beta_{-,-,-}}$ is normalized. Taking $x=1$, we obtain using  (\ref{betatilefxyz}) and the commutativity of the diagrams in (\ref{diagramasdaunidade}), that  $	\widetilde{\beta_{1,y,z}}(ru_y\stackrel{\T}{\circ}u_z)=(ru_y\stackrel{\T}{\circ}u_z),$
	for all  $r \in R$, $u_x\in \T_x$ and $u_y \in \T_y$. Putting $r=1$ and following the above argument, we conclude that   $\widetilde{\beta_{1,y,z}}=1_y1_{yz}$. Analogously, we obtain 
	
	\begin{equation*}
		\widetilde{\beta_{x,1,z}}=1_x1_{xz} \ \ \mbox{e} \ \ \widetilde{\beta_{x,y,1}}=1_x1_{xy}, \ \mbox{for all} \ x,y,z \in G.
	\end{equation*}
	Therefore,  $\widetilde{\beta_{-,-,-}} \in   Z^3_\T(G,\al,\Z) $ is  normalized.

	We now   verify that  $\bar{f}^\T=\{\bar{f}_{x,y}^\T, \ x,y \in G \}$ is a factor set  for  $\T$.  Let $x,y,z \in G$ and $u_x \in \T_x,u_y \in \T_y, u_z\in \T_z .$ Then,
	\begin{eqnarray*}
		\bar{f}_{x,yz}^\T(u_x\ot \bar{f}_{y,z}^\T(u_y\ot u_z)) & = & \sigma_{x,yz}(u_x \stackrel{\T}{\circ}\sigma_{y,z}(u_y\stackrel{\T}{\circ}u_z))\\
		& = & \sigma_{x,yz}(u_x\sigma_{y,z} \stackrel{\T}{\circ}(u_y\stackrel{\T}{\circ}u_z))\\
		& \stackrel{(\ref{uxercomutam})}{=}&  \sigma_{x,yz}(\al_{x}(\sigma_{y,z}1_{x^\m})u_x \stackrel{\T}{\circ}(u_y\stackrel{\T}{\circ}u_z))\\
		& =&  \sigma_{x,yz}\al_{x}(\sigma_{y,z}1_{x^\m})(u_x \stackrel{\T}{\circ}(u_y\stackrel{\T}{\circ}u_z))\\
		& = & \sigma_{x,yz}\al_{x}(\sigma_{y,z}1_{x^\m})\widetilde{\beta_{x,y,z}}((u_x \stackrel{\T}{\circ}u_y)\stackrel{\T}{\circ}u_z)\\
		& \stackrel{(\ref{betaedelta2desigma})}{=} & \sigma_{xy,z}{\sigma_{x,y}}((u_x \stackrel{\T}{\circ}u_y)\stackrel{\T}{\circ}u_z)\\
		& = & \sigma_{xy,z}({\sigma_{x,y}}(u_x \stackrel{\T}{\circ}u_y)\stackrel{\T}{\circ}u_z)\\
		& =& \bar{f}_{xy,z}^\T(\bar{f}_{x,y}^\T(u_x\ot u_y)\ot u_z).
	\end{eqnarray*}
	Consequently,  $\bar{f}^\T=\{\bar{f}_{x,y}^\T, \ x,y\in G\}$ is a factor set for   $\T$.
\end{dem}

The next fact  follows from the proof of  Proposition~\ref{3cocicloquecorrigeaassociatividade} and (\ref{mercomutam}).

\begin{cor}\label{obs3cocicloparasimilares} Let $\Gm:G\longrightarrow \Pics(R)$ be a unital partial representation with  $\Gm_x\ot_R\Gm_{x^\m}\simeq R1_x$ and $\Gm_x|\T_x$, for all  $x\in G$. Let, furthermore,  $f^\Gm=\{f_{x,y}^\Gm:\Gm_x\ot \Gm_y\longrightarrow1_x\Gm_{xy}, \ \in x,y \in G\}$ be a family of $R$-bimodule isomorphisms which satisfies the commutative diagrams in (\ref{diagramasdaunidade}) and let $\beta^{\Gm}_{x,y,z}:1_x1_{xy}\T_{xyz}\longrightarrow 1_x1_{xy}\T_{xyz}$ be an $R$-bimodule isomorphism such that 
	\begin{equation*}
		\beta^\Gm_{x,y,z}\circ f^\Gm_{x,yz}\circ ({\Gm_x}\ot f^\Gm_{y,z})=f^\Gm_{xy,z}\circ (f^\Gm_{x,y}\ot {\Gm_z}), \ \forall \ x,y,z\in G. 
	\end{equation*}
	Then $\widetilde{\beta_{-,-,-}^\Gm}\in Z_{\T}^3(G,\al,\Z)$.
\end{cor}


The next result establishes  the uniqueness of the class in $H^3_\T(G,\al,\Z)$ of the  $3$-cocycle given in  Corollary~\ref{obs3cocicloparasimilares}.

\begin{pro}\label{proparaboadefinicaodedelta}
Let $\Gm:G\longrightarrow \Pics(R)$ be a unital partial representation with  $\Gm_x\ot_R\Gm_{x^\m}\simeq  R1_x$  and $\Gm_x|\T_x$, for all  $x\in G .$ Let, furthermore,  $f^\Gm=\{f_{x,y}^\Gm:\Gm_x\ot_R\Gm_y\longrightarrow1_x\Gm_{xy},\ x,y\in G\}$
be a  family of $R$-bimodule isomorphisms which satisfies  the commutative diagrams in  (\ref{diagramasdaunidade}).  { Then the class   $[\widetilde{\beta_{-,-,-}^\Gm}]$ in $H^3_\T(G,\al,\Z)$ does not depend  on the choice of the representatives in   $[\Gm_x],$ nor on the choice of   $R$-bimodule isomorphisms  $f^\Gm$.}
\end{pro}
\begin{dem} Let $\Om_x \in [\Gm_x].$ Then there exists an $R$-bimodule isomorphism  $a_x:\Om_x\longrightarrow \Gm_x$, for each  $x \in G$. Let 
\begin{equation*}
	f^\Om=\{f_{x,y}^\Om:\Om_x\ot_R\Om_y\longrightarrow1_x\Om_{xy}, \ x,y \in G\}, \ \ f^\Gm=\{f_{x,y}^\Gm:\Gm_x\ot_R\Gm_y\longrightarrow1_x\Gm_{xy},\ x,y\in G\}
\end{equation*}
be  families of $R$-module isomorphisms which satisfy the commutative diagrams in  (\ref{diagramasdaunidade}), and let     $\widetilde{\beta_{-,-,-}^\Om}$ and $\widetilde{\beta_{-,-,-}^\Gm}$ be the  $3$-cocycles  associated with  $f^\Om$ and $f^\Gm$, respectively. By definition, we have that  
\begin{equation*}
	\beta_{x,y,z}^\Om \circ f_{x,yz}^\Om\circ (\Om_x\ot f_{y,z}^\Om)=f_{xy,z}^\Om\circ (f_{x,y}^\Om\ot \Om_z),
\end{equation*}
\begin{equation*}
	\beta_{x,y,z}^\Gm \circ f_{x,yz}^\Gm\circ (\Gm_x\ot f_{y,z}^\Gm)=f_{xy,z}^\Gm\circ (f_{x,y}^\Gm\ot \Gm_z).
\end{equation*}
Denote by  $b_{x,y}:1_x\Om_{xy}\longrightarrow1_x\Om_{xy}$  the $R$-bimodule isomorphism  defined by the commutative diagram 
\begin{equation*}
	\xymatrix{  \Om_x\ot_R\Om_y\ar[rr]^{f_{x,y}^\Om}\ar[d]_{a_x\ot a_y} & & 1_x\Om_{xy}\ar@{-->}[d]^{b_{x,y}}\\
		\Gm_x\ot_R\Gm_y\ar[rd]_{f_{x,y}^\Gm} & & 1_x\Om_{xy}\\
		& 1_x\Gm_{xy} \ar[ru]_{a_{xy}^{\m}}&   }
\end{equation*}
that is, 
\begin{equation*}
	a_{xy}\circ b_{x,y}\circ f_{x,y}^\Om=f_{x,y}^\Gm\circ (a_x\ot a_y), \ \ \forall \ x,y \in G.
\end{equation*}
Let $\widetilde{b_{x,y}}\in \U(\Z1_x1_{xy})$ be the image of  $b_{x,y}$ under the isomorphism of  Lemma~\ref{gammatil}. By (\ref{relacaodefcomftil}) we obtain that  $\widetilde{b_{x,y}}t_{xy}=b_{x,y}(1_xt_{xy})$, for all  $t_{xy}\in \Om_{xy}$. Writing $f_{x,y}^\Om(t_x\ot t_y)=(t_x\stackrel{\Om}{\circ} t_y)$ and $f_{x,y}^\Gm(v_x\ot v_y)=(v_x\stackrel{\Gm}{\circ}v_y)$, we have:
\begin{equation}
	\widetilde{\beta_{x,y,z}^\Om}(t_x\stackrel{\Om}{\circ}(t_y\stackrel{\Om}{\circ}t_z))=((t_x\stackrel{\Om}{\circ}t_y)\stackrel{\Om}{\circ}t_z), \label{betaOmega}
\end{equation}
\begin{equation}
	\widetilde{\beta_{x,y,z}^\Gm}(v_x\stackrel{\Gm}{\circ}({v_y}\stackrel{\Gm}{\circ}v_z))=((v_x\stackrel{\Gm}{\circ}v_y)\stackrel{\Gm}{\circ}v_z), \label{betaGamma}
\end{equation}
\begin{equation}
	\widetilde{b_{x,y}}a_{xy}(t_x\stackrel{\Om}{\circ}t_y)=(a_x(t_x)\stackrel{\Gm}{\circ}a_y(t_y)). \label{bxytilcomaxy}
\end{equation}

Given $x,y, z\in G$ and $t_x \in \Om_x,t_y\in \Om_y, t_z\in \Om_z$ we compute:
\begin{eqnarray*}
	a_{xyz}((t_x\stackrel{\Om}{\circ}t_y)\stackrel{\Om}{\circ}t_z) & \stackrel{(\ref{bxytilcomaxy})}{=} & \widetilde{b_{xy,z}^\m}(a_{xy}(t_x\stackrel{\Om}{\circ}t_y) \stackrel{\Gm}{\circ}a_z(t_z))\\
	& \stackrel{(\ref{bxytilcomaxy})}{=} & \widetilde{b_{xy,z}^\m}\widetilde{b_{x,y}^\m}((a_{x}(t_x)\stackrel{\Gm}{\circ}a_y(t_y)) \stackrel{\Gm}{\circ}a_z(t_z))\\
	& \stackrel{(\ref{betaGamma})}{=} & \widetilde{b_{xy,z}^\m}\widetilde{b_{x,y}^\m}\widetilde{\beta_{x,y,z}^\Gm}(a_{x}(t_x)\stackrel{\Gm}{\circ}(a_y(t_y) \stackrel{\Gm}{\circ}a_z(t_z)))\\
	& \stackrel{(\ref{bxytilcomaxy})}{=} & \widetilde{b_{xy,z}^\m}\widetilde{b_{x,y}^\m}\widetilde{\beta_{x,y,z}^\Gm}(a_{x}(t_x)\stackrel{\Gm}{\circ}\widetilde{b_{y,z}}a_{yz}(t_y\stackrel{\Om}{\circ}t_z))\\
	& \stackrel{(\ref{uxercomutam})}{=} & \widetilde{b_{xy,z}^\m}\widetilde{b_{x,y}^\m}\widetilde{\beta_{x,y,z}^\Gm}\al_x(\widetilde{b_{y,z}}1_{x^\m})(a_{x}(t_x)\stackrel{\Gm}{\circ}a_{yz}(t_y\stackrel{\Om}{\circ}t_z))\\
	& \stackrel{(\ref{bxytilcomaxy})}{=} & \widetilde{b_{xy,z}^\m}\widetilde{b_{x,y}^\m}\widetilde{\beta_{x,y,z}^\Gm}\al_x(\widetilde{b_{y,z}}1_{x^\m})\widetilde{b_{x,yz}}a_{xyz}(t_x\stackrel{\Om}{\circ}(t_y\stackrel{\Om}{\circ}t_z))\\
	& \stackrel{(\ref{betaOmega})}{=} & \widetilde{b_{xy,z}^\m}\widetilde{b_{x,y}^\m}\widetilde{\beta_{x,y,z}^\Gm}\al_x(\widetilde{b_{y,z}}1_{x^\m})\widetilde{b_{x,yz}}\widetilde{\beta_{x,y,z}^\Om}^\m a_{xyz}((t_x\stackrel{\Om}{\circ}t_y)\stackrel{\Om}{\circ}t_z)\\
	& = &  a_{xyz}(\widetilde{b_{xy,z}^\m}\widetilde{b_{x,y}^\m}\widetilde{\beta_{x,y,z}^\Gm}\al_x(\widetilde{b_{y,z}}1_{x^\m})\widetilde{b_{x,yz}}\widetilde{\beta_{x,y,z}^\Om}^\m((t_x\stackrel{\Om}{\circ}t_y)\stackrel{\Om}{\circ}t_z)).
\end{eqnarray*} 
Since $a_{xyz}$ is an isomorphism, it follows that 
\begin{equation*}
	((t_x\stackrel{\Om}{\circ}t_y)\stackrel{\Om}{\circ}t_z)=\widetilde{b_{xy,z}^\m}\widetilde{b_{x,y}^\m}\widetilde{\beta_{x,y,z}^\Gm}\al_x(\widetilde{b_{y,z}}1_{x^\m})\widetilde{b_{x,yz}}\widetilde{\beta_{x,y,z}^\Om}^\m((t_x\stackrel{\Om}{\circ}t_y)\stackrel{\Om}{\circ}t_z),
\end{equation*}
for all   $t_x \in \Om_x,t_y\in \Om_y, t_z\in \Om_z$. { Using the same argument of \cite[Lemma 4.3]{DoRo}}, we conclude:
\begin{equation*}
	\widetilde{b_{xy,z}^\m}\widetilde{b_{x,y}^\m}\widetilde{\beta_{x,y,z}^\Gm}\al_x(\widetilde{b_{y,z}}1_{x^\m})\widetilde{b_{x,yz}}\widetilde{\beta_{x,y,z}^\Om}^\m=1_x1_{xy}1_{xyz}.
\end{equation*}
Hence, 
\begin{eqnarray*}
	\widetilde{\beta_{x,y,z}^\Om} & = & \widetilde{\beta_{x,y,z}^\Gm}\widetilde{b_{xy,z}^\m}\widetilde{b_{x,y}^\m}\al_x(\widetilde{b_{y,z}}1_{x^\m})\widetilde{b_{x,yz}}\\
	& = & \widetilde{\beta_{x,y,z}^\Gm}(\delta^2\widetilde{b_{-,-}})(x,y,z),
\end{eqnarray*}
where 
\begin{equation*}
	\begin{array}{c c c l}
		\widetilde{b_{-,-}}:& G\times G& \longrightarrow&  \U(\Z),\\
		& (x,y) & \longmapsto & \widetilde{b_{x,y}} \in \U(\Z1_x1_{xy}).
	\end{array}
\end{equation*}
Consequently, $[\widetilde{\beta_{-,-,-}^\Om}]=[\widetilde{\beta_{-,-,-}^\Gm}]$ in $H^3_{\T}(G,\al, \Z).$
\end{dem}

We proceed with the following technical fact.

\begin{lem}\label{conjuntodefatores}  Let $\T,\Om,\Gm:G\longrightarrow \Pics(R)$ be unital partial representations with  $\Gm_x|\T_x$, $\Om_x|\T_x$ and $\Om_x\ot_R\Om_{x^\m}\simeq \Gm_x\ot_R \Gm_{x^\m}\simeq \T_x\ot_R \T_{x^\m}\simeq R1_x$, for all  $x \in G$. Let, furthermore,  $f^\T=\{f_{x,y}^\T:\T_x\ot_R \T_y\longrightarrow 1_x\T_{xy}, \ x,y \in G\}$ be a factor set for  $\T ,$ let 
	\begin{equation*}
		f^\Gm=\{f^\Gm_{x,y}:\Gm_x\ot_R \Gm_y\longrightarrow 1_x\Gm_{xy}, \ x,y \in G\}, \ f^\Om=\{f^\Om_{x,y}:\Om_x\ot_R \Om_y\longrightarrow 1_x\Om_{xy}, \ x,y\in G \}
	\end{equation*}
	be families of $R$-bimodule isomorphisms which satisfy the commutative diagrams in  (\ref{diagramasdaunidade}) 
	and let  $\widetilde{\beta^\Gm_{-,-,-}},$ and    $\widetilde{\beta^\Om_{-,-,-}}$ be the  $3$-cocycles in $Z^3_\T(G,\al,\Z)$ associated to $f^\Gm$ and $f^\Om$, respectively. Then, 
	$$\begin{array}{c c c l}
		\Lambda: & G & \longrightarrow & \Pics(R),\\
		& x  & \longmapsto & [\Gamma_x\ot_R\T_{x^\m}\ot_R\Om_x]
	\end{array}$$
	is a unital partial  representation with   $\Lambda_x\ot_R \Lambda_{x^\m}\simeq R1_x$ and $\Lambda_x|\T_x$, for all  $x\in G$. The $3$-cocycle  associated with the family of $R$-bimodule isomorphisms  
	$f^\Lambda=\{f_{x,y}^\Lambda:\Lambda_x\ot_R \Lambda_y\longrightarrow 1_x\Lambda_{xy}, \ x,y\in G \}$  defined by  
	\begin{equation*}
		\xymatrix@C=1.0cm{  \Gamma_x\ot_R \T_{x^\m}\ot_R \Omega_x\ot_R \Gamma_y\ot_R \T_{y^\m}\ot_R \Omega_y \ar[rd]^{\ \ \ \ \ \  \ \ \ \ \Gamma_x\ot T_{\T_{x^\m}\ot\Omega_x,\Gamma_y\ot \T_{y^\m}}\ot \Om_y}\ar@/_1.7cm/@{-->}[rdddd]_{f^{\Lambda}_{x,y}} & & \\
			&  \Gamma_x\ot_R\Gamma_y\ot_R \T_{y^\m} \ot_R \T_{x^\m}\ot_R \Omega_x \ot_R\Omega_y\ar[dd]|{f_{x,y}^\Gamma\ot f_{y^{\m},x^\m}^\T\ot f_{x,y}^\Omega} \\
			&  \\
			& 1_x\Gamma_{xy}\ot_R 1_{y^\m}\T_{(xy)^\m}\ot_R 1_x\Omega_{xy}\ar@{=}[d]  \\
			&  1_x\Gamma_{xy}\ot_R \T_{(xy)^\m}\ot_R \Omega_{xy},  }
	\end{equation*}
	is given by  $\widetilde{\beta^\Lambda_{-,-,-}}=\widetilde{\beta^\Gm_{-,-,-}}\widetilde{\beta^\Om_{-,-,-}},$
	where $T_{-,-}$ is the   isomorphism from Proposition~\ref{isomorfismoT}. 
	{If $f^\Gm$ and $f^\Om$ are factor sets for  $\Gm$ and $\Om$, respectively, then  $f^\Lambda$ is a factor set for  $\Lambda$ and    $\D(\Lambda)$ is a generalized partial crossed product. }
\end{lem}

\begin{dem} Clearly $\Lambda_1=R$. By the above diagram there is an $R$-bimodule isomorphism  $\Lambda_x\ot_R\Lambda_y\simeq R1_x\ot_R\Lambda_{xy}$, where $x,y\in G$. In particular, $\Lambda_x\ot_R\Lambda_{x^\m}\simeq R1_x$. Moreover, it is easy { to see} that $\Lambda_x\ot_R R1_{x^\m}\simeq \Lambda_x\simeq R1_x\ot_R\Lambda_x$ and $\Lambda_x\ot_RR1_y\simeq R1_{xy}\Lambda_x$, for all $x,y \in G$. Then,
	$$	\Lambda_{x^\m}\ot_R \Lambda_x\ot_R \Lambda_y \simeq \Lambda_{x^\m}\ot_R R1_x\ot_R\Lambda_{xy}\simeq \Lambda_{x^\m}\ot_{R} \Lambda_{xy}.$$
	On the other hand,
	$$\Lambda_x\ot_R \Lambda_y\ot_R \Lambda_{y^\m}  \simeq R1_x\ot_R \Lambda_{xy}\ot_R \Lambda_{y^\m}
		\simeq \Lambda_{xy}\ot_R R1_{y^\m}\ot_R\Lambda_{y^\m}
		\simeq  \Lambda_{xy}\ot_R\Lambda_{y^\m}.$$
	
	Thus, $\Lambda$ is a unital partial representation.
	Since $\Gm_x|\T_x$ and $\Om_x|\T_x$, then $\Lambda_x\simeq \Gm_x\ot \T_{x^\m}\ot \Om_x|\T_x$, for all $x \in G$.
%
%
{		By definition,
		\begin{equation}
		\widetilde{\beta_{x,y,z}^\Gm}(v_x\stackrel{\Gm}{\circ}(v_y\stackrel{\Gm}{\circ}v_z))=((v_x\stackrel{\Gm}{\circ}v_y)\stackrel{\Gm}{\circ}v_z), \ \forall \ v_i \in \Gm_i,\label{betaGm}
	\end{equation}
	\begin{equation}
		\widetilde{\beta_{x,y,z}^\Om}(w_x\stackrel{\Om}{\circ}(w_y\stackrel{\Om}{\circ}w_z))=((w_x\stackrel{\Om}{\circ}w_y)\stackrel{\Om}{\circ}w_z), \ \forall \ w_i \in \Om_i.\label{betaOm}
	\end{equation}}

{	Let us  show that  $\widetilde{\beta^\Gm_{x,y,z}}\widetilde{\beta^\Om_{x,y,z}}(\lambda_x\stackrel{\Lambda}{\circ} (\lambda_y\stackrel{\Lambda}{\circ}\lambda_z))=((\lambda_x\stackrel{\Lambda}{\circ} \lambda_y)\stackrel{\Lambda}{\circ}\lambda_z)$, for all  $\lambda_i \in \Lambda_i$.}
	Consider the $R$-bimodule isomorphisms $T_1=T_{\T_{x^\m}\ot\Om_x,\Gm_y\ot \T_{y^\m}}$ and $T_2=T_{\T_{(xy)^\m}\ot \Om_{xy},\Gm_z\ot\T_{z^\m}} $. Note that $f_{xy,z}^\Lambda\circ (f_{x,y}^\Lambda\ot \Lambda_z)$ is given by:
	\raisetag{2cm}
	{\small 	\begin{equation}
			\xymatrix@C=-2cm{  \Gamma_x\ot \T_{x^\m}\ot \Omega_x\ot\Gamma_y\ot \T_{y^\m}\ot\Omega_y\ot \Gamma_z\ot \T_{z^\m}\ot \Omega_z\ar@{-->}@/_2.5cm/[rdddd]\ar[rd]^-{\Gm_x\ot T_1\ot \Om_y\ot \Lambda_z}& \\
				& \Gamma_x\ot\Gamma_y\ot \T_{y^\m}\ot \T_{x^\m}\ot \Omega_x\ot\Omega_y\ot \Gamma_z\ot \T_{z^\m}\ot \Omega_z\ar[d]|-{f_{x,y}^\Gm\ot f_{y^\m,x^\m}^\T\ot f_{x,y}^\Om\ot  \Lambda_z}\\
				&  1_x\Gamma_{xy}\ot \T_{(xy)^\m}\ot \Omega_{xy}\ot \Gamma_z\ot \T_{z^\m}\ot \Omega_z\ar[d]^{1_x\Gm_{xy}\ot T_2\ot \Om_z} \\
				& 1_x\Gamma_{xy}\ot \Gamma_z\ot \T_{z^\m}\ot \T_{(xy)^\m}\ot \Omega_{xy}\ot \Omega_z\ar[d]|-{f_{xy,z}^\Gm\ot f_{z^\m,(xy)^\m}^\T\ot f_{xy,z}^\Om}\\
				&  1_x1_{xy}\Gamma_{xyz}\ot\T_{(xyz)^\m}\ot \Omega_{xyz}   } \label{d1}
	\end{equation}}
	
	On the other hand, let $T_3=T_{\T_{y^\m\ot \Om_y,\Gm_z\ot\T_{z^\m}}}$ and $T_4=T_{\T_{x^\m}\ot\Om_x,\Gm_{xy}\ot \T_{(yz)^\m}}$, then $f_{x,yz}^\Lambda\circ (\Lambda_x\ot f_{y,z}^\Lambda)$ is given by
	{\small 		\begin{equation}
			\xymatrix@C=-2.5cm{ \Gamma_x\ot \T_{x^\m}\ot \Omega_x\ot\Gamma_y\ot \T_{y^\m}\ot\Omega_y\ot \Gamma_z\ot \T_{z^\m}\ot \Omega_z\ar@{-->}@/_2.5cm/[rdddd]\ar[rd]^-{\Lambda_x\ot \Gm_y\ot T_3\ot \Om_z} & \\
				& 	\Gamma_x\ot \T_{x^\m}\ot \Omega_x\ot\Gamma_y\ot \Gamma_z\ot \T_{z^\m}\ot \T_{y^\m}\ot\Omega_y\ot \Omega_z\ar[d]|-{\Lambda_x\ot f_{y,z}^\Gm\ot f_{z^\m,y^\m}^\T\ot f_{y,z}^\Om}\\
				& \Gamma_x\ot \T_{x^\m}\ot \Omega_x\ot1_y\Gamma_{yz}\ot \T_{(yz)^\m}\ot\Omega_{yz}\ar[d]^{\Gm_x\ot T_4\ot \Om_{yz}}\\
				&  1_{xy}\Gamma_x\ot\Gamma_{yz}\ot \T_{(yz)^\m}\ot \T_{x^\m}\ot \Omega_x\ot\Omega_{yz}\ar[d]|-{f_{x,yz}^\Gm\ot f_{(yz)^\m,x^\m}^\T\ot f_{x,yz}^\Om}\\
				& 1_x1_{xy}\Gamma_{xyz}\ot\T_{(xyz)^\m}\ot \Omega_{xyz} , } \label{d2}
	\end{equation}}
	
	We will use  Proposition \ref{isomorfismoT} to build the isomorphisms $T_i$,  $i=1,2,3,4$. To this end consider the $R$-bilinear maps $\varphi_i:\T_{x^\m}\ot_R \Omega_x\rightarrow R$ and $\psi_i:R\rightarrow \T_{x^\m}\ot_R \Omega_x$, with $i=1,2,...,n$, such that $\dsum_{i=1}^n\psi_i\varphi_i=Id_{\T_{x^\m}\ot \Omega_x},$ and the $R$-bilinear maps $\varphi_j':\Gamma_z\ot_R \T_{z^\m}\rightarrow R$ and $\psi_j':R\rightarrow \Gamma_z\ot_R \T_{z^\m}$, with $j=1,2,...,m$,  such that $\dsum_{j=1}^m\psi_j'\varphi_j'=Id_{\Gamma_z\ot \T_{z^\m}}$. 	Denote
	$$\psi_i(1)=\dsum_{l}u_{x^{-1}}^{i,l}\ot \omega_x^{i,l} \ \ \mbox{and} \ \ \ \psi'_j(1)=\dsum_{k}\widetilde{v}_z^{j,k}\ot \widetilde{u}_{z^\m}^{j,k}.$$
	 By Proposition~\ref{isomorfismoT}, we have
	$$\begin{array}{c c c l}
		T_1:& \T_{x^\m}\ot_R \Omega_x\ot_R \Gamma_y\ot_R \T_{y^\m}& \rightarrow&  \Gamma_y\ot_R \T_{y^\m}\ot_R\T_{x^\m}\ot_R \Omega_x\\
		& u_{x^\m}\ot \omega_x\ot v_y\ot u_{y^\m} & \mapsto & \dsum_{i,l}\varphi_i(u_{x^\m}\ot \omega_x)v_y\ot u_{y^\m}\ot u_{x^{-1}}^{i,l}\ot \omega_x^{i,l},
	\end{array}$$
	$$\begin{array}{c c c l}
		T_2:& \T_{(xy)^\m}\ot_R \Omega_{xy}\ot_R\Gamma_z\ot_R \T_{z^\m}& \rightarrow&  \Gamma_z\ot_R \T_{z^\m}\ot_R\T_{(xy)^\m}\ot_R \Omega_{xy}\\
		& u_{(xy)^\m}\ot \omega_{xy}\ot v_{z}\ot u_{z^\m} & \mapsto & \dsum_{j,k}\widetilde{v}_z^{j,k}\ot \widetilde{u}_{z^\m}^{j,k}\ot u_{(xy)^\m}\ot \omega_{xy}\varphi_j'(v_z\ot u_{z^\m}),
	\end{array}$$
	$$\begin{array}{c c c l}
		T_3:& \T_{y^\m}\ot_R \Omega_y\ot_R\Gamma_z\ot_R \T_{z^\m}& \rightarrow&  \Gamma_z\ot_R \T_{z^\m}\ot_R\T_{y^\m}\ot_R \Omega_y\\
		& u_{y^\m}\ot \omega_{y}\ot v_{z}\ot u_{z^\m} & \mapsto & \dsum_{j,k}\widetilde{v}_z^{j,k}\ot \widetilde{u}_{z^\m}^{j,k}\ot u_{y^\m}\ot \omega_{y}\varphi_j'(v_z\ot u_{z^\m}),
	\end{array}$$
	$$\begin{array}{c c c l}
		T_4:& \T_{x^\m}\ot_R \Omega_x\ot_R \Gamma_{yz}\ot_R \T_{(yz)^\m}& \rightarrow&  \Gamma_{yz}\ot_R \T_{(yz)^\m}\ot_R\T_{x^\m}\ot_R \Omega_x\\
		& u_{x^\m}\ot \omega_x\ot v_{yz}\ot u_{(yz)^\m} & \mapsto & \dsum_{i,l}\varphi_i(u_{x^\m}\ot \omega_x)v_{yz}\ot u_{(yz)^\m}\ot u_{x^{-1}}^{i,l}\ot \omega_x^{i,l}.
	\end{array}$$

	In diagram (\ref{d1}) the composition $f^\Lambda_{x,yz}\circ (\Lambda_x\ot f_{y,z}^\Lambda)$ is given by:
	\begin{eqnarray*}
		& & v_x\ot u_{x^\m}\ot \omega_x\ot v_y\ot u_{y^\m}\ot \omega_y\ot v_z\ot u_{z^\m}\ot \omega_z\\
		&\stackrel{T_3}{\mapsto}& \dsum_{j,k}v_x\ot u_{x^\m}\ot \omega_x\ot v_y\ot\widetilde{v}_z^{j,k}\ot \widetilde{u}_{z^\m}^{j,k}\ot u_{y^\m}\ot \omega_{y}\varphi_j'(v_z\ot u_{z^\m})\ot \omega_z\\
		& \mapsto & \dsum_{j,k}v_x\ot u_{x^\m}\ot \omega_x\ot ( v_y\stackrel{\Gm}{\circ}\widetilde{v}_z^{j,k})\ot (\widetilde{u}_{z^\m}^{j,k}\stackrel{\T}{\circ}u_{y^\m})\ot (\omega_{y}\varphi_j'(v_z\ot u_{z^\m})\stackrel{\Om}{\circ} \omega_z)\\
		& \stackrel{T_4}{\mapsto} & \dsum_{i,j,k,l}v_x\ot \varphi_i(u_{x^\m}\ot \omega_x)( v_y\stackrel{\Gm}{\circ}\widetilde{v}_z^{j,k})\ot (\widetilde{u}_{z^\m}^{j,k}\stackrel{\T}{\circ}u_{y^\m})\ot u_{x^{-1}}^{i,l}\ot \omega_x^{i,l}\ot  (\omega_{y}\varphi_j'(v_z\ot u_{z^\m})\stackrel{\Om}{\circ} \omega_z)\\
		& \mapsto & \dsum_{i,j,k,l}(v_x\stackrel{\Gm}{\circ} (\varphi_i(u_{x^\m}\ot \omega_x) v_y\stackrel{\Gm}{\circ}\widetilde{v}_z^{j,k}))\ot ((\widetilde{u}_{z^\m}^{j,k}\stackrel{\T}{\circ}u_{y^\m})\stackrel{\T}{\circ} u_{x^{-1}}^{i,l})\ot (\omega_x^{i,l}\stackrel{\Om}{\circ} (\omega_{y}\varphi_j'(v_z\ot u_{z^\m})\stackrel{\Om}{\circ} \omega_z))
	\end{eqnarray*}

	On the other hand, in (\ref{d2}) the composition 
	$f_{xy,z}^\Lambda\circ(f_{x,y}^\Lambda\ot \Lambda_z)$ is given by:
	\begin{eqnarray*}
		& & v_x\ot u_{x^\m}\ot \omega_x\ot v_y\ot u_{y^\m}\ot \omega_y\ot v_z\ot u_{z^\m}\ot \omega_z\\
		&\stackrel{T_1}{\mapsto} & \dsum_{i,l} v_x\ot \varphi_i(u_{x^\m}\ot \omega_x)v_y\ot u_{y^\m}\ot u_{x^{-1}}^{i,l}\ot \omega_x^{i,l}\ot \omega_y\ot v_z\ot u_{z^\m}\ot \omega_z\\
		&\mapsto & \dsum_{i,l} (v_x\stackrel{\Gm}{\circ} \varphi_i(u_{x^\m}\ot \omega_x)v_y)\ot (u_{y^\m}\stackrel{\T}{\circ}u_{x^{-1}}^{i,l})\ot (\omega_x^{i,l}\stackrel{\Om}{
				\circ} \omega_y)\ot v_z\ot u_{z^\m}\ot \omega_z\\
		&\stackrel{T_2}{\mapsto} & \dsum_{i,j,k,l} (v_x\stackrel{\Gm}{\circ} \varphi_i(u_{x^\m}\ot \omega_x)v_y)\ot \widetilde{v}_z^{j,k}\ot \widetilde{u}_{z^\m}^{j,k}\ot (u_{y^\m}\stackrel{\T}{\circ}u_{x^{-1}}^{i,l})\ot (\omega_x^{i,l}\stackrel{\Om}{
				\circ} \omega_y)\varphi_j'( v_z\ot u_{z^\m})\ot \omega_z\\
		&\mapsto & \dsum_{i,j,k,l} ((v_x\stackrel{\Gm}{\circ} \varphi_i(u_{x^\m}\ot \omega_x)v_y)\stackrel{\Gm}{\circ} \widetilde{v}_z^{j,k})\ot (\widetilde{u}_{z^\m}^{j,k}\stackrel{\T}{\circ} (u_{y^\m}\stackrel{\T}{\circ}u_{x^{-1}}^{i,l})) \ot ((\omega_x^{i,l}\stackrel{\Om}{
				\circ} \omega_y)\varphi_j'( v_z\ot u_{z^\m})\stackrel{\Om}{
				\circ} \omega_z)\\
			& \stackrel{(\ref{betaGm}), (\ref{betaOm})}{=} & {\dsum_{i,j,k,l}	\widetilde{\beta_{x,y,z}^\Gm}(v_x\stackrel{\Gm}{\circ} (\varphi_i(u_{x^\m}\ot \omega_x) v_y\stackrel{\Gm}{\circ}\widetilde{v}_z^{j,k}))\ot (\widetilde{u}_{z^\m}^{j,k}\stackrel{\T}{\circ} (u_{y^\m}\stackrel{\T}{\circ}u_{x^{-1}}^{i,l})) \ot 	\widetilde{\beta_{x,y,z}^\Om}(\omega_x^{i,l}\stackrel{\Om}{\circ} (\omega_{y}\varphi_j'(v_z\ot u_{z^\m})\stackrel{\Om}{\circ} \omega_z))}\\  
				& = & {\dsum_{i,j,k,l}	\widetilde{\beta_{x,y,z}^\Gm}\widetilde{\beta_{x,y,z}^\Om}(v_x\stackrel{\Gm}{\circ} (\varphi_i(u_{x^\m}\ot \omega_x) v_y\stackrel{\Gm}{\circ}\widetilde{v}_z^{j,k}))\ot ((\widetilde{u}_{z^\m}^{j,k}\stackrel{\T}{\circ}u_{y^\m})\stackrel{\T}{\circ} u_{x^{-1}}^{i,l})\ot 	(\omega_x^{i,l}\stackrel{\Om}{\circ} (\omega_{y}\varphi_j'(v_z\ot u_{z^\m})\stackrel{\Om}{\circ} \omega_z)),}
	\end{eqnarray*}
\noindent {where the last equality follows by (\ref{mercomutam}) and associativity of $\D(\T)$.}	Thus we conclude that
	\begin{equation*}
	{\beta_{x,y,z}^\Gm\beta_{x,y,z}^\Om\circ	f_{x,yz}^\Lambda\circ (\Lambda_x\ot f_{y,z}^\Lambda)=f_{xy,z}^\Lambda\circ (f_{x,y}^\Lambda\ot \Lambda_z), \ x,y,z \in G. }
	\end{equation*}
	{In particular, if $f^\Gm$ is a factor set for $\Gm$ and $f^\Om$ is a factor set for $\Om$, then $\beta^\Gm$ and $\beta^\Om$ are trivial, thus $\beta^\Gm\beta^\Om$ is trivial too.} Therefore $f^\Lambda$ is a factor set for $\Lambda.$
	
\end{dem}

\begin{defi} Let $\D(\T)$ and $\D(\Gm)$ be  partial generalized crossed products with $R$-bimodule and ring isomorphism   $i:R\longrightarrow \T_1$ and $j:R\longrightarrow \Gm_1$, respectively, as in Remark \ref{unidadeediagramasemDT}. A morphism of partial generalized crossed products  $F:\D(\T)\longrightarrow \D(\Gm)$ is a set of $R$-bimolude morphisms  $\{F_x:\T_x\longrightarrow \Gm_x, x \in G \}$ such that $F_1\circ i=j$ and  the following  commutative diagram is satisfied:
	\begin{equation}
	\xymatrix{ \T_x\ot_R\T_y\ar[rr]^{f_{x,y}^\T}\ar[dd]_{F_x\ot F_y} & & 1_x\T_{xy}\ar[dd]^{F_{xy}}\\
		& & \\
		\Gm_x\ot_R\Gm_y\ar[rr]_{f_{x,y}^\Gm} & & 1_x\Gm_{xy} } \label{morfismodepcgp}
	\end{equation}
	A morphism $F$ of partial generalized crossed products is called an isomorphism if each  morphism $F_x:\T_x\longrightarrow \Gm_x$ is an $R$-bimodule isomorphism.
\end{defi}

{We denote  by  $[\D(\Omega)]$ the isomorphism  class of the partial generalized crossed product  $\D(\Omega)$.}
\begin{obs}
If  $F=\{F_x:\T_x\longrightarrow \Gm_x, \ x\in G \}$  is a family of $R$-bimodule isomorphism satisfying the commutative diagram (\ref{morfismodepcgp}), then $F$ is an isomorphism of partial generalized crossed products. Indeed, it is enough to show that  $F_1\circ i=j$. It follows from Remark~\ref{unidadeediagramasemDT} and the commutative diagram (\ref{morfismodepcgp}) that { the} following diagram is commutative:
\begin{equation*}
\xymatrix{ R\ot_R \T_x\ar[rrdd]_{\simeq}\ar[rr]^{i\ot \T_x}&  & \T_1\ot_R \T_x\ar[rr]^{F_1\ot F_x}\ar[dd]_{f_{1,x}^\T} & & \Gm_1\ot_R \Gm_x\ar[dd]^{f_{1,x}^\Gm}& & R\ot_R \Gm_x\ar[lldd]^{\simeq}\ar[ll]_{j\ot \Gm_x}\\
	& 	& & & & & \\
	& 	& \T_x\ar[rr]_{F_x} & & \Gm_x  &  &  }
\end{equation*}
Given $v_x \in \Gm_x$, there is $u_x \in \T_x$ such that $F_x(u_x)=v_x$. Then
$(F_1(i(1))\stackrel{\Gm}{\circ} F_x(u_x))=F_x(u_x),$ that is,
$$(F_1(i(1))\stackrel{\Gm}{\circ}v_x)=v_x, \ \ \mbox{for all } \ \ v_x \in \Gm_x.$$
On the other hand, the commutative diagram
\begin{equation*}
\xymatrix{ \T_x\ot_R R\ar[rrdd]_{\simeq}\ar[rr]^{\T_x\ot i}&  & \T_x\ot_R \T_1\ar[rr]^{F_x\ot F_1}\ar[dd]_{f_{x,1}^\T} & & \Gm_x\ot_R \Gm_1\ar[dd]^{f_{x,1}^\Gm}& & \Gm_x\ot_R R\ar[lldd]^{\simeq}\ar[ll]_{\Gm_x\ot j}\\
	& 	& & & & & \\
	& 	& \T_x\ar[rr]_{F_x} & & \Gm_x  &  &  }
\end{equation*}
implies that $(v_x\stackrel{\Gm}{\circ}F_1(i(1)))=v_x$, for all $v_x \in \Gm_x.$ Thus, $F(i(1))=j(1)$ is the unity of $\D(\Gm)$. Therefore,  $F\circ i=j$.

\end{obs}

\begin{obs}\label{obsinversodefxxm} Let $\T$ be a partial representation and 
 $f^\T=\{f_{x,y}^\T:\T_x\ot_R\T_y\longrightarrow 1_x\T_{xy}, \ x,y \in G\}$ a factor set for $\T$.
  For each $x \in G$, 
we have an  $R$-bimodule isomorphism $$f_{x,x\m}^\T:\T_x\ot_R\T_{x\m}\longrightarrow R 1_x,$$ whose   inverse is
	\begin{equation*}
	\begin{array}{c c c l}
	(f_{x,x^\m}^\T)^\m: &  R1_x & \longrightarrow & \T_x\ot_R \T_{x^\m}\\
	& r1_x & \longrightarrow & \dsum_{(x)}r\om_x\ot \om_{x^\m},
	\end{array}
	\end{equation*}
where $\dsum_{(x)}(\om_x\stackrel{\T}{\circ} \om_{x^\m})=1_x$. Restricting, we obtain an $R$-bimodule isomorphism  $1_xf_{xy,(xy)^\m}^\T:1_x\T_{xy}\ot_R\T_{(xy)^\m}\longrightarrow R1_x1_{xy}.$ Denote
	\begin{equation*}
	1_x=\dsum_{(x)}(\om_x\stackrel{\T}{\circ}\om_{x^\m}) \ \ \mbox{and} \ \ 1_y=\dsum_{(y)}(\om_y\stackrel{\T}{\circ}\om_{y^\m}).
	\end{equation*}
Then, by associativity, we have
	\begin{eqnarray*}
		\dsum_{(x),(y)}\left( (\om_x\stackrel{\T}{\circ}\om_y)\stackrel{\T}{\circ}(\om_{y^\m}\stackrel{\T}{\circ}\om_{x^\m})\right)& = & \dsum_{(x),(y)}\left( \om_x\stackrel{\T}{\circ}(\om_y\stackrel{\T}{\circ}\om_{y^\m})\stackrel{\T}{\circ}\om_{x^\m}\right)\\
		& = & \dsum_{(x)}(\om_x1_y\stackrel{\T}{\circ}\om_{x^\m})=\dsum_{(x)}1_{xy}(\om_x\stackrel{\T}{\circ}\om_{x^\m})\\
		& = & 1_{xy}1_x.
	\end{eqnarray*}
Thus,
	\begin{equation*}
	\begin{array}{c c c l}
	(1_xf_{xy,y^\m x^\m}^\T)^\m: & R1_x1_{xy} & \longrightarrow & 1_x\T_{xy}\ot_R \T_{(xy)^\m}\\
	& r1_x1_{{xy}} & \longmapsto & \dsum_{(x),(y)} r((\om_x\stackrel{\T}{\circ}\om_y)\ot(\om_{y^\m}\stackrel{\T}{\circ}\om_{x^\m})).
	\end{array}
	\end{equation*}
Indeed,  let $r \in R1_x1_{xy}$ and $\om_{xy}\in \T_{xy}, \om_{(xy)^\m}\in \T_{(xy)^\m}$  be such that $\dsum_{(xy)}\om_{xy}\stackrel{\T}{\circ} \om_{(xy)^\m}=1_{xy}$. Then
	\begin{eqnarray*}
		\dsum_{(xy)}r\om_{xy}\ot \om_{(xy)^\m} & = & \dsum_{(xy)}r1_x1_{xy}\om_{xy}\ot \om_{(xy)^\m}\\
		& = & \dsum_{(x),(y),(xy)}r\left( (\om_x\stackrel{\T}{\circ}\om_y)\stackrel{\T}{\circ}(\om_{y^\m}\stackrel{\T}{\circ}\om_{x^\m})\right)\om_{xy}\ot \om_{(xy)^\m}\\
		& = & \dsum_{(x),(y),(xy)}r\left( (\om_x\stackrel{\T}{\circ}\om_y)\stackrel{\T}{\circ}(\om_{y^\m}\stackrel{\T}{\circ}\om_{x^\m})\stackrel{\T}{\circ}\om_{xy}\right)\ot \om_{(xy)^\m}\\
		& = & \dsum_{(x),(y),(xy)}r(\om_x\stackrel{\T}{\circ}\om_y)\underbrace{\left( (\om_{y^\m}\stackrel{\T}{\circ}\om_{x^\m})\stackrel{\T}{\circ}\om_{xy}\right)}_{\in R}\ot \om_{(xy)^\m}\\
		& = & \dsum_{(x),(y),(xy)}r(\om_x\stackrel{\T}{\circ}\om_y)\ot \left( (\om_{y^\m}\stackrel{\T}{\circ}\om_{x^\m})\stackrel{\T}{\circ}\om_{xy}\right)\om_{(xy)^\m}\\
		& = & \dsum_{(x),(y),(xy)}r(\om_x\stackrel{\T}{\circ}\om_y)\ot (\om_{y^\m}\stackrel{\T}{\circ}\om_{x^\m})( \om_{xy}\stackrel{\T}{\circ}\om_{(xy)^\m})\\
		& = & \dsum_{(x),(y)}r(\om_x\stackrel{\T}{\circ}\om_y)\ot (\om_{y^\m}\stackrel{\T}{\circ}\om_{x^\m})1_{xy}\\
		& = & \dsum_{(x),(y)}r(\om_x\stackrel{\T}{\circ}\om_y)\ot (\om_{y^\m}\stackrel{\T}{\circ}1_{y}\om_{x^\m})\\
		& = & \dsum_{(x),(y)}r(\om_x\stackrel{\T}{\circ}\om_y)\ot (\om_{y^\m}\stackrel{\T}{\circ}\om_{x^\m}).\\
	\end{eqnarray*}
	
\end{obs}

\subsection{Unital partial representations $G\to \mathcal{S}_R(S)$.} 
\label{subsection: SRS}	
As we have seen in Section~\ref{sec: pcgp}, if   $\T:G\longrightarrow \Pics(R)$ is a unital partial representation such that we have a generalized partial crossed product $\D(\T)$, then $R\subseteq \D(\T)$ is  a ring extension with the same unity.  In this section we shall see how  to obtain a generalized partial crossed  product from   a ring extension and a certain partial representation.


Let $R\subseteq S$ be a ring extension with the same unity element. Denote by  $\mathcal{S}_R(S)$ the set of the   $R$-subbimodules of $S$. Given  $M,N\in \mathcal{S}_R(S),$ define the  product: 
$$MN=\left\lbrace  \dsum_{i=1}^k m_in_i; \ m_i \in M, n_i \in N\right\rbrace .$$
With the above operation $\mathcal{S}_R(S)$ is a monoid with neutral element $R.$ 

A partial  representation  
$$\begin{array}{c c c l}
	\Theta: & G & \longrightarrow & \mathcal{S}_R(S),\\
	& x & \longmapsto & \T_x
\end{array}$$
will be  called  \textit{unital} if  $\e_x=\T_x\T_{x^{-1}}=R1_x$, where $1_x$ is a central idempotent in $R,$ for each 
$x\in G$.
The equality $\T_x\T_{x^\m}\T_x=\T_x$,  $(x\in G)$, implies  that  $R1_x\T_x=\T_x$ and $\T_xR1_{x^\m}=\T_x .$ Hence $1_xu_x=u_x$ and $u_x1_{x^\m}=u_x$, for all  $u_x \in \T_x$ and $x\in G$. { Furthermore, if $M$ is an  $S$-bimodule then $\T_xM=1_xM=\{m \in M, 1_xm=m\}$. Indeed, obviously, $\T_xM\subseteq 1_xM$. If $m \in 1_xM$, then
	\begin{equation*}
		m=1_xm=\dsum_{(x)}(\om_x\om_{x^\m})m=\dsum_{(x)}\om_x(\om_{x^\m}m) \in \T_xM.
	\end{equation*}
	Analogously, $M\T_{x}=M1_{x^\m}$, for each  $x \in G$.}

\begin{pro}\label{SRemPicS}
	Let $\T:G\longrightarrow \mathcal{S}_R(S)$ be a unital partial representation. Then,
	\begin{itemize}
		 \item[(i)]   $[\T_x]\in \Pics(R)$.
		\item[(ii)] If $M$ is an  $S$-bimodule, then 	$$\begin{array}{c c c}
			\begin{array}{c c c l}
				m_l:& \T_x\ot_R M & \longrightarrow & \T_x M\\
				& u_x\ot m & \longmapsto & u_xm 
			\end{array} & and & \begin{array}{c c c l}
				m_r: & M\ot_R \T_{x} & \longrightarrow &M\T_x \\
				& m\ot u_x & \longmapsto & mu_x
			\end{array}
		\end{array}$$
		are  $R$-$S$-bimodule and  $S$-$R$-bimodule isomorphisms, respectively. 
		\item[(iii)] Let $M$ be an $S$-bimodule and let $N$ be an $R$-subbimodule of $M$. Then we have the following 
		$R$-bimodule isomorphisms: 
		$$\begin{array}{c c c}
			\begin{array}{ c c l}
				\T_x\ot_R N & \longrightarrow & \T_x N\\
				u_x\ot n & \longmapsto & u_xn 
			\end{array} & and & \begin{array}{ c c l}
				N\ot_R \T_{x} & \longrightarrow &N\T_x \\
				n\ot u_{x^\m} & \longmapsto & nu_{x^\m}
			\end{array}
		\end{array}.$$
		\end{itemize}
	
\end{pro}
\begin{dem}
	{(i)} Let $\om_x^i\in \T_x$ and $\om_{x\m}^i \in \T_{x^\m}$ be such that  $\dsum_{i=1}^m\om_x^i\om_{x^\m}^i=1_x$.
	Define,
	\begin{equation*}
		\begin{array}{c c c l}
			f_i: & \T_x & \longrightarrow & R\\
			& u_x & \longmapsto & \om_{x^\m}^iu_x
		\end{array}.
	\end{equation*}
	Then $f_i$ is right  $R$-linear and 
	\begin{equation*}
		\dsum_{i=1}^n\om_x^if_i(u_x)=\dsum_{i=1}^n\om_x^i\om_{x^\m}^iu_x=1_xu_x=u_x,
	\end{equation*}
	for all   $u_x \in \T_x$. Consequently, $\{\om_x^i, f_i\}$ is a dual basis for the right $R$-module $\T_x .$
	
	Analogously, let $\overline{\om}_x^j\in \T_x$ and $\overline{\om}_{x^\m}^j \in \T_{x^\m}$,  $j=1,2,...,m$, be such that  $\dsum_{j=1}^m\overline{\om}_{x^\m}^j\overline{\om}_x^j=1_{x^\m}$. Define
	\begin{equation*}
		\begin{array}{c c c l}
			g_j: & \T_x & \longrightarrow & R,\\
			& u_x & \longrightarrow & u_x\overline{\om}_{x^\m}^j. 
		\end{array}
	\end{equation*}
	Then $g_j$ is left  $R$-linear and
	\begin{equation*}
		\dsum_{j=1}^mg_j(u_x)\overline{\om}_x^j=\dsum_{j=1}^m u_x\overline{\om}_{x^\m}^j\overline{\om}_x^j=u_x1_{x^\m}=u_x,
	\end{equation*}
	for all  $u_x \in \T_x$. Therefore, $\{ \overline{\om}_x^j,g_j\}$ is a dual basis for the left $R$-module $\T_x .$ 
	
	Let $\varphi:\T_x\longrightarrow \T_x$ be a right $R$-linear map.  Then,
	$\widetilde{\varphi}=\dsum_{i=1}^n\varphi(\om_x^i)\om_{x^\m}^i\in R1_{x}$ and
	\begin{equation*}
		\widetilde{\varphi}u_x=\dsum_{i=1}^n\varphi(\om_x^i)\om_{x^\m}^iu_x=\dsum_{i=1}^n\varphi(\om_x^i\om_{x^\m}^iu_x)=\varphi(1_xu_x)=\varphi(u_x),
	\end{equation*}
	for all  $u_x \in \T_x$. Then the map  $R\longrightarrow \End({\T_x}_R),$ defined by  $r\longmapsto(u_x \mapsto ru_x),$ is surjective.  
	
	Similarly, if $\psi:\T_x\longrightarrow \T_x$ is a left  $R$-linear map, then $\widetilde{\psi}=\dsum_{j=1}^m\overline{\om}_{x^\m}^j\psi(\overline{\om}_x^j)\in R1_{x^\m}$ and  
	\begin{equation*}
		u_x\widetilde{\psi}=\dsum_{j=1}^mu_x\overline{\om}_{x^\m}^j\psi(\overline{\om}_x^j)=\dsum_{j=1}^m\psi(u_x\overline{\om}_{x^\m}^j\overline{\om}_x^j)=\psi(u_x1_{x^\m})=\psi(u_{x}),
	\end{equation*}
	for all  $u_x \in \T_x$. Hence, the map  $R\longrightarrow \End(_R\T_x),$ defined by  $r \longmapsto (u_x\mapsto u_xr),$ is surjective.   Thus, $[\T_x]\in \Pics(R)$. 
	
	{(ii)} 	Clearly,  $m_l$ is a well-defined  $R$-$S$-bilinear map. Its inverse is given by 
	$$\begin{array}{c c c l}
		m_l^\m: & \T_xM & \longrightarrow & \T_x\ot_R M,\\
		& m & \longmapsto & \dsum_{(x)}\om_x\ot \om_{x^\m}m,
	\end{array}$$
	where $1_x=\dsum_{(x)}\om_x\om_{x^\m}$, with $\om_x \in \T_x$ and $\om_{x^\m}\in \T_{x^\m}$.
	Indeed, let us check that  $m_l^\m$ is well-defined. Let $1_x=\dsum_{\widetilde{x}}\widetilde{\om}_x\widetilde{\om}_{x^\m}$ be another decomposition of  $1_x$. Then,
	\begin{eqnarray*}
		m_l^\m(m) & = & \dsum_{(x)}\om_x\ot \om_{x^\m}m = \dsum_{(x)}\om_x\ot \om_{x^\m}1_xm =  \dsum_{(x),(\widetilde{x})}\om_x\ot \om_{x^\m}\widetilde{\om}_x\widetilde{\om}_{x^\m}m\\
		& = & \dsum_{(x),(\widetilde{x})}\om_x\om_{x^\m}\widetilde{\om}_x\ot \widetilde{\om}_{x^\m}m=\dsum_{(\widetilde{x})}1_x\widetilde{\om}_x\ot \widetilde{\om}_{x^\m}m= \dsum_{(\widetilde{x})}\widetilde{\om}_x\ot \widetilde{\om}_{x^\m}m.
	\end{eqnarray*}
	Thus, given  $m \in \T_xM$, we have:
	\begin{equation*}
		(m_l\circ m_l^\m)(m)=\dsum_{(x)}\om_x\om_{x^\m}m=1_xm=m.
	\end{equation*}
	On the other hand, given  $u_x \in \T_x$ and $m \in M$, we obtain: 
	\begin{eqnarray*}
		(m_l^\m\circ m_l)(u_x\ot m) & = & \dsum_{(x)}\om_x\ot \om_{x^\m}u_xm= \dsum_{(x)}\om_x\om_{x^\m}u_x\ot m\\
		& = & 1_xu_x\ot m= u_x\ot m.
	\end{eqnarray*}
	Consequently, $m_l$ is an $R$-$S$-bimodule isomorphism.
	Similarly, the inverse of $m_r$ is given by 
	$$\begin{array}{c c c l}
		m_r^\m: & M\T_x & \longrightarrow & M\ot_R\T_{x},\\
		& m & \longmapsto & \dsum_{(x)} m\overline{\om}_{x^\m}\ot \overline{\om}_{x},
	\end{array}$$
	where $1_{x^\m}=\dsum_{(\overline{x})}\overline{\om}_{x^\m}\overline{\om}_x$.

	{(iii)} Note that $\T_xN\subseteq M$ is an  $R$-subbimodule of $M$, such that { $1_xn=n$ and $u_{x^\m}n\in N$, for all  $n \in \T_xN$ and $u_{x^\m}\in \T_{x^\m}$. }
	Indeed,  if $n \in \T_xN$, then  $n=\dsum_{i=1}^nu_x^in_i$, with $u_x^i \in \T_x$ and $n_i \in N$. Since $N\subseteq M$ is an $R$-subbimodule, we have: 
	\begin{equation*}
		u_{x^\m}n=\dsum_{i=1}^nu_{x^\m}(u_x^in_i)=\dsum_{i=1}^n(\underbrace{u_{x^\m}u_x^i}_{\in R1_{x^\m}})n_i \in N.
	\end{equation*}
	
	Analogously, $N\T_{x}$ is an  $R$-subbimodule of $M$ which  satisfies {$n'1_{x^\m}=n$ and $n'u_{x^\m}\in N$, for all  $n'\in N\T_{x}$ and $u_{x^\m}\in \T_{x^\m}$.}
	Hence, the isomorphism follows as in   item (ii).
\end{dem}

\begin{cor}\label{TxotTyisoTxTy} Let $N \in \mathcal{S}_R(S).$ Then $\T_x\ot_RN\simeq \T_xN$, for each  $x \in G$. In particular, $\T_x\ot_R\T_y\simeq\T_{x}\T_y$, for all  $x,y \in G$.
\end{cor}

\begin{obs}\label{pcgpviaSRS}
	Let 
	\begin{equation*}
		\begin{array}{c c c l}
			\T: & G & \longrightarrow & \mathcal{S}_R(S)\\
			& x& \longmapsto & \T_x
		\end{array}
	\end{equation*}
	be a unital partial representations with  $\T_x\T_{x^\m}=R1_x$. By Proposition~\ref{SRemPicS} and  Corolary~\ref{TxotTyisoTxTy} the map 
	\begin{equation*}
		\begin{array}{c c c l}
			\underline{\T}: & G & \longrightarrow & \Pics(R),\\
			& x & \longmapsto & [\T_x]
		\end{array}
	\end{equation*}
	is a unital partial representation with  ${ \T_x  \ot _R  \T_{x^\m}   }   \simeq R1_x,$   ($x \in G)$. Let $f^\T$ be the family of $R$-bimodule isomorpshims whose members  $$\begin{array}{c c c l}
		f^\T_{x,y}:& \T_x\ot_R \T_y& \longrightarrow & \T_x\T_y=1_x\T_{xy},\\
		& u_x \ot u_y & \longmapsto & u_xu_y
	\end{array}$$ 
	are  induced by the multiplication in  $S$. Then $f^\T$ is factor set for  $\underline{\T}$. Therefore, we have a generalized partial crossed product  $\D(\T),$ where each  $\T_x\subseteq S$ is an  $R$-subbimodule.
	
	Conversely, if   $\T:G\longrightarrow \Pics(R)$ is a unital partial representation with a factor set   $f^\T$, then the generalized  partial crossed product  $\D(\T)$ is an extension of $R$ with the same unity. In this case, each   $\T_x\subseteq \D(\T)$ is an  $R$-subbimodule. By Corollary~\ref{TxotTyisoTxTy}, the map   $\underline{\T}:G\longrightarrow \mathcal{S}_R(\D(\T))$ with $\underline{\T}(x)=\T_x$, is a unital partial representation. 
\end{obs}

\begin{exe}(Partial crossed product) Let $\al=(\{D_x\}_{x\in G}, \{\al_x\}_{x\in G}, \{\om_{x,y}\}_{x,y \in G})$ be a unital  twisted  partial action of $G$ on $R$ with $D_x=R1_x$, for each  $x\in G$   (see \cite[Definition 2.1]{dokuchaev2008crossed}).  Consider the partial crossed product  $R\rtimes_{\al,\om}G=\bigoplus_{x\in G}D_x\delta_x,$ whose multiplication  is defined by   
	$$(u_x\delta_x)(u_y\delta_y)=u_x\al_x(u_y1_{x^\m})\om_{x,y}\delta_{xy}.$$
	{ By \cite[Theorem 2.4]{dokuchaev2008crossed},} $R\rtimes_{\al,\om}G$ is an associative unital ring.   Let
	\begin{equation*}
		\begin{array}{c c c l}
			\T: & G & \longrightarrow & \mathcal{S}_R(R\rtimes_{\al,\om}G)\\
			& x & \longmapsto & D_x\delta_x 
		\end{array}.
	\end{equation*} 
	It is easy to see that  $(D_x\delta_x)(D_{x^\m}\delta_{x^\m})(D_x\delta_x)\subseteq D_x\delta_x$, for each $x\in G$. 
	On the other hand, if  $u_x\delta_x \in D_x\delta_x ,$ then 
	$$(u_x\delta_x)(1_{x^\m}\delta_{x^\m})(\om_{x,x^\m}^\m\delta_x)=(u_x\om_{x,x^\m}\delta_1)(\om_{x,x^\m}^\m\delta_x)=u_x\om_{x,x^\m}\om_{x,x^\m}^\m\delta_x=u_x1_x\delta_x=u_x\delta_x.$$
	Hence,  $(D_x\delta_x)(D_{x^\m}\delta_{x^\m})(D_x\delta_x)= D_x\delta_x$, for each  $x\in G$. By 
	\cite[Lemma 5.3]{dokuchaev2008crossed}, we have that 
	$$(D_x\delta_x)(D_y\delta_y)(D_{y^\m}\delta_{y^\m})=(D_{xy}\delta_{xy})(D_{y^\m}\delta_{y^\m}) \ \ \mbox{and} \ \ (D_x^\m\delta_x)(D_x\delta_x)(D_{y}\delta_{y})=(D_{x^\m}\delta_{x^\m})(D_{xy}\delta_{xy}),$$
	for all  $x,y \in G$. Consequently, $\T$ is a partial representation. Observe now that  $(D_x\delta_x)(D_{x^\m}\delta_{x^\m})=D_x$, for each  $x\in G$. Indeed,  the inclusion   $(D_x\delta_x)(D_{x^\m}\delta_{x^\m})\subseteq D_x$ is immediate.  Given $u_x \in D_x$, we have
	$(u_x\om_{x,x^\m}^\m\delta_x)(1_{x^\m}\delta_x)=u_x\om_{x,x^\m}^\m\om_{x,x^\m}\delta_1=u_x1_x\delta_1=u_x\delta_1\in D_x$.
	Therefore, $\T$ is a unital partial representation.  By Remark~\ref{pcgpviaSRS}, we obtain that 
	$$\begin{array}{c c c l}
		f_{x,y}^\T: & (D_x\delta_x)\ot (D_y\delta_y) & \longrightarrow & D_xD_{xy}\delta_{xy}\\
		& u_x\delta_x\ot u_y\delta_y & \longmapsto & u_x\al_x(u_y1_{x^\m})\om_{x,y}\delta_{xy}
	\end{array},$$
	is a factor set for  $\T$. Hence, $\D(\T)=\bigoplus_{x\in G}D_x\delta_x=R\rtimes_{\al,\om}G$ is a generalized partial crossed product.
\end{exe}


\section{The group $\C(\T/R)$ and partial cohomology}\label{sec:groupC}

 {{\it In all what follows the unadorned $\otimes $ will stand for $\otimes _R,$ unless otherwise stated.}} 
 
  Let 
\begin{equation*}
\begin{array}{c c c c}
\T: & G & \longrightarrow & \Pics(R)\\
& x & \longmapsto & [\T_x]
\end{array}
\end{equation*}
be a unital partial representation with  $\varepsilon_x=\T_x\ot\T_{x^\m}\simeq R1_x$, for all $x\in G$ and let $f^\T=\{ f_{x,y}^\T:\T_x\ot\T_y\longrightarrow 1_x\T_{xy},\ x,y \in G  \}$ be a factor set for $\T$ and $\D(\T)$  { the} partial generalized crossed product with factor set $f^\T=\{f_{x,y}^\T:\T_x\ot\T_y\rightarrow 1_x\T_{xy}, x,y\in G\}$.

\begin{teo}\label{grupoC} The set
	\begin{equation*}
	\C(\T/R)=\{[\D(\Gamma)]; \ \Gamma_x|\T_x  \ \mbox{and} \ \Gamma_x\ot\Gamma_{x^\m}\simeq R1_x, \ \mbox{for all }x \in G\}
	\end{equation*}
	is an abelian group with multiplication defined by
	\begin{equation*}
	[\D(\Omega)][\D(\Gamma)]=\left[\bigoplus_{x \in G}\Omega_x\ot\T_{x^{-1}}\ot\Gamma_x \right],
	\end{equation*}
	where the factor set 
	$$f^{\Om\Gm}=\{f_{x,y}^{\Omega\Gamma}:\Gamma_x\ot \T_{x^\m}\ot\Omega_x\ot \Gamma_y\ot \T_{y^\m}\ot \Omega_y\longrightarrow 1_x\Gamma_{xy}\ot \T_{(xy)^\m}\ot \Omega_{xy}\}$$
consists of  isomorphisms defined in Lemma \ref{conjuntodefatores}.
\end{teo}
\begin{dem} 
	By Lemma \ref{conjuntodefatores} we have $$[\D(\Omega)][\D(\Gamma)]=\left[\bigoplus_{x \in G}\Omega_x\ot\T_{x^{-1}}\ot\Gamma_x \right] \in \C(\T/R).$$

First, we need to show that this operation is well defined, i.e. it  does not depend on the choice of the representative of the  class. Let $[\D(\Gm)]=[\D(\Sigma)]$ and $[\D(\Om)]=[\D(\Lambda)]$ in $\C(\T/R)$, then there is $R$-bimodule isomorphisms  $a_x:\Gm_x\longrightarrow \Sigma_x$ and $b_x:\Om_x\longrightarrow \Lambda_x$ and commutative diagrams
\begin{center}
	\begin{tabular}{c c c}
		$$
		\xymatrix{ \Gm_x\ot \Gm_y\ar[dd]_{a_x\ot a_y}\ar[rr]^{f_{x,y}^\Gm} & &
			1_x\Gm_{xy}\ar[dd]^{a_{xy}}\\
			& & \\
			\Sigma_x\ot \Sigma_y\ar[rr]_{f_{x,y}^\Sigma}& & 1_x\Sigma_{xy}  }  
		$$ &   \ \ \ \ \ \ &$$
		\xymatrix{ \Om_x\ot \Om_y\ar[dd]_{b_x\ot b_y}\ar[rr]^{f_{x,y}^\Om} & & 1_x\Om_{xy}\ar[dd]^{b_{xy}}\\
			& & \\
			\Lambda_x\ot \Lambda_y\ar[rr]_{f_{x,y}^\Lambda}& & 1_x\Lambda_{xy}  }  
		$$
	\end{tabular}
\end{center}
Thus, we have
\begin{equation}
	a_{xy}(v_x\stackrel{\Gm}{\circ}v_y)=(a_x(v_x)\stackrel{\Sigma}{\circ}a_y(v_y))\ \ \mbox{and} \ \ b_{xy}(\om_x\stackrel{\Om}{\circ}\om_y)=(b_x(\om_x)\stackrel{\Lambda}{\circ}b_y(\om_y)), \label{comutatividadedosdiagramasdaigualdade}
\end{equation}
for all $v_x \in \Gm_x,v_y \in \Gm_y,\om_x\in \Om_x $ and $\om_y \in \Om_y$.
By definition 
\begin{equation*}
	[\D(\Gm)][\D(\Om)]
	=\left[ \bigoplus_{x\in G}\Gm_x\ot \T_{x^\m}\ot \Om_x\right]\  \ \mbox{and} \ \  [\D(\Sigma)][\D(\Lambda)]
	=\left[ \bigoplus_{x\in G}\Sigma_x\ot \T_{x^\m}\ot \Lambda_x\right] \end{equation*}
Define the   $R$-bimodule isomorphism:
$$\begin{array}{c c c l}
	d_{xy}: & \Gm_x\ot \T_{x^\m}\ot \Om_x & \longrightarrow & \Sigma_x\ot \T_{x^\m}\ot \Lambda_x,\\
	& v_x\ot u_{x^\m}\ot \om_x & \longmapsto & a_x(v_x)\ot u_{x^\m}\ot b_x(\om_x).
\end{array}$$
We shall see that the following diagram is commutative:
\begin{equation*}
	\xymatrix@C=1.5cm{ \Gm_x\ot \T_{x^\m}\ot \Om_{x}\ot \Gm_y\ot \T_{y^\m}\ot \Om_{y}\ar[dd]_{d_x\ot d_y}\ar[rr]^{f_{x,y}^{\Gm\Om}} & & 1_x\Gm_{xy}\ot \T_{(xy)^\m}\ot \Om_{xy}\ar[dd]^{d_{xy}} \\
		& 	& \\
		\Sigma_x\ot \T_{x^\m}\ot \Lambda_{x}\ot\Sigma_y\ot \T_{y^\m}\ot \Lambda_{y}\ar[rr]_{f_{x,y}^{\Sigma\Lambda}} & & 1_x\Sigma_{xy}\ot \T_{(xy)^\m}\ot \Lambda_{xy}	  }
\end{equation*}
On  one hand, $$d_{xy}\circ f_{x,y}^{\Gm\Om}=d_{xy}\circ(f_{x,y}^{\Gm}\ot f_{y^\m,x^\m}^\T\ot f_{x,y}^\Om)\circ (\Gm_x\ot T_1\ot \Om_y),$$ where
$T_1=T_{\T_{x^\m}\ot \Om_x,\Gm_y\ot \T_{y^\m}}$.
On the other hand, 
$$f_{x,y}^{\Sigma\Lambda}\circ(d_x\ot d_y)=(f_{x,y}^\Sigma\ot f_{y^\m,x^\m}^\T\ot f_{x,y}^\Lambda)\circ (\Sigma_x\ot T_2\ot \Lambda_y)\circ (d_x\ot d_y),$$
where $T_2=T_{\T_{x^\m}\ot \Lambda_x,\Sigma_y\ot \T_{y^\m}}$. We use   Proposition ~\ref{isomorfismoT} to construct the isomorphisms $T_1$ and $T_2$. To this end, consider the $R$-bilinear maps $f_i:\T_{x^\m}\ot \Om_x\longrightarrow R$  and $g_i: R\longrightarrow \T_{x^\m}\ot \Om_x$,  $i=1,2,...,n$, such that $\dsum_{i=1}^ng_i\circ f_i=Id_{\T_{x^\m}\ot \Om_x}$. We define $\overline{f}_i:\T_{x^\m}\ot \Lambda_x\rightarrow R$ and $\overline{g}_i:R\rightarrow \T_{x^\m}\ot\Lambda_x$ by $\overline{f}_i=f_i\circ (\T_{x^\m}\ot b_x^{-1})$ and $\overline{g}_i:(\T_{x^\m}\ot b_x)\circ g_i$, $i=1,2,...,n$. Then, $\overline{f}_i$ and $\overline{g}_i$ are $R$-bilinear and we have $\dsum_{i=1}^n\overline{g}_i\overline{f}_i=Id_{\T_{x^\m}\ot\Lambda_x}$.
Denote $g_i(1)=\dsum_{l}u_{x^\m}^{i,l}\ot \om_x^{i,l},$ then $\overline{g}_i(1)=\dsum_lu_{x^\m}^{i,l}\ot b_x(\om_{x}^{i,l})$. Thus,
\begin{equation*}
	\begin{array}{c c c l}
		T_1: & \T_{x^\m}\ot \Om_x\ot\Gm_y\ot \T_{y^\m}& \longrightarrow & \Gm_y\ot \T_{y^\m}\ot\T_{x^\m}\ot \Om_x,\\
		& u_{x^\m}\ot \om_x\ot v_y\ot u_{y^\m} & \longmapsto & \dsum_{i,l}f_i(u_{x^\m}\ot \om_x)v_y\ot u_{y^\m}\ot u_{x^\m}^{i,l}\ot \om_{x}^{i,l},
	\end{array}
\end{equation*}
\begin{equation*}
	\begin{array}{c c c l}
		T_2:& \T_{x^\m}\ot \Lambda_x\ot \Sigma_y\ot \T_{y^\m}&\longrightarrow & \Sigma_y\ot \T_{y^\m}\ot\T_{x^\m}\ot \Lambda_x,\\
		& u_{x^\m}\ot l_x\ot t_y\ot u_{y^\m} & \longmapsto & \dsum_{i,l}f_i(u_{x^\m}\ot b_x^\m(l_x))t_y\ot u_{y^\m}\ot u_{x^\m}^{i,l}\ot b_x(\om_x^{i,l}).
	\end{array}
\end{equation*}
Therefore,
\begin{eqnarray*}
	&  & 	f_{x,y}^{\Sigma\Lambda}\circ(d_x\ot d_y)(v_x\ot u_{x^\m}\ot \om_x\ot v_y\ot u_{y^\m}\ot \om_y)\\
	& = & (f_{x,y}^\Sigma\ot f_{y^\m,x^\m}^\T\ot f_{x,y}^\Lambda)\circ (\Sigma_x\ot T_2\ot \Lambda_y)\circ (d_x\ot d_y)(v_x\ot u_{x^\m}\ot \om_x\ot v_y\ot u_{y^\m}\ot \om_y)\\
	& = & (f_{x,y}^\Sigma\ot f_{y^\m,x^\m}^\T\ot f_{x,y}^\Lambda)\circ (\Sigma_x\ot T_2\ot \Lambda_y)(a_x(v_x)\ot u_{x^\m}\ot b_x(\om_x)\ot a_y(v_y)\ot u_{y^\m}\ot b_y(\om_y))\\
	& = & \dsum_{i,l} (f_{x,y}^\Sigma\ot f_{y^\m,x^\m}^\T\ot f_{x,y}^\Lambda)( a_x(v_x)\ot f_i(u_{x^\m}\ot \om_x)a_y(v_y)\ot u_{y^\m}\ot u_{x^\m}^{i,l}\ot b_x(\om_x^{i,l})\ot b_y(\om_y))\\
	& = & \dsum_{i,l}(a_x(v_x)\stackrel{\Sigma}{\circ} f_i(u_{x^\m}\ot \om_x)a_y(v_y))\ot (u_{y^\m}\stackrel{\T}{\circ} u_{x^\m}^{i,l})\ot (b_x(\om_x^{i,l})\stackrel{\Lambda}{\circ} b_y(\om_y)).
\end{eqnarray*}
On  the other hand, 
\begin{eqnarray*}
	& & d_{xy}\circ f_{x,y}^{\Gm\Om}( v_x\ot u_{x^\m}\ot \om_x\ot v_y\ot u_{y^\m}\ot \om_y)\\
	& = & d_{xy}\circ(f_{x,y}^{\Gm}\ot f_{y^\m,x^\m}^\T\ot f_{x,y}^\Om)\circ (\Gm_x\ot T_1\ot \Om_y)(v_x\ot u_{x^\m}\ot \om_x\ot v_y\ot u_{y^\m}\ot \om_y)\\
	& = & \dsum_{i,l}d_{xy}\circ(f_{x,y}^{\Gm}\ot f_{y^\m,x^\m}^\T\ot f_{x,y}^\Om)( v_x\ot f_i(u_{x^\m}\ot \om_x)v_y\ot u_{y^\m}\ot u_{x^\m}^{i,l}\ot \om_x^{i,l}\ot \om_y)\\
	& = & \dsum_{i,l}d_{xy}((v_x\stackrel{\Gm}{\circ} f_i(u_{x^\m}\ot \om_x)v_y)\ot (u_{y^\m}\stackrel{\T}{\circ}u_{x^\m}^{i,l})\ot (\om_x^{i,l}\stackrel{\Om}{\circ} \om_y))\\
	& =  &\dsum_{i,l} a_{xy}(v_x\stackrel{\Gm}{\circ}f_i(u_{x^\m}\ot \om_x)v_y)\ot (u_{y^\m}\stackrel{\T}{\circ} u_{x^\m}^{i,l})\ot b_{xy}(\om_x^{i,l}\stackrel{\Om}{\circ} \om_y)\\
	& = & \dsum_{i,l}(a_x(v_x)\stackrel{\Sigma}{\circ} f_i(u_{x^\m}\ot \om_x)a_y(v_y))\ot (u_{y^\m}\stackrel{\T}{\circ} u_{x^\m}^{i,l})\ot (b_x(\om_x^{i,l})\stackrel{\Lambda}{\circ} b_y(\om_y)), 
\end{eqnarray*}
where the last equality follows from (\ref{comutatividadedosdiagramasdaigualdade}).
Consequently, $d_{xy}\circ f_{x,y}^{\Gm\Om}=f_{x,y}^{\Sigma\Lambda}\circ (d_x\ot d_y)$, for all $x,y \in G$. This implies that $[\D(\Gm)][\D(\Om)]=[\D(\Sigma)][\D(\Lambda)]$ in $\C(\T/R)$
and the assertion follows. 
	
	Let us show that with this operation $\C(\T/R)$ is, in fact, a group.

\underline{Associativity:} Take $[\D(\Gamma)],[\Delta(\Omega)],[\D(\Sigma)] \in \C(\T/R)$. Then,
\begin{equation*}
	([\D(\Gamma)][\D(\Omega)])[\D(\Sigma)]=\left[\bigoplus_{x \in G}\Gamma_x\ot \T_{x^\m}\ot \Om_x\ot \T_{x^\m} \ot \Sigma_{x} \right] 
\end{equation*}
\begin{equation*}
	[\D(\Gamma)]([\D(\Omega)][\D(\Sigma)])=\left[\bigoplus_{x \in G}\Gamma_x\ot \T_{x^\m}\ot \Om_x\ot \T_{x^\m} \ot \Sigma_{x} \right] 
\end{equation*}
We will check that  $f_{x,y}^{(\Gamma\Om)\Sigma}=f_{x,y}^{\Gamma(\Om\Sigma)}$, for all $x,y \in G$. Consider the $R$-bimodule isomorphisms
$$T_1=T_{\T_{x^\m}\ot\Sigma_x,\Gm_y\ot \T_{y^\m}\ot\Om_y\ot \T_{y^\m}}, T_2=T_{\T_{x^\m}\ot \Om_x,\Gm_y \ot\T_{y^\m}}, T_3=T_{\T_{x^\m}\ot \Om_x\ot \T_{x^\m}\ot \Sigma_x,\Gm_y\ot \T_{y^\m}}$$ and  $T_4=T_{\T_{x^\m}\ot \Sigma_x,\Om_y\ot \T_{y^\m}}.$ On  one hand, $f_{x,y}^{(\Gamma\Om)\Sigma}$ is given by 
\begin{equation}
	{\small \xymatrix@C=-5cm{  \Gamma_x\ot \T_{x^\m}\ot \Om_x\ot \T_{x^\m} \ot \Sigma_{x}\ot \Gamma_y\ot \T_{y^\m}\ot \Om_y\ot \T_{y^\m} \ot \Sigma_{y} \ar[rdd]^-{\Gm_x\ot \T_{x^\m}\ot \Om_x\ot T_1\ot \Sigma_y}\ar@/_4.2cm/@{-->}[rdddddd]_*-{f_{x,y}^{(\Gamma\Om)\Sigma}} & 
			\\
			& \\
			&  \Gamma_x\ot \T_{x^\m}\ot \Om_x\ot \Gamma_y\ot \T_{y^\m}\ot \Om_y\ot \T_{y^\m} \ot \T_{x^\m} \ot \Sigma_{x}\ot \Sigma_{y}\ar[dd]|-{\Gm_x\ot T_2\ot  \Om_y\ot \T_{y^\m} \ot \T_{x^\m} \ot \Sigma_{x}\ot \Sigma_{y}}\\
			& \\
			&  \Gamma_x\ot \Gamma_y\ot \T_{y^\m}\ot \T_{x^\m}\ot \Om_x\ot \Om_y\ot \T_{y^\m} \ot \T_{x^\m} \ot \Sigma_{x}\ot \Sigma_{y}\ar[dd]|{f_{x,y}^{\Gamma}\ot f_{y^\m,x^\m}^\T\ot f_{x,y}^\Om\ot f_{y^\m,x^\m}^\T\ot f_{x,y}^\Sigma}\\
			& \\
			& 1_x\Gamma_{xy}\ot \T_{(xy)^\m}\ot \Om_{xy}\ot \T_{(xy)^\m} \ot \Sigma_{xy}}} \label{associatividade1}
\end{equation}
On the other hand, $f_{x,y}^{\Gamma(\Om\Sigma)}$ is given by
\begin{equation}
	{\small 	\xymatrix@C=-5cm{  \Gamma_x\ot \T_{x^\m}\ot \Om_x\ot \T_{x^\m} \ot \Sigma_{x}\ot \Gamma_y\ot \T_{y^\m}\ot \Om_y\ot \T_{y^\m} \ot \Sigma_{y} \ar[rdd]^{\Gm_x\ot T_3\ot \Om_y\ot \T_{y^\m}\ot \Sigma_y}\ar@/_4.2cm/@{-->}[rdddddd]_{f_{x,y}^{\Gamma(\Om\Sigma)}} &  \\ 
			& \\
			& \Gamma_x\ot \Gamma_y\ot \T_{y^\m}\ot \T_{x^\m}\ot \Om_x\ot \T_{x^\m} \ot \Sigma_{x}\ot \Om_y\ot \T_{y^\m} \ot \Sigma_{y} \ar[dd]|-{\Gamma_x\ot \Gamma_y\ot \T_{y^\m}\ot \T_{x^\m}\ot \Om_x\ot T_4\ot \Sigma_y}  \\
			& \\
			& \Gamma_x\ot \Gamma_y\ot \T_{y^\m}\ot \T_{x^\m}\ot \Om_x\ot \Om_y\ot \T_{y^\m}\ot \T_{x^\m} \ot \Sigma_{x} \ot \Sigma_{y} \ar[dd]|{f_{x,y}^{\Gamma}\ot f_{y^\m,x^\m}^\T\ot f_{x,y}^\Om\ot f_{y^\m,x^\m}^\T\ot f_{x,y}^\Sigma}  \\
			& \\
			& 1_x\Gamma_{xy}\ot \T_{(xy)^\m}\ot \Om_{xy}\ot \T_{(xy)^\m} \ot \Sigma_{xy}}} \label{associatividade2}
\end{equation}
We will construct the isomorphisms $T_1,T_2,T_3$ and $T_4.$ Let $f_i:\T_{x^{-1}}\ot \Sigma_x\longrightarrow R$ and $g_i:R\longrightarrow\T_{x^{-1}}\ot \Sigma_x$, $i=1,2,...,n$, be $R$-bilinear maps satisfying  $\dsum_{i=1}^ng_if_i=Id_{\T_{x^{-1}}\ot \Sigma_x}$. Denote $g_i(1)=\dsum_{l}u_{x^\m}^{i,l}\ot t_x^{i,l}.$ Then
$${\small \begin{array}{c c c l}
	T_1:& \T_{x^{-1}}\ot \Sigma_x\ot \Gamma_y\ot \T_{y^{-1}}\ot \Om_y\ot\T_{y^{-1}}& \longrightarrow& \Gamma_y\ot \T_{y^{-1}}\ot \Om_y\ot\T_{y^{-1}}\ot \T_{x^{-1}}\ot \Sigma_x\\
	& u_{x^\m}\ot t_x\ot v_y\ot u_{y^\m}\ot \om_y\ot u_{y^\m}' & \longmapsto & \dsum_{i,l}f_i(u_{x^\m}\ot t_x)v_y\ot u_{y^\m}\ot \om_y\ot u_{y^\m}'\ot u_{x^\m}^{i,l}\ot t_x^{i,l},
\end{array}}$$
and 
$$\begin{array}{c c c l}
	T_4:& \T_{x^{-1}}\ot \Sigma_x\ot\Om_y\ot \T_{y^{-1}}& \longrightarrow & \Om_y\ot \T_{y^{-1}}\ot\T_{x^{-1}}\ot \Sigma_x\\
	& u_{x^\m}\ot t_x\ot \om_y\ot u_{y^\m} & \longmapsto & \dsum_{i,l}f_i(u_{x^\m}\ot t_x)\om_y\ot u_{y^\m}\ot u_{x^\m}^{i,l}\ot t_x^{i,l}. 
\end{array}$$
Consider now  $R$-bilinear maps $f_j':\Gamma_y\ot \T_{y^\m}\longrightarrow R $ and $g_j':R\longrightarrow \Gamma_y\ot \T_{y^\m}$,  $j=1,2,...,m$, such that $\dsum_{j=1}^mg_j'f_j'=Id_{\Gamma_y\ot \T_{y^\m}}$. Write $g_j(1)=\dsum_{k}\widetilde{v}_y^{j,k}\ot \widetilde{u}_{y^\m}^{j,k}.$ Then
$$\begin{array}{c c c l}
	T_2:&  \T_{x^\m}\ot \Om_x\ot \Gamma_y\ot \T_{y^{-1}}& \longrightarrow& \Gamma_y\ot \T_{y^{-1}}\ot\T_{x^\m}\ot \Om_x,\\
	& u_{x^\m}\ot \om_x\ot v_y\ot u_{y^\m} & \longmapsto & \dsum_{j,k} \widetilde{v}_y^{j,k}\ot \widetilde{u}_{y^\m}^{j,k} \ot u_{x^\m}\ot \om_xf_j'(v_y\ot u_{y^\m}),
\end{array}$$
and 
$${\small \begin{array}{c c c l}
	T_3:& \T_{x^\m}\ot \Om_x\ot \T_{x^\m} \ot \Sigma_{x}\ot \Gamma_y\ot \T_{y^{-1}}&\longrightarrow & \Gamma_y\ot \T_{y^{-1}}\ot\T_{x^\m}\ot \Om_x\ot \T_{x^\m} \ot \Sigma_{x},\\
	& u_{x^\m}\ot \om_x\ot u'_{x^\m}\ot t_x\ot v_y\ot u_{y^\m} & \longmapsto & \dsum_{j,k}\widetilde{v}_y^{j,k}\ot \widetilde{u}_{y^\m}^{j,k}\ot u_{x^\m}\ot \om_x\ot u'_{x^\m}\ot t_x f_j'(v_y\ot u_{y^\m}).
\end{array}}$$

Thus, by diagram (\ref{associatividade1}) and using the above notation, we have that $f_{x,y}^{(\Gamma\Om)\Sigma}$ is given by:
\begin{eqnarray}
	&  &v_x\ot u_{x^\m}\ot \om_x\ot u'_{x^\m}\ot t_x\ot v_y\ot u_{y^\m}\ot \om_y\ot u_{y^\m}'\ot t_y\nonumber\\
	& \stackrel{T_1}{\mapsto} &\dsum_{i,l} v_x\ot u_{x^\m}\ot \om_x\ot f_i(u'_{x^\m}\ot t_x)v_y\ot u_{y^\m}\ot \om_y\ot u_{y^\m}' \ot u_{x^\m}^{i,l}\ot t_x^{i,l}\ot  t_y\nonumber\\
	& \stackrel{T_2}{\mapsto} &\dsum_{i,j,l,k} v_x\ot\widetilde{v}_y^{j,k}\ot \widetilde{u}_{y^\m}^{j,k}\ot  u_{x^\m}\ot \om_xf_j'( f_i(u'_{x^\m}\ot t_x)v_y\ot u_{y^\m}) \ot \om_y\ot u_{y^\m}' \ot u_{x^\m}^{i,l}\ot t_x^{i,l}\ot  t_y\nonumber\\
	& \mapsto &\dsum_{i,j,l,k} (v_x\stackrel{\Gamma}{\circ}\widetilde{v}_y^{j,k})\ot (\widetilde{u}_{y^\m}^{j,k}\stackrel{\T}{\circ}  u_{x^\m})\ot  (\om_xf_j'( f_i(u'_{x^\m}\ot t_x)v_y\ot u_{y^\m}) \stackrel{\Om}{\circ} \om_y)\ot (u_{y^\m}' \stackrel{\T}{\circ}u_{x^\m}^{i,l})\ot (t_x^{i,l}\stackrel{\Sigma}{\circ}  t_y)\nonumber.
\end{eqnarray}

As $f_i$  and $f_j'$ are $R$-bilinear, then $f_i(u_x'\ot t_xf_j'(v_y \ot u_{y^\m}))=f_i(u_x'\ot t_x)f_j'(v_y \ot u_{y^\m})=f_j'(f_i(u_x'\ot t_x)v_y \ot u_{y^\m})$.  Hence, by diagram (\ref{associatividade2}) we see that  $f_{x,y}^{\Gamma(\Om\Sigma)}$  is given by:  
\begin{eqnarray}
	&  &v_x\ot u_{x^\m}\ot \om_x\ot u'_{x^\m}\ot t_x\ot v_y\ot u_{y^\m}\ot \om_y\ot u_{y^\m}'\ot t_y\nonumber\\
	& \stackrel{T_3}{\mapsto} & \dsum_{j,k} v_x\ot \widetilde{v}_y^{j,k}\ot \widetilde{u}_{y^\m}^{j,k} \ot u_{x^\m}\ot \om_x\ot u'_{x^\m}\ot t_xf_j'(v_y\ot u_{y^\m})\ot \om_y\ot u_{y^\m}'\ot t_y\nonumber\\
	& \stackrel{T_4}{\mapsto} & \dsum_{i,j,k,l} v_x\ot \widetilde{v}_y^{j,k}\ot \widetilde{u}_{y^\m}^{j,k} \ot u_{x^\m}\ot \om_x\ot f_i(u'_{x^\m}\ot t_xf_j'(v_y\ot u_{y^\m}))\om_y\ot u_{y^\m}'\ot u_{x^\m}^{i,l}\ot t_x^{i,l}\ot  t_y\nonumber\\
	& = &  \dsum_{i,j,k,l} v_x\ot \widetilde{v}_y^{j,k}\ot \widetilde{u}_{y^\m}^{j,k} \ot u_{x^\m}\ot \om_x\ot f_j'(f_i(u'_{x^\m}\ot t_x)v_y\ot u_{y^\m})\om_y\ot u_{y^\m}'\ot u_{x^\m}^{i,l}\ot t_x^{i,l}\ot  t_y\nonumber\\
	& \mapsto &\dsum_{i,j,l,k} (v_x\stackrel{\Gamma}{\circ}\widetilde{v}_y^{j,k})\ot (\widetilde{u}_{y^\m}^{j,k}\stackrel{\T}{\circ}  u_{x^\m})\ot  (\om_xf_j'( f_i(u'_{x^\m}\ot t_x)v_y\ot u_{y^\m}) \stackrel{\Om}{\circ} \om_y)\ot (u_{y^\m}' \stackrel{\T}{\circ}u_{x^\m}^{i,l})\ot (t_x^{i,l}\stackrel{\Sigma}{\circ}  t_y)\nonumber.
\end{eqnarray}

Thus, $f_{x,y}^{(\Gm\Om)\Sigma}=f_{x,y}^{\Gm(\Om\Sigma)}$. Therefore, this operation is associative.

	\underline{Neutral element:} We will check that $[\D(\T)]$ is the neutral element. Given $[\D(\Omega)]\in \C(\T/R),$ we have
$$[\D(\Omega)][\D(\T)]=\left[ \bigoplus_{x\in G}\Omega_x\ot
\T_{x^\m}\ot\T_x\right].$$
For each $x \in G$, consider the $R$-bimodule isomorphism
$$\begin{array}{c c c l}
	\varphi_x: & \Om_x\ot\T_x\ot\T_{x^\m}& \longrightarrow & \Om_x\\
	&\om_x\ot u_x\ot u_{x^\m}&  \longmapsto & \om_x(u_x\stackrel{\T}{\circ}u_{x^\m}).
\end{array}$$
Then, the following diagram is commutative
$$\xymatrix@C=3cm{ \Om_x\ot\T_{x^\m}\ot\T_x\ot\Om_y\ot\T_{y^\m}\ot\T_y \ar[dd]_{\varphi_x\ot \varphi_y}\ar[r]^{f_{x,y}^{\Om\T}}  & 1_x\Om_{xy}\ot\T_{(xy)^\m}\ot\T_{xy}\ar[dd]^{\varphi_{xy}}\\
	& \\
	\Omega_x\ot\Om_y\ar[r]_{f_{x,y}^{\Om}}& 1_x\Om_{xy}     }$$
Indeed, on  one hand,
\begin{equation}
	f_{x,y}^\Om\circ(\varphi_x\ot \varphi_y)(\om_x\ot u_{x^\m}\ot u_x\ot \om_y\ot u_{y^\m}\ot u_y)=(\om_x(u_{x^\m}\stackrel{\T}{\circ} u_x)\stackrel{\Om}{\circ} \om_y( u_{y^\m}\stackrel{\T}{\circ} u_y)). \label{unidadeemClado1}
\end{equation}
On the other hand, observe that $\varphi_{xy}\circ f_{x,y}^{\Om\T}$ is given by
\begin{equation}
	{\small 	\xymatrix@C=1cm{ \Om_x\ot\T_{x^\m}\ot\T_x\ot\Om_y\ot\T_{y^\m}\ot\T_y\ar[rd]^-{\Om_x\ot T\ot \T_y}\ar@/_1.5cm/@{-->}[rdddd]_{\varphi_{xy}\circ f_{x,y}^{\Om\T}} & \\
			&\Om_x\ot\Om_y\ot\T_{y^\m}\ot\T_{x^\m}\ot\T_x\ot\T_y\ar[dd]|{f_{x,y}^\Om\ot f_{y^\m,x^\m}^\T\ot f_{x,y}^\T} \\
			& \\
			& 1_x\Om_{xy}\ot\T_{(xy)^\m}\ot\T_{xy}\ar[d]^{\varphi_{xy}}\\
			& 1_x\Om_{xy}}} \label{inverso2}
\end{equation}
where $T=T_{\T_{x^\m}\ot_R\T_x,\Om_y\ot_R\T_{y^\m}}$. By Example \ref{execomRe2}, we have
\begin{equation*}
	\begin{array}{c c c l}
		T: & \T_{x^\m}\ot\T_x\ot\Om_y\ot\T_{y^\m} & \longrightarrow & \Om_y\ot\T_{y^\m}\ot\T_{x^\m}\ot\T_x\\
		& u_{x^\m}\ot u_x\ot \om_y\ot u_{y^\m} & \longmapsto & \dsum_{(x)}(u_{x^\m}\stackrel{\T}{\circ}u_x)\om_y\ot u_{y^\m}\ot \widetilde{u}_{x^\m}\ot \widetilde{u}_x,
	\end{array}
\end{equation*}
where $\dsum_{(x)}(\widetilde{u}_{x^\m}\stackrel{\T}{\circ} \widetilde{u}_x)=1_{x^\m}$. Thus, by diagram (\ref{inverso2})  we have that $\varphi_{xy}\circ f_{x,y}^{\Om\T}$ is given by
	\begin{eqnarray*}
& &	\varphi_{xy}\circ f_{x,y}^{\Om\T}(\om_x\ot u_{x^\m}\ot u_x\ot \om_y\ot u_{y^\m}\ot u_y)\\
& = & \varphi_{xy}\circ (f_{x,y}^\Om\ot f_{y^\m,x^\m}^\T\ot f_{x,y}^\T)\circ (\Om_x\ot T\ot \T_y) (\om_x\ot u_{x^\m}\ot u_x\ot \om_y\ot u_{y^\m}\ot u_y)\\
&=  & \dsum_{(x)}\varphi_{xy}\circ (f_{x,y}^\Om\ot f_{y^\m,x^\m}^\T\ot f_{x,y}^\T)(\om_x\ot (u_{x^\m}\stackrel{\T}{\circ}u_x)\om_y\ot u_{y^\m}\ot \widetilde{u}_{x^\m}\ot \widetilde{u}_x \ot u_y)\\
	& = &  \dsum_{(x)}\varphi_{xy} ((\om_x\stackrel{\Om}{\circ} (u_{x^\m}\stackrel{\T}{\circ}u_x)\om_y)\ot (u_{y^\m}\stackrel{\T}{\circ} \widetilde{u}_{x^\m})\ot (\widetilde{u}_x\stackrel{\T}{\circ} u_y))\\
	& = & \dsum_{(x)}  (\om_x\stackrel{\Om}{\circ} (u_{x^\m}\stackrel{\T}{\circ}u_x)\om_y(u_{y^\m}\stackrel{\T}{\circ} (\widetilde{u}_{x^\m}\stackrel{\T}{\circ} \widetilde{u}_x)\stackrel{\T}{\circ} u_y))\\
	& = &  (\om_x(u_{x^\m}\stackrel{\T}{\circ}u_x)\stackrel{\Om}{\circ}\om_y(u_{y^\m}1_{x^\m}\stackrel{\T}{\circ} u_y))\\
	& = &  (\om_x(u_{x^\m}\stackrel{\T}{\circ}u_x)\stackrel{\Om}{\circ} \om_y(u_{y^\m}\stackrel{\T}{\circ} u_y)),
\end{eqnarray*}
where the last equality follows from (\ref{uxercomutam}) and (\ref{mercomutam}). Hence, $f_{x,y}^\Om\circ (\varphi_x\ot \varphi_y)=\varphi_{xy}\circ f_{x,y}^{\Om\T}$ and the assertion follows. The commutativity of the diagram  yields  $[\D(\Omega)][\D(\T)]=[\D(\Om)]$ in $\C(\T/R).$
The commutativity of the operation will be proved later and will imply that  $[\Delta(\T)]$ is the neutral element.

\underline{Inverse element:} Next we check that the inverse element of the class $[\D(\Om)] \in \C(\T/R)$ is given by
$$[\D(\Om)]^{\m}=\left[ \bigoplus_{x\in G}\T_x\ot
\Om_{x^\m}\ot
\T_x\right]. $$
By Lemma \ref{conjuntodefatores} we have $[\D(\Om)]^{\m}\in \C(\T/R)$. By definition
\begin{equation*}
	[\D(\Om)][\D(\Om)]^{\m}=\left[ \bigoplus_{x \in G}\Om_x\ot\T_{x^\m}\ot\T_x\ot\Om_{x^\m}\ot\T_x \right]. 
\end{equation*}
Consider the $R$-bimodule isomorphism $\psi_x:\Om_x\ot\T_{x^\m}\ot\T_x\ot\Om_{x^\m}\ot\T_x \rightarrow \T_x$  defined via
$$\xymatrix{ \Om_x\ot\T_{x^\m}\ot\T_x\ot\Om_{x^\m}\ot\T_x\ar[r]\ar@/_1cm/@{-->}[rddd]_{\psi_x} &  \Om_x\ot R1_{x^\m}\ot\Om_{x^\m}\ot\T_x\ar[d]\\
	&\Om_x\ot\Om_{x^\m}\ot\T_x\ar[d]\\
	& R1_x\ot\T_x\ar[d]\\
	& \T_x }$$
that is,
$$\psi_x(\om_x\ot u_{x^\m}\ot u_x\ot \om_{x^\m}\ot u'_x)=(\om_x(u_{x^\m}\stackrel{\T}{\circ}u_x)\stackrel{\Om}{\circ}\om_{x^\m})u'_x.$$
Write, for simplicity of notation, $\Pi_x=\Om_x\ot\T_{x^\m}\ot\T_x\ot\Om_{x^\m}\ot\T_x$, for $x \in G$.  We check that the following diagram 
\begin{equation}
	\xymatrix@C=2cm{ \Pi_x\ot\Pi_y\ar[rr]^{f_{x,y}^{\Om\Om^{\m}}}\ar[dd]_{\psi_x\ot \psi_y} & &  1_x\Pi_{xy}\ar[dd]^{\psi_{xy}}\\
		& & \\
		\T_x\ot\T_y\ar[rr]_{f_{x,y}^\T} & & 1_x\T_{xx}  } \label{diagramaelementoneutro}
\end{equation} 
is commutative. On  one hand, we have
\begin{eqnarray}
	& & (f_{x,y}^\T\circ (\psi_x\ot \psi_y))(\om_x\ot u_{x^\m}\ot u_x\ot \om_{x^\m}\ot u'_x\ot \om_y\ot u_{y^\m}\ot u_y\ot \om_{y^\m}\ot u'_y) \nonumber\\
	& = & \left( (\om_x(u_{x^\m}\stackrel{\T}{\circ}u_x)\stackrel{\Om}{\circ}\om_{x^\m})u'_x\stackrel{\T}{\circ}(\om_y(u_{y^\m}\stackrel{\T}{\circ}u_y)\stackrel{\Om}{\circ}\om_{y^\m})u'_y\right) .\label{diagramainversolado1}
\end{eqnarray}

On the other hand, let $T_1=T_{\T_{x^\m}\ot \T_{x}\ot \Om_{x^\m}\ot \T_x,\Om_y\ot \T_{y^\m}}$ and $T_2=T_{\Om_{x^\m\ot \T_x,\T_y\ot \Om_{y^\m}}}$, then $\psi_{xy}\circ f^{\Om\Om^\m}_{x,y}$ is given by the composition: 
\begin{equation}
	{\small 	\xymatrix@C=-5cm{ \Om_x\ot\T_{x^\m}\ot\T_x\ot\Om_{x^\m}\ot\T_x\ot\Om_y\ot\T_{y^\m}\ot\T_y\ot\Om_{y^\m}\ot\T_y\ar[rd]^-{\Om_x\ot T_1\ot \T_y\ot\Om_{y^\m}\ot\T_y}\ar@/_3cm/@{-->}[rdddddd]_{\psi_{xy}\circ f_{x,y}^{\Om\Om^{\m}}} & \\
			& \Om_x\ot\Om_y\ot\T_{y^\m}\ot\T_{x^\m}\ot\T_x\ot\Om_{x^\m}\ot\T_x\ot\T_y\ot\Om_{y^\m}\ot\T_y\ar[dd]|-{\Om_x\ot\Om_y\ot\T_{y^\m}\ot\T_{x^\m}\ot\T_x\ot T_2\ot \T_y} \\
			& \\
			& \Om_x\ot\Om_y\ot\T_{y^\m}\ot\T_{x^\m}\ot\T_x\ot\T_y\ot\Om_{y^\m}\ot\Om_{x^\m}\ot\T_x\ot\T_y\ar[dd]|{f_{x,y}^\Om\ot f_{y^\m,x^\m}^\T\ot f_{x,y}^\T\ot f_{y^\m,x^\m}^\Om\ot f_{x,y}^\T} \\
			& \\
			& 1_x\Om_{xy}\ot\T_{(xy)^\m}\ot\T_{xy}\ot\Om_{(xy)^\m}\ot\T_{xy}\ar[d]^{\psi_{xy}}\\
			& 1_x\T_{xy}}}\label{diagramainverso2}
\end{equation}
Let $f_i:\Om_{x^\m}\ot\T_x\rightarrow R$ and $g_i:R\rightarrow \Om_{x^\m}\ot\T_x$,  $i=1,2,...,n$, be $R$-bilinear maps such that $\dsum_{i=1}^n g_if_i=Id_{\Om_{x^\m}\ot_R\T_x}.$ Denote $g_i(1)=\dsum_{l}\om_{x^\m}^{i,l}\ot u_x^{i,l}.$ Then,  by Proposition~\ref{isomorfismoT},
$$\begin{array}{c c c l}
	T_2: & \Om_{x^\m}\ot\T_x\ot\T_y\ot\Om_{y^\m} & \longrightarrow & \T_y\ot\Om_{y^\m}\ot\Om_{x^\m}\ot\T_x \\
	& \om_{x^\m}\ot u_x\ot u_y\ot \om_{y^\m} & \longmapsto & \dsum_{i,l}f_i(\om_{x^\m}\ot u_x)u_y\ot \om_{y^\m} \ot \om^{i,l}_{x^\m}\ot u_x^{i,l}.
\end{array}$$ 
For $T_1$ consider the $R$-bilinear maps $$\bar{f}_i:\T_{x^\m}\ot\T_x\ot\Om_{x^\m}\ot\T_x\longrightarrow R1_{x^\m}\ot\Om_{x^\m}\ot\T_x\longrightarrow \Om_{x^\m}\ot\T_x\stackrel{f_i}{\longrightarrow} R,$$
$$\bar{g}_i:R\stackrel{g_i}{\longrightarrow} \Om_{x^\m}\ot\T_x\longrightarrow R1_{x^\m}\ot\Om_{x^\m}\ot\T_x\longrightarrow\T_{x^\m}\ot \T_x\ot\Om_{x^\m}\ot\T_x.$$
Then, $$\bar{f}_i(u_{x^\m}\ot u_x\ot \om_{x^\m}\ot u_x')=f_i((u_{x^\m}\stackrel{\T}{\circ}u_x)\om_{x^\m}\ot u_x')$$ $$\bar{g}_i(1)=\dsum_{j,l}\bar{u}_{x^\m}^j\ot \bar{u}_x^j\ot\om_{x^\m}^{i,l}\ot u^{i,l}_x,$$
where  $\dsum_{j=1}^m(\bar{u}_{x^\m}^j\stackrel{\T}{\circ} \bar{u}_x^j) =1_{x^\m}$. It is easy to see that $\dsum_{i=1}^n\bar{g}_i\bar{f}_i=Id_{\T_{x^\m}\ot \T_x\ot_R\Om_{x^\m}\ot_R\T_x}$. Then we can write $T_1$ depending on $\bar{f}_i$ and $\bar{g}_i$ as follows
	\begin{equation*}
		\begin{array}{c c c l}
			T_1:&  \T_{x^\m}\ot\T_x\ot\Om_{x^\m}\ot\T_x\ot\Om_y\ot\T_{y^\m}& \rightarrow & \Om_y\ot\T_{y^\m}\ot\T_{x^\m}\ot\T_x\ot\Om_{x^\m}\ot\T_x \\
			& u_{x^\m}\ot u_x \ot \om_{x^\m} \ot u'_{x^\m}\ot \om_y\ot u_{y^\m} & \mapsto & \dsum_{i,j,l} f_i((u_{x^\m}\stackrel{\T}{\circ}u_x)\om_{x^\m}\ot u_x')\om_y\ot u_{y^\m}\ot \bar{g}_i(1).
		\end{array}
\end{equation*}

Applying  $T_1$ and  $T_2$ in $\om_x\ot u_{x^\m}\ot u_x\ot \om_{x^\m}\ot u'_x\ot \om_y\ot u_{y^\m}\ot u_y\ot \om_{y^\m}\ot u'_y$, we have:
	\begin{eqnarray*}
		& &\om_x\ot u_{x^\m}\ot u_x\ot \om_{x^\m}\ot u'_x\ot \om_y\ot u_{y^\m}\ot u_y\ot \om_{y^\m}\ot u'_y \\
		& \mapsto & \dsum_{i,j,k,l} \om_x\ot f_k((u_{x^\m}\stackrel{\T}{\circ} u_x)\om_{x^\m}\ot u'_x)\om_y\ot u_{y^\m}\ot \bar{u}^j_{x^\m}\ot \bar{u}^j_x \ot f_i(g_k(1))u_y\ot \om_{y^\m}\ot g_i(1)\ot u'_y
	\end{eqnarray*}
	As $f_i$ and $g_i$ are $R$-bilinear, then $f_i(g_k(1)) \in \Z$, for all $i,k=1,2,...,n$. By (\ref{mercomutam}) follows that
	\begin{equation*}
		\om_y\ot u_{y^\m}\ot \bar{u}^j_{x^\m}\ot \bar{u}^j_x \ot f_k(g_i(1))u_y=f_k(g_i(1))\om_y\ot u_{y^\m}\ot \bar{u}^j_{x^\m}\ot \bar{u}^j_x \ot u_y
	\end{equation*}
	Thus,
	\begin{eqnarray*}
		& & \dsum_{i,j,k,l} \om_x\ot f_i((u_{x^\m}\stackrel{\T}{\circ} u_x)\om_{x^\m}\ot u'_x)\om_y\ot u_{y^\m}\ot \bar{u}^j_{x^\m}\ot \bar{u}^j_x \ot f_k(g_i(1))u_y\ot \om_{y^\m}\ot g_k(1)\ot u'_y\\
		& = & \dsum_{i,j,k,l} \om_x\ot f_i((u_{x^\m}\stackrel{\T}{\circ} u_x)\om_{x^\m}\ot u'_x)f_k(g_i(1))\om_y\ot u_{y^\m}\ot \bar{u}^j_{x^\m}\ot \bar{u}^j_x \ot u_y\ot \om_{y^\m}\ot g_k(1)\ot u'_y\\
		& = & \dsum_{i,j,k,l} \om_x\ot f_k(g_i(f_i((u_{x^\m}\stackrel{\T}{\circ} u_x)\om_{x^\m}\ot u'_x)))\om_y\ot u_{y^\m}\ot \bar{u}^j_{x^\m}\ot \bar{u}^j_x \ot u_y\ot \om_{y^\m}\ot g_k(1)\ot u'_y\\
		& = & \dsum_{k,j,l} \om_x\ot f_k((u_{x^\m}\stackrel{\T}{\circ} u_x)\om_{x^\m}\ot u'_x)\om_y\ot u_{y^\m}\ot \bar{u}^j_{x^\m}\ot \bar{u}^j_x \ot u_y\ot \om_{y^\m}\ot \omega_{x^\m}^{k,l}\ot u_x^{k,l}\ot u'_y.
	\end{eqnarray*}
	Therefore,\begin{eqnarray*}
		& &\psi_{xy}\circ f^{\Om\Om^\m}_{x,y}( \om_x\ot u_{x^\m}\ot u_x\ot \om_{x^\m}\ot u'_x\ot \om_y\ot u_{y^\m}\ot u_y\ot \om_{y^\m}\ot u'_y)\nonumber\\
		& = & \dsum_{j,k,l} \left( \om_x \stackrel{\Om}{\circ} f_k((u_{x^\m}\stackrel{\T}{\circ}u_x)\om_{x^\m}\ot u_x')\om_y (u_{y^\m} \stackrel{\T}{\circ}\bar{u}^j_{x^\m} \stackrel{\T}{\circ}\bar{u}^j_x\stackrel{\T}{\circ} u_y)\stackrel{\Om}{\circ} \om_{y^\m} \stackrel{\Om}{\circ} \om^{k,l}_{x^\m}\right)(u^{k,l}_x\stackrel{\T}{\circ} u'_y)\nonumber\\
		& = & \dsum_{k,l} \left( \om_x(u_{x^\m}\stackrel{\T}{\circ}u_x) \stackrel{\Om}{\circ} f_k(\om_{x^\m}\ot u_x')\om_y (u_{y^\m} 1_{x^\m}\stackrel{\T}{\circ} u_y)\stackrel{\Om}{\circ} \om_{y^\m} \stackrel{\Om}{\circ} \om^{k,l}_{x^\m}\right)(u^{k,l}_x\stackrel{\T}{\circ} u'_y).\nonumber
	\end{eqnarray*}
	Denote $a=f_k(\om_{x^\m}\ot u_x')(\om_y (u_{y^\m} 1_{x^\m}\stackrel{\T}{\circ} u_y)\stackrel{\Om}{\circ} \om_{y^\m}) \in R$. By (\ref{uxe1y}), we have $a=f_k(\om_{x^\m}\ot u_x')(\om_y (u_{y^\m}\stackrel{\T}{\circ} u_y)\stackrel{\Om}{\circ} \om_{y^\m})$  and 
	\begin{eqnarray*}
		\dsum_{k,l}a\om^{k,l}_{x^\m}\ot u^{k,l}_x & = & \dsum_kag_k(1)= \dsum_kg_k(1)a= \dsum_kg_k(1)f_k(\om_{x^\m}\ot u_x')(\om_y (u_{y^\m}\stackrel{\T}{\circ} u_y)\stackrel{\Om}{\circ} \om_{y^\m})\\
		& = & \dsum_kg_k(f_k(\om_{x^\m}\ot u_x'))(\om_y (u_{y^\m}\stackrel{\T}{\circ} u_y)\stackrel{\Om}{\circ} \om_{y^\m})= \om_{x^\m}\ot u_x'(\om_y (u_{y^\m}\stackrel{\T}{\circ} u_y)\stackrel{\Om}{\circ} \om_{y^\m}).
	\end{eqnarray*}
	Then, 
	\begin{eqnarray*}
		& &\psi_{xy}\circ f^{\Om\Om^\m}_{x,y}( \om_x\ot u_{x^\m}\ot u_x\ot \om_{x^\m}\ot u'_x\ot \om_y\ot u_{y^\m}\ot u_y\ot \om_{y^\m}\ot u'_y)\nonumber\\
		& = & \dsum_{k,l} \left( \om_x(u_{x^\m}\stackrel{\T}{\circ}u_x) \stackrel{\Om}{\circ} a\om^{k,l}_{x^\m}\right)(u^{k,l}_x\stackrel{\T}{\circ} u'_y)\nonumber\\
		& = & \dsum_{k,l} (f_{x,y}^\T\circ (f_{x,x^{\m}}^\Om\ot \T_x\ot \T_y))(\om_x(u_{x^{-1}}\stackrel{\T}{\circ}u_x)\ot a\om_{x^\m}^{k,l}\ot u_x^{k,l}\ot u_y')\nonumber\\
		& = & (f_{x,y}^\T\circ (f_{x,x^{\m}}^\Om\ot \T_x\ot \T_y))(\om_x(u_{x^{-1}}\stackrel{\T}{\circ}u_x)\ot \om_{x^\m}\ot u_x'(\om_y (u_{y^\m}\stackrel{\T}{\circ} u_y)\stackrel{\Om}{\circ} \om_{y^\m})\ot u_y')\nonumber\\
		& = & (f_{x,y}^\T\circ (f_{x,x^{\m}}^\Om\ot \T_x\ot \T_y))(\om_x(u_{x^{-1}}\stackrel{\T}{\circ}u_x)\ot \om_{x^\m}\ot u_x'\ot (\om_y (u_{y^\m}\stackrel{\T}{\circ} u_y)\stackrel{\Om}{\circ} \om_{y^\m})u_y')\nonumber\\
		& = & ((\om_x(u_{x^{-1}}\stackrel{\T}{\circ}u_x)\stackrel{\Om}{\circ} \om_{x^\m})u_x'\stackrel{\T}{\circ} (\om_y (u_{y^\m}\stackrel{\T}{\circ} u_y)\stackrel{\Om}{\circ} \om_{y^\m})u_y').
\end{eqnarray*}
Comparing with (\ref{diagramainversolado1}) it follows that the diagram (\ref{diagramaelementoneutro}) is commutative. Therefore, $[\D(\Om)][\D(\Om)]^{\m}=[\D(\T)] $.

\underline{Commutativity:} Given $[\D(\Gamma)], [\D(\Om)]\in \C(\T/R)$, we have
	$$[\D(\Gamma)][\D(\Om)]=\left[ \bigoplus_{x \in G}\Gamma_x\ot\T_{x^\m}\ot\Om_x \right] \ \ \mbox{and} \ \ [\D(\Om)][\D(\Gamma)]=\left[ \bigoplus_{x \in G}\Om_x\ot\T_{x^\m}\ot\Gamma_x \right]. $$
For each $x\in G$, consider the isomorphism $\lambda_x$ defined by
	$$\xymatrix{ \Gamma_x\ot\T_{x^\m}\ot\Om_x\ar[r]\ar@/_1.5cm/@{-->}[rrddd]_{\lambda_x}  & \Gamma_x\ot\T_{x^\m}\ot\Om_x\ot R1_{x^\m}\ar[r] & \Gamma_x\ot\T_{x^\m}\ot\Om_x\ot\T_{x^\m}\ot\T_{x}\ar[d]\\
		& & \Om_x\ot\T_{x^\m}\ot\Gamma_x\ot\T_{x^\m}\ot\T_{x}\ar[d]\\
		& & \Om_x\ot\T_{x^\m}\ot\Gamma_x\ot R1_{x^\m}\ar[d]\\
		& & \Om_x\ot\T_{x^\m}\ot\Gamma_x}$$
	that is,
	$$\lambda_x=(\Om_x\ot \T_{x^\m}\ot \Gm_x\ot f_{x^\m,x}^\T)\circ (T_1\ot \T_x)\circ (\Gm_x\ot \T_{x^\m} \ot \Om_x\ot (f_{x^\m,x}^\T)^\m),$$
	where $T_1=T_{\Gamma_x\ot\T_{x^\m},\Om_x\ot\T_{x^\m}}$.
	We need to show that  the following diagram is commutative:
	\begin{equation}
\xymatrix@C=1cm{ \Gamma_x\ot\T_{x^\m}\ot\Om_x\ot\Gamma_y\ot\T_{y^\m}\ot\Om_y\ar[rr]^{\ \ \ \ \ \ f_{x,y}^{\Gamma\Om}}\ar[ddd]_{\lambda_x\ot \lambda_y} & &  1_x\Gamma_{xy}\ot\T_{(xy)^\m}\ot\Om_{xy}\ar[ddd]^{1_x\lambda_{xy}}\\
	& & \\
	& & \\
	\Om_x\ot\T_{x^\m}\ot\Gamma_x\ot\Om_y\ot\T_{y^\m}\ot\Gamma_y\ar[rr]_{\ \ \ \ \ f_{x,y}^{\Om\Gamma}} & & 1_x\Om_{xy}\ot\T_{(xy)^\m}\ot\Gamma_{xy}.	  } \label{comutatividadediagrama}
	\end{equation}
Consider the isomorphisms $T_2=T_{\Gamma_y\ot \T_{y^\m},\Om_y\ot \T_{y^\m}}, T_3=T_{\T_{x^\m}\ot\Gamma_x,\Om_y\ot\T_{y^\m}},$  $T_4=T_{\T_{x^\m}\ot\Om_x,\Gamma_y\ot \T_{y^\m}},$ and $T_5=T_{\Gamma_{xy}\ot \T_{(xy)^\m},1_x\Om_{xy}\ot 1_{y^\m}\T_{(xy)^\m}}$.
Then
 $f_{x,y}^{\Om\Lambda}\circ (\lambda_x\ot \lambda_y)$ is given by the sequence of isomorphisms:
	\begin{equation}
{\small 	\xymatrix@C=-2.5cm{ & \Gamma_x\ot\T_{x^\m}\ot\Om_x\ot\Gamma_y\ot\T_{y^\m}\ot\Om_y\ar[ld]\ar@{-->}@/^3.5cm/[ldddddd]^{f_{x,y}^{\Om\Gm}\circ(\lambda_x\ot \lambda_y)}\\  
		\Gamma_x\ot\T_{x^\m}\ot\Om_x\ot\T_{x^{-1}}\ot\T_x\ot \Gamma_y\ot\T_{y^\m}\ot\Om_y\ot \T_{y^{-1}}\ot\T_y \ar[d]_{T_1\ot\T_x\ot T_2\ot \T_y} & \\
		\Om_x\ot\T_{x^{-1}}\ot\Gamma_x\ot\T_{x^\m}\ot\T_x\ot \Om_y\ot \T_{y^{-1}}\ot\Gamma_y\ot\T_{y^\m}\ot\T_y\ar[d]& \\
		\Om_x\ot\T_{x^{-1}}\ot\Gamma_x\ot\Om_y\ot \T_{y^{-1}}\ot\Gamma_y\ar[d]_{\Om_x\ot T_3\ot \Gm_y}& \\
		\Om_x\ot\Om_y\ot \T_{y^{-1}}\ot\T_{x^{-1}}\ot\Gamma_x\ot\Gamma_y\ar[dd]|{f_{x,y}^\Om\ot f_{y^\m,x^\m}^\T\ot f_{x,y}^\Gm}  & \\
		& \\
		1_x\Om_{xy}\ot\T_{(xy)^\m}\ot\Gamma_{xy} & }} \label{diagramacomutatividade1}
	\end{equation}
	%

On the other hand, $1_x\lambda_{xy}\circ f_{x,y}^{\Gamma\Om}$ is given by the sequence of isomorphisms:
	
	\begin{equation}
{\small 	\xymatrix{  \Gamma_x\ot\T_{x^\m}\ot\Om_x\ot\Gamma_y\ot\T_{y^\m}\ot\Om_y\ar[rd]^{\Gm_x\ot T_4\ot \Om_y}\ar@{-->}@/_4.5cm/[rdddddddd]_{1_x\lambda_{xy}\circ f_{x,y}^{\Gamma\Om}}& \\
		&\Gamma_x\ot\Gamma_y\ot\T_{y^\m}\ot\T_{x^\m}\ot\Om_x\ot\Om_y\ar[dd]|{f_{x,y}^{\Gm}\ot f_{y^\m,x^\m}^\T\ot f_{x,y}^\Om}\\
		& \\
		&  1_x\Gamma_{xy}\ot 1_{y^\m}\T_{(xy)^\m}\ot 1_x\Om_{xy}\ar@{=}[d]\\
		&  \Gamma_{xy}\ot \T_{(xy)^\m}\ot1_x\Om_{xy}1_{y^\m}\ar[d]\\
		&  \Gamma_{xy}\ot \T_{(xy)^\m}\ot 1_{x}\Om_{xy}\ot R1_{y^\m}1_{(xy)^\m}\ar[d]\\
		& \Gamma_{xy}\ot \T_{(xy)^\m}\ot 1_{x}\Om_{xy}\ot 1_{y^\m}\T_{(xy)^{-1}} \ot \T_{xy}\ar[d]^{T_5\ot \T_{xy}}\\
		&  1_x\Om_{xy}\ot 1_{y^\m}\T_{(xy)^{-1}} \ot \Gamma_{xy}\ot \T_{(xy)^\m}\ot\T_{xy}\ar[d]\\
		& 1_x\Om_{xy}\ot\T_{(xy)^\m}\ot\Gamma_{xy}}} \label{diagramacomutatividade2}
	\end{equation}

		Consider the $R$-bilinear maps $f_l:\Gamma_y\ot\Om_{y^\m}\longrightarrow R$ and $g_l:R\longrightarrow \Gamma_y\ot\Om_{y^\m}$, $l=1,2,...,k$, such that $\dsum_{l=1}^kg_lf_l=Id_{\Gamma_y\ot\Om_{y^\m}}$. Denote $g_l(1)=\dsum_{t}v_y^{l,t}\ot u_{y^\m}^{l,t}$. Then, by Proposition \ref{isomorfismoT},
	$$\begin{array}{c c c l}
		T_2: & \Gamma_y\ot\T_{y^\m}\ot\Om_y\ot\T_{y^\m} &\rightarrow&  \Om_y\ot\T_{y^\m}\ot\Gamma_y\ot\T_{y^\m},\\
		& v_y\ot u_{y^\m} \ot \om_y\ot u'_{y^\m} & \mapsto & \dsum_{l,t} f_l(v_y\ot u_{y^\m})\om_y\ot u'_{y^\m}\ot v_y^{l,t}\ot u_{y^\m}^{l,t},
	\end{array}$$
	and, by Remark \ref{remarkT},
	$$\begin{array}{c c c l}
		T_4: & \T_{x^\m}\ot\Om_x\ot\Gamma_y\ot\T_{y^\m} & \rightarrow &  \Gamma_y\ot\T_{y^\m}\ot\T_{x^\m}\ot\Om_x,\\
		& u_{x^\m}\ot \om_x\ot v_y\ot u_{y^\m} & \mapsto & \dsum_{l,t}v_y^{l,t}\ot u_{y^\m}^{l,t}\ot u_{x^\m}\ot \om_xf_l( v_y\ot u_{y^\m}).
	\end{array}$$
	Now, for  $T_1$ and $T_3$, consider the $R$-bilinear maps $f_i^x:\Om_x\ot\T_{x^\m}\longrightarrow R$ , $g_i^x:R\longrightarrow \Om_x\ot\T_{x^\m}$, $f_j^y:\Om_y\ot\T_{y^\m}\longrightarrow R$ and $g_j^y:R\longrightarrow \Om_y\ot\T_{y^\m} $, with $i=1,2,...,n$ and $j=1,2,...,m$, such that $\dsum_{i=1}^ng_i^xf_i^x=Id_{\Om_x\ot\T_{x^\m}} $ and $\dsum_{j=1}^mg_j^yf_j^y=Id_{\Om_y\ot\T_{y^\m}}$. Denote,
	$$g_i^x(1)=\dsum_{i,p}\widetilde{\om}_x^{i,p}\ot\widetilde{u}_{x^\m}^{i,p} \ \ \mbox{e} \ \ g_j^y(1)=\dsum_{j,q}\widehat{\om}_y^{j,q}\ot \widehat{u}_{y^\m}^{j,q}.$$
	Then, by Remark \ref{remarkT}, 
	$$\begin{array}{c c c l}
	T_1: & \Gamma_x\ot\T_{x^\m}\ot\Om_x\ot\T_{x^\m} & \rightarrow &  \Om_x\ot\T_{x^\m}\ot\Gamma_x\ot\T_{x^\m},\\
	& v_x\ot u_{x^\m}\ot \om_x\ot u'_{x^\m} &\mapsto& \dsum_{i,p} \widetilde{\om}_x^{i,p}\ot\widetilde{u}_{x^\m}^{i,p}\ot v_x\ot u_{x^\m}f_i^x(\om_x\ot u'_{x^\m}).
\end{array}$$
$$\begin{array}{c c c l}
	T_3: &  \T_{x^\m}\ot\Gamma_x\ot\Om_y\ot\T_{y^\m} & \rightarrow &  \Om_y\ot\T_{y^\m}\ot\T_{x^\m}\ot\Gamma_x ,\\
	& u_{x^\m} \ot v_x\ot \om_y\ot u'_{y^\m} & \mapsto & \dsum_{j,q} \widehat{\om}_y^{j,q}\ot \widehat{u}_{y^\m}^{j,q}\ot u_{x^\m}\ot v_xf_j^y(\om_y\ot u'_{y^\m}). 
\end{array}$$
	Finally, for $T_5$, consider $f_{i,j}^{xy}:1_x\Om_x\ot 1_{y^\m}\T_{(xy)^\m}\rightarrow R$ defined by
	$$\xymatrix@C=4cm{1_x\Om_{xy}\ot 1_{y^\m}\T_{(xy)^{-1}}\ar[r]^{(f_{x,y}^\Om)^\m\ot (f_{x,y}^\T)^{\m}}\ar@{-->}@/_1cm/[rdd]_{f_{i,j}^{xy}} & \Om_x\ot_R\Om_y\ot \T_{y^{-1}}\ot_R\T_{x^{-1}}\ar[d]^{\Om_x\ot f_j^y\ot \T_{x^\m}}\\
		& \Om_x\ot \T_{x^{-1}}\ar[d]^{f_i^x}\\
		& R} $$
	Then, $f_{i,j}^{xy}$ is $R$-bilinear and 
	\begin{equation*}
	f_{i,j}^{xy}(\om_{xy}\ot u_{(xy)^{-1}})=\dsum f_i^x(\om_x^gf_j^y(\om_y^g\ot u_{y^\m}^h)\ot u_{x^\m}^h), \label{fijparacomutatividadeemC}
	\end{equation*} 
	where $\dsum_{g} (\om_x^g\stackrel{\Om}{\circ} \om_y^g)=\om_{xy}$ and $\dsum_h (u_{y^\m}^h\stackrel{\T}{\circ} u_{x^\m}^h)=u_{(xy)^\m}$. 
	Analogously, consider $g_{i,j}^{xy}:R \rightarrow 1_x\Om_{xy}\ot 1_{y^\m}\T_{(xy)^{-1}}$ defined by
	$$\xymatrix{ R\ar[r]^{g_i^x\ \ \  \ }\ar@{-->}@/_1cm/[rrrdd]_{g_{i,j}^{xy}} & \Om_x\ot\T_{x^{-1}} \ar[r] & \Om_x\ot R\ot \T_{x^\m}\ar[r]^{g_j^y\ \ \  \ \ \ \ \ \ \ \ } & \Om_x\ot_R\Om_y\ot_R\T_{y^{-1}}\ot \T_{x^{-1}}\ar[dd]|{f_{x,y}^\Om\ot f_{x,y}^\T}\\
		& & & \\
		& 	& & 1_x\Om_{xy}\ot 1_{y^\m}\T_{(xy)^{-1}}}$$
	Then, $g_{i,j}^{xy}$ is $R$-bilinear and 
	$$g_{i,j}^{xy}(1)=\dsum_{i,j,p,q} (\widetilde{\om}_x^{i,p}\stackrel{\Om}{\circ}\widehat{\om}_y^{j,q})\ot (\widehat{u}_{y^\m}^{j,q}\stackrel{\T}{\circ}\widetilde{u}_{x^\m}^{i,p}).$$
	It is easy to see that $\dsum_{i,j}^{n,m}g_{i,j}^{xy}f_{i,j}^{xy}=Id_{1_x\Om_{xy}\ot 1_{y^\m}\T_{(xy)^\m}}$. Then,
	
	$$\begin{array}{c c c l}
	T_5: & \Gamma_{xy}\ot\T_{(xy)^\m}\ot 1_x\Om_{xy}\ot 1_{y^\m}\T_{(xy)^\m}& \rightarrow & 1_x\Om_{xy}\ot 1_{y^\m}\T_{(xy)^\m}\ot\Gamma_{xy}\ot\T_{(xy)^\m},\\
	& v_{xy}\ot u_{(xy)^\m}\ot \om_{xy}\ot u'_{(xy)^\m} & \mapsto & \dsum g_{i,j}^{xy}(1)\ot v_{xy}\ot u_{(xy)^\m}f_{i,j}^{xy}(\om_{xy}\ot u_{(xy)^\m}').
	\end{array}$$
	We denote
	$$\dsum_{(x)}(\bar{u}_{x^\m}\stackrel{\T}{\circ} \bar{u}_x)=1_{x^\m} \ \ \mbox{and} \ \ \dsum_{(y)}(\bar{u}_{y^\m}\stackrel{\T}{\circ} \bar{u}_y)=1_{y^\m}.$$

	Let us check that  diagram (\ref{comutatividadediagrama}) is commutative. First, by diagram (\ref{diagramacomutatividade1}), we have that $f^{\Om\Lambda}\circ (\lambda_x\ot \lambda_y)$ is: 
	\begin{eqnarray*}
		& &  v_x\ot u_{x^\m}\ot \om_x\ot v_y\ot u_{y^\m}\ot \om_y\\
		&\mapsto &  v_x\ot u_{x^\m}\ot \om_x\ot 1_{x^\m}\ot v_y\ot u_{y^\m}\ot \om_y\ot 1_{y^\m}\\
		& \mapsto & \dsum v_x\ot u_{x^\m}\ot \om_x\ot \bar{u}_{x^\m}\ot \bar{u}_x\ot v_y\ot u_{y^\m}\ot \om_y\ot \bar{u}_{y^\m}\ot \bar{u}_y\\
		& \stackrel{T_1,T_2}{\mapsto} & \dsum \widetilde{\om}_x^{i,p}\ot\widetilde{u}_{x^\m}^{i,p}\ot v_x\ot u_{x^\m}f_i^x(\om_x\ot \bar{u}_{x^\m})\ot \bar{u}_x\ot f_l(v_y\ot u_{y^\m})\om_y\ot \bar{u}_{y^\m}\ot v_y^{l,t}\ot u_{y^\m}^{l,t}\ot  \bar{u}_y\\
		&\mapsto & \dsum  \widetilde{\om}_x^{i,p}\ot\widetilde{u}_{x^\m}^{i,p}\ot v_x(u_{x^\m}f_i^x(\om_x\ot \bar{u}_{x^\m})\stackrel{\T}{\circ} \bar{u}_x)\ot f_l(v_y\ot u_{y^\m})\om_y\ot \bar{u}_{y^\m}\ot v_y^{l,t}(u_{y^\m}^{l,t}\stackrel{\T}{\circ}  \bar{u}_y)\\
		&\stackrel{T_3}{\mapsto} & \dsum \widetilde{\om}_x^{i,p}\ot\widehat{\om}_y^{j,q}\ot \widehat{u}_{y^\m}^{j,q}\ot \widetilde{u}_{x^\m}^{i,p}\ot v_x(u_{x^\m}f_i^x(\om_x\ot \bar{u}_{x^\m})\stackrel{\T}{\circ} \bar{u}_x)\underbrace{f_j^y(f_l(v_y\ot u_{y^\m})\om_y\ot \bar{u}_{y^\m})}_{b_{l,j}\in R}\ot v_y^{l,t} (u_{y^\m}^{l,t}\stackrel{\T}{\circ}  \bar{u}_y).	\\
		& = &  \dsum \widetilde{\om}_x^{i,p}\ot\widehat{\om}_y^{j,q}\ot \widehat{u}_{y^\m}^{j,q}\ot \widetilde{u}_{x^\m}^{i,p}\ot v_x(u_{x^\m}f_i^x(\om_x\ot \bar{u}_{x^\m})\stackrel{\T}{\circ} \bar{u}_x)b_{l,j}\ot v_y^{l,t} (u_{y^\m}^{l,t}\stackrel{\T}{\circ}  \bar{u}_y)\\
		& \mapsto & \dsum (\widetilde{\om}_x^{i,p}\stackrel{\Om}{\circ}\widehat{\om}_y^{j,q})\ot (\widehat{u}_{y^\m}^{j,q}\stackrel{\T}{\circ} \widetilde{u}_{x^\m}^{i,p})\ot \left( v_x(u_{x^\m}f_i^x(\om_x\ot \bar{u}_{x^\m})\stackrel{\T}{\circ} \bar{u}_x)b_{l,j}\stackrel{\Gm}{\circ} v_y^{l,t} (u_{y^\m}^{l,t}\stackrel{\T}{\circ}  \bar{u}_y)\right),
	\end{eqnarray*}
where $b_{l,j}=f_j^y(f_l(v_y\ot u_{y^\m})\om_y\ot \bar{u}_{y^\m}) \in R$.

Now, by diagram (\ref{diagramacomutatividade2})	and using  Remark~\ref{obsinversodefxxm}, we have that $1_x\lambda_{xy}\circ f_{y,x}^{\Gamma\Omega}$ is given by:
	\begin{eqnarray*}
		& & v_x\ot u_{x^\m}\ot \om_x\ot v_y\ot u_{y^\m}\ot \om_y \\
		& \stackrel{T_4}{\mapsto}&  \dsum_{l,t} v_x\ot v_y^{l,t}\ot u_{y^\m}^{l,t} \ot u_{x^\m}\ot \om_xf_l(v_y\ot u_{y^\m})\ot \om_y\\
		& \mapsto & \dsum_{l,t} (v_x\stackrel{\Gamma}{\circ} v_y^{l,t})\ot  (u_{y^\m}^{l,t} \stackrel{\T}{\circ} u_{x^\m})\ot (\om_xf_l(v_y\ot u_{y^\m})\stackrel{\Om}{\circ} \om_y)\\
		& \mapsto & \dsum_{l,t} (v_x\stackrel{\Gamma}{\circ} v_y^{l,t})\ot  (u_{y^\m}^{l,t} \stackrel{\T}{\circ} u_{x^\m})\ot (\om_xf_l(v_y\ot u_{y^\m})\stackrel{\Om}{\circ} \om_y)\ot 1_{y^\m}1_{(xy)^\m}\\
		& \mapsto & \dsum_{l,t} (v_x\stackrel{\Gamma}{\circ} v_y^{l,t})\ot  (u_{y^\m}^{l,t} \stackrel{\T}{\circ} u_{x^\m})\ot (\om_xf_l(v_y\ot u_{y^\m})\stackrel{\Om}{\circ} \om_y)\ot (\bar{u}_{y^\m}\stackrel{\T}{\circ} \bar{u}_{x^\m})\ot (\bar{u}_{x}\stackrel{\T}{\circ} \bar{u}_{y})\\
		%
		%
		%
		& \stackrel{T_5}{\mapsto} & \dsum (\widetilde{\om}_x^{i,p}\stackrel{\Om}{\circ}\widehat{\om}_y^{j,q})\ot (\widehat{u}_{y^\m}^{j,q}\stackrel{\T}{\circ} \widetilde{u}_{x^\m}^{i,p})\ot (v_x\stackrel{\Gamma}{\circ} v_y^{l,t})\ot  (u_{y^\m}^{l,t} \stackrel{\T}{\circ} u_{x^\m})f_i^x(\om_x\underbrace{f_l(v_y\ot u_{y^\m})f_j^y(\om_y\ot \bar{u}_{y^\m})}_{=b_{l,j} \in R}\ot \bar{u}_{x^\m})\\
		& & \ot (\bar{u}_{x}\stackrel{\T}{\circ} \bar{u}_{y})\\ 
		& = & \dsum (\widetilde{\om}_x^{i,p}\stackrel{\Om}{\circ}\widehat{\om}_y^{j,q})\ot (\widehat{u}_{y^\m}^{j,q}\stackrel{\T}{\circ} \widetilde{u}_{x^\m}^{i,p})\ot (v_x\stackrel{\Gamma}{\circ} v_y^{l,t})\ot  (u_{y^\m}^{l,t} \stackrel{\T}{\circ} u_{x^\m})f_i^x(\om_xb_{l,j}\ot \bar{u}_{x^\m}) \ot (\bar{u}_{x}\stackrel{\T}{\circ} \bar{u}_{y})\\
		& = & \dsum (\widetilde{\om}_x^{i,p}\stackrel{\Om}{\circ}\widehat{\om}_y^{j,q})\ot (\widehat{u}_{y^\m}^{j,q}\stackrel{\T}{\circ} \widetilde{u}_{x^\m}^{i,p})\ot (v_x\stackrel{\Gamma}{\circ} v_y^{l,t})\ot  (u_{y^\m}^{l,t} \stackrel{\T}{\circ} u_{x^\m}f_i^x(\om_xb_{l,j}\ot \bar{u}_{x^\m})) \ot (\bar{u}_{x}\stackrel{\T}{\circ} \bar{u}_{y})\\
		& \mapsto & \dsum (\widetilde{\om}_x^{i,p}\stackrel{\Om}{\circ}\widehat{\om}_y^{j,q})\ot (\widehat{u}_{y^\m}^{j,q}\stackrel{\T}{\circ} \widetilde{u}_{x^\m}^{i,p})\ot (v_x\stackrel{\Gamma}{\circ} v_y^{l,t})\left( (u_{y^\m}^{l,t} \stackrel{\T}{\circ} u_{x^\m}f_i^x(\om_xb_{l,j}\ot \bar{u}_{x^\m})) \stackrel{\T}{\circ} (\bar{u}_{x}\stackrel{\T}{\circ} \bar{u}_{y})\right) \\
		& = & \dsum (\widetilde{\om}_x^{i,p}\stackrel{\Om}{\circ}\widehat{\om}_y^{j,q})\ot (\widehat{u}_{y^\m}^{j,q}\stackrel{\T}{\circ} \widetilde{u}_{x^\m}^{i,p})\ot (v_x\stackrel{\Gamma}{\circ} v_y^{l,t})\left( u_{y^\m}^{l,t} \stackrel{\T}{\circ} \underbrace{(u_{x^\m}f_i^x(\om_xb_{l,j}\ot \bar{u}_{x^\m}) \stackrel{\T}{\circ} \bar{u}_{x})}_{=c_{i,j,l }\in R}\bar{u}_{y}\right)\\
		& = & \dsum (\widetilde{\om}_x^{i,p}\stackrel{\Om}{\circ}\widehat{\om}_y^{j,q})\ot (\widehat{u}_{y^\m}^{j,q}\stackrel{\T}{\circ} \widetilde{u}_{x^\m}^{i,p})\ot (v_x\stackrel{\Gamma}{\circ} v_y^{l,t})\left( u_{y^\m}^{l,t} \stackrel{\T}{\circ} c_{i,j,l}\bar{u}_{y}\right)\\
		%
		& = & \dsum (\widetilde{\om}_x^{i,p}\stackrel{\Om}{\circ}\widehat{\om}_y^{j,q})\ot (\widehat{u}_{y^\m}^{j,q}\stackrel{\T}{\circ} \widetilde{u}_{x^\m}^{i,p})\ot \underbrace{\left( v_x\stackrel{\Gamma}{\circ} v_y^{l,t}(u_{y^\m}^{l,t} \stackrel{\T}{\circ} c_{i,j,l}\bar{u}_{y})\right)}_{(*)}.\\
	\end{eqnarray*}

Note that since $g_l$ is $R$-bilinear, then $r\dsum_{t}v_y^{l,t}\ot u_{y^\m}^{l,t}=\dsum_{t}v_y^{l,t}\ot u_{y^\m}^{l,t}r$, for all $r \in R$. Thus, $(*)$ can be written as:
\begin{eqnarray}
	\dsum_{t}\left( v_x\stackrel{\Gamma}{\circ} v_y^{l,t}(u_{y^\m}^{l,t} \stackrel{\T}{\circ} c_{i,j,l}\bar{u}_{y})\right) & = & \dsum_{t} (f_{x,y}^\Gm\circ(\Gm_x\ot \Gm_y\ot f_{y^\m,y}^\T))(v_x\ot v_y^{l,t}\ot u_{y^\m}^{l,t}c_{i,j,l}\ot \bar{u}_y) \nonumber\\
	& = & \dsum_{t} (f_{x,y}^\Gm\circ(\Gm_x\ot \Gm_y\ot f_{y^\m,y}^\T))(v_x\ot c_{i,j,l}v_y^{l,t}\ot u_{y^\m}^{l,t}\ot \bar{u}_y)\nonumber\\
	& = & \dsum_{t}\left( v_x\stackrel{\Gamma}{\circ}c_{i,j,l} v_y^{l,t}(u_{y^\m}^{l,t} \stackrel{\T}{\circ} \bar{u}_{y})\right). \label{cijldooutrolado}
\end{eqnarray}
Analogously, since $\dsum_{(x)}(\bar{u}_{x^\m}\stackrel{\T}{\circ}\bar{u}_x)=1_{x^\m}\in \Z$, then 
\begin{eqnarray}
c_{i,j,l}& = & \dsum_{(x)} f^\T_{x^\m,x}\circ (\T_x\ot f_i^x\ot\T_{x^\m})(u_{x^\m}\ot \om_xb_{l,j}\ot \bar{u}_{x^\m}\ot \bar{u}_x) \nonumber\\
& = & \dsum_{(x)} f^\T_{x^\m,x}\circ (\T_x\ot f_i^x\ot\T_{x^\m})(u_{x^\m}\ot \om_x\ot  \bar{u}_{x^\m}\ot \bar{u}_xb_{l,j})\nonumber\\
& = & \dsum_{(x)}(u_{x^\m}f_i^x(\om_x\ot \bar{u}_{x^\m}) \stackrel{\T}{\circ} \bar{u}_{x}b_{l,j}). \label{cijl}
\end{eqnarray}
Thus, by (\ref{cijldooutrolado}) and (\ref{cijl}) we have
\begin{eqnarray*}
& & (1_x\lambda_{xy}\circ f_{y,x}^{\Gamma\Omega})(v_x\ot u_{x^\m}\ot \om_x\ot v_y\ot u_{y^\m}\ot \om_y)\nonumber\\
& = & \dsum (\widetilde{\om}_x^{i,p}\stackrel{\Om}{\circ}\widehat{\om}_y^{j,q})\ot (\widehat{u}_{y^\m}^{j,q}\stackrel{\T}{\circ} \widetilde{u}_{x^\m}^{i,p})\ot \left( v_x\stackrel{\Gamma}{\circ} c_{i,j,l}v_y^{l,t}(u_{y^\m}^{l,t} \stackrel{\T}{\circ} \bar{u}_{y})\right)\nonumber\\
& = & \dsum (\widetilde{\om}_x^{i,p}\stackrel{\Om}{\circ}\widehat{\om}_y^{j,q})\ot (\widehat{u}_{y^\m}^{j,q}\stackrel{\T}{\circ} \widetilde{u}_{x^\m}^{i,p})\ot \left( v_x(u_{x^\m}f_i^x(\om_x\ot \bar{u}_{x^\m})\stackrel{\T}{\circ}\bar{u}_x)b_{l,j}\stackrel{\Gamma}{\circ}v_y^{l,t}(u_{y^\m}^{l,t} \stackrel{\T}{\circ} \bar{u}_{y})\right)\label{comutatvidadelado1}
\end{eqnarray*}

Therefore, $1_x\lambda_{xy}\circ f_{y,x}^{\Gamma\Omega}=f_{x,y}^{\Om\Gamma}\circ(\lambda_x\ot \lambda_y)$. 
 Consequently, the diagram (\ref{comutatividadediagrama}) is commutative and
	$$[\D(\Om)][\D(\Gm)]=[\D(\Gm)][\D(\Om)] \ \ \mbox{in} \ \C(\T/R).$$
	\end{dem}

\begin{pro}\label{pro:subgrC_0} The set
	$$\C_0(\T/R)=\{[\D(\Gamma)]; \Gamma_x\simeq \T_x, \ \mbox{for all } \ x\in G\},$$ 
	is a subgroup of $\C(\T/R)$.
\end{pro}
\begin{dem}
	If $[\D(\Gamma)] \in \C_0(\T/R)$, then $\Gamma_x \simeq \T_x$, for all $x \in G$. In particular, $\Gamma_x|\T_x$ and 
	\begin{equation*}
		\Gamma_x\ot\Gamma_{x^\m}\simeq \T_x\ot\T_{x^{-1}}\simeq R1_x, \ \ \mbox{for all} \ x \in G.
	\end{equation*}
	Then, $[\D(\Gamma)]\in \C(\T/R)$. Given $[\D(\Gamma)],[\D(\Om)]\in \C_0(\T/R)$, we have
	\begin{equation*}
		\Gamma_x\ot\T_{x^{-1}}\ot\Om_x\simeq \T_x\ot\T_{x^{-1}}\ot\T_x\simeq \T_x, \ \ \mbox{for all} \ x \in G.
	\end{equation*}
	Thus, $[\D(\Gamma)][\D(\Om)]=\left[ \bigoplus_{x \in G}\Gamma_x\ot\T_{x^{-1}}\ot\Om_x\right] \in \C_0(\T/R)$. Analogously,  $[\D(\Gamma)]\in \C_0(\T/R)$, then
	\begin{equation*}
		\T_x\ot\Gamma_{x^\m}\ot\T_x\simeq \T_x\ot\T_{x^\m}\ot\T_x\simeq \T_x, \ \mbox{for all}\  x \in G.
	\end{equation*}
	Then, $[\D(\Gamma)]^{-1}=\left[ \bigoplus_{x\in G}\T_x\ot\Gamma_{x^\m}\ot\T_x\right] \in \C_0(\T/R)$. Therefore, $\C_0(\T/R)$ is a subgroup of $\C(\T/R)$.
	
\end{dem}


\begin{lem}\label{tautilehcociclo} Let $[\D(\Gm)]\in \C_0(\T/R)$ and  $a_x:\Gamma_x\longrightarrow \T_x$ be an $R$-bimodule isomorphism, for each $x \in G$. Let  $\tau_{x,y}:1_x\T_{xy}\longrightarrow 1_x\T_{xy}$ be  the $R$-bimodule isomorphism defined by the commutative diagram
	$$\xymatrix{  \Gamma_x\ot\Gamma_y\ar[rr]^{f_{x,y}^\Gamma}\ar[dd]_{a_x\ot a_y}  & & 1_x\Gamma_{xy}\ar[dd]^{a_{xy}}\\
		& & \\
		\T_x\ot\T_y\ar[rd]_{f_{x,y}^\T} & & 1_x\T_{xy}\\
		& 1_x\T_{xy}\ar@{-->}[ru]_{\tau_{x,y}}& 	  }$$
	that is,
	\begin{equation*}
		\tau_{x,y}\circ f_{x,y}^\T\circ (a_x\ot a_y)= a_{xy}\circ f_{x,y}^\Gm,  \ x,y \in G.
	\end{equation*}
	Then, $\widetilde{\tau}_{-,-}$ is a normalized element in $Z^2_{\T}(G,\al,\Z)$, where $\widetilde{\tau}_{x,y}$ is defined in Lemma \ref{gammatil}. 
\end{lem}
\begin{dem}	By Lemma \ref{gammatil} we have $\widetilde{\tau}_{x,y}\in \U(\Z1_x1_{xy})$ and
	\begin{equation}
	\widetilde{\tau}_{x,y}(a_x(v_x)\stackrel{\T}{\circ}a_y(v_y))=a_{xy}(v_x\stackrel{\Gamma}{\circ}v_y), \label{tautilcomfTefGamma}
	\end{equation}
	for all $v_x \in \Gm_x$ and $v_y \in \Gm_y$.
	
	For $x,y,z \in G$,  $v_x \in \Gamma_x, v_y \in \Gamma_y$ and $v_z\in \Gamma_z,$ we obtain
	\begin{eqnarray*}
		a_{xyz}((v_x\stackrel{\Gamma}{\circ}v_y)\stackrel{\Gamma}{\circ}v_z) & \stackrel{(\ref{tautilcomfTefGamma})}{=} & \widetilde{\tau}_{xy,z}(a_{xy}(v_x\stackrel{\Gamma}{\circ}v_y)\stackrel{\T}{\circ}a_z(v_z))\\
		& \stackrel{(\ref{tautilcomfTefGamma})}{=} & \widetilde{\tau}_{xy,z}\widetilde{\tau}_{x,y}((a_{x}(v_x)\stackrel{\T}{\circ}a_y(v_y))\stackrel{\T}{\circ}a_z(v_z))\\
		& = & \widetilde{\tau}_{xy,z}\widetilde{\tau}_{x,y}(a_{x}(v_x)\stackrel{\T}{\circ}(a_y(v_y)\stackrel{\T}{\circ}a_z(v_z)))\\
		& \stackrel{(\ref{tautilcomfTefGamma})}{=} &\widetilde{\tau}_{xy,z}\widetilde{\tau}_{x,y}(a_{x}(v_x)\stackrel{\T}{\circ}\widetilde{\tau}_{y,z}^\m a_{yz}(v_y\stackrel{\Gamma}{\circ}v_z))\\
		& = & \widetilde{\tau}_{xy,z}\widetilde{\tau}_{x,y}(a_{x}(v_x)\widetilde{\tau}_{y,z}^\m \stackrel{\T}{\circ}a_{yz}(v_y\stackrel{\Gamma}{\circ}v_z))\\
		& \stackrel{(\ref{uxercomutam})}{=}& \widetilde{\tau}_{xy,z}\widetilde{\tau}_{x,y}\al_{x}(\widetilde{\tau}_{y,z}^\m1_{x^\m})(a_{x}(v_x)\stackrel{\T}{\circ}a_{yz}(v_y\stackrel{\Gamma}{\circ}v_z))\\
		& \stackrel{(\ref{tautilcomfTefGamma})}{=} & \widetilde{\tau}_{xy,z}\widetilde{\tau}_{x,y}\al_{x}(\widetilde{\tau}_{y,z}^\m1_{x^\m})\widetilde{\tau}_{x,yz}^\m a_{xyz}(v_x\stackrel{\Gamma}{\circ}(v_y\stackrel{\Gamma}{\circ}v_z))\\
		& =& \widetilde{\tau}_{xy,z}\widetilde{\tau}_{x,y}\al_{x}(\widetilde{\tau}_{y,z}^\m1_{x^\m})\widetilde{\tau}_{x,yz}^\m a_{xyz}((v_x\stackrel{\Gamma}{\circ}v_y)\stackrel{\Gamma}{\circ}v_z)\\
		& =&  a_{xyz}(\widetilde{\tau}_{xy,z}\widetilde{\tau}_{x,y}\al_{x}(\widetilde{\tau}_{y,z}^\m1_{x^\m})\widetilde{\tau}_{x,yz}^\m(v_x\stackrel{\Gamma}{\circ}v_y)\stackrel{\Gamma}{\circ}v_z).
	\end{eqnarray*}
	Since $a_{xyz}$ is an $R$-bimodule isomorphism, we see that 
	\begin{equation}
	((v_x\stackrel{\Gamma}{\circ}v_y)\stackrel{\Gamma}{\circ}v_z) = \widetilde{\tau}_{xy,z}\widetilde{\tau}_{x,y}\al_{x}(\widetilde{\tau}_{y,z}^\m1_{x^\m})\widetilde{\tau}_{x,yz}^\m((v_x\stackrel{\Gamma}{\circ}v_y)\stackrel{\Gamma}{\circ}v_z), \label{vxvyvztaus}
	\end{equation}
for all $v_x \in \Gamma_x, v_y \in \Gamma_y$ and $v_z\in \Gamma_z$. 
{	Using the same argument as in Proposition \ref{3cocicloquecorrigeaassociatividade} we obtain}
		\begin{equation*}
		 \widetilde{\tau}_{xy,z}\widetilde{\tau}_{x,y}\al_{x}(\widetilde{\tau}_{y,z}^\m1_{x^\m})\widetilde{\tau}_{x,yz}^\m1_x=1_x1_{xy}1_{xyz}.
		\end{equation*}
	This implies that $\widetilde{\tau}^\m_{-,-} \in Z^2_\T(G,\al,\Z)$ and consequently  $\widetilde{\tau}_{-,-} \in Z^2_\T(G,\al,\Z)$. For $x=1$, by (\ref{produtoporr}) and (\ref{tautilcomfTefGamma}) we have $\widetilde{\tau}_{1,y}(ra_y(v_y))=a_y(rv_y)$ for all $r \in R, v_y \in \Gm_y$. Take $r=1.$ Since $a_y$ is an $R$-bimodule isomorphism, 
then $\widetilde{\tau}_{1,y}v_y=v_y,$ for all $v_y \in \Gm_y.$
	The same argument as above implies that $\widetilde{\tau}_{1,y}=1_y.$ Analogously, $\widetilde{\tau}_{x,1}=1_x$. Therefore, $\widetilde{\tau}_{-,-}$ is normalized. 
	 \end{dem}

\begin{teo}\label{C0isomorfoH2} The map
	$$\begin{array}{c c c l}
	\zeta: & \C_0(\T/R) & \longrightarrow & H^2_\T(G,\al,\Z)\\
	& [\D(\Gm)] & \longmapsto & [\widetilde{\tau}_{-,-}],
	\end{array}$$
	where $\widetilde{\tau}_{-,-} \in Z^2(G,\al,\Z)$ is defined in Lemma \ref{tautilehcociclo}, is a group isomorphism.
\end{teo}
\begin{dem} 
	Let $[\D(\Gm)],[\D(\Om)] \in \C_0(\T/R)$ and take $R$-bimodule isomorphisms
	$a_x:\Gamma_x\longrightarrow\T_x$ and $b_x:\Om_x\longrightarrow \T_x$, for $x\in G$. Denote,
	\begin{equation}
	\tau_{x,y} \circ f_{x,y}^\T\circ (a_x\ot a_y)=a_{xy}\circ f_{x,y}^\Gm, \ \ \forall \ x,y\in G, \label{definicaodetau}
	\end{equation}
	\begin{equation}
	\gamma_{x,y} \circ f_{x,y}^\T\circ (b_x\ot b_y)=b_{xy}\circ f_{x,y}^\Om, \ \ \forall \ x,y\in G. \label{definicaodegamma}
	\end{equation}
	
	Let us first show    that  $\zeta$ is well-defined. Suppose that  $[\D(\Gm)]=[\D(\Om)]$ in $\C_0(\T/R)$. Then there exist  $R$-bimodule isomorphisms $\xi_x:\Gm_x\longrightarrow \Om_x$,  $x \in  G$, such that the following diagram is commutative:
	\begin{equation}
	\xymatrix{ \Gm_x\ot \Gm_y\ar[rr]^{f_{x,y}^\Gm}\ar[dd]_{\xi_x\ot \xi_{y}} & & 1_x\Gm_{xy}\ar[dd]^{\xi_{xy}}\\
		&  & \\
		\Om_x\ot\Om_y\ar[rr]_{f_{x,y}^\Om} & & 1_x\Om_{xy}. }\label{diagramaisomorfismoentreDGammaeDOmega}
	\end{equation}
	For each $x \in G$ consider the  $R$-bimodule isomorphism $\beta_x:\T_x\rightarrow\T_x$ defined by
	\begin{equation}
	\beta_x=b_x\circ \xi_x\circ a_{x}^\m, \ \ \forall \ x\in G. \label{definicaodebetaemC0isoH2}
	\end{equation}
	\begin{afr}
		$\beta_{xy}\circ \tau_{x,y}\circ f_{x,y}^\T=\gamma_{x,y}\circ f_{x,y}^\T\circ (\beta_x \ot \beta_y)$, for all $x,y \in G$. 
	\end{afr} 	
	Indeed, by (\ref{definicaodetau}), (\ref{definicaodegamma}) and by the commutative diagram (\ref{diagramaisomorfismoentreDGammaeDOmega}), we obtain
	\begin{eqnarray*}
		\gamma_{x,y}\circ f_{x,y}^\T\circ (\beta_x \ot \beta_y) & = & \gamma_{x,y}\circ f_{x,y}^\T\circ ((b_x\circ \xi_x\circ a_{x}^\m) \ot (b_y\circ \xi_y\circ a_{y}^\m))\\
		& = & \gamma_{x,y}\circ f_{x,y}^\T\circ (b_x\ot b_y)\circ (\xi_x\ot\xi_y)\circ  (a_{x}^\m \ot  a_{y}^\m)\\
		& \stackrel{(\ref{definicaodegamma})}{=} & b_{xy}\circ f_{x,y}^{\Om}\circ (\xi_x\ot\xi_y)\circ  (a_{x}^\m \ot  a_{y}^\m) \\
		& = & b_{xy}\circ \xi_{xy}\circ f_{x,y}^\Gm\circ (a_x^\m\ot a_y^\m)\\
		& \stackrel{(\ref{definicaodetau})}{=} &  b_{xy}\circ \xi_{xy}\circ a_{xy}^\m\circ \tau_{x,y}\circ f_{x,y}^\T\\
		& = &  \beta_{xy}\circ \tau_{x,y}\circ f_{x,y}^\T.
	\end{eqnarray*}
	Then, for each $u_x \in \T_x$, $u_y \in \T_y$ we have
	\begin{equation*}
	\beta_{xy}(\tau_{x,y}(u_x\stackrel{\T}{\circ} u_{y}))=\gamma_{x,y}(\beta_x(u_x)\stackrel{\T}{\circ} \beta_y(u_{y})).
	\end{equation*}
	On  one hand, by (\ref{relacaodefcomftil}), we see that
	\begin{equation*}
	\beta_{xy}(\tau_{x,y}(u_x\stackrel{\T}{\circ} u))=\beta_{xy}(\widetilde{\tau}_{x,y}(u_x\stackrel{\T}{\circ} u_{y}))=\widetilde{\beta}_{xy}\widetilde{\tau}_{x,y}(u_x\stackrel{\T}{\circ} u_{y}).
	\end{equation*}
	On the other hand, using (\ref{relacaodefcomftil}) and (\ref{uxercomutam}), we obtain
	\begin{eqnarray*}
		\gamma_{x,y}(\beta_x(u_x)\stackrel{\T}{\circ} \beta_y(u_{y}))& = &\widetilde{\gamma}_{x,y}(\beta_x(u_x)\stackrel{\T}{\circ} \beta_y(u_{y}))= \widetilde{\gamma}_{x,y}(\widetilde{\beta}_xu_x\stackrel{\T}{\circ} \widetilde{\beta}_yu_{y})\\
		& = &\widetilde{\gamma}_{x,y}(\widetilde{\beta}_xu_x\widetilde{\beta}_y\stackrel{\T}{\circ} u_{y}) =\widetilde{\gamma}_{x,y}(\widetilde{\beta}_x\al_x(\widetilde{\beta}_y1_{x^\m})u_x\stackrel{\T}{\circ} u_{y})\\
		& = & \widetilde{\gamma}_{x,y}\widetilde{\beta}_x\al_x(\widetilde{\beta}_y1_{x^\m})(u_x\stackrel{\T}{\circ}u_{y}).
	\end{eqnarray*}
	Then,
	\begin{equation*}
	\widetilde{\beta}_{xy}\widetilde{\tau}_{x,y}(u_x\stackrel{\T}{\circ} u_{y})=\widetilde{\gamma}_{x,y}\widetilde{\beta}_x\al_x(\widetilde{\beta}_y1_{x^\m})(u_x\stackrel{\T}{\circ}u_{y}), \label{betatau}
	\end{equation*}
	for all $u_x \in \T_x$ and $u_y \in \T_y$. 
Using the argument in Lemma \ref{tautilehcociclo}, we have
	\begin{equation}
	\widetilde{\beta}_{xy}\widetilde{\tau}_{x,y}=\widetilde{\gamma}_{x,y}\widetilde{\beta}_x\al_x(\widetilde{\beta}_y1_{x^\m}), \ \ \mbox{for all } \ x,y \in G. \label{c0isoH2}
	\end{equation}
	
%
%
	Let
	$$\begin{array}{c c c l}
	h: & G & \longrightarrow & \Z\\
	& x & \longmapsto & \widetilde{\beta}_x.
	\end{array}$$
	Then $h_x \in \U(\Z1_x)$, for all $x\in G$ and we get
	\begin{equation*}
	\widetilde{\tau}_{x,y}=\widetilde{\gamma}_{x,y}\al_x(h_y1_{x^\m})h_{xy}^\m h_x=\widetilde{\gamma}_{x,y}(\delta^1h)(x,y), \ \ \forall \ x,y \in G.
	\end{equation*}
	This implies that $[\widetilde{\tau}_{-,-}]=[\widetilde{\gamma}_{-,-}]$ in $H^2_\T(G,\al,\Z),$ proving that   $\zeta$ is well-defined.

	To show that $\zeta$ is injective, suppose that $\zeta([\D(\Gm)])=\zeta([\D(\Om)])$ in $H^2_\T(G,\al,\Z)$. Then there exist $h:G\longrightarrow \Z$, with $h(x)=h_x \in \U(\Z1_x)$, for all $x\in G$, and
	\begin{equation*}
	\widetilde{\tau}_{x,y}=\widetilde{\gamma}_{x,y}\al_x(h_y1_{x^\m})h_{xy}^\m h_x, \ \ \forall \ x,y \in G.
	\end{equation*}  
	Consider the map
	$$\begin{array}{c c c l}
	\beta_x: & \T_x & \longrightarrow & \T_x\\
	& u_x & \longmapsto & h_xu_x.
	\end{array}$$
	Then $\beta_x$ { is} an $R$-bimodule isomorphism and $\widetilde{\beta}_x=h_x$, for all $x\in G$ (see Lemma \ref{gammatil}).
	Hence, 
	\begin{equation}
	\widetilde{\tau}_{x,y}\widetilde{\beta}_{xy}=\widetilde{\gamma}_{x,y}\al_x(\widetilde{\beta}_y1_{x^\m})\widetilde{\beta}_y, \ \ \forall \ x,y \in G. \label{hemC0isoH2}
	\end{equation}
	Let $\lambda_x:\Gamma_x\rightarrow \Om_x$, $(x \in G )$, be  the isomorphism defined by 
	\begin{equation*}
	\lambda_x=b_x^\m\circ \beta_x\circ a_x.
	\end{equation*}
	We will verify  that following diagram is commutative:
	\begin{equation}
	\xymatrix{ \Gm_x\ot_R \Gm_y\ar[rr]^{f_{x,y}^\Gm}\ar[dd]_{\lambda_x\ot \lambda_y} & & 1_x\Gm_{xy}\ar[dd]^{\lambda_{xy}}\\
		&  & \\
		\Om_x\ot_R\Om_y\ar[rr]_{f_{x,y}^\Om} & & 1_x\Om_{xy}. } \label{diagramaC0isomorfo}
	\end{equation}
	Firstly, we check that
	\begin{equation}
		\beta_{xy}\circ \tau_{x,y}\circ f_{x,y}^\T=\gamma_{x,y}\circ f_{x,y}^\T\circ (\beta_x\ot \beta_y), \ \mbox{for all}\  x,y \in G.  \label{obs2}
	\end{equation} 
	Indeed, for $u_x \in \T_x$ and $u_y \in \T_y$,  by (\ref{relacaodefcomftil}) we see that 
	\begin{equation*}
	\beta_{xy}(\tau_{x,y}(u_x\stackrel{\T}{\circ}u_y))= \widetilde{\beta}_{x,y}\widetilde{\tau}_{x,y}(u_x\stackrel{\T}{\circ}u_y). \label{qualquercoisa}
	\end{equation*}
Analogously, by (\ref{relacaodefcomftil}), (\ref{uxercomutam}) and (\ref{hemC0isoH2}) it follows that
	\begin{eqnarray*}
		\gamma_{x,y}(\beta_x(u_x)\stackrel{\T}{\circ}\beta_y(u_y)) & = & \widetilde{\gamma}_{x,y}(\widetilde{\beta}_xu_x\stackrel{\T}{\circ}\widetilde{\beta}_yu_y)= \widetilde{\gamma}_{x,y}(\widetilde{\beta}_xu_x\widetilde{\beta}_y\stackrel{\T}{\circ}u_y)= \widetilde{\gamma}_{x,y}(\widetilde{\beta}_x\al_x(\widetilde{\beta}_y1_{x^\m})u_x\stackrel{\T}{\circ}u_y)\\
		& = & \widetilde{\gamma}_{x,y}\widetilde{\beta}_x\al_x(\widetilde{\beta}_y1_{x^\m})(u_x\stackrel{\T}{\circ}u_y)=\widetilde{\beta}_{xy}\widetilde{\tau}_{x,y}(u_x\stackrel{\T}{\circ}u_y)
	\end{eqnarray*}
Consequently, (\ref{obs2}) holds. 	
	
	Using (\ref{definicaodetau}), (\ref{definicaodegamma}) and 
	 (\ref{obs2}),  we have
	\begin{eqnarray*}
		\lambda_{xy}\circ f_{x,y}^\Gm & = & b_{xy}^\m\circ \beta_{xy}\circ a_{xy}\circ f_{x,y}^\Gm\\
		& \stackrel{(\ref{definicaodetau})}{=} &  b_{xy}^\m\circ \beta_{xy}\circ \tau_{x,y}\circ f_{x,y}^\T\circ (a_x\ot a_y)\\
		& \stackrel{(\ref{obs2})}{=} & b_{xy}^\m\circ \gamma_{xy}\circ f_{x,y}^\T\circ (\beta_x\ot \beta_y)\circ (a_x\ot a_y)\\
		& \stackrel{(\ref{definicaodegamma})}{=} & b_{xy}^\m\circ b_{xy}\circ f_{x,y}^\Om\circ (b_x^\m\ot b_y^\m)\circ (\beta_x\ot \beta_y)\circ (a_x\ot a_y)\\
		& = & f_{x,y}^\Om\ot (\lambda_x\ot \lambda_y).
	\end{eqnarray*}
This implies that the diagram (\ref{diagramaC0isomorfo}) is commutative. Thus, we have a partial generalized crossed product isomorphism. Therefore,  $[\D(\Gm)]=[\D(\Om)]$ in $\C_0(\T/R)$ and  $\zeta$ is injective. 
	
	Let $\sigma:G\times G\longrightarrow \Z$ be a normalized $2$-cocycle in $Z^2_\T(G,\al,\Z) .$ Then 
	\begin{equation}
	\al_x(\sigma_{y,z}1_{x^\m})\sigma_{x,yz}=\sigma_{xy,z}\sigma_{x,y}, \ \forall \ x,y,z \in G.\label{sigmacociclo}
	\end{equation}
	Consider the $R$-bimodule isomorphism
	\begin{equation*}
	\begin{array}{c c c l}
	\rho_x: & \T_x & \longrightarrow & \T_x\\
	& u_x & \longmapsto & \sigma_{x,x}u_x,
	\end{array}
	\end{equation*}
	and  define $\Sigma_x=\T_x$, as $R$-bimodules. Then, $\rho_x:\Sigma_x\longrightarrow \T_x$ is  an $R$-bimodule isomorphism. Let 
	\begin{equation*}
	\begin{array}{c c c l}
	f_{x,y}^\Sigma: & \Sigma_x\ot_R \Sigma_y & \longrightarrow & 1_x\Sigma_{xy}\\
	& u_x\ot u_y & \longmapsto & \sigma_{x,y}(u_x\stackrel{\T}{\circ}u_y).
	\end{array}
	\end{equation*}
Given $u_x \in \Sigma_x,u_y\in \Sigma_y$ and $u_z\in \Sigma_z$, by (\ref{sigmacociclo}) we obtain that
	\begin{eqnarray*}
		f_{x,yz}^{\Sigma}(u_x\ot f_{y,z}^\Sigma(u_y\ot u_z)) & = & \s_{x,yz}(u_x\stackrel{\T}{\circ}\s_{y,z}(u_y\stackrel{\T}{\circ}u_z))\\
		& =&  \s_{x,yz}\al_x(\s_{y,z}1_{x^\m})(u_x\stackrel{\T}{\circ}(u_y\stackrel{\T}{\circ}u_z))\\
		& = &  \s_{xy,z}\s_{x,y}((u_x\stackrel{\T}{\circ}u_y)\stackrel{\T}{\circ}u_z)\\
		& = &  \s_{xy,z}(\s_{x,y}(u_x\stackrel{\T}{\circ}u_y)\stackrel{\T}{\circ}u_z)\\
		& = &  f_{xy,z}^\Sigma(f_{x,y}^\Sigma(u_x\ot u_y)\ot u_z).
	\end{eqnarray*}
	Thus, $f^\Sigma=\{f_{x,y}^\Sigma:\sigma_x\ot \Sigma_y\rightarrow 1_x\Sigma_{xy}, \ x,y \in G\}$ is a factor set for $\Sigma$ and $\D(\Sigma)=\bigoplus_{x}\Sigma_x$ is a partial generalized crossed product with $[\D(\Sigma)]\in \C_0(\T/R)$. Let $\lambda \in C^1(G,\al,\Z)$ be defined by $\lambda_x=\sigma_{x,x}^\m$. Consider the $R$-bimodule isomorphism
	$$\begin{array}{c c c l}
	\Upsilon: & 1_x\T_{xy} & \rightarrow & 1_x\T_{xy}\\
	 & u_x & \mapsto & \sigma_{x,y}(\delta^1\lambda)(x,y)u_{xy}.
	\end{array}$$
Then, $\tilde{\Upsilon}=\sigma_{x,y}(\delta^1\lambda)(x,y)$ and 
	\begin{eqnarray*}
		\rho_{xy}(f_{x,y}^\Sigma(u_x\ot u_y)) & = & \rho_{xy}(\sigma_{x,y}(u_x\stackrel{\T}{\circ} u_y))\\
		& = & \sigma_{x,y}\rho_{xy}(u_x\stackrel{\T}{\circ} u_y)\\
		& = & \sigma_{x,y}\sigma_{xy,xy}(u_x\stackrel{\T}{\circ} u_y)\\
		& = & \sigma_{x,y}\sigma_{xy,xy}(\sigma_{x,x}^\m\rho_x(u_x)\stackrel{\T}{\circ} \sigma_{y,y}^\m\rho_{y}(u_y))\\
		& = & \sigma_{x,y}\sigma_{xy,xy}(\sigma_{x,x}^\m\rho_x(u_x)\sigma_{y,y}^\m\stackrel{\T}{\circ} \rho_{y}(u_y))\\
		& = & \sigma_{x,y}\sigma_{xy,xy}(\sigma_{x,x}^\m\al_x(\sigma_{y,y}^\m1_{x^\m})\rho_x(u_x)\stackrel{\T}{\circ} \rho_{y}(u_y))\\
		& = & \sigma_{x,y}\sigma_{xy,xy}\sigma_{x,x}^\m\al_x(\sigma_{y,y}^\m1_{x^\m})(\rho_x(u_x)\stackrel{\T}{\circ} \rho_{y}(u_y))\\
		& = & \sigma_{x,y}(\delta^1\lambda)(x,y)f_{x,y}^\T(\rho_x(u_x)\ot \rho_{y}(u_y)).
	\end{eqnarray*}
This implies that $\rho_{x,y}\circ f_{x,y}^\Sigma=\Upsilon\circ f_{x,y}^\T\circ (\rho_x\ot \rho_y)$. Thus, $\zeta([\D(\Sigma)])=[\s(\delta^1\lambda)]=[\sigma_{-,-}]$ in $H^2_\T(G,\al,\Z),$ proving that   $\zeta$ is onto.
	
	Finally, we will show that $\zeta$ is a group homomorphism. Since $[\D(\Gm)][\D(\Om)]=\left[ \bigoplus_{x \in G}\Gm_x\ot \T_{x^{-1}}\ot \Om_x\right] $ in $\C_0(\T/R)$, there exists an  $R$-module homomorphism 
	\begin{equation*}
	\begin{array}{c c c l}
	d_x: & \Gm_x\ot\T_{x^\m}\ot\Om_x & \longrightarrow & \T_x \\
	& v_x\ot u_{x^\m}\ot \om_x & \longmapsto & a_x\left( v_x(u_{x^\m}\stackrel{\T}{\circ} b_x(\om_x))\right) .
	\end{array}
	\end{equation*}
	We will check that $d_{xy}\circ f_{x,y}^\Gm\Om=\sigma_{x,y}\gamma_{x,y}\circ f_{x,y}^\T\circ (d_x\ot d_y)$, for all $x,y \in G$.
	Observe that $d_{xy}\circ f_{x,y}^{\Gamma\Omega}=d_{xy}\circ (f_{x,y}^\Gamma\ot f_{y^\m,x^\m}^\T\ot f_{xy}^\Om)\circ (\Gamma_x\ot T\ot \Om_y)$, where $T=T_{\T_{x^{-1}}\ot \Om_x,\Gm_y\ot \T_{y^{-1}}}$.
	Consider the $R$-bimodule isomorphisms $$f_x:\T_{x^\m}\ot \Om_x\stackrel{\T_{x^\m}\ot b_x}{\longrightarrow} \T_{x^{-1}}\ot \T_x\longrightarrow R1_{x^\m}\stackrel{i}{\hookrightarrow} R,$$   
	$$g_x:R\stackrel{\pi}{\longrightarrow} R1_{x^\m}\longrightarrow\T_{x^\m}\ot \T_x\stackrel{\T_{x^\m}\ot b_x^\m}{\longrightarrow} \T_{x^\m}\ot \Om_x,$$
	where $i$ and $\pi$ are the canonical projection and injection, respectively. It is easy to see that $g_xf_x=Id_{\T_{x^\m}\ot \Om_x}$. Since
	$f_x(u_{x^\m}\ot \om_x)=(u_{x^\m}\stackrel{\T}{\circ}b_x(\om_x))$ and $g_x(1)=\dsum_{(x)} \bar{u}_{x^\m}\ot b_{x}^{\m}(\bar{u}_x)$, where $\dsum_{(x)}(\bar{u}_{x^\m}\stackrel{\T}{\circ} \bar{u}_x)=1_{x^\m},$ then
	\begin{equation*}
	T(u_{x^\m}\ot \om_x\ot v_y\ot u_{y^\m})=\dsum_{(x)}(u_{x^\m}\stackrel{\T}{\circ} b_x(\om_x))v_y\ot u_{y^\m}\ot \bar{u}_{x^\m}\ot b_{x}^{\m}(\bar{u}_x),
	\end{equation*}
	Thus, $d_{xy}\circ f_{x,y}^{\Gm\Om}$:
	\begin{eqnarray*}
		& & v_x\ot u_{x^\m}\ot \om_x\ot v_y\ot u_{y^\m}\ot \om_{y}\\
		& \stackrel{T}{\mapsto} & \dsum_{(x)}v_x \ot (u_{x^\m}\stackrel{\T}{\circ} b_x(\om_x))v_y\ot u_{y^\m}\ot \bar{u}_{x^\m}\ot b_{x}^{\m}(\bar{u}_x)\ot \om_y\\
		& \mapsto& \dsum_{(x)} (v_x(u_{x^\m}\stackrel{\T}{\circ} b_x(\om_x)) \stackrel{\Gm}{\circ} v_y)\ot (u_{y^\m}\stackrel{\T}{\circ} \bar{u}_{x^\m})\ot (b_{x}^{\m}(\bar{u}_x)\stackrel{\Om}{\circ} \om_y)\\
		& \stackrel{d_{xy}}{\mapsto}& \dsum_{(x)} a_{xy}\left(( v_x(u_{x^\m}\stackrel{\T}{\circ} b_x(\om_x)) \stackrel{\Gm}{\circ} v_y)((u_{y^\m}\stackrel{\T}{\circ} \bar{u}_{x^\m})\stackrel{\T}{\circ} b_{xy}(b_{x}^{\m}(\bar{u}_x)\stackrel{\Om}{\circ} \om_y))\right) \\
		& =& \dsum_{(x)} a_{xy}\left( v_x(u_{x^\m}\stackrel{\T}{\circ} b_x(\om_x)) \stackrel{\Gm}{\circ} v_y((u_{y^\m}\stackrel{\T}{\circ} \bar{u}_{x^\m})\stackrel{\T}{\circ} b_{xy}(b_{x}^{\m}(\bar{u}_x)\stackrel{\Om}{\circ} \om_y))\right) \\
		& =& \dsum_{(x)}\widetilde{\tau}_{x,y}\left( a_x\left( v_x(u_{x^\m}\stackrel{\T}{\circ} b_x(\om_x))\right)  \stackrel{\T}{\circ} a_y\left( v_y((u_{y^\m}\stackrel{\T}{\circ} \bar{u}_{x^\m})\stackrel{\T}{\circ} \widetilde{\gamma}_{x,y}(\bar{u}_x\stackrel{\T}{\circ} b_y(\om_y)))\right) \right) \\
		& =& \dsum_{(x)}\widetilde{\tau}_{x,y}\widetilde{\gamma}_{x,y}\left( a_x\left( v_x(u_{x^\m}\stackrel{\T}{\circ} b_x(\om_x))\right)  \stackrel{\T}{\circ} a_y\left( v_y((u_{y^\m}\stackrel{\T}{\circ} \bar{u}_{x^\m})\stackrel{\T}{\circ} (\bar{u}_x\stackrel{\T}{\circ} b_y(\om_y)))\right) \right) \\
		& =& \dsum_{(x)}\widetilde{\tau}_{x,y}\widetilde{\gamma}_{x,y}\left( a_x\left( v_x(u_{x^\m}\stackrel{\T}{\circ} b_x(\om_x))\right)  \stackrel{\T}{\circ} a_y\left( v_y(u_{y^\m}\stackrel{\T}{\circ} (\bar{u}_{x^\m}\stackrel{\T}{\circ} \bar{u}_x)\stackrel{\T}{\circ} b_y(\om_y))\right) \right) \\
		& =& \widetilde{\tau}_{x,y}\widetilde{\gamma}_{x,y}\left( a_x\left( v_x(u_{x^\m}\stackrel{\T}{\circ} b_x(\om_x))\right)  \stackrel{\T}{\circ} a_y\left( v_y(u_{y^\m}1_{x^\m}\stackrel{\T}{\circ} b_y(\om_y))\right) \right) \\
		& =& \widetilde{\tau}_{x,y}\widetilde{\gamma}_{x,y}\left( a_x\left( v_x(u_{x^\m}\stackrel{\T}{\circ} b_x(\om_x1_{x^\m}))\right)  \stackrel{\T}{\circ} a_y\left( v_y(u_{y^\m}\stackrel{\T}{\circ} b_y(\om_y))\right) \right) \\
		& =& \widetilde{\tau}_{x,y}\widetilde{\gamma}_{x,y}\left( a_x\left( v_x(u_{x^\m}\stackrel{\T}{\circ} b_x(\om_x))\right)  \stackrel{\T}{\circ} a_y\left( v_y(u_{y^\m}\stackrel{\T}{\circ} b_y(\om_y))\right) \right) \\
		& = & \widetilde{\tau}_{x,y}\widetilde{\gamma}_{x,y}\left( d_x(v_x\ot u_{x^\m}\ot  \om_x) \stackrel{\T}{\circ} d_y(v_y\ot u_{y^\m}\ot \om_y)\right). 
	\end{eqnarray*}
	Hence, 	$\tau_{x,y}\gamma_{x,y}\circ f_{x,y}^\T\circ(d_x\ot d_y)=d_{xy}\circ f_{x,y}^{\Gm\Om},$ for all $x,y \in G,$
	and we conclude that   $$\zeta([\D(\Gm)][\D(\Om)])=[\widetilde{\tau}_{-,-}\widetilde{\gamma}_{-,-}]=\zeta([\D(\Gm)])\zeta([\D(\Om)]),$$ completing the proof of the theorem.
\end{dem}

\section{The seven term exact sequence}
\label{cap: sequencia}
In this section we shall construct an exact sequence which generalizes   the Miyashita's sequence for non-commutative unital rings  \cite[Teorema 2.12]{miyashita1973exact}.
To this end fix a ring extension $R\subseteq S$ with the same unity and a unital partial representation  
$$\begin{array}{c c c l}
	\Theta: & G & \longrightarrow & \mathcal{S}_R(S),\\
	& x & \longmapsto & \T_x ,
\end{array}$$
with $\e_x=\T_x\T_{x^{-1}}=R1_x$, for each  $x\in G$. Let $\D(\T)$ be the generalized partial crossed product  provided by  Remark~\ref{pcgpviaSRS} with  factor set $f^\T=\{f_{x,y}^\T:\T_x\ot \T_y\longrightarrow 1_x\T_{xy}, \ x,y \in G\}$    and the ring and $R$-bimodule isomorphism   $\iota:R\longrightarrow \T_1$  as in \eqref{diagramasdaunidade}. 


\subsection{The first exact sequence}
\label{seq1}

{ In this subsection we establish the initial three term part of the final sequence.} Let $\p_\Z(S/R)$ be the group  defined in Section~\ref{sec: pSR}. Denote
\begin{equation*}
	\p_\Z(S/R)^{(G)}=\{\xymatrix@C=1.2cm{ [P]\ar@{=>}[r]|{[\phi]} & [X]}\in \p_\Z(S/R); \ \T_x\phi(P)=\phi(P)\T_x \ \mbox{for all } \ x\in G\}.
\end{equation*}
\begin{obs} Let $\xymatrix@C=1.2cm{ [P]\ar@{=>}[r]|{[\phi]} & [X]}\in \p_\Z(S/R)$.  Then,  
	\begin{equation}
		\T_x\phi(P)=\phi(P)\T_x \ \mbox{ if and only if } \  \T_x\phi(P)\T_{x^{-1}}=\phi(P)1_x, \ \mbox{for each} \ x \in G.\label{obsPSRG}
	\end{equation}
\end{obs}
Indeed,  if  $\T_x\phi(P)=\phi(P)\T_x$, for each  $x\in G$, then  $$\T_x\phi(P)\T_{x^{-1}}=\phi(P)\T_x\T_{x^\m}=\phi(P)1_x, \ \mbox{for each }\ x \in G.$$ 
On the other hand, { since $P$ is a central  $\Z$-bimodule, then   $P1_x=1_xP$, for each  $x\in G$. Thus, by the  $R$-bilinearity of  $\phi$ we have  $1_x\phi(P)=\phi(P)1_x$, for each  $x\in G$. Then} 
\begin{eqnarray*}
	\T_x\phi(P)\T_{x^\m}=\phi(P)1_x & \R & \T_x\phi(P)\T_{x^{-1}}\T_x=\phi(P)1_x\T_x\\
	& \R & \T_x\phi(P)1_{x^\m}=\phi(P)\T_x\\
	& \R & \T_x1_{x^\m}\phi(P)=\phi(P)\T_x\\
	& \R & \T_x\phi(P)=\phi(P)\T_x.
\end{eqnarray*}

\begin{lem}\label{lemma:SubgrP^G} $\p_\Z(S/R)^{(G)}$ is a  subgroup of $\p_\Z(S/R)$.
\end{lem}
\begin{dem} It is easy to see that  $\p_\Z(S/R)^{(G)}$ is closed under the multiplication. Let us verify that it is also closed with respect to taking  inverses.  Let $\xymatrix@C=1.2cm{ [P]\ar@{=>}[r]|{[\phi]} & [X]} \in \p_\Z(S/R)^{(G)}$. { As we know from  Section~\ref{sec:Group P},}
	its inverse in $\p_\Z(S/R)$ is given by  $\xymatrix@C=1.2cm{ [P^*]\ar@{=>}[r]|{[\phi^*]} & [X^*]}$. { Since $\phi$ is injective,  we may assume  $P\subseteq X$ and $P^{*}\subseteq X^*$.} Observe that $P^*=\{f \in X^*; \ f(P)\subseteq R\}$. It follows that $P^*\cdot  1_x=\{f\in X^*; \ f(P)\subseteq R1_x\}$. Indeed,   if $f \in X^*$ is such that  $f(P)\subseteq R1_x$, then $f \in P^*$ and 
	\begin{equation*}
		(f\cdot 1_x)(p)=f(p)1_x=f(p).
	\end{equation*}
	Hence, $f \in P^*\cdot 1_x$. On the other hand, if  $f \in P^*\cdot 1_x$, then $f=f'\cdot 1_x$, for some   $f' \in P^*$. Then,
	\begin{equation*}
		f(p)=(f'\cdot 1_x)(p)=f'(p)1_x \in R1_x, \ \ \mbox{for all} \ \ p \in P,
	\end{equation*}
	and the claimed equality follows. 
	Let us now  see  that  $\T_x\cdot P^*\cdot \T_{x^\m}=P^*\cdot 1_x$, for each  $x\in G$.
	We have:
	\begin{eqnarray*}
		(\T_x\cdot P^*\cdot \T_{x^{-1}})(P) & = & (P^*\cdot \T_{x^\m})(P\T_x)=(P^*\cdot \T_{x^\m})(\T_xP)\\
		& = & \T_x(P^*\cdot \T_{x^\m})(P)=\T_x[P^*(P)]\T_{x^\m}\\
		& \subseteq&  \T_xR\T_{x^{-1}}=R1_x.
	\end{eqnarray*}
	Therefore, $\T_x\cdot P^*\cdot \T_{x^{-1}}\subseteq P^*\cdot R1_x,$ for each  $x\in G$. Now, observe that 
	\begin{eqnarray*}
		P^*\cdot R1_x & = & R1_x\cdot P^*\cdot R1_x = (\T_x\T_{x^{-1}})\cdot P^*\cdot (\T_x\T_{x^\m})\\
		& = & \T_x\cdot (\T_{x^\m}\cdot P^*\cdot \T_x)\cdot \T_{x^\m}\\
		& \subseteq & \T_x\cdot (P^*\cdot R1_{x^\m})\cdot\T_{x^\m}\\
		& = & \T_x\cdot P^*\cdot (R1_{x^\m} \T_{x^\m})\\
		& = & \T_x\cdot P^*\cdot \T_{x^\m}. 
	\end{eqnarray*}
	Consequently, the desired equality follows, and we conclude that  $\xymatrix@C=1.2cm{ [P^*]\ar@{=>}[r]|{[\phi^*]} & [X^*]}\in \p_\Z(S/R)^{(G)}$. 	
\end{dem}

\begin{pro}\label{osegundomorfismo} The map  
	$$\begin{array}{c c c l}
		\varphi_2: &  \p_\Z(S/R)^{(G)} &\longrightarrow& \Pic_\Z(R)\cap \Pics_\Z(R)^{\al^*},\\
		& \xymatrix@C=1.2cm{ [P]\ar@{=>}[r]|{[\phi]} & [X]} & \longmapsto & [P]
	\end{array}$$
	is a well-defined group homomorphism. 
\end{pro}
\begin{dem}
	Clearly,  $\varphi_2$ respects the group operations. Thus, it is enough to show that if   $\xymatrix@C=1.2cm{ [P]\ar@{=>}[r]|{[\phi]} & [X]}\in \p_\Z(S/R)^{(G)}$, then $[P]\in \Pics_\Z(R)^{\al^*}$. By  (\ref{obsPSRG}) we have that  $\T_x\phi(P)\T_{x^{-1}}=\phi(PR1_x)$, for each  $x\in G$.
	Since $\phi$ is $R$-bilinear, then  $\phi(P)$ is an  $R$-subbimodule of $X$.  By Proposition~\ref{SRemPicS} and since   $\phi$ is injective, we obtain that   $f_x$ is defined by 
	$$\xymatrix{\T_x\ot P\ot \T_{x^{-1}}\ar[r]\ar@/_1cm/@{-->}[rrdd]_{f_x} & \T_x\ot \phi(P)\ot \T_{x^{-1}}\ar[r] & \T_x\phi(P)\T_{x^{-1}}=\phi(P1_x)\ar[d]\\
		& &   P1_x\ar[d]\\
		& & P\ot R1_x  }$$
	that is,
	\begin{equation*}
		f_x(u_x\ot p \ot u_{x^\m})=p'\ot 1_x, 
	\end{equation*}
	where $\phi(p')=u_x\phi(p)u_{x^\m}$ is an $R$-bimodule  isomorphism. Hence $[P]\in \Pics_\Z(R)^{\al^*}$.
\end{dem}

Denote 
\begin{equation*}
	\Aut_{R\mbox{-rings}}(S)^{(G)}=\{f\in \Aut_{R\mbox{-rings}}(S); \ f(\T_x)=\T_x, \ \forall \ x\in G\}.
\end{equation*}
Obviously, $\Aut_{R\mbox{-rings}}(S)^{(G)}$ is a subgroup of $\Aut_{R\mbox{-rings}}(S)$.

\begin{lem}\label{diagramaparaoprimieromorfismo} The  following diagram is commutative with exact rows:
	$$\xymatrix{ \U(\Z)\ar[r]^{\F\ \ \ \ \ \ }\ar@{=}[dd] & \Aut_{R\mbox{-rings}}(S)\ar[r]^{\mathcal{E}} & \p(S/R)\ar[r]^{\vartheta} & \Pic(R)\\
		& & & \\
		\U(\Z)\ar[r] &  \Aut_{R\mbox{-rings}}(S)^{(G)}\ar[r]^{{\ \ \ \varphi_1}}\ar@{_(->}[uu] & \p_\Z(S/R)^{(G)}\ar[r]^{\varphi_2\ \ \ \ \ \ \ }\ar@{_(->}[uu] & \Pic_\Z(R)\cap \Pics_\Z(R)^{\al^*}\ar@{_(->}[uu]     }$$ 
	
\end{lem}
\begin{dem}
	The exactness of the first row is given by  Theorem~\ref{sequenciaparaodiagramadoprimeiromorfismo}. The restriction
	of   $\vartheta$ to $\p_\Z(S/R)^{(G)}$ is the homomorphism  $\varphi_2$ from Proposition~\ref{osegundomorfismo}. If $r \in \U(\Z)$, then since  $\T_x$ is an  $R$-bimodule, we have that  $\F(r) \in \Aut_{R\mbox{-rings}}(S)^{(G)}$. If $f \in \Aut_{R\mbox{-rings}}(S)^{(G)}$, then
	\begin{equation*}
		i_f(R)\cdot \T_x=Rf(\T_x)=R\T_x=\T_xR=\T_x\cdot i_f(R), \ \mbox{for each } \ x \in G.
	\end{equation*}
	Since $R$ is a central   $\Z$-bimodule, then $\mathcal{E}(f)=(\xymatrix@C=1.2cm{ [R]\ar@{=>}[r]|{[\iota_f]} & [S_f]}) \in \p_\Z(S/R)^{(G)}$.  Consequently, the homomorphisms of the second row are well-defined.  Let us verify that the second row is exact.  
	
	The exactness at the first term and the inclusion  $\mathcal{E}(\Aut_{R\mbox{-Rings}}(S)^{(G)})\subseteq \ker(\varphi_2)$ directly follow from the exactness of the first row. Let   $\xymatrix@C=1.2cm{ [P]\ar@{=>}[r]|{[\phi]} & [X]} \in \ker(\varphi_2)$. Then there exists   $f \in \Aut_{R\mbox{-rings}}(S)$ such that  $\mathcal{E}(f)=(\xymatrix@C=1.2cm{ [R]\ar@{=>}[r]|{[\iota_f]} & [S_f]})=(\xymatrix@C=1.2cm{ [P]\ar@{=>}[r]|{[\phi]} & [X]})$ in $\p(S/R)$. { More specifically, by the proof of  Theorem~\ref{sequenciaparaodiagramadoprimeiromorfismo} } we have that  if $\lambda:R\longrightarrow P$ is an $R$-bimodule isomorphism, then, defining  $\al$ and $\beta$ by
	\begin{equation*}
		\al: S\longrightarrow R\ot_RS\longrightarrow P\ot_RS\stackrel{\bar{\phi}_r}{\longrightarrow} X,
	\end{equation*}
	\begin{equation*}
		\beta: S\longrightarrow S\ot_RR\longrightarrow S\ot_RP\stackrel{\bar{\phi}_l}{\longrightarrow} X,
	\end{equation*}
	{ we may assume that }  $f=\beta^\m\circ \al ,$ and the  diagram
	$$\xymatrix{ R\ar[rrr]^{i_f}\ar[dd]_{\lambda} & & & S_f\ar[dd]^{\beta}\\
		& & & \\
		P\ar[rrr]_{\phi}& &  & X   }$$
	is commutative. Since $\xymatrix@C=1.2cm{ [P]\ar@{=>}[r]|{[\phi]} & [X]}\in \p_\Z(S/R)^{(G)}$, then $\phi(P)\T_x=\T_x\phi(P)$, for each $x \in G$. We shall verify that  $f(\T_x)=\T_x$, for all  $x\in G$.

	Given $u_x \in \T_x$, we have that  $\phi(\lambda(1))u_x\in \phi(P)\T_x=\T_x\phi(P)$, and, consequently,   there exist  $u_x^{i}\in \T_x$ and $p_i \in P,$  $(i=1,2,...,n)$, such that  $\phi(\lambda(1))u_x=\dsum_{i=1}^nu_x^i\phi(p_i)$. As $\lambda$ is an isomorphism, there exists   $r_i \in R$ such that  $\lambda(r_i)=p_i$, for all  $i=1,2,...,n$. Hence, 
	\begin{eqnarray*}
		\al(u_x) & = & \phi(\lambda(1))u_x = \dsum_{i=1}^nu_x^{i}\phi(p_i)\\
		& = & \dsum_{i=1}^nu_x^{i}\phi(\lambda(r_i))=\dsum_{i=1}^nu_x^{i}\beta(i_f(r_i))\\
		& = & \dsum_{i=1}^nu_x^{i}\beta(r_i)  =\dsum_{i=1}^n\beta(u_x^{i}r_i).
	\end{eqnarray*}
	Thus, 
	$$f(u_x)=(\beta^{-1}\circ \al)(u_x)=\dsum_{i=1}^nu_x^{i}r_i \in \T_x.$$
	Consequently, $f(\T_x)\subseteq \T_x$, for all  $x \in G$. 
	
	On the other hand, given  $u_x \in \T_x,$ we have $u_x\phi(\lambda(1))\in \T_x\phi(P),$ and thus there exist   $v_x^j\in \T_x$ and $p_j\in P$ with $j=1,...,m$, such that  $u_x\phi(\lambda(1))=\dsum_{j=1}^m\phi(p_j)v_x^j$. Using again that   $\lambda$ is an isomorphism, we get that there exist  elements $r_j\in R$ such that  $\lambda(r_j)=p_j$ for $j=1,...,m$. Then
	\begin{eqnarray*}
		\beta(u_x) & = & u_x\phi(\lambda(1)) = \dsum_{j=1}^m\phi(p_j)v_x^j= \dsum_{j=1}^m\phi(\lambda(r_j))v_x^j \\
		& = & \dsum_{j=1}^m\phi(\lambda(1))r_jv_x^j=\al\left( \dsum_{j=1}^mr_jv_x^j\right). 
	\end{eqnarray*}
	Therefore, 
	$$u_x=(\beta^{-1}\circ \al)\left( \dsum_{j=1}^mr_jv_x^j\right)=f\left( \dsum_{j=1}^mr_jv_x^j\right) \in f(\T_x).$$
	Hence, $\T_x\subseteq f(\T_x)$ and, consequently, we obtain the desired equality.
	It follows that  $f\in \Aut_{R\mbox{-rings}}(S)^{(G)}$ and, thus,   $\xymatrix@C=1.2cm{ [P]\ar@{=>}[r]|{[\phi]} & [X]} \in \mathcal{E}(\Aut_{R\mbox{-rings}}(S)^{(G)}).$
\end{dem}

Since $R\subseteq \D(\T)$ is a ring extension with the same unity,   Lemma~\ref{diagramaparaoprimieromorfismo}
is applicable:

\begin{teo}\label{primeiraseqexata}
	The following sequence of { group homomorphisms} is exact: 
	\begin{equation*}
		\xymatrix{ 1\ar[r] &   H^1_\T(G,\al, \Z)\ar[r]^{\varphi_1\ \ \  } &  \p_\Z(\D(\T)/R)^{(G)}\ar[r]^{\varphi_2\ \ \ \ \ } &\Pic_\Z(R)\cap \Pics_\Z(R)^{\al^*} .}
	\end{equation*}
\end{teo}
\begin{dem} By Lemma~\ref{diagramaparaoprimieromorfismo}, the sequence  
	\begin{equation*}
		\xymatrix{ \U(\Z) \ar[r]^{\mathcal{F}\ \ \ \ \ \ \ \ \ } & \Aut_{R\mbox{-rings}}(\Delta(\T))^{(G)}\ar[r] & \p_\Z(\Delta(\T)/R)^{(G)}\ar[r]^{\varphi_2\ \ \ \ } & \Pic_\Z(R)\cap\Pics_\Z(R)^{\al^*} }
	\end{equation*}
	is exact. Hence, 
	\begin{equation*}
		\xymatrix{1\ar[r]  & \dfrac{ \Aut_{R\mbox{-rings}}(\D(\T))^{(G)}}{\mbox{Im}(\mathcal{F})}\ar[r] & \p_\Z(\Delta(\T)/R)^{(G)}\ar[r]^{\varphi_2\ \ \ \ \ \ } & \Pic_\Z(R)\cap\Pics_\Z(R)^{\al^*} }
	\end{equation*}
	is exact. Thus, it suffices to show that there exists a group isomorphism  
	$ \dfrac{\Aut_{R\mbox{-rings}}(\D(\T))^{(G)}}{\mbox{Im}(\mathcal{F})}\simeq H^1_\T(G,\al,\Z)$.
	
	Let $f \in \Aut_{R\mbox{-rings}}(\D(\T))^{(G)}$. Since $f(\T_x)=\T_x$ for all  $x\in G$ and $f$ fixes each element of  $R$,  restricting  $f$ to $\T_x$, we get an $R$-bimodule isomorphism  $f_x:\T_x\longrightarrow\T_x$, for each  $x \in G$. Let $\widetilde{f}_x \in \U(\Z1_x)$ be as in   Lemma~\ref{gammatil}.
	As $f$ is a ring automorphism, we have  that $f(u_x\stackrel{\T}{\circ}u_y)=f(u_x)\stackrel{\T}{\circ} f(u_y)$, for all   $u_x \in \T_x$ and $u_y \in \T_y$. Then 
	\begin{equation*}
		f_{xy}(u_x\stackrel{\T}{\circ}u_y)=f_x(u_x)\stackrel{\T}{\circ} f_y(u_y), \ u_x \in \T_x, u_y \in \T_y.
	\end{equation*}
	By (\ref{relacaodefcomftil}) we obtain
	\begin{equation*}
		f_{xy}(u_x\stackrel{\T}{\circ}u_y)=\widetilde{f_{xy}}(u_x\stackrel{\T}{\circ}u_y),\label{fxyuxuy}
	\end{equation*}
	\begin{eqnarray*}
		f_x(u_x)\stackrel{\T}{\circ} f_y(u_y) & =& \widetilde{f_x}u_x\stackrel{\T}{\circ}\widetilde{f_y}u_y=\widetilde{f_x}u_x\widetilde{f_y}\stackrel{\T}{\circ}u_y\nonumber\\
		&\stackrel{(\ref{uxercomutam})}{=} & \widetilde{f_x}\al_x(\widetilde{f_y}1_{x^\m})(u_x\stackrel{\T}{\circ}u_y).
	\end{eqnarray*}
	Consequently, $\widetilde{f_{xy}}(u_x\stackrel{\T}{\circ}u_y)=\widetilde{f_x}\al_x(\widetilde{f_y}1_{x^\m})(u_x\stackrel{\T}{\circ}u_y),$ for all  $u_x \in \T_x$ and $u_y \in \T_y$. { {	Using the same argument as in Proposition \ref{3cocicloquecorrigeaassociatividade} we obtain}}
	\begin{equation*}
		\widetilde{f_x}\al_x(\widetilde{f_y}1_{x^\m})=\widetilde{f_{xy}}1_x, \ \mbox{for all} \  x,y \in G.
	\end{equation*} Then
	\begin{equation*}
		\begin{array}{c c c l}
			\widetilde{f}: & G & \longrightarrow & \Z ,\\
			& x & \longmapsto & \widetilde{f_x},
		\end{array}
	\end{equation*}
	belongs to  $Z_\T^1(G,\al,\Z)$. Since $\widetilde{f_1}=f_1(1)=1,$ it follows that    $\widetilde{f}$ is a normalized element in   $Z_\T^1(G,\al,\Z)$. We have the following homomorphism of groups:
	\begin{equation*}
		\begin{array}{c c c l}
			\Psi:&	\Aut_{R\mbox{-rings}}(\D(\T))^{(G)} & \longrightarrow & Z^1_\T(G,\al,\Z)\\
			&	f & \longmapsto & \widetilde{f}
		\end{array}.
	\end{equation*}
	
	Conversely, if $\sigma \in Z^1_\T(G,\al,\Z)$ is a normalized $1$-cocycle, then  
	\begin{equation}
		\al_x(\s_{y}{1_{x^\m}})\s_x=\s_{xy}1_x, \ \ \mbox{for all} \ x,y \in G. \label{sigmae1cociclonoprimeiromorfismo}
	\end{equation}
	For each  $x\in G$ define
	\begin{equation*}
		\begin{array}{c c c l}
			g_x:& \T_x & \longrightarrow & \T_x , \\
			& u_x & \longmapsto & \s_xu_x. 
		\end{array}
	\end{equation*}
	Since $\sigma_x\in \U(\Z1_x)$, it follows by  Lemma~\ref{gammatil} that  $g_x$ is an $R$-bilinear isomorphism  and  $\widetilde{g_x}=\sigma_x$, for each  $x\in G$.
	Consider the map
	\begin{equation*}
		\begin{array}{ c c c l}
			g:=\bigoplus_{x \in G}g_x: & \D(\T) & \longrightarrow & \D(\T)\\
			& u_x & \longrightarrow & g_x(u_x)
		\end{array}.
	\end{equation*}
	Then $g_x(\T_x)=\T_x$, for all  $x\in G$. Since $\s$ is normalized, we have:  
	\begin{equation*}
		g(r)=g_1(r)=\s_1r=1r=r, \ \ \forall \ r\in R.
	\end{equation*}
	Hence, $g$ fixes each element of   $R$. Now, given  $u_x\in \T_x$ and $u_y\in \T_y$, we compute:
	\begin{eqnarray*}
		g(u_x\stackrel{\T}{\circ}u_y) & = & \sigma_{xy}(u_x\stackrel{\T}{\circ}u_y) \stackrel{(\ref{sigmae1cociclonoprimeiromorfismo})}{=}\s_x\al_x(\s_y1_{x^\m})(u_x\stackrel{\T}{\circ}u_y)\\
		& = &( \s_x\al_x(\s_y1_{x^\m})u_x\stackrel{\T}{\circ}u_y)\stackrel{(\ref{uxercomutam})}{=}( \s_xu_x\sigma_y\stackrel{\T}{\circ}u_y)\\
		& = & ( \s_xu_x\stackrel{\T}{\circ}\sigma_yu_y)= ( g_x(u_x)\stackrel{\T}{\circ}g_y(u_y)).
	\end{eqnarray*}
	Therefore, $g \in 	\Aut_{R\mbox{-rings}}(\D(\T))^{(G)}$ and we have the map  
	\begin{equation*}
		\begin{array}{c c c l}
			\Phi:& Z^1_\T(G,\al,\Z) & \longrightarrow & 	\Aut_{R\mbox{-rings}}(\D(\T)/R)^{(G)},\\
			& \sigma & \longmapsto & g.
		\end{array}
	\end{equation*}
	Let us check that   $\Phi$ is the  inverse of $\Psi$. If $f \in \Aut_{R\mbox{-rings}}(\D(\T))^{(G)}$, then 
	\begin{equation*}
		\Phi(\Psi(f))=\Phi(\widetilde{f})=\bigoplus_{x \in G}g_x,
	\end{equation*}
	where $g_x(u_x)=\widetilde{f_x}u_x=f_x(u_x)$, for all   $u_x \in \T_x$, $x\in G$. Hence, $\Psi(\Phi(f))=f$, for all 
	$f \in \Aut_{R\mbox{-rings}}(\D(\T))^{(G)}$.
	
	On the other hand, if  $\s \in Z^1_\T(G,\al,\Z)$, then
	\begin{equation*}
		\Psi(\Phi(\s))=\Psi(g)=\widetilde{g},
	\end{equation*}
	where $g_x(u_x)=\sigma_xu_x$, for each  $x\in G$. Thus, $\widetilde{g_x}=\sigma_x$, for all  $x\in G$. Consequently, $\Psi(\Phi(\s))=\s, $ for all  $\s\in Z^1_\T(G,\al,\Z)$. 
	Therefore, $\Psi:	\Aut_{R\mbox{-rings}}(\D(\T))^{(G)} \longrightarrow  Z^1_\T(G,\al,\Z)$ is a group isomorphism.
	{Now}, if $r \in \U(\Z)$, then
	\begin{equation*}
		\mathcal{F}(r)(u_x)={r u_xr^\m=r\al_x(r^\m1_{x^\m})u_x=(\delta^0r^\m)(x)u_x}, \ \mbox{for all } \ u_x \in \T_x.
	\end{equation*}
	Hence, $\Psi(\mathcal{F}(r))=\delta^0r^\m\in B^1_\T(G,\al,\Z)$. Conversely, if  $\s \in B^1_\T(G,\al,\Z)$, then there exists  $r \in \U(\Z)$ such that  $\delta^{0}r=\s$, that is, 
	\begin{equation*}
		\sigma(x)=\al_x(r1_{x^\m})r^\m.
	\end{equation*}
	Then, $\mathcal{F}({{r^\m}})(u_x)=r^\m u_xr=r^\m\al_x(r1_{x^\m})u_x$, for all   $u_x \in \T_x.$
	Thus, $\widetilde{\mathcal{F}(r)}_x={\s(x)}$, for each  $x\in G$. Consequently, $\Psi(\mathcal{F}(r))=\s,$. {Thus 	$\Psi(\mbox{Im}(\mathcal{F}))=B^1_\T(G,\al,\Z)$ and we have the isomorphism
		\begin{equation*}
			\dfrac{ \Aut_{R\mbox{-rings}}(\D(\T))^{(G)}}{\mbox{Im}(\mathcal{F})}\simeq H^1_\T(G,\al,,\Z), \label{isomorfismodeH1}
		\end{equation*} 
		as desired.}
\end{dem}

\subsection{The second exact sequence}\label{sec:secondExSeq}
Let 
\begin{equation*}
	\Pic_\Z(R)^{(G)}=\{[P]\in \Pic_\Z(R); P\ot\T_x\ot P^{-1}|\T_x, \ \mbox{for all} \ x \in G\}.
\end{equation*}

\begin{lem}\label{lemma:PicZ^G} $\Pic_\Z(R)^{(G)}$ is a  subgroup of $\Pic_\Z(R)$ which contains  $\Pic_\Z(R)\cap\Pics_\Z(R)^{\al^*}.$ 
\end{lem}
\begin{dem} Let $[P],[Q]\in \Pic_\Z(R)^{(G)}$. Since $Q\ot\T_x\ot Q^{\m}|\T_x ,$ it follows by compatibility with the tensor product that   $P\ot Q\ot\T_x\ot Q^{\m}\ot P^{\m}|P\ot \T_x\ot P^{\m}$. Therefore, by transitivity, we obtain that  $$P\ot Q\ot \T_x\ot Q^{\m}\ot P^{\m}|\T_x, \ \ \mbox{for all } x \in G.$$
	Hence, $[P\ot Q]\in \Pic_\Z(R)^{(G)}$.  Let us verify that  $\Pic_\Z(R)^{(G)}$ is also closed with respect to inverses. Take  $[P]\in \Pic_\Z(R)^{(G)}$. Then
	\begin{eqnarray*}
		P^\m\ot \T_x\ot P & \simeq & P^\m\ot \T_x\ot R1_{x^\m}\ot P \simeq P^\m\ot \T_x\ot P\ot R1_{x^\m}\\
		&\simeq & P^\m\ot \T_x\ot P\ot \T_{x^\m}\ot \T_x \simeq P^\m\ot \T_x\ot P\ot \T_{x^\m}\ot R\ot \T_x \\
		& \simeq & P^\m\ot \T_x\ot P\ot \T_{x^\m}\ot P^\m \ot P\ot \T_x. 
	\end{eqnarray*}
	Thus, since $P\ot \T_x\ot P^\m|\T_x$, for each  $x\in G$, we have 
	\begin{equation*}
		P^\m\ot \T_x\ot P\ot \T_{x^\m}\ot P^\m \ot P\ot \T_x |P^\m \ot \T_x\ot \T_{x^\m}\ot P\ot \T_x.
	\end{equation*}
	But,
	\begin{equation*}
		P^\m \ot \T_x\ot \T_{x^\m}\ot P\ot \T_x\simeq P^\m\ot R1_x\ot P\ot \T_x\simeq P^\m\ot P\ot R1_x\ot \T_x\simeq \T_x.
	\end{equation*}
	Consequently, $P^\m\ot \T_x\ot P|\T_x$, for each $x\in G ,$ showing that  $[P^\m]\in \Pic_\Z(R)^{(G)}$. Therefore, $\Pic_\Z(R)^{(G)}$ is a subgroup of  $\Pic_\Z(R)$.
	
	Let $[P]\in\Pic_\Z(R)\cap\Pics_\Z(R)^{\al^*}$. By (\ref{invariantesdePicsZ}) we have that  $\T_x\ot P\simeq P\ot\T_x$, for each  $x \in G$. As $[P]\in \Pic_\Z(R)$, it follows that 
	\begin{equation}
		P\ot\T_x\ot P^{\m}\simeq \T_x, \ \forall \ x\in G, \ [P]\in \Pic_\Z(R)\cap\Pics_\Z(R)^{\al^*}.\label{PicRePicsRinvarinates}
	\end{equation} 
	In particular, $P\ot\T_x\ot P^{\m}|\T_x$, for all  $x \in G$. Hence, $\Pic_\Z(R)\cap\Pics_\Z(R)^{\al^*}\subseteq \Pic_\Z(R)^{(G)}$.
\end{dem}

Our next purpose is to construct a generalized partial crossed product starting with an element $[P]\in \Pic_\Z(R)^{(G)}.$
We begin by producing  a unital partial representation.

\begin{lem} Let $[P]\in \Pic_\Z(R)^{(G)}$ and  denote $\Om_x^{P}=P\ot\T_x\ot P^{^\m}$,  $(x \in G)$. Then 
	\begin{equation*}
		\begin{array}{c c c l}
			\Om^P: & G & \longrightarrow & \Pics(R),\\
			& x & \longmapsto & [\Om_x^P]
		\end{array}
	\end{equation*}
	is a unital partial representation with  $\Om^P_x\ot \Om^P_{x^\m}\simeq R1_x$ and $\Om_x^P|\T_x$, for all  $x\in G$.
\end{lem}
\begin{dem} Obviously, $[\Om^P_1]=[R]$ { and $\Om_x^P|\T_x$, for all  $x\in G$}. Given $x,y \in G ,$ we have:
	\begin{eqnarray*}
		\Om_{x}^P\ot\Om_y^P\ot\Om_{y^\m}^P & = & P\ot\T_x\ot P^{\m}\ot P\ot\T_y\ot P^{\m}\ot P\ot\T_{y^\m}\ot P^{\m}\\
		& \simeq & P\ot\T_x\ot\T_y\ot\T_{y^\m}\ot P^{\m}\\
		& \simeq & P\ot\T_{xy}\ot\T_{y^\m}\ot P^{\m}\\
		& \simeq & P\ot\T_{xy}\ot P^\m\ot P\ot\T_{y^\m}\ot P^{\m}\\
		& = & \Om_{xy}^P\ot\Om_{y^\m}^P.
	\end{eqnarray*}
	Similarly,  we obtain that  $\Om_{x^\m}^P\ot\Om_x^P\ot\Om_y^P\simeq \Om_{x^\m}^P\ot\Om_{xy}^P$, for all  $x,y \in G$.
	Moreover, 
	\begin{eqnarray*}
		\Om_x^P\ot\Om_{x^\m}^P& = & P\ot\T_x\ot P^{^\m}\ot P\ot\T_{x^\m}\ot P^{^\m}\\
		&\simeq & P\ot\T_x\ot\T_{x^\m}\ot P^{^\m}\\
		& \simeq & P\ot R1_{x}\ot P^{^\m}\\
		& \simeq &  R1_x\ot P\ot P^{\m}\\
		& \simeq & R1_x\ot R\simeq R1_x, \ \ \forall  \ x\in G.
	\end{eqnarray*}
	Therefore, $\Om^P$ is a unital partial representation with  $\Om_x^P\ot\Om_{x^\m}^P\simeq R1_x$, for all $x \in G$.
\end{dem}

Given $[P]\in \Pic_\Z(R)^{(G)}$ let  $P^{-1}\ot P\stackrel{\Rr}{\longrightarrow} R \stackrel{\Ll}{\longleftarrow}P\ot P^{-1}$ be $R$-bimodule isomorphisms. Define the $R$-bimodule  isomorphisms   $f_{x,y}^P:\Om_x^P\ot\Om_y^P\longrightarrow 1_x\Om_{xy}^P$ via 
$$\xymatrix{ P\ot\T_x\ot P^{^\m}\ot P\ot\T_{y}\ot P^{^\m}\ar[r]\ar@/_1cm/@{-->}[rdd]_{f_{x,y}^P} & P\ot\T_x\ot\T_{y}\ot P^{^\m}\ar[d]^{P\ot f_{x,y}^\T\ot P^{-1}} \\
	& P\ot 1_x\T_{xy}\ot P{^\m}\ar[d]\\
	& 1_xP\ot\T_{xy}\ot P{^\m}, }$$
where the first isomorphism is induced by  $\Rr$, that is,
\begin{equation*}
	f_{x,y}^P(p_1\ot u_x\ot \bar{p}_1\ot p_2\ot u_y \ot \bar{p}_2)=p_1\ot (u_x\Rr(\bar{p}_1\ot p_2)\stackrel{\T}{\circ}u_y )\ot \bar{p}_2,
\end{equation*}
for $u_x \in \T_x,u_y\in \T_y$, $p_1,p_2 \in P$ and $\bar{p}_1,\bar{p}_2\in P^{-1}$. 

\begin{afr}
	$f^P=\{f_{x,y}^P:\Om_x^P\ot\Om_y^P\longrightarrow 1_x\Om_{xy}^P, \ x,y \in G\}$ is a factor set for  $\Om^P$.
\end{afr}  Indeed, given $x,y,z \in G$, let us check that the following  diagram is  commutative:
\begin{equation*}
	\xymatrix{ P\ot\T_x\ot P^{-1}\ot P\ot\T_y\ot P^{-1}\ot P\ot\T_z\ot P^{-1}\ar[r]\ar[dd] &  1_xP\ot\T_{xy}\ot P^{-1}\ot P\ot \T_z\ot P^{-1}\ar[dd]\\
		& & \\
		P\ot\T_x\ot P^{-1}\ot 1_yP\ot \T_{yz}\ot P^{-1}\ar[r]&  1_x1_{xy}P\ot\T_{xyz}\ot P^{-1}.   } \label{diagramaconjuntodefatoresparaP}
\end{equation*}
Take $u_x\in\T_x, u_y \in \T_y, u_z\in \T_z$, $p_1,p_2,p_3\in P$ and $\bar{p}_1,\bar{p}_2,\bar{p}_3\in P^{-1}.$  Then  
\begin{eqnarray*}
	& & [f_{xy,z}^P\circ (f_{x,y}^P\ot \Om_z^P)](p_1\ot u_x\ot\bar{p}_1 \ot p_2\ot u_y\ot \bar{p}_2\ot p_3\ot u_z\ot \bar{p}_3)\\
	& = & f_{xy,z}^P(p_1\ot (u_x\Rr(\bar{p}_1 \ot p_2)\stackrel{\T}{\circ} u_y)\ot \bar{p}_2\ot p_3\ot u_z\ot \bar{p}_3)\\
	& = & p_1\ot \left( (u_x\Rr(\bar{p}_1 \ot p_2)\stackrel{\T}{\circ} u_y)\Rr( \bar{p}_2\ot p_3)\stackrel{\T}{\circ} u_z\right) \ot \bar{p}_3\\
	& = & p_1\ot \left( u_x\Rr(\bar{p}_1 \ot p_2)\stackrel{\T}{\circ} (u_y\Rr( \bar{p}_2\ot p_3)\stackrel{\T}{\circ} u_z)\right) \ot \bar{p}_3\\
	& = & f_{x,yz}^{P}(p_1\ot u_x\ot \bar{p}_1 \ot p_2 \ot  (u_y\Rr( \bar{p}_2\ot p_3)\stackrel{\T}{\circ} u_z)\ot \bar{p}_3)\\
	& = & [f_{x,yz}^{P}\circ (\Om_x^{P}\ot f_{y,z}^{P})](p_1\ot u_x\ot \bar{p}_1 \ot p_2 \ot  u_y\ot  \bar{p}_2\ot p_3\ot  u_z\ot \bar{p}_3).\\
\end{eqnarray*}
Consequently, $f_{xy,z}^P\circ (f_{x,y}^P\ot \Om_z^P)=f_{x,yz}^{P}\circ (\Om_x^{P}\ot f_{y,z}^{P}),$ which shows the   commutativity of the above   diagram.  Therefore,   $f^P=\{f_{x,y}^P:\Om_x^P\ot\Om_y^P\longrightarrow 1_x\Om_{xy}^P, \ x,y \in G\}$ is a factor set for  $\Om^P,$ proving the claim.

Thanks to the above claim we have  a generalized partial crossed product 
\begin{equation*}
	\D(\Om^{P})=\bigoplus_{x\in G}P\ot\T_x\ot P^\m,
\end{equation*}
such that  $P\ot\T_x\ot P^\m|\T_x$ and $\Om_x^P\ot\Om_{x^\m}^P\simeq R1_x$, for each  $x \in G$. Hence, $[\D(\Om^{P})]\in \C(\T/R)$, for all  $[P]\in \Pic_\Z(R)^{(G)}$.

\begin{obs} It is important to observe that the isomorphism class of  $\D(\Om^P)$ in $\C(\T/R)$ does not  depend on the choice of the isomorphism  $\Rr$ used to define the isomorphism  $f_{x,y}^P$ in $\Om^P$.
\end{obs}
Indeed, let us consider another
$R$-bimodule  isomorphism  $\widetilde{\Rr}:{P^\m\ot P}\longrightarrow R$, and denote by   $\widetilde{f}_{x,y}^P$ the  isomorphism  determined by  $\widetilde{\Rr}$ and by  $\D(\widetilde{\Om}^P)$ the generalized partial crossed product with the factor set  $\widetilde{f}^P$. Since $\widetilde{\Rr}\circ \Rr^\m \in \Aut_{R-R}(R)\simeq \U(\Z)$, there exists   $a \in \U(\Z)$ such that  $\widetilde{\Rr}=a\Rr$. Consider the $R$-bimodule  isomorphism 
\begin{equation*}
	\begin{array}{c c c l}
		\varphi_x: & P\ot \T_x\ot P^\m & \longrightarrow & P\ot\T_{x}\ot P^\m\\
		& p\ot u_x\ot {\bar{p}} & \longmapsto & p\ot a^\m u_x\ot {\bar{p}}.
	\end{array}
\end{equation*}
Then the diagram below is commutative:
\begin{equation*}
	\xymatrix{ P\ot\T_x\ot P^\m\ot P\ot \T_y\ot P^\m\ar[rr]^{\ \ \  \ \ \ \  f_{x,y}^P}\ar[dd]_{\varphi_x\ot \varphi_y} & & 1_xP\ot \T_{xy}\ot P^\m\ar[dd]^{\varphi_{xy}} \\
		& & \\
		P\ot\T_x\ot P^\m\ot P\ot \T_y\ot P^\m\ar[rr]_{\ \ \ \ \ \ \ \widetilde{f}_{x,y}^P} & &   1_xP\ot \T_{xy}\ot P^\m . }
\end{equation*}
Indeed, let  $p_1,p_2\in P, \bp_1,\bp_2 \in P^\m, u_x\in \T_x$ and $u_y \in \T_y$. On one hand,
\begin{equation*}
	\varphi_{xy}\circ f_{x,y}^P(p_1\ot u_x\ot \bp_1\ot p_2\ot u_y\ot \bp_2)=p_1\ot (a^\m u_x \Rr(\bp_1\ot p_2)\stackrel{\T}{\circ} u_y)\ot \bp_2.
\end{equation*}
On the other hand, 
\begin{eqnarray*}
	\widetilde{f}_{x,y}^P\circ (\varphi_x\ot \varphi_y)(p_1\ot u_x\ot \bp_1\ot p_2\ot u_y\ot \bp_2) & = & \widetilde{f}_{x,y}^P(p_1\ot a^\m u_x\ot \bp_1\ot p_2\ot a^\m u_y\ot \bp_2)\\
	& = & p_1\ot (a^\m u_x\widetilde{\Rr}(\bp_1\ot p_2)\stackrel{\T}{\circ} a^\m u_y)\ot \bp_2\\
	& = & p_1\ot (a^\m u_x a\Rr(\bp_1\ot p_2)\stackrel{\T}{\circ} a^\m u_y)\ot \bp_2\\
	& = & p_1\ot (a^\m u_x \Rr(\bp_1\ot p_2)\stackrel{\T}{\circ} aa^\m u_y)\ot \bp_2\\
	& = & p_1\ot (a^\m u_x \Rr(\bp_1\ot p_2)\stackrel{\T}{\circ} u_y)\ot \bp_2,
\end{eqnarray*}
which shows the commutativity of  the diagram. Therefore, $[\D(\Om^P)]=[\D(\widetilde{\Om}^P)]$ in $\C(\T/R)$.


{ With respect to the generalized partial crossed product  $\D(\Om^P)$ we know that 
	there exists an $R$-bimodule and a ring  isomorphism  $\nu_P:R\longrightarrow P\ot \T_1\ot P^\m$ with commutative diagrams  like in  (\ref{diagramasdaunidade}).  We shall specify this isomorphism in the next remark, which will be used later.}

\begin{obs}\label{morfismonuparaOmP} Let $\iota:R\longrightarrow \T_1$ be the { $R$-bimodule  and ring isomorphism} which satisfies the commutative diagrams in (\ref{diagramasdaunidade}). Let $\nu_P:R\longrightarrow P\ot\T_1\ot P^\m$ be the isomorphism  defined by 
	\begin{equation*}
		\nu_P(r)=\dsum_{k=1}^nrp_k\ot \iota(1)\ot \bp_k,
	\end{equation*}
	where $\dsum_{k=1}^n\Ll(p_k\ot \bp_k)=1$. Then, $\nu_P$ satisfies the commutative diagrams  
	\begin{equation*}
		\begin{tabular}{c c}
			
			\xymatrix{ R\ot \Om_x^P\ar[rr]^{\simeq}\ar[rdd]_{\nu_P} & & \Om_x^P\\
				& & \\
				& \Om_1^P\ot \Om_x^P\ar[uur]_{f_{1,x}^P}  }
			& 	
			\xymatrix{ \Om_x^P\ot R\ar[rr]^{\simeq}\ar[rdd]_{\nu_P} & & \Om_x^P\\
				& & \\
				&  \Om_x^P\ot \Om_1^P\ar[uur]_{f_{x,1}^P}  }
			
		\end{tabular}
	\end{equation*}
\end{obs}
Indeed, given $r \in R$, $p\in P, u_x\in \T_x$ and $\bp\in P^\m,$ we have 
\begin{eqnarray*}
	f_{1,x}^P(\nu_P(r\ot p\ot u_x\ot \bp))& = & \dsum_{k=1}^nf_{1,x}^P(rp_k\ot \iota(1)\ot \bp_k\ot p\ot u_x\ot \bp)\\
	& = & \dsum_{k=1}^nrp_k\ot {(\iota(1)\stackrel{\T}{\circ} \Rr(\bp_k\ot p)u_x)}\ot \bp\\
	& = & \dsum_{k=1}^nrp_k\ot \Rr(\bp_k\ot p)u_x\ot \bp\\
	& = & \dsum_{k=1}^nrp_k\Rr(\bp_k\ot p)\ot u_x\ot \bp\\
	& = & \dsum_{k=1}^nr\Ll(p_k\ot\bp_k) p\ot u_x\ot \bp\\
	& = & rp\ot u_x\ot \bp.
\end{eqnarray*}
Analogously, we obtain the commutativity of the second diagram.

\begin{teo}\label{teo:homL} The map 
	$$\begin{matrix}
		\mathcal{L}: & \Pic_\Z(R)^{(G)} & \longrightarrow & \C(\T/R),\\
		& [P] & \longmapsto & [\D(\Om^{P})]
	\end{matrix}$$
	is a group homomorphism.
\end{teo}
\begin{dem}
	We first verify  that   $\mathcal{L}$ is well-defined.  If $[P]=[Q]$ in $\Pic_\Z(R)^{(G)}$, then there exist $R$-bimodule isomorphisms   $\varphi:Q\longrightarrow P$ and $\overline{\varphi}:Q^\m\longrightarrow P^\m$. Let $\Rr:P^\m\ot P\longrightarrow R$ be an $R$-bimodule isomorphism and let   
	$\Rr':Q^\m\ot Q\longrightarrow R$ be  the $R$-bimodule isomorphism defined by  $\Rr'(\bar{q}\ot q)=\Rr(\overline{\varphi}(\bar{q})\ot \varphi(q))$.
	
	For each  $x \in G$ consider the $R$-bimodule isomorphism  
	$$\begin{array}{c c c l}
		\varphi_x: & Q\ot \T_x\ot Q^\m & \longrightarrow & P\ot \T_x\ot P^\m ,\\
		& q\ot u_x\ot \bar{q} & \longmapsto &\varphi(q) \ot u_x \ot \overline{\varphi}(\bar{q}).
	\end{array}$$
	Then the   diagram 
	$$\xymatrix@C=3cm{ Q\ot\T_x\ot Q^{\m}\ot Q\ot\T_y\ot Q^{\m} \ar[r]^{f_{x,y}^Q}\ar[dd]_{\varphi_x\ot \varphi_y} & 1_xQ\ot\T_{xy}\ot Q^{\m}\ar[dd]^{\varphi_{xy}}    \\
		& \\
		P\ot\T_x\ot P^{\m}\ot P\ot\T_y\ot P^{\m}\ar[r]_{f_{x,y}^P} & 1_x{P}\ot\T_{xy}\ot P^{\m}  }$$
	
	is commutative. Indeed, if $u_x \in \T_x,u_y \in \T_y,q_1,q_2 \in Q$ and $\bar{q}_1,\bar{q}_2\in Q^\m$, then
	\begin{eqnarray*}
		(\varphi_{xy}\circ f_{x,y}^Q)(q_1\ot u_x\ot \bar{q}_1\ot q_2\ot u_y\ot \bar{q}_2) & = & \varphi_{xy}(q_1\ot (u_x\Rr'(\bar{q}_1\ot q_2)\stackrel{\T}{\circ}u_y)\ot \bar{q}_2)\\
		& = & \varphi_{xy}(q_1\ot (u_x\Rr(\overline{\varphi}(\bar{q}_1)\ot \varphi(q_2))\stackrel{\T}{\circ}u_y)\ot \bar{q}_2)\\
		& = & \varphi(q_1)\ot (u_x\Rr(\overline{\varphi}(\bar{q}_1)\ot \varphi(q_2))\stackrel{\T}{\circ}u_y)\ot \overline{\varphi}(\bar{q}_2)\\
		& = & f_{x,y}^P(\varphi(q_1)\ot u_x\ot \overline{\varphi}(\bar{q}_1)\ot \varphi(q_2)\ot u_y\ot \overline{\varphi}(\bar{q}_2))\\
		& = & (f_{x,y}^P\circ (\varphi_x\ot \varphi_y))(q_1\ot u_x\ot \bar{q}_1\ot q_2\ot u_y\ot \bar{q}_2).
	\end{eqnarray*}
	Therefore, $[\D(\Om^P)]=[\D(\Om^Q)]$ in $\C(\T/R)$ and, consequently,   $\mathcal{L}$ is well-defined. 
	
	Let $[P],[Q]\in \Pic_\Z(R)^{(G)}$ { and let, furthermore,} 
	\begin{equation*}
		P\ot P^{-1}\stackrel{\Ll}{\longrightarrow} R \stackrel{\Rr}{\longleftarrow} P^{-1}\ot P \ \ \ \ \mbox{and} \ \ \ \ Q\ot Q^{-1}\stackrel{\Ll'}{\longrightarrow} R \stackrel{\Rr'}{\longleftarrow} Q^{-1}\ot Q
	\end{equation*}
	{ be  $R$-bimodule isomorphisms.}
	Then, 
	$$\mathcal{L}([P\ot Q])=\left[ \bigoplus_{x\in G}P\ot Q\ot \T_{x}\ot Q^{-1}\ot P^{-1}\right].$$
	On the other hand, 
	$$\mathcal{L}([P])\mathcal{L}([Q])=\left[ \bigoplus_{x\in G} P\ot\T_{x}\ot P^{-1}\ot\T_{x^{-1}}\ot Q\ot\T_x\ot Q^{-1}\right]. $$
	
	Since $Q\ot\T_x\ot Q^{-1}|\T_x$, then $\T_{x^{-1}}\ot Q\ot \T_x\ot Q^{-1}|R$. As $P^{-1}$ is  a central  $\Z$-bimodule, { by Corollary~\ref{TcomZbomodcentral}} { we have an $R$-bimodule isomorphism  $P^\m\ot\T_{x^{-1}}\ot Q\ot\T_x\ot Q^{-1} \simeq  \T_{x^{-1}}\ot Q\ot\T_x\ot Q^{-1}\ot P^\m$. } 
	Similarly, since  $Q$ is a central  $\Z$-bimodule and $R1_x|R$, we have an $R$-bimodule isomorphism  {$R1_x\ot Q\longrightarrow Q\ot R1_x$, for each $x \in G.$ }
	
	Define the $R$-bimodule isomorphism  $F_x:\Om_x^P\ot \T_{x^\m}\ot \Om_x^Q\longrightarrow \Om_x^{P\ot Q}$ by:
	{	\begin{equation*}
			\xymatrix@C=1.5cm{ P\ot\T_{x}\ot P^{-1}\ot \T_{x^{-1}}\ot Q\ot \T_x\ot Q^{-1} \ar[r]^{P\ot \T_x\ot T}\ar@/_1.5cm/@{-->}[rddddd]_{F_x} & P\ot \T_{x}\ot \T_{x^{-1}}\ot Q\ot \T_x\ot Q^{-1}\ot P^{-1}\ar[dd]|{P\ot f_{x,x^\m}^\T\ot Q\ot \T_x\ot Q^{-1}\ot P^{-1}}\\
				& \\
				&  P\ot R1_{x}\ot Q\ot \T_x\ot Q^{-1}\ot P^{-1}\ar[dd]|{P\ot Q\ot T'\ot Q^\m\ot P^\m}\\
				& \\
				& P\ot Q\ot R1_{x}\ot\T_x\ot Q^{-1}\ot P^{-1}\ar[d]|{\simeq}\\
				& P\ot Q\ot \T_x\ot Q^{-1}\ot P^{-1}. }
	\end{equation*}}
	Let us check  that the following  diagram is commutative
	\begin{equation*}
		\xymatrix@C=3cm{ \Om_x^P\ot\T_{x^\m}\ot \Om_x^Q\ot \Om_y^P\ot\T_{y^\m}\ot\Om_y^Q\ar[r]^{f_{x,y}^{PQ}}\ar[dd]_{F_x\ot F_y} &  1_x\Om_{xy}^P\ot \T_{(xy)^\m}\ot \Om_{xy}^Q\ar[dd]^{F_{xy}}\\
			& \\
			\Om_x^{P\ot Q}\ot \Om_y^{P\ot Q}\ar[r]_{\ \ \  f_{x,y}^{P\ot Q}} & 1_x\Om_{xy}^{P\ot Q}
		}
	\end{equation*}
	

	{ Consider the $R$-bimodule isomorphisms
		$$T_1:\T_{x^{-1}}\ot Q\ot \T_x\ot Q^{^\m}\ot P\ot \T_y\ot P^{\m}\ot \T_{y^{-1}}\longrightarrow P\ot \T_y\ot P^{\m}\ot \T_{y^{-1}}\ot \T_{x^{-1}}\ot Q\ot \T_x\ot Q^{^\m},$$
		$$T_2:P^{\m}\ot 1_{y^\m}\T_{(xy)^{-1}}\ot Q\ot 1_x\T_{xy}\ot Q^{\m}\longrightarrow  {1_{y^\m}\T_{(xy)^{-1}}\ot Q\ot 1_x\T_{xy}\ot Q^{\m}\ot P^{\m}}$$
		$$T_3:R1_{xy}\ot Q\longrightarrow Q\ot R1_{xy}.$$
		$$T_4:P^{-1}\ot \T_{x^\m}\ot Q\ot \T_x\ot Q^{^\m}\longrightarrow \T_{x^\m}\ot Q\ot \T_x\ot Q^{^\m}\ot P^{-1} , $$
		$$T_5:P^{-1}\ot \T_{y^\m}\ot Q\ot \T_y\ot Q^{^\m}\longrightarrow \T_{y^\m}\ot Q\ot \T_y\ot Q^{^\m}\ot P^{-1}$$
		and
		$$T_6:R1_x\ot Q\longrightarrow Q\ot R1_x, \ \ T_7:R1_y\ot Q\longrightarrow Q\ot R1_y.$$}
	
	{Then},  $F_{xy}\circ f_{x,y}^{PQ}$ is given by the sequence of isomorphisms: 
	
	%
	{	\begin{equation*}
			\xymatrix@C=-1.5cm{   \Om_x^P\ot\T_{x^{-1}}\ot \Om_x^Q\ot \Om_y^P\ot\T_{y^{-1}}\ot \Om_y^Q \ar[rd]^{\Om_x^P\ot T_1\ot \Om_y^Q}\ar@{-->}@/_5cm/[rddddddd]_{F_{xy}\circ f_{x,y}^{PQ}}
				&   \\
				& \Om_x^P\ot \Om_y^P\ot\T_{y^{-1}}\ot\T_{x^{-1}}\ot \Om_x^Q\ot \Om_y^Q\ar[d]^{f_{x,y}^P\ot f_{y^\m,x^\m}^\T\ot f_{x,y}^Q}  \\
				& {R1_x\ot}  P\ot \T_{xy}\ot P^\m\ot 1_{y^\m}\T_{(xy)^\m}\ot Q\ot 1_x\T_{xy}\ot Q^\m\ar[d]|-{P\ot\T_{xy}\ot T_2}\\
				&  {R1_x\ot}P\ot \T_{xy}\ot 1_{y^\m}\T_{(xy)^\m}\ot Q\ot 1_x\T_{xy}\ot Q^\m\ot P^\m\ar@{=}[d] \\ 
				& {R1_x\ot} P\ot \T_{xy}\ot\T_{(xy)^\m}\ot Q\ot \T_{xy}\ot Q^\m\ot P^\m\ar[d]\\
				& {R1_x\ot} P\ot R1_{xy}\ot Q\ot \T_{xy}\ot Q^\m\ot P^\m\ar[d]^{P\ot T_3\ot \T_{xy}\ot Q^\m\ot P^\m}\\
				& {R1_x\ot} P\ot Q\ot R1_{xy}\ot \T_{xy}\ot Q^\m\ot P^\m\ar[d]\\
				& {R1_x\ot} P\ot Q\ot \T_{xy}\ot Q^\m\ot P^\m  }
	\end{equation*} }
	
	On the other hand,  $f_{x,y}^{P\ot Q}\circ (F_x\ot F_y)$ is given as follows 
	
	$$\xymatrix@C=-3cm{ &  \Om_x^P\ot\T_{x^{-1}}\ot \Om_x^Q\ot \Om_y^P\ot\T_{y^{-1}}\ot \Om_y^Q \ar[ldd]_-{P\ot \T_x\ot T_4\ot P\ot \T_y \ot T_5} \ar@/^3.8cm/@{-->}[ldddddd]^{f_{x,y}^{P\ot Q}\circ (F_x\ot F_y)}   \\
		& \\
		P\ot R1_{x}\ot Q\ot \T_x\ot Q^{\m}\ot P^{\m}\ot P\ot R1_{y}\ot Q\ot \T_y\ot Q^{\m}\ot P^{\m}\ar[dd]|-{P\ot T_6\ot \T_x\ot Q^{\m}\ot P^{\m}\ot P\ot T_7\ot \T_y\ot Q^{\m}\ot P^{\m}} &  \\ 
		& \\
		P\ot Q\ot R1_{x}\ot \T_x\ot Q^{\m}\ot P^{\m}\ot P\ot Q\ot R1_{y}\ot \T_y\ot Q^{\m}\ot P^{\m}\ar[d] &  \\  
		P\ot Q\ot \T_x\ot Q^{\m}\ot P^{\m}\ot P\ot Q\ot \T_y\ot Q^{\m}\ot P^{\m}\ar[d] &  \\ 
		1_xP\ot Q\ot \T_{(xy)^{-1}}\ot Q^{\m}\ot P^{\m} &  }$$

	
	We next  construct the isomorphisms. Consider the $R$-bilinear maps  
	$$f_i^x:\T_{x^\m}\ot Q\ot \T_x\ot Q^{\m}\longrightarrow R \ \mbox{e} \ g_i^x:R\longrightarrow \T_{x^\m}\ot Q\ot \T_x\ot Q^{\m}, \ \mbox{with} \ i=1,2,...,n , $$
	$$f_j^y:\T_{y^\m}\ot Q\ot \T_y\ot Q^{\m}\longrightarrow R \ \mbox{e} \ g_j^y:R\longrightarrow \T_{y^\m}\ot Q\ot \T_y\ot Q^{\m}, \ \mbox{with} \ j=1,2,...,m,$$
	such that  $\dsum_{i=1}^ng_i^xf_i^x=Id_{\T_{x^\m}\ot Q\ot \T_x\ot Q^{\m}}$ and $\dsum_{j=1}^mg_j^yf_j^y=Id_{\T_{y^\m}\ot Q\ot \T_y\ot Q^{\m}}$. 
	Denote,
	$$g_i^x(1)=\dsum_{l}u_{x^\m}^{i,l}\ot q_x^{i,l}\ot u_{x}^{i,l}\ot q_{x^\m}^{i,l}  \ \mbox{and} \ \ g_j^y(1)=\dsum_{k}\widetilde{u}_{y^\m}^{j,k}\ot \widetilde{q}_{y}^{j,k}\ot \widetilde{u}_{y}^{j,k}\ot \widetilde{q}_{y^\m}^{j,k}.$$
	Then
	$$\begin{array}{c c c l}
		T_4:& P^{-1}\ot \T_{x^\m}\ot Q\ot \T_x\ot Q^{^\m}& \longrightarrow & \T_{x^\m}\ot Q\ot \T_x\ot Q^{^\m}\ot P^{-1},\\
		& \bar{p}\ot u_{x^\m}\ot q\ot u_x\ot \bar{q} & \longmapsto & \dsum u_{x^\m}^{i,l}\ot q_{x}^{i,l}\ot u_{x}^{i,l}\ot q_{x^\m}^{i,l}\ot \bar{p}f_i^x(u_{x^\m}\ot q\ot u_x\ot \bar{q}),
	\end{array}$$
	
	$$\begin{array}{c c c l}
		T_5:& P^{-1}\ot \T_{y^\m}\ot Q\ot \T_y\ot Q^{^\m}& \longrightarrow & \T_{y^\m}\ot Q\ot \T_y\ot Q^{^\m}\ot P^{-1},\\
		& \bar{p}\ot u_{y^\m}\ot q\ot u_y\ot \bar{q} & \longmapsto & \dsum \widetilde{u}_{y^\m}^{j,k}\ot \widetilde{q}_{y}^{j,k}\ot \widetilde{u}_{y}^{j,k}\ot \widetilde{q}_{y^\m}^{j,k}\ot \bar{p}f_j^y(u_{y^\m}\ot q\ot u_y\ot \bar{q}),
	\end{array}$$
	and
	$${\small \begin{array}{ c c c l}
			T_1:&  \T_{x^{-1}}\ot Q\ot \T_x\ot Q^{^\m}\ot P\ot \T_y\ot P^{\m}\ot \T_{y^{-1}}& \rightarrow & P\ot \T_y\ot P^{\m}\ot \T_{y^{-1}}\ot \T_{x^{-1}}\ot Q\ot \T_x\ot Q^{^\m},\\
			&  u_{x^\m}\ot q\ot u_x\ot \bar{q}\ot p \ot u_y\ot \bar{p}\ot u_{y^\m} & \mapsto & \dsum f_i^x(u_{x^\m}\ot q\ot u_x\ot \bar{q})p \ot u_y\ot \bar{p}\ot u_{y^\m}\ot g_i^x(1).
	\end{array}}$$
	
	In order to construct   $T_2$,  consider the $R$-bilinear maps   $f_{i,j}^{xy}:1_{y^\m}\T_{{(xy)^\m}}\ot Q\ot 1_{x}\T_{xy}\ot Q^{-1} \longrightarrow R$ and $g_{i,j}^{xy}:R\longrightarrow 1_{y^\m}\T_{xy}\ot Q\ot 1_{x}\T_{xy}\ot Q^{-1}$ 
	{ defined as follows:}
	$$\xymatrix@C=4cm{ 1_{y^\m}\T_{{(xy)^\m}}\ot Q\ot 1_{x}\T_{xy}\ot Q^{-1}\ar[r]^{(f_{y^\m,x^\m}^\T)^\m\ot Q\ot (f_{x,y}^\T)^\m\ot Q^\m}\ar@{-->}@/_1cm/[rddd]_{f_{i,j}^{xy}} & \T_{y^\m}\ot \T_{x^{-1}}\ot Q\ot \T_x\ot \T_y\ot Q^{\m}\ar[d]\\
		& \T_{y^\m}\ot \T_{x^{-1}}\ot Q\ot \T_x\ot{ Q^{-1}\ot Q}\ot  \T_y\ot Q^{\m}\ar[d]^{\T_{y^\m}\ot f_i^x\ot Q\ot \T_y\ot Q}\\
		&  \T_{y^\m}\ot Q^{-1}\ot  \T_y\ot Q^{\m}\ar[d]^{f_j^y}\\
		& R }$$ 
	{ and}
	\begin{equation*}
		\xymatrix{ R\ar[r]^{g_j^y\ \ \ \ \ \ \ \ \ \ \ \ \ }\ar@{-->}@/_1.5cm/[rrddd]_{g_{i,j}^{xy}} & \T_{y^{-1}} \ot Q\ot \T_y\ot Q^{\m} \ar[r] &  \T_{y^{-1}}\ot R\ot Q\ot \T_y\ot Q^{\m}\ar[d]^{\T_{y^\m}\ot g_i^x\ot Q\ot \T_y\ot Q^\m}\\
			& &  \T_{y^{-1}} \ot \T_{x^{-1}} \ot Q\ot \T_x\ot Q^{\m} \ot Q\ot \T_y\ot Q^{\m}\ar[d]\\
			& & \T_{y^{-1}}\ot \T_{x^{-1}}\ot Q\ot \T_x\ot \T_y\ot Q^{\m}\ar[d]^{f_{y^\m,x^\m}^\T\ot Q\ot f_{x,y}^\T\ot Q^\m}\\
			& & 1_{y^\m}\T_{(xy)^{-1}}\ot Q\ot 1_x\T_{xy}\ot Q^{\m}}
	\end{equation*}
	
	Then, 
	\begin{equation}
		f_{i,j}^{xy}(u_{(xy)^{-1}}\ot q\ot u_{xy}\ot \bar{q})=\dsum f_j^y(u_{y^\m}f_i^x(u_{x^\m}\ot q \ot u_x\ot {\bar{q}_t})\ot {q_t}\ot u_{y}\ot \bar{q}),
	\end{equation}  
	where $\dsum (u_{y^\m}\stackrel{\T}{\circ} u_{x^\m})=u_{(xy)^{-1}}$, $\dsum (u_x\stackrel{\T}{\circ}u_y)=u_{xy}$ and $\dsum_{t}\Rr'(\bar{q}_t\ot q_t)=1,$ { and} 
	
		\begin{equation}
			g_{i,j}^{xy}(1)=\dsum (\widetilde{u}_{y^\m}^{j,k}\stackrel{\T}{\circ} {u_{x^\m}^{i,l}})\ot{ q_{x}^{i,l}}\ot ({u_{x}^{i,l}}\Rr'( {q_{x^\m}^{i,l}}\ot  \widetilde{q}_{y}^{j,k})\stackrel{\T}{\circ} \widetilde{u}_{y}^{j,k})\ot \widetilde{q}_{y^\m}^{j,k}.
		\end{equation}
		It is easy to see that  $\dsum_{i,j}g_{i,j}^{xy}f_{i,j}^{xy}=Id_{1_{y^\m}\T_{(xy)^{-1}}\ot Q\ot 1_x\T_{xy}\ot Q^{\m}}$. Then the isomorphism   $T_2$ is given by 
		$$\begin{array}{c c c l}
			T_2:& P^{\m}\ot \T_{(xy)^{-1}}\ot Q\ot \T_{xy}\ot Q^{\m}& \longrightarrow&   \T_{(xy)^{-1}}\ot Q\ot \T_{xy}\ot Q^{\m}\ot P^{\m},\\
			& \bar{p}\ot u_{(xy)^{\m}}\ot q \ot u_{xy}\ot \bar{q} & \longmapsto & \dsum g_{i,j}^{xy}(1)\ot \bar{p}f_{i,j}^{xy}(u_{(xy)^{\m}}\ot q \ot u_{xy}\ot \bar{q}).
		\end{array}$$
		
		By Example~\ref{execomRe1}, the   isomorphisms  $T_3, T_6$ and $T_7$ are  defined by:
		%
		\begin{equation*}
			\begin{tabular}{c c c}
				$\begin{array}{c c l}
					R1_{xy}\ot Q & \stackrel{T_3}{\longrightarrow}&  Q\ot R1_{xy}\\
					r \ot q & \longmapsto & rq\ot 1_{xy}
				\end{array}$,	& $\begin{array}{ c c l}
					R1_x\ot Q & \stackrel{T_6}{\longrightarrow}&  Q\ot R1_x\\
					r \ot q & \longmapsto & rq\ot 1_x
				\end{array}$,  & $\begin{array}{c c c l}
					R1_y\ot Q & \stackrel{T_7}{\longrightarrow}&  Q\ot R1_y\\
					r \ot q & \longmapsto & rq\ot 1_y.
				\end{array}$
			\end{tabular}
		\end{equation*}
		We first  compute  $f_{x,y}^{P\ot Q}\circ (F_{x}\ot F_{y}):$ {Given $p_1,p_2 \in P, q_1,q_2\in Q,\bar{p_1},\bar{p_2}\in P^{-1},\bar{q_1},\bar{q_2}\in Q^{-1}, u_x,u_x'\in \T_x$ and $u_y,u_y'\in \T_y$, denote $a_i^x=f_i^x(u_{x^\m}\ot q_1\ot u_x' \ot \bar{q_1})$ and $a_y^y=f_j^y(u_{y^\m}\ot q_2\ot u_y' \ot \bar{q_2})$.} {  Then}

		\begin{eqnarray*}
			& & p_1\ot u_x\ot \bar{p}_1\ot u_{x^\m}\ot q_1\ot u_x'\ot \bar{q}_1\ot p_2\ot u_y\ot \bar{p}_2\ot u_{y^\m}\ot q_2\ot u_y'\ot \bar{q}_2\\
			&\stackrel{T_4,T_5}{\mapsto} & \dsum_{i,l,j,k} p_1\ot u_x \ot u_{x^\m}^{i,l}\ot q_{x}^{i,l}\ot u_{x}^{i,l}\ot q_{x^\m}^{i,l}\ot \bar{p}_1{a_i^x}\ot p_2\ot u_y \ot  \widetilde{u}_{y^\m}^{j,k}\ot \widetilde{q}_{y}^{j,k}\ot \widetilde{u}_{y}^{j,k}\ot \widetilde{q}_{y^\m}^{j,k}\ot \bar{p}_2{a_j^y}\\
			& \mapsto & \dsum_{i,l,j,k} p_1\ot (u_x \stackrel{\T}{\circ} u_{x^\m}^{i,l}) \ot q_{x}^{i,l}\ot u_{x}^{i,l}\ot q_{x^\m}^{i,l}\ot \bar{p}_1a_i^x\ot p_2\ot  (u_y\stackrel{\T}{\circ}\widetilde{u}_{y^\m}^{j,k})\ot \widetilde{q}_{y}^{j,k}\ot \widetilde{u}_{y}^{j,k}\ot \widetilde{q}_{y^\m}^{j,k}\ot \bar{p}_2a_j^y\\
			& \stackrel{T_6,T_7}{\mapsto} & \dsum_{i,l,j,k} p_1\ot (u_x \stackrel{\T}{\circ} u_{x^\m}^{i,l})q_{x}^{i,l}\ot 1_x\ot u_{x}^{i,l}\ot q_{x^\m}^{i,l}\ot \bar{p}_1a_i^x\ot p_2\ot (u_y\stackrel{\T}{\circ}\widetilde{u}_{y^\m}^{j,k})\widetilde{q}_{y}^{j,k}\ot 1_y\ot \widetilde{u}_{y}^{j,k}\ot \widetilde{q}_{y^\m}^{j,k}\ot \bar{p}_2a_j^y\\
			& \mapsto & \dsum_{i,l,j,k} p_1\ot (u_x \stackrel{\T}{\circ} u_{x^\m}^{i,l})q_{x}^{i,l}\ot u_{x}^{i,l}\ot q_{x^\m}^{i,l}\ot \bar{p}_1a_i^x\ot p_2\ot (u_y\stackrel{\T}{\circ}\widetilde{u}_{y^\m}^{j,k})\widetilde{q}_{y}^{j,k}\ot \widetilde{u}_{y}^{j,k}\ot \widetilde{q}_{y^\m}^{j,k}\ot \bar{p}_2a_j^y\\
			& \mapsto & \dsum_{i,l,j,k} p_1\ot (u_x \stackrel{\T}{\circ} u_{x^\m}^{i,l})q_{x}^{i,l} \ot \left( u_{x}^{i,l}\Rr'(\widetilde{q}_{x^\m}^{i,l}\Rr(\bar{p}_1a_i^x\ot p_2)\ot (u_y\stackrel{\T}{\circ}\widetilde{u}_{y^\m}^{j,k})\widetilde{q}_{y}^{j,k})\stackrel{\T}{\circ} \widetilde{u}_{y}^{j,k}\right) \ot \widetilde{q}_{y^\m}^{j,k}\ot \bar{p}_2a_j^y\\
			& = & \dsum_{i,l,j,k} p_1\ot (u_x \stackrel{\T}{\circ} u_{x^\m}^{i,l})q_{x}^{i,l}
			\ot \left( u_{x}^{i,l}a \stackrel{\T}{\circ} \widetilde{u}_{y}^{j,k}\right) \ot \widetilde{q}_{y^\m}^{j,k}\ot \bar{p}_2f_j^y(u_{y^\m}\ot q_2\ot u'_y\ot \bar{q}_2),
		\end{eqnarray*}
		{	{ where} $a=\Rr'(\widetilde{q}_{x^\m}^{i,l}\Rr(\bar{p}_1a_i^x\ot p_2)\ot (u_y\stackrel{\T}{\circ}\widetilde{u}_{y^\m}^{j,k})\widetilde{q}_{y}^{j,k}).$}
		
		For $F_{xy}\circ f_{x,y}^{PQ}$, we have:
		
		\begin{eqnarray}
			\nonumber
			& & p_1\ot u_x\ot \bar{p}_1\ot u_{x^\m}\ot q_1\ot u_x'\ot \bar{q}_1\ot p_2\ot u_y\ot \bar{p}_2\ot u_{y^\m}\ot q_2\ot u_y'\ot \bar{q}_2\\ \nonumber
			& \stackrel{T_1}{\mapsto} &  \dsum_{h,l'} p_1\ot u_x\ot \bar{p}_1\ot   {a_h^x}p_2 \ot u_y\ot \bar{p}_2\ot u_{y^\m}  \ot u_{x^\m}^{h,l'}\ot q_{x}^{h,l'} \ot u_{x}^{h,l'}\ot q_{x^\m}^{h,l'}  \ot    q_2\ot u_y'\ot \bar{q}_2\\ \nonumber
			& \mapsto & \dsum_{h,l'}  p_1\ot (u_x\Rr(\bar{p}_1{a_h^x}\ot p_2) \stackrel{\T}{\circ} u_y)\ot \bar{p}_2\ot (u_{y^\m}  \stackrel{\T}{\circ} u_{x^\m}^{h,l'} )\ot q_{x}^{h,l'} \ot (u_{x}^{h,l'}\Rr' (q_{x^\m}^{h,l'}  \ot    q_2)\stackrel{\T}{\circ} u_y')\ot \bar{q}_2\\ \nonumber
			& \stackrel{T_2}{\mapsto} & \dsum_{h,l',i,l,j,k,t}    p_1\ot (u_x \Rr(\bar{p}_1{a_h^x}\ot   p_2) \stackrel{\T}{\circ} u_y)\ot  (\widetilde{u}_{y^\m}^{j,k}\stackrel{\T}{\circ} u_{x^\m}^{i,l})\ot q_{x}^{i,l}\ot (u_{x}^{i,l}\Rr'( q_{x^\m}^{i,l}\ot  \widetilde{q}_{y}^{j,k})\stackrel{\T}{\circ} \widetilde{u}_{y}^{j,k})\ot \widetilde{q}_{y^\m}^{j,k}\\ \nonumber
			& & \ot\bar{p}_2 f_j^y(u_{y^\m}f_i^x(u_{x^\m}^{h,l'} \ot q_{x}^{h,l'}\ot u_{x}^{h,l'}\Rr' (q_{x^\m}^{h,l'}  \ot    q_2)\ot \bar{q}_t)\ot q_t\ot  u_y'\ot \bar{q}_2 ) \\ \nonumber
			& \mapsto & \dsum_{h,l',i,l,j,k,t}    p_1\ot (u_x \Rr(\bar{p}_1{a_h^x}\ot   p_2) \stackrel{\T}{\circ} u_y){\stackrel{\T}{\circ} } (\widetilde{u}_{y^\m}^{j,k}\stackrel{\T}{\circ} u_{x^\m}^{i,l})\ot q_{x}^{i,l}\ot (u_{x}^{i,l}\Rr'( q_{x^\m}^{i,l}\ot  \widetilde{q}_{y}^{j,k})\stackrel{\T}{\circ} \widetilde{u}_{y}^{j,k})\ot \widetilde{q}_{y^\m}^{j,k}\\ \nonumber
			& & \ot\bar{p}_2 f_j^y(u_{y^\m}f_i^x(u_{x^\m}^{h,l'} \ot q_{x}^{h,l'}\ot u_{x}^{h,l'}\Rr' (q_{x^\m}^{h,l'}  \ot    q_2)\ot \bar{q}_t)\ot q_t\ot  u_y'\ot \bar{q}_2 ) \\ \nonumber
						\end{eqnarray}
		\begin{eqnarray}
			& \stackrel{T_3}{\mapsto} & \dsum_{h,l',i,l,j,k,t}    p_1\ot (u_x \Rr(\bar{p}_1{a_h^x}\ot   p_2) \stackrel{\T}{\circ} u_y){\stackrel{\T}{\circ} } (\widetilde{u}_{y^\m}^{j,k}\stackrel{\T}{\circ} u_{x^\m}^{i,l}) q_{x}^{i,l}\ot 1_{xy}\ot (u_{x}^{i,l}\Rr'( q_{x^\m}^{i,l}\ot  \widetilde{q}_{y}^{j,k})\stackrel{\T}{\circ} \widetilde{u}_{y}^{j,k})\\ \nonumber
			& & \ot \widetilde{q}_{y^\m}^{j,k}\ot\bar{p}_2 f_j^y(u_{y^\m}f_i^x(u_{x^\m}^{h,l'} \ot q_{x}^{h,l'}\ot u_{x}^{h,l'}\Rr' (q_{x^\m}^{h,l'}  \ot    q_2)\ot \bar{q}_t)\ot q_t\ot  u_y'\ot \bar{q}_2 ) \\ \nonumber
			& \mapsto & \dsum_{h,l',i,l,j,k,t}    p_1\ot (u_x \Rr(\bar{p}_1{a_h^x}\ot   p_2) \stackrel{\T}{\circ} (u_y{\stackrel{\T}{\circ} } \widetilde{u}_{y^\m}^{j,k})\stackrel{\T}{\circ} u_{x^\m}^{i,l}) q_{x}^{i,l}\ot (u_{x}^{i,l}\Rr'( q_{x^\m}^{i,l}\ot  \widetilde{q}_{y}^{j,k})\stackrel{\T}{\circ} \widetilde{u}_{y}^{j,k})\\ 
			& & \ot \widetilde{q}_{y^\m}^{j,k}\ot\bar{p}_2 f_j^y(u_{y^\m}f_i^x(u_{x^\m}^{h,l'} \ot q_{x}^{h,l'}\ot u_{x}^{h,l'}\Rr' (q_{x^\m}^{h,l'}  \ot    q_2)\ot \bar{q}_t)\ot q_t\ot  u_y'\ot \bar{q}_2 )  \label{fxyPotQFxFyparte1}
		\end{eqnarray}
		
		{	Notice that 
			\begin{eqnarray}
				\nonumber
				& & \dsum_{t}\bar{p}_2 f_j^y(u_{y^\m}f_i^x(u_{x^\m}^{h,l'} \ot q_{x}^{h,l'}\ot u_{x}^{h,l'}\Rr' (q_{x^\m}^{h,l'}  \ot    q_2)\ot \bar{q}_t)\ot q_t\ot  u_y'\ot \bar{q}_2 )\\ \nonumber
				&=& \dsum_{t}\bar{p}_2 f_j^y(u_{y^\m}f_i^x(u_{x^\m}^{h,l'} \ot q_{x}^{h,l'}\ot u_{x}^{h,l'}\ot q_{x^\m}^{h,l'}  \Ll'   ( q_2\ot\bar{q}_t))\ot q_t\ot  u_y'\ot \bar{q}_2 )\\ \nonumber
				&=&\dsum_{t} \bar{p}_2 f_j^y(u_{y^\m}f_i^x(u_{x^\m}^{h,l'} \ot q_{x}^{h,l'}\ot u_{x}^{h,l'}\ot q_{x^\m}^{h,l'}  )\ot   q_2\Rr'(\bar{q}_t\ot q_t)\ot  u_y'\ot \bar{q}_2 )\\ \nonumber
				&=&\dsum_{t} \bar{p}_2 f_j^y(u_{y^\m}f_i^x(g_h^x(1) )\ot  q_2\Rr'(\bar{q}_t\ot q_t)\ot  u_y'\ot \bar{q}_2 )\\
				&=& \bar{p}_2 f_j^y(u_{y^\m}f_i^x(g_h^x(1) )\ot  q_2\ot  u_y'\ot \bar{q}_2 ). \label{fixghx}
		\end{eqnarray}}
		
		{ Moreover, we have
			\begin{eqnarray}
				\dsum_{h,l'}a_h^xf_i^x(g_h^x(1))=a_i^x. \label{aix}
			\end{eqnarray}
			Indeed,
			\begin{eqnarray*}
				\dsum_{h,l'}a_h^xf_i^x(g_h^x(1))& = & \dsum_{h,l'} f_i^x(g_h^x(a_h^x))=\dsum_{h,l'}f_i^x(g_h^x(f_h^x(u_{x^\m}\ot q_1\ot u_x' \ot \bar{q_1})))\\
				& = & f_i^x(u_{x^\m}\ot q_1\ot u_x' \ot \bar{q_1}) = a_i^x.
		\end{eqnarray*}}
		
		{  { Substituting (\ref{fixghx}) in (\ref{fxyPotQFxFyparte1}) and using (\ref{uxercomutam}) and (\ref{aix}), we obtain that}
			
			\begin{eqnarray*}
				& & F_{xy}\circ f_{x,y}^{PQ}(p_1\ot u_x\ot \bar{p}_1\ot u_{x^\m}\ot q_1\ot u_x'\ot \bar{q}_1\ot p_2\ot u_y\ot \bar{p}_2\ot u_{y^\m}\ot q_2\ot u_y'\ot \bar{q}_2)\\ 
				& = & 	\dsum_{h,l',i,l,j,k,t}    p_1\ot (u_x \Rr(\bar{p}_1{a_h^x}\ot   p_2) \stackrel{\T}{\circ} (u_y{\stackrel{\T}{\circ} } \widetilde{u}_{y^\m}^{j,k})\stackrel{\T}{\circ} u_{x^\m}^{i,l}) q_{x}^{i,l}\ot (u_{x}^{i,l}\Rr'( q_{x^\m}^{i,l}\ot  \widetilde{q}_{y}^{j,k})\stackrel{\T}{\circ} \widetilde{u}_{y}^{j,k})\\ 
				& & \ot \widetilde{q}_{y^\m}^{j,k}\ot\bar{p}_2 f_j^y(u_{y^\m}f_i^x(g_h^x(1) )\ot  q_2\ot  u_y'\ot \bar{q}_2 )\\
				& \stackrel{(\ref{uxercomutam})}{=} & \dsum_{h,l',i,l,j,k,t}    p_1\ot (u_x \Rr(\bar{p}_1{a_h^x}f_i^x(g_h^x(1) )\ot   p_2) \stackrel{\T}{\circ} (u_y{\stackrel{\T}{\circ} } \widetilde{u}_{y^\m}^{j,k})\stackrel{\T}{\circ} u_{x^\m}^{i,l}) q_{x}^{i,l}\ot (u_{x}^{i,l}\Rr'( q_{x^\m}^{i,l}\ot  \widetilde{q}_{y}^{j,k})\stackrel{\T}{\circ} \widetilde{u}_{y}^{j,k})\\ 
				& & \ot \widetilde{q}_{y^\m}^{j,k}\ot\bar{p}_2 f_j^y(u_{y^\m}\ot  q_2\ot  u_y'\ot \bar{q}_2 )\\
				& \stackrel{(\ref{aix})}{=} & \dsum_{i,l,j,k,t}    p_1\ot (u_x \Rr(\bar{p}_1{a_i^x}\ot   p_2) \stackrel{\T}{\circ} (u_y{\stackrel{\T}{\circ} } \widetilde{u}_{y^\m}^{j,k})\stackrel{\T}{\circ} u_{x^\m}^{i,l}) q_{x}^{i,l}\ot (u_{x}^{i,l}\Rr'( q_{x^\m}^{i,l}\ot  \widetilde{q}_{y}^{j,k})\stackrel{\T}{\circ} \widetilde{u}_{y}^{j,k})\\ 
				& & \ot \widetilde{q}_{y^\m}^{j,k}\ot\bar{p}_2 f_j^y(u_{y^\m}\ot  q_2\ot  u_y'\ot \bar{q}_2 )\\
				& \stackrel{(\ref{aix})}{=} & \dsum_{i,l,j,k,t}    p_1\ot (u_x b\stackrel{\T}{\circ} u_{x^\m}^{i,l}) q_{x}^{i,l}\ot (u_{x}^{i,l}\Rr'( q_{x^\m}^{i,l}\ot  \widetilde{q}_{y}^{j,k})\stackrel{\T}{\circ} \widetilde{u}_{y}^{j,k}) \ot \widetilde{q}_{y^\m}^{j,k}\ot\bar{p}_2 f_j^y(u_{y^\m}\ot  q_2\ot  u_y'\ot \bar{q}_2 ),
			\end{eqnarray*}
			where $b=\Rr(\bar{p}_1{a_i^x}\ot   p_2)(u_y{\stackrel{\T}{\circ} } \widetilde{u}_{y^\m}^{j,k})\in R.$}

		{	{ Since}  $g_i^x=\dsum_{i,l}u_{x^\m}^{i,l}\ot q_{x}^{i,l}\ot u_{x}^{i,l}\ot q_{x^\m}^{i,l}\in C_{\T_{x^\m}\ot Q\ot \T_x\ot Q^{\m}}(R)$, {then} }
		\begin{eqnarray*}
			& & (u_x  \stackrel{\T}{\circ} b u_{x^\m}^{i,l})q_{x}^{i,l}\ot (u_{x}^{i,l}\Rr'( q_{x^\m}^{i,l}\ot  \widetilde{q}_{y}^{j,k})\stackrel{\T}{\circ} \widetilde{u}_{y}^{j,k})\\
			& = &  (Q\ot f_{x,y}^\T)\circ(f_{x,x^{\m}}^\T\ot Q \ot \T_x {\ot} \Rr'\ot \T_y)(u_x \ot  bu_{x^\m}^{i,l}\ot q_{x}^{i,l}\ot u_{x}^{i,l}\ot q_{x^\m}^{i,l}\ot  \widetilde{q}_{y}^{j,k}\ot  \widetilde{u}_{y}^{j,k})\\
			& = &  (Q\ot f_{x,y}^\T)\circ(f_{x,x^{\m}}^\T\ot Q \ot \T_x {\ot}\Rr'\ot \T_y)(u_x \ot u_{x^\m}^{i,l}\ot q_{x}^{i,l}\ot u_{x}^{i,l}\ot q_{x^\m}^{i,l}b\ot  \widetilde{q}_{y}^{j,k}\ot  \widetilde{u}_{y}^{j,k})\\
			& = &   (u_x  \stackrel{\T}{\circ} u_{x^\m}^{i,l})q_{x}^{i,l}\ot (u_{x}^{i,l}\Rr'( q_{x^\m}^{i,l}\Rr(\bar{p}_1a_i^x\ot   p_2)(u_y\stackrel{\T}{\circ}  \widetilde{u}_{y^\m}^{j,k})\ot  \widetilde{q}_{y}^{j,k})\stackrel{\T}{\circ} \widetilde{u}_{y}^{j,k})\\
			& = &   (u_x  \stackrel{\T}{\circ} u_{x^\m}^{i,l})q_{x}^{i,l}\ot (u_{x}^{i,l}\Rr'( q_{x^\m}^{i,l}\Rr(\bar{p}_1a_i^x\ot   p_2)\ot  (u_y\stackrel{\T}{\circ}  \widetilde{u}_{y^\m}^{j,k})\widetilde{q}_{y}^{j,k})\stackrel{\T}{\circ} \widetilde{u}_{y}^{j,k})\\
			& = & (u_x  \stackrel{\T}{\circ} u_{x^\m}^{i,l})q_{x}^{i,l}\ot (u_{x}^{i,l}a\stackrel{\T}{\circ} \widetilde{u}_{y}^{j,k})
		\end{eqnarray*}
		Hence,  
		\begin{eqnarray*}
			& & (F_{xy}\circ f_{x,y}^{PQ})(p_1\ot u_x\ot \bar{p}_1\ot u_{x^\m}\ot q_1\ot u_x'\ot \bar{q}_1\ot p_2\ot u_y\ot \bar{p}_2\ot u_{y^\m}\ot q_2\ot u_y'\ot \bar{q}_2)\\
			& = & \dsum_{i,l,j,k}  p_1\ot (u_x  \stackrel{\T}{\circ} u_{x^\m}^{i,l})q_{x}^{i,l}\ot (u_{x}^{i,l}a\stackrel{\T}{\circ} \widetilde{u}_{y}^{j,k})  \ot \widetilde{q}_{y^\m}^{j,k}\ot\bar{p}_2 f_j^y(u_{y^\m}\ot q_2\ot  u_y'\ot \bar{q}_2 ).
		\end{eqnarray*}
		Therefore,  $\mathcal{L}$ is a group homomorphism.  
	\end{dem}

	By (\ref{PicRePicsRinvarinates}), if $[P]\in \Pic_\Z(R)\cap \Pics_\Z(R)^{\al^*}$, then $[\D(\Om^{P})]\in \C_0(\T/R)$. We define the homomorphism  $\varphi_3:\Pic_\Z(R)\cap\Pics_\Z(R)^{\al^*}\longrightarrow \C_0(\T/R)$ as a restriction of  $\mathcal{L},$ that is,
	\begin{equation*}
		\begin{array}{c c c l}
			\varphi_3: & \Pic_\Z(R)\cap \Pics_\Z(R)^{\al^*} & \longrightarrow & \C_0(\T/R)\\
			& [P] & \longmapsto & [\D(\Om^P)]
		\end{array}.
	\end{equation*}

{ For reader's convenience we separate in the next  lemma a fact observed in 
	\cite[p. 161]{el2010extended}. }

\begin{lem}\label{estruturadedeltathetabimodulo}
	Let $R\subseteq S$ be a ring extension  with the same unity and let $P$ be an  $R$-bimodule. Denote by $\mu$ the multiplication in $S$. 
	\begin{itemize}
		\item[(i)] Let $\rho:S\ot_RP\longrightarrow P\ot_RS$  be an $R$-bilinear map such that
		\begin{equation}
			(P\ot \mu)\circ (\rho\ot S)\circ (S\ot \rho)=\rho \circ (\mu \ot P), \label{eq5.1}
		\end{equation}
		\begin{equation}
			\rho(1\ot p)=p\ot 1, \ \ \mbox{for all} \ p\in P. \label{eq5.2}
		\end{equation}
		If $M$ is a unital left  $S$-module, then $P\ot_RM$ possesses a structure of a left  $S$-module, given by 
		\begin{equation*}
			s*(p\ot m)=\dsum_{i=1}^np_i\ot s_i m,
		\end{equation*}
		where $\rho(s\ot p)=\dsum_{i=1}^np_i\ot s_i$. Moreover, if   $M$ is a unital  $S$-bimodule, then $P\ot_RM$ is
		a  unital   $S$-bimodule, where the structure of a right  $S$-module is given by  the structure of  $M.$
	 	\item[(ii)] Let $\rho':P\ot_RS\longrightarrow S\ot_RP$  be an $R$-bilinear map such that 
		\begin{equation}
			(\mu \ot P)\circ (S\ot \rho1')\circ (\rho'\ot S)=\rho' \circ (P\ot \mu), \label{eq5.1linha}
		\end{equation}
		\begin{equation}
			\rho(p\ot 1)=1\ot p, \ \ \mbox{for all} \ p\in P. \label{eq5.2linha}
		\end{equation}
		If $M$ is a right   $S$-module, then   $M\ot_RRP$ possesses a structure of a right  $S$-module given by 
		\begin{equation*}
			(m\ot p)*s=\dsum_{i=1}^nms_i\ot p_i,
		\end{equation*}
		where $\rho'(p\ot s)=\dsum_{i=1}^ns_i\ot p_i$. Moreover, if  $M$ is a unital  $S$-bimodule, then $M\ot_RP$ 
		is a unital  $S$-bimodule, where the structure of the right  $S$-module is given by the structure of   $M$.
		\end{itemize}
\end{lem}

\begin{exe}\label{acaosobrePZM}
let $P$ be a central  $\Z$-bimodule. Since $R$ is a central  $\Z$-bimodule too, we have  a  $\Z$-bimodule isomorphism  $\rho:R\ot_\Z P\longrightarrow P\ot_\Z R$ defined by  $\rho(r\ot p)=p\ot r$. Let us verify that $\rho$ satisfies the   conditions  of  Lemma~\ref{estruturadedeltathetabimodulo}.
If $r, r' \in R$ and $p \in P$, then
\begin{eqnarray*}
	(P\ot m)\circ (\rho\ot R)\circ (R\ot \rho)(r\ot r' \ot p) & = & (P\ot m)\circ (\rho\ot R)(r\ot p \ot r')\\
	& = & (P\ot m)(p\ot r\ot r')\\
	& = & p\ot rr'= \rho(rr'\ot p)\\
	& = & \rho\circ (m\ot P)(r\ot r'\ot p).
\end{eqnarray*} 
Therefore, $\rho$ satisfies (\ref{eq5.1}). Condition  (\ref{eq5.2}) directly follows from the definition of  $\rho$.
For any   $R$-bimodule  $M,$ by   Lemma~\ref{estruturadedeltathetabimodulo}, the tensor product   $P\ot_\Z M$ has a structure of an   $R$-bimodule defined by  
$$r_1*(p\ot m)*r_2=p\ot r_1mr_2,$$
for all  $p \in P, m \in M$ and $r_1,r_2\in R.$

\end{exe}

\begin{lem}\label{lemaPSMisoPM} Let $R\subseteq S$ be a ring extension with the same unity,  $P$ a unital 
$R$-bimodule and  $\rho:S\ot_RP\longrightarrow P\ot_RS$ an $R$-bimodule homomorphism which satisfies conditions 
(\ref{eq5.1}) and (\ref{eq5.2}). If $M$ is an  $S$-bimodule, then there exists an $S$-bimodule  isomorphism 
\begin{equation*}
	\begin{array}{c c c l}
		\eta: & (P\ot_RS)\ot_SM & \longrightarrow & P\ot_RM\\
		& (p\ot_R s)\ot_Sm & \longmapsto & p\ot_Rs m
	\end{array},
\end{equation*}
where $P\ot_RS$ and $P\ot_RM$ are $S$-bimodules with the actions given in  Lemma~\ref{estruturadedeltathetabimodulo}.
\end{lem}
\begin{dem} Clearly, $\eta$ is a well-defined right $S$-linear map.   Given,  $s,s'\in S$, $p\in P$ and $m\in M$, denote $\rho(s'\ot_R p)=\dsum_{i=1}^np_i\ot_R s_i'$. Then 
\begin{eqnarray*}
	\eta(s'*((p\ot_R s)\ot_S m))&=&\eta\left( \dsum_{i=1}^np_i\ot_R s_i's\ot_S m\right)=\dsum_{i=1}^np_i\ot_R (s_i's) m \\
	& = & \dsum_{i=1}^np_i\ot_R s_i'(s m)= s'*(p\ot_R sm)\\
	& = & s' \eta((p\ot_R s)\ot_S  m).
\end{eqnarray*}
Hence $\eta$ is $S$-bilinear. Its inverse is given by:
\begin{equation*}
	\begin{array}{c c c l}
		\eta^\m: & P\ot_RM & \longrightarrow & (P\ot_RS)\ot_SM,\\
		& p\ot_R m & \longmapsto & p\ot_R 1\ot_S m.
	\end{array}
\end{equation*}
Indeed, if  $p\in P,m\in M$ and $s\in S$, then 
\begin{equation*}
	\eta^\m(\eta(p\ot_Rs\ot_Sm))=\eta^\m (p\ot_Rs m)= p\ot_R1\ot_Ss m= p\ot_Rs\ot_S m.
\end{equation*}
On the other hand, 
\begin{equation*}
	\eta(\eta^\m(p\ot_Rm))=\eta(p\ot_R1\ot_Sm)=p\ot_R1 m=p\ot_Rm.
\end{equation*}
Consequently, $\eta$ is an $S$-bimodule isomorphism. 
\end{dem}

\begin{lem}\label{lemacomoverlinerho} Let $R\subseteq S$ be an extension of rings with the same unity,  let $[P]\in \Pic(R)$ and suppose that $\rho: S\ot_RP\longrightarrow P\ot_RS$ is an $R$-bimodule isomorphism, which satisfies the conditions of  Lemma~\ref{estruturadedeltathetabimodulo}. Denote by $[P^\m]$ the  inverse of $[P]$ in $\Pic(R)$ and let  $P\ot_RP^\m\stackrel{\Ll}{\longrightarrow} R\stackrel{\Rr}{\longleftarrow}P^\m\ot_RP$ be isomorphisms  of $R$-bimodules. 
\begin{itemize}
	\item[(i)]   Let $\overline{\rho}:S\ot_RP^\m\longrightarrow P^\m\ot_RS$ be the
	$R$-bimodule isomorphism  induced by  $\overline{\rho}=P^\m\ot \rho^\m\ot P^\m$, that is,  
	\begin{equation*}
		\xymatrix{S\ot_RP^\m\ar[r]\ar@{-->}@/_0.5cm/[rrdd]_{\overline{\rho}} & R\ot_RS\ot_R P^\m\ar[r]&  P^\m\ot_RP\ot_RS\ot_RP^\m\ar[d]^{P^\m\ot \rho^\m \ot P^\m} \\
			& &  P^\m\ot_RS\ot_RP\ot_RP^\m\ar[d]\\
			& &  P^\m\ot_R S. }
	\end{equation*}
Suppose that   $\overline{\rho}$ satisfies the conditions of  Lemma~\ref{estruturadedeltathetabimodulo} and endow   $P^\m\ot_RS$ with the structure of an  $S$-bimodule determined by $\overline{\rho}$ as in  Lemma~\ref{estruturadedeltathetabimodulo}. 
Let $\Psi:P\ot_R(P^\m\ot_RS)\longrightarrow S$ be the  isomorphism defined by 
\begin{equation*}
	\Psi: P\ot_R(P^\m\ot_RS)\stackrel{P\ot \overline{\rho}^\m}{\longrightarrow} P\ot_R(S\ot_RP^\m)\stackrel{\rho^\m\ot P^\m}{\longrightarrow} S\ot_RP\ot P^\m\stackrel{S\ot\Ll}{\longrightarrow} S,
\end{equation*} 
that is, $\Psi=(S\ot \Ll)\circ (\rho^\m\ot P^\m)\circ (P\ot \overline{\rho}^\m).$ Then 
$$\Psi(p\ot \bar{p}\ot s)=\Ll(p\ot \bar{p})s,$$
for all  $s \in S$, $p\in P$ and $\bp\in P^\m$.

\item[(ii)] Let $s \in S$, $p\in P$ and $\bp \in P^\m$. Write  $\dsum_{k=1}^l\Ll(p_k\ot \bp_k)=1,$ and denote $\rho(s\ot p)=\dsum_{i=1}^np_i\ot s_i$ and $\rho^\m(p_k\ot s_i)=\dsum_{j=1}^ms_{i,j}\ot \bp_{k,j}$. Then 
\begin{equation}
	s\Ll(p\ot \bp)=\dsum_{i,j,k}\Ll(p_i\ot \bp_k)s_{i,j}\Ll(p_{k,j}\ot \bp). \label{slpbp}
\end{equation}
\item[(iii)] $\Psi$ is an isomorphism of  $S$-bimodules. 
\item[(iv)] $[P\ot_RS]\in \Pic(S)$.
\end{itemize}
\end{lem}
\begin{dem}
{ Let} $s\in S,$  $\bp \in P^\m$ 
and take $p_l \in P,$  $\bp_l \in P^\m$, $l=1,2,...,n$, such that  $\dsum_{l=1}^n\Rr(\bp_l\ot p_l)=1$. Denoting,  $\rho^\m(p_l\ot s)=\dsum_{i=1}^ms_i\ot p_{l,i}$, we have 
\begin{equation*}
\overline{\rho}(s\ot \bp)=\dsum_{l,i}\bp_l\ot s_i\Ll(p_{l,i}\ot \bp).
\end{equation*}
Its  inverse is given by  
\begin{equation*}
\overline{\rho}^\m(\bp\ot s)=\dsum_{k,j}\Rr(\bp\ot p_{k,j})s_j\ot \bp_k,
\end{equation*}
where $\dsum_{k=1}^l\Ll(p_k\ot \bp_k)=1$ and $\rho(s\ot p_k)=\dsum_{j=1}^mp_{k,j}\ot s_j$. 

{ (i)} Let $s \in S$, $p\in P$ and  $\bp \in P^\m$. Write $\dsum_{k=1}^n\Ll(p_k\ot \bp_k)=1$ and denote $\rho(s\ot p_k)=\dsum_{j=1}^mp_{k,j}\ot s_j$. Then 
\begin{eqnarray*}
p\ot \bp \ot s & \stackrel{P\ot \overline{\rho}^\m}{\longmapsto} & \dsum_{k,j} p\ot \Rr(\bp\ot p_{k,j})s_j\ot \bp_k\\
&\stackrel{\rho^\m\ot P^\m}{\longmapsto} & \dsum_{k,j} \rho^\m(p\Rr(\bp\ot p_{k,j})\ot s_j)\ot \bp_k\\
& = & \dsum_{k,j} \rho^\m(\Ll(p\ot\bp)p_{k,j}\ot s_j)\ot \bp_k\\
& = & \dsum_{k,j} \Ll(p\ot\bp)\rho^\m(p_{k,j}\ot s_j)\ot \bp_k\\
& = & \dsum_{k} \Ll(p\ot\bp)s\ot p_k\ot \bp_k\\
& \stackrel{S\ot \Ll}{\longmapsto} & \dsum_{k} \Ll(p\ot\bp)s\Ll(p_k\ot \bp_k)\\
& = &  \Ll(p\ot\bp)s.
\end{eqnarray*}
Hence, 
\begin{equation*}
\Psi(p\ot \bp\ot s)=(S\ot \Ll)\circ (\rho^\m\ot P^\m)\circ (P\ot \overline{\rho}^\m)(p\ot \bp \ot s)=\Ll(p\ot \bp)s.
\end{equation*}

{ (ii)} Let us  compute  $\Psi^\m(s\Ll(p\ot \bp))$: observe that $\Psi^\m=(P\ot \overline{\rho})\circ (\rho\ot P^\m)\circ (S\ot \Ll^\m)$. If $\dsum_{k=1}^n\Ll(p_k\ot \bp_k)=1$, then  
\begin{eqnarray*}
s\Ll(p\ot \bp) & \stackrel{S\ot\Ll^\m}{\longmapsto}& \dsum_{k=1}^l s\Ll(p\ot \bp)\ot p_k\ot \bp_k = \dsum_{k=1}^l s\ot \Ll(p\ot \bp)p_k\ot \bp_k\\
& = &   \dsum_{k=1}^l s\ot p\Rr (\bp\ot p_k)\ot \bp_k=\dsum_{k=1}^l s\ot p\ot\Rr (\bp\ot p_k) \bp_k\\
& = & \dsum_{k=1}^l s\ot p\ot \bp\Ll(p_k\ot \bp_k)=s\ot p\ot \bp\\
& \stackrel{\rho \ot P^\m}{\longmapsto} & \dsum_{i=1}^np_i\ot s_i \ot \bp \stackrel{P\ot \overline{\rho}}{\longmapsto} \dsum_{i,j,k} p_i\ot \bp_k\ot s_{i,j}\Ll(p_{k,j}\ot \bp),
\end{eqnarray*}
that is, $\Psi^\m(s\Ll(p\ot \bp))= \dsum_{i,j,k} p_i\ot \bp_k\ot s_{i,j}\Ll(p_{k,j}\ot \bp).$
Therefore, 
\begin{eqnarray*}
s\Ll(p\ot \bp)& = & \Psi\Psi^\m(s\Ll(p\ot \bp))= \Psi\left( \dsum_{i,j,k} p_i\ot \bp_k\ot s_{i,j}\Ll(p_{k,j}\ot \bp)\right) \\
& = & \dsum_{i,j,k} \Ll(p_i\ot \bp_k) s_{i,j}\Ll(p_{k,j}\ot \bp).
\end{eqnarray*}

{ (iii)} It is enough to  verify that  $\Psi$ is left  $S$-linear. Take $s,s' \in S$, $p\in P$ and $\bp \in P^\m$. Let $p_k \in P$
and   $\bp_k \in P^\m$, be such that  $\dsum_{k=1}^l\Ll(p_k\ot \bp_k)=1$ and denote $\rho(s'\ot p)=\dsum_{i=1}^np_i\ot s_i'$ e $\rho^\m(p_k\ot s_i')=\dsum_{j=1}^l s_{i,j}'\ot p_{k,j}$. Then 
$\overline{\rho}(s_i'\ot \bp)=\dsum_{k,l}\bp_k\ot s_{i,j}'\Ll(p_{k,j}\ot \bp)$. Thus,
\begin{eqnarray*}
\Psi(s'*(p\ot (\bp\ot s)))& = & \Psi\left( \dsum_{i=1}^np_i\ot s_i'*(\bp\ot s)\right)\\
& = & \Psi\left( \dsum_{i,j,k}p_i\ot \bp_k\ot s_{i,j}'\Ll(p_{k,j}\ot \bp)s\right) \\
& = & \dsum_{i,j,k}\Ll(p_i\ot \bp_k)s_{i,j}'\Ll(p_{k,j}\ot \bp)s\\
& \stackrel{(\ref{slpbp})}{=} & s'\Ll(p\ot \bp)s=s'\Psi(p\ot \bp \ot s). 
\end{eqnarray*}

{ (iv)} By Lemma~\ref{lemaPSMisoPM} and by   item (iii), we have the $R$-bimodule  isomorphisms 
$$(P\ot_RS)\ot_S(P^\m\ot_RS)\simeq P\ot_R(P^\m\ot_RS)\simeq S.$$
Similarly, $(P^\m\ot_RS)\ot_S(P\ot_RS)\simeq S ,$ as $S$-bimodules. Consequently, $[P\ot_RS]\in \Pic(S).$
\end{dem}

\begin{teo}\label{segundaseqexata}
The sequence 
\begin{equation*}
\xymatrix@C=1.2cm{ \p_\Z(\D(\T)/R)^{(G)}\ar[r]^-{\varphi_2}  &  \Pic_\Z(R)\cap\Pics_\Z(R)^{\al^*}\ar[r]^-{\varphi_3}  & \C_0(\T/R)  }
\end{equation*}
is exact.
\end{teo}
\begin{dem}
If $[Q]\in \mbox{Im}(\varphi_2)$, then there exists  $\xymatrix@C=1.2cm{ [Q]\ar@{=>}[r]|-{[\phi]} & [X]}\in \p_\Z(\D(\T)/R)^{(G)}$. By the proof of  Proposition~\ref{osegundomorfismo} there exists an $R$-bimodule isomorphism  $g_x:\T_x\ot Q\ot \T_{x^{-1}}\longrightarrow Q\ot R1_x$ determined  by 
\begin{equation*}
g_x(u_x\ot q\ot u_{x^\m})=q'\ot 1_x, \ \mbox{where} \ \phi(q')=u_x\phi(q)u_{x^\m}.
\end{equation*} 
%

Let $h_x:Q\ot \T_x\longrightarrow \T_x\ot Q$ be an $R$-bimodule isomorphism, such that  the following diagram is commutative: 
\begin{equation*}
\xymatrix@C=2cm{ R1_x\ot Q\ot \T_x\ar[r] & \T_x\ot \T_{x^{-1}}\ot Q\ot \T_x\ar[r]^{\T_x \ot g_{x^\m}} & \T_x\ot Q\ot R1_{x^\m}\ar[d]\\
	Q\ot R1_{x}\ot \T_{x}\ar[u] & & \T_x\ot R1_{x^\m}\ot Q\ar[d]\\
	Q\ot \T_x\ar[u]\ar@{-->}[rr]_{h_x} & & \T_x\ot Q   } 
\end{equation*}	
that is, 
\begin{equation*}
h_x(q\ot u_x)= \dsum_{(x)} \om_x\ot g_{x^\m}(\om_{x^\m}\ot q \ot u_x),
\end{equation*}
where $1_x=\dsum_{(x)}\om_x\stackrel{\T}{\circ}\om_{x^\m}$, with  $\om_x\in \T_x ,$  $\om_{x^\m}\in \T_{x^{-1}}$, { and arbitrary} $q\in Q,$  $u_x \in \T_x$.

\begin{afr} The isomorphisms  $h_x$ satisfy the commutative   diagram:
\end{afr}
\begin{equation}
\xymatrix@C=2cm{ & Q\ot \T_x\ot\T_y\ar[ld]_{h_x\ot \T_y}\ar[rd]^{Q\ot f_{x,y}^\T} & \\
	\T_x\ot Q\ot\T_y\ar[dd]_{\T_x\ot h_y} & & Q\ot 1_x\T_{xy}\ar[dd]^{1_xh_{xy}}\\
	& & \\
	\T_x\ot \T_y \ot Q\ar[rr]_{f_{x,y}^\T\ot Q} & & 1_x\T_{xy}\ot Q.}\label{diagramafxyTfxfy}
\end{equation}
Indeed, let $1_x=\dsum_{(x)}(\om_x\stackrel{\T}{\circ}\om_{x^\m})$ and $1_y=\dsum_{(y)}(\om_y\stackrel{\T}{\circ}\om_{y^\m})$. By Remark~\ref{obsinversodefxxm}, { the restriction }   $1_xh_{xy}:Q\ot 1_x\T_{xy}\longrightarrow 1_x\T_{xy}\ot Q$ is given by   
\begin{equation*}
1_xh_{xy}(q\ot u_{xy})=\dsum_{(x),(y)}(\om_x\stackrel{\T}{\circ}\om_y)\ot g_{(xy)^\m}((\om_{y^\m}\stackrel{\T}{\circ}\om_{x^\m})\ot q \ot u_{xy}),
\end{equation*}
for all  $u_{xy}\in 1_x\T_{xy}$. 
Thus, given $q \in Q$, $u_x \in \T_x$ and $u_y \in \T_y$, we have:
\begin{eqnarray}
1_xh_{xy}\circ(Q\ot f_{x,y}^\T)(q \ot u_x\ot u_y)  & = & 1_xh_{xy}( q\ot (u_x\stackrel{\T}{\circ} u_y))\nonumber \\
& = & \dsum_{(x),(y)}(\om_x\stackrel{\T}{\circ }\om_y)\ot g_{(xy)^\m}((\om_{y^\m}\stackrel{\T}{\circ}\om_{x^\m})\ot q \ot (u_x\stackrel{\T}{\circ} u_y))\nonumber\\
& = & \dsum_{(x),(y)}(\om_x\stackrel{\T}{\circ }\om_y)\ot q_{x,y},\label{lado1defxyTefxfy}
\end{eqnarray}
where $\phi(q_{x,y})=\dsum_{(x),(y)}(\om_{y^\m}\stackrel{\T}{\circ}\om_{x^\m})\phi(q) (u_x\stackrel{\T}{\circ} u_y)$.

On the other hand, applying  $(f_{x,y}^\T\ot Q)\circ (\T_x\ot h_{y})\ot (h_x\ot \T_y) , $ we obtain:
\begin{eqnarray}
q\ot u_x\ot u_y & \stackrel{h_x\ot \T_y}{\mapsto}& \dsum_{(x)} \om_x\ot g_x(\om_{x^\m}\ot q \ot u_x)\ot u_y\nonumber\\
& = & \dsum_{(x)} \om_x\ot q'_x\ot u_y, \ \mbox{{ where}} \ \phi(q'_x)=\dsum_{(x)}\om_{x^\m}\phi(q)u_x,\nonumber\\
& \stackrel{\T_x\ot h_y}{\mapsto} & \dsum_{(x),(y)} \om_x\ot \om_y\ot g_y(\om_{y^\m}\ot {q_x'}\ot u_y)\nonumber\\
& = &  \dsum_{(x),(y)} \om_x\ot \om_y\ot q'_{x,y},\nonumber\\
& \stackrel{f_{x,y}^\T\ot Q}{\mapsto} &  \dsum_{(x),(y)} (\om_x\stackrel{\T}{\circ} \om_y)\ot q'_{x,y},\label{lado2defxyTefxfy}
\end{eqnarray}
where 
\begin{eqnarray*}
\phi(q_{x,y}')& = & \dsum_{(x),(y)}\om_{y^\m}\phi(q_x')u_y=\dsum_{(x),(y)}\om_{y^\m}(\om_{x^\m}\phi(q)u_x)u_y\\
& = & \dsum_{(x),(y)}(\om_{y^\m}\stackrel{\T}{\circ}\om_{x^\m})\phi(q)(u_x\stackrel{\T}{\circ}u_y)=\phi(q_{x,y}).
\end{eqnarray*}
Since $\phi$ is injective, we get   $q_{x,y}'=q_{x,y}$. Consequently,  (\ref{lado1defxyTefxfy}) and (\ref{lado2defxyTefxfy}) imply that the  diagram is commutative. 

The commutativity of  (\ref{diagramafxyTfxfy}) may be expressed  as follows: given $q\in Q, u_x \in \T_x$ and $u_y\in \T_y$, denote
\begin{equation*}
h_x(q\ot u_x)=\dsum_{i=1}^n u_x^i\ot q_i \ \ \mbox{and} \ \ h_y(q_i\ot u_y)=\dsum_{j=1}^mu_y^{i,j}\ot q_{i,j}.
\end{equation*}
Then, 
\begin{eqnarray}
h_{xy}(q \ot (u_x\stackrel{\T}{\circ}u_y)) & = & \dsum_{i,j}(u_x^i\stackrel{\T}{\circ} u_{y}^{i,j})\ot q_{i,j}. \label{comutatividadedodiagramafxyTfxfy}
\end{eqnarray}

Consider now the $R$-bimodule isomorphism  $F_x:Q\ot \T_x\ot Q^\m\longrightarrow\T_x$ defined  by  $F_x:=h_x\ot Q^\m$, that is, given  $q\in Q, u_x \in \T_x$ and $\bar{q}\in Q^\m$, if we denote   $h_x(q\ot u_x)=\dsum_{i=1}^nu_x^i\ot q_i$, then 
\begin{equation*}
F_x(q\ot u_x\ot \bar{q})=\dsum_{i=1}^nu_x^{i}\Ll(q_i\ot \bar{q}) \in Q,
\end{equation*}
where $Q\ot_RQ^\m\stackrel{\Ll}{\longrightarrow} R\stackrel{\Rr}{\longleftarrow}Q^\m\ot_RQ$ are $R$-bimodule  isomorphisms. 
\begin{afr}
$\varphi_3([Q])=[\D(\T)]$ in $\C_0(\T/R)$.
\end{afr}
Indeed, it suffices to  verify that the following  diagram is commutative
\begin{equation*}
\xymatrix@C=2cm{ Q\ot\T_x\ot Q^\m\ot Q\ot \T_y\ot Q^\m\ar[rr]^{f_{x,y}^Q}\ar[dd]_{F_x\ot F_y} & &  R1_x\ot Q\ot \T_{xy}\ot Q^\m\ar[dd]^{F_{xy}}\\
	& & \\
	\T_x\ot \T_y\ar[rr]_{f_{x,y}^\T} & &  R1_x\ot \T_{xy}    .      }
\end{equation*}

Let $q_1,q_2\in Q, \bar{q}_1,\bar{q}_2\in Q^\m, u_x\in \T_x$ e $u_y \in \T_y$. Write
\begin{equation*}
h_x(q_1\ot u_x)=\dsum_{i=1}^n u_x^i\ot q_i \ \ \mbox{and} \  \ h_y(q_i\ot \Rr(\bar{q}_1\ot q_2)u_y)=\dsum_{j=1}^mu_y^{i,j}\ot q_{i,j}.
\end{equation*}
Then, by (\ref{comutatividadedodiagramafxyTfxfy}) 
\begin{equation*}
h_{xy}(q_i\ot (u_x\stackrel{\T}{\circ}\Rr(\bar{q}_1\ot q_2)u_y))=\dsum_{i,j}(u_x^i\stackrel{\T}{\circ}u_y^{i,j})\ot q_{i,j}.
\end{equation*}
Thus, 
\begin{eqnarray*}
(F_{xy}\circ f_{x,y}^Q)(q_1\ot u_x\ot \bar{q}_1\ot q_2\ot u_y\ot  \bar{q}_2) & = & F_{xy}(q_1\ot (u_x\stackrel{\T}{\circ}\Rr(\bar{q}_1\ot q_2)u_y) \ot \bar{q}_2)\nonumber\\
& = & \dsum_{i,j}(u_x^i\stackrel{\T}{\circ}u_y^{i,j})\Ll(q_{i,j}\ot \bar{q}_2).\label{lado1defQfT}
\end{eqnarray*}
On the other hand, applying  $f_{x,y}^\T\circ (F_x\ot F_y)$ to $q_1\ot u_x\ot \bar{q}_1\ot q_2\ot u_y\ot  \bar{q}_2$, we obtain
\begin{eqnarray*}
& & (q_1\ot u_x\ot \bar{q}_1\ot q_2\ot u_y \ot  \bar{q}_2)  \mapsto  \dsum_{i=1}^n u_x^i\Ll(q_i\ot \bar{q}_1)\ot h_y(q_2\ot u_y)\ot \bar{q}_2\nonumber\\
& = & \dsum_{i=1}^nu_x^i\ot \Ll(q_i\ot \bar{q}_1) {h_y}(q_2\ot u_y)\ot \bar{q}_2 = \dsum_{i=1}^nu_x^i\ot  h_y(\Ll(q_i\ot \bar{q}_1)q_2\ot u_y)\ot \bar{q}_2\nonumber\\
& = & \dsum_{i=1}^nu_x^i\ot  h_y(q_i\Rr(\bar{q}_1\ot q_2)\ot u_y)\ot \bar{q}_2 = \dsum_{i=1}^nu_x^i\ot  h_y(q_i\ot \Rr(\bar{q}_1\ot q_2) u_y)\ot \bar{q}_2\nonumber\\
& = & \dsum_{i,j}u_x^i\ot u_y^{i,j}\ot q_{i,j}\ot \bar{q}_2 {\mapsto} \dsum_{i,j}u_x^i\ot u_y^{i,j}\Ll(q_{i,j}\ot \bar{q}_2)\nonumber\\
&\mapsto & \dsum_{i,j}(u_x^i\stackrel{\T}{\circ}u_y^{i,j})\Ll(q_{i,j}\ot \bar{q}_2).\label{lado2defQfT}
\end{eqnarray*}
Hence, $F_x\circ f_{x,y}^Q=f_{x,y}^\T\circ (F_x\ot F_y)$, for all $x,y \in G$. 
Therefore, $\D(\Om^Q)\simeq \D(\T)$ as generalized partial crossed products. This yields that  $\varphi_3([Q])=[\D(\Om^Q)]=[\D(\T)]$ in $\C_0(\T/R)$ and, consequently,   $[Q]\in \ker(\varphi_3)$.

If $[P]\in \ker(\varphi_3)$, then there  exists an isomorphism of generalized partial crossed products  $j:\D(\Om^P)\longrightarrow \D(\T).$ This means that,  for each  $x \in G$, there exists an $R$-bimodule isomorphism  $j_x:P\ot\T_x\ot P^\m\longrightarrow \T_x$ such that the  diagram 
\begin{equation}
\xymatrix{ P\ot\T_x\ot P^\m\ot P\ot \T_y\ot P^\m\ar[rr]^{f_{x,y}^P}\ar[dd]_{j_x\ot j_y} & &  R1_x\ot P\ot \T_{xy}\ot P^\m\ar[dd]^{j_{xy}}\\
	& & \\
	\T_x\ot \T_y\ar[rr]_{f_{x,y}^\T} & &  R1_x\ot \T_{xy}          } \label{diagramaj}
\end{equation}
is commutative. Given $u_x\in \T_x$ e $u_y\in \T_y$, denote
\begin{equation}
j_x^\m(u_x)=\dsum_{i=1}^n p_i\ot u_x^i\ot\bar{p}_i \ \ \mbox{e} \ \ j_y^\m(u_y)=\dsum_{j=1}^mp_j\ot u_y^j\ot\bar{p}_j. \label{jxinversoejyinverso}
\end{equation}
Then, by the  commutativity of the diagram, we obtain 
\begin{equation}
j_{xy}^\m(u_x\stackrel{\T}{\circ}u_y)=\dsum_{i,j}p_i\ot (u_x^i\Rr(\bar{p}_i\ot p_j)\stackrel{\T}{\circ}u_y^j)\ot \bar{p}_j.\label{comutatividadedodiagramacomj}
\end{equation}

Let $i_x:\T_x\ot P\longrightarrow P\ot \T_x$ be defined by  
\begin{equation*}
i_x: \T_x\ot P \stackrel{j_x^\m\ot P}{\longrightarrow } P\ot \T_x\ot P^\m\ot P\longrightarrow P\ot \T_x,
\end{equation*}
that is, given $p\in P$ and $u_x\in \T_x$, if $j_x^\m(u_x)=\dsum_{i=1}^np_i\ot u_x^i\ot \bar{p}_i$, then
\begin{equation*}
i_x(u_x\ot p)=\dsum_{i=1}^np_i\ot u_x^i\Rr(\bar{p}_i\ot p).
\end{equation*}

Put $i:=\bigoplus_{x\in G}i_x:\D(\T)\ot P\longrightarrow P\ot \D(\T)$. Let us verify that  $i$ satisfies  conditions  (\ref{eq5.1}) and (\ref{eq5.2}) of Lemma~\ref{estruturadedeltathetabimodulo}.  The first condition is:
\begin{equation*}
(P\ot f_{x,y}^\T)\circ (i_x\ot \T_y)\circ (\T_x\ot i_y)=i_{xy}\circ (f_{x,y}^\T\ot P) \ \ \forall \ x,y \in G. 
\end{equation*}
Let $p \in P, u_x \in \T_x, u_y \in \T_y$ in the  notation of  (\ref{jxinversoejyinverso}), then $i_{xy}\circ (f_{x,y}^\T\ot P)$ is given by:
\begin{eqnarray*}
u_x\ot u_y\ot p & \mapsto & (u_x\stackrel{\T}{\circ}u_y)\ot p \mapsto j_{xy}^\m(u_x\stackrel{\T}{\circ}u_y)\ot p\\
& \stackrel{(\ref{comutatividadedodiagramacomj})}{=} & \dsum_{i,j}p_i\ot \left( u_x^i\Rr(\bar{p}_i\ot p_j)\stackrel{\T}{\circ}u_y^j\right) \ot \bar{p}_j\ot p\\
& \mapsto & \dsum_{i,j}p_i\ot \left( u_x^i\Rr(\bar{p}_i\ot p_j)\stackrel{\T}{\circ}u_y^j\right) \Rr(\bar{p}_j\ot p).
\end{eqnarray*}
On the other hand,  $(P\ot f_{x,y}^\T)\circ (i_x\ot \T_y)\circ (\T_x\ot i_y)$ is given by:
\begin{eqnarray*}
u_x\ot u_y\ot p & \stackrel{\T_x\ot i_y}{\longmapsto} & \dsum_{j=1}^mu_x\ot p_j\ot u_y^j\Rr(\bar{p}_j\ot p)\\
& \stackrel{i_x\ot \T_y}{\longmapsto} & \dsum_{i,j}p_i\ot u_x^i\Rr(\bar{p}_i\ot p_j)\ot u_y^j\Rr(\bar{p}_j\ot p)\\
& \stackrel{P\ot f_{x,y}^\T}{\longmapsto} & \dsum_{i,j}p_i\ot \left( u_x^i\Rr(\bar{p}_i\ot p_j)\stackrel{\T}{\circ}u_y^j\Rr(\bar{p}_j\ot p)\right) \\
&=& \dsum_{i,j}p_i\ot \left( u_x^i\Rr(\bar{p}_i\ot p_j)\stackrel{\T}{\circ} u_y^j\right) \Rr(\bar{p}_j\ot p).
\end{eqnarray*}
Consequently, $i$ satisfies  condition (\ref{eq5.1}). Since $j$ is an  isomorphism of generalized partial crossed products,  then  $j_1\circ \nu_P=\iota$, where $\nu_P$ is given in  Remark~\ref{morfismonuparaOmP}. Then  $$j_1^\m(1)=j_1^\m(\iota(1))=\nu_P(1)=\dsum_{k=1}^np_k \ot \iota(1)\ot \bp_k,$$ where $\dsum_{k=1}^n\Ll(p_k\ot \bar{p}_k)=1$, for all  $r \in R$. Hence,
\begin{eqnarray*}
i_1(1\ot p)  & = & \dsum_{k=1}^np_k\ot \Rr(\bar{p}_k\ot p)=\dsum_{k=1}^np_k\Rr(\bar{p}_k\ot p)\ot 1\\
& = & \dsum_{k=1}^n\Ll(p_k\ot\bar{p}_k)p\ot 1=p\ot 1.
\end{eqnarray*}
Therefore, $i$  satisfies also condition  (\ref{eq5.2}). It follows from  Lemma~\ref{estruturadedeltathetabimodulo} that $P\ot \D(\T)$ is an  $\D(\T)$-bimodule via:
{	\begin{equation*}
	(p\ot u_x)*u_y=p\ot (u_x\stackrel{\T}{\circ}u_y) \ \mbox{and} \ u_y*(p\ot u_x)= \dsum_{i=1}^np_i\ot (u_y^i\stackrel{\T}{\circ}u_x),
\end{equation*}
where $i_y(u_y\ot p)=\dsum_{i=1}^np_i\ot u_y^i.$	}
We   construct the  isomorphism  $\overline{i_x}:\T_x\ot P^\m\longrightarrow P^\m\ot \T_x$ by:
\begin{equation*}
\overline{i_x}:\T_x\ot P^\m\longrightarrow P^\m\ot P\ot_R\T_x\ot P^\m\stackrel{j_x}{\longrightarrow} P^\m\ot \T_x,
\end{equation*}
that is,
\begin{equation*}
\overline{i_x}(u_x\ot \bp)=\dsum_{k=1}^n\bp_k\ot j_x(p_k\ot u_x\ot \bp),
\end{equation*}
{ where $\dsum_{k=1}^n\Rr(\bp_k\ot p_k)=1$}. Let $\overline{i}=\bigoplus_{x\in G}\overline{i_x}:\D(\T)\ot {P^\m}\longrightarrow {P^\m}\ot_R\D(\T)$. We shall verify that   $\overline{i}$ also   satisfies the conditions of 
Lemma~\ref{estruturadedeltathetabimodulo}. Let $u_x \in \T_x,u_y \in \T_y$ and $\bp \in P^\m$. { Take any two decompositions  $\dsum_{k=1}^n\Rr(\bp_k\ot p_k)=1$ and
$\dsum_{l=1}^m\Rr(\bp_l\ot p_l)=1$ }. Then, $(P\ot f_{x,y}^\T)\circ (\overline{i_x}\ot \T_y)\circ (\T_x\ot \overline{i_y})$ is given by:
\begin{eqnarray*}
u_x\ot u_y \ot \bp & \stackrel{\T_x\ot \overline{i_y}}{\longmapsto}& \dsum_{l=1}^mu_x\ot \bp_l\ot j_y(p_l\ot u_y\ot \bp)\\
& \stackrel{\overline{i_x}\ot \T_y}{\longmapsto} & \dsum_{k,l} \bp_k\ot j_x(p_k\ot u_x\ot \bp_l)\ot j_y(p_l\ot u_y\ot \bp)\\
& \stackrel{P\ot f_{x,y}^\T}{\longmapsto} & \dsum_{k,l} \bp_k\ot \left( j_x(p_k\ot u_x\ot \bp_l)\stackrel{\T}{\circ} j_y(p_l\ot u_y\ot \bp)\right) \\
& \stackrel{(*)}{=} & \dsum_{k,l}  \bp_k\ot j_{xy}(p_k\ot (u_x\Rr (\bp_l\ot p_l)\stackrel{\T}{\circ} u_y)\ot \bp)\\
& = & \dsum_{k=1}^n  \bp_k\ot j_{xy}(p_k\ot (u_x\stackrel{\T}{\circ} u_y)\ot \bp)\\
& = & i_{xy}((u_x\stackrel{\T}{\circ}u_y)\ot \bp)\\
& = & i_{xy}\circ (f_{x,y}^\T\ot P^\m)(u_x\ot u_y \ot \bp),
\end{eqnarray*}
where the equality   $(*)$ comes from the commutativity of diagram  (\ref{diagramaj}). Consequently,    $(P\ot f_{x,y}^\T)\circ (\overline{i_x}\ot \T_y)\circ (\T_x\ot \overline{i_y})=\overline{i_{xy}}\circ (f_{x,y}^\T\ot P^\m)$. For the second condition we have 
\begin{eqnarray*}
\overline{i_1}(1\ot \bp) & = & \dsum_{k=1}^n\bp_k\ot j_1(p_k\ot 1 \ot \bp) = \dsum_{k=1}^n\bp_k\ot \Ll(p_k\ot \bp)\\
& = & \dsum_{k=1}^n\bp_k\Ll(p_k\ot \bp)\ot 1=\dsum_{k=1}^n\Rr(\bp_k\ot p_k)\bp\ot 1\\
& = & \bp\ot 1.
\end{eqnarray*}
Lemma~\ref{estruturadedeltathetabimodulo} implies that  $P^\m\ot_R\D(\T)$ possesses a   $\D(\T)$-bimodule structure  given by  $\overline{i}$. 
\begin{afr}
The isomorphism  $\overline{i}$ is given by   $\overline{i}=P^\m\ot i^\m\ot P^\m$, where $P^\m\ot i^\m\ot P^\m:\D(\T)\ot P^\m\longrightarrow P^\m\ot \D(\T)$ is defined as in  item (i) of Lemma~\ref{lemacomoverlinerho}.
\end{afr}
Indeed, we first observe that  $i_x^\m: P\ot \T_x\longrightarrow \T_x\ot P$ is given by 
$$i_x^\m(p\ot u_x)=\dsum_{l=1}^mj_x(p\ot u_x\ot \bp_l)\ot p_l,$$
where $\dsum_{l=1}^m\Rr(\bp_l\ot p_l)=1$. Now, applying  $P^\m\ot i^\m\ot P^\m$ for $u_x \ot \bp$ we obtain: 

\begin{eqnarray*}
u_x\ot \bp & \mapsto & \dsum_{k=1}^n \bp_k\ot p_k\ot u_x\ot \bp \mapsto \dsum_{k=1}^n \bp_k\ot i^\m(p_k\ot u_x)\ot \bp\\
& = & \dsum_{k,l}\bp_k\ot j_x(p_k\ot u_x\ot \bp_l)\ot p_l\ot \bp\\
& \mapsto &  \dsum_{k,l}\bp_k\ot j_x(p_k\ot u_x\ot \bp_l)\Ll( p_l\ot \bp)\\
& = &  \dsum_{k,l}\bp_k\ot j_x(p_k\ot u_x\ot \bp_l\Ll( p_l\ot \bp))\\
& = &  \dsum_{k,l}\bp_k\ot j_x(p_k\ot u_x\ot\Rr( \bp_l\ot p_l)\bp)\\
& = &  \dsum_{k}\bp_k\ot j_x(p_k\ot u_x\ot\bp)\\
& = & \overline{i_x}(u_x\ot \bp). 
\end{eqnarray*}

By item (iv) of Lemma~\ref{lemacomoverlinerho}, we get that  $[P\ot\Delta(\T)]\in \Pic(\Delta(\T))$.

Let us verify that  $\xymatrix@C=1.2cm{ [P]\ar@{=>}[r]|-{[\phi]} & [P\ot\D(\T)]}\in \p_\Z(\Delta(\T)/R)^{(G)},$
where 
$$\begin{array}{c c c l}
\phi: & P & \longrightarrow & P\ot\Delta(\T) , \\
& p & \longmapsto &p\ot 1.
\end{array}$$
Note that 
$$\begin{array}{c c c l}
\bar{\phi}_r: & P\ot_R\D(\T) & \longrightarrow & P\ot\D(\T),\\
& p\ot u_x & \longrightarrow & \phi(p)*u_x
\end{array}$$
is the identity, because   
\begin{equation*}
\phi(p)*u_x=(p\ot 1)*u_x=p\ot u_x=p\ot u_x.
\end{equation*}
Thus, $\bar{\phi}_r$ is an   $R$-$\D(\T)$-bimodule isomorphism. By Remark~\ref{obsphirouphil},  $\xymatrix@C=1.2cm{ [P]\ar@{=>}[r]|-{[\phi]} & [P\ot\D(\T)]}\in \p_\Z(\Delta(\T)/R)$. Let us now  verify that   $\phi(P)*\T_x=\T_x*\phi(P)$, for all  $x\in G$. Let $p \in P$ and $u_x\in \T_x$. Denoting  $i_x(u_x\ot p)=\dsum_{i=1}^np_i\ot u_x^{i}$, we compute 
\begin{eqnarray*}
u_x*\phi(p) & = & u_x*(p\ot 1)= \dsum_{i=1}^np_i\ot (u_x^i\stackrel{\T}{\circ} 1)\\
& = & \dsum_{i=1}^np_i\ot u_x^{i}=\dsum_{i=1}^n(p_i\ot 1)*u_x^i\\
& = & \dsum_{i=1}^n\phi(p_i)*u_x^{i} \in \phi(P)*\T_x.
\end{eqnarray*}
Consequently, $\T_x*\phi(P)\subseteq \phi(P)*\T_x$. On the other hand, since  $p\ot u_x \in P\ot_R\Delta(\T)$ and $i:\Delta(\T)\ot_RP\longrightarrow P\ot_R\Delta(\T)$ is an  isomorphism, there exist  $u_x^{j}\in \T_x$ and $p_j\in P$,  $j=1,2,...,m$, such that 
$$i_x\left( \dsum_{j=1}^m u_x^{j}\ot p_j\right)=p\ot u_x.$$
Then, 
\begin{eqnarray*}
\phi(p)*u_x & = & p\ot u_x= i_x\left( \dsum_{j=1}^mu_x^j\ot p_j\right) \\
& = & \dsum_{j=1}^mu_x^j*(p_j\ot 1)=\dsum_{j=1}^mu_x^j*\phi(p_j)  \in \T_x*\phi(P).
\end{eqnarray*}
Thus, $\phi(P){*}\T_x\subseteq \T_x{*}\phi(P) $ and, consequently,  we obtain the desired equality.  Therefore,  
$$\xymatrix@C=1.2cm{ [P]\ar@{=>}[r]|-{[\phi]} & [P\ot_R\D(\T)]}\in \p_\Z(\Delta(\T)/R)^{(G)}.$$
Obviously,  $\varphi_2(\xymatrix@C=1.2cm{ [P]\ar@{=>}[r]|-{[\phi]} & [P\ot\D(\T)]})=[P]$ and, hence, $[P]\in \mbox{Im}(\varphi_2)$.
\end{dem}

\subsection{The group $\mathcal{B}(\T/R)$ and the third exact sequence}\label{sec:B}

Define the  group $\mathcal{B}(\T/R)$ by  the exact sequence  
\begin{equation*}
	\xymatrix{ \Pic_\Z(R)^{(G)}\ar[r]^-{\mathcal{L}} & \C(\T/R)\ar[r]& \mathcal{B}(\T/R)\ar[r] & 1 },
\end{equation*}
that is,  
\begin{equation*}
	\mathcal{B}(\T/R)=\dfrac{\C(\T/R)}{\mbox{Im}(\mathcal{L})}.
\end{equation*}

\textbf{Notation:} If $f:G_1\longrightarrow G_2$ is a homomorphism of groups, we denote by  $f^c$ the co-kernel of $f$, that is, $f^c:G_2\longrightarrow \dfrac{G_2}{\mbox{Im}(f)}$.
\begin{pro}\label{terceiraseqexata}
	The sequence
	\begin{equation*}
		\xymatrix@C=1.2cm{ \Pic_\Z(R)\cap\Pics_\Z(R)^{\al^*}\ar[r]^-{\varphi_3}  & \C_0(\T/R)\ar[r]^-{\varphi_4}  & \mathcal{B}(\T/R)  , }
	\end{equation*}
	where $\varphi_4$ is $\mathcal{L}^c$ restricted to  $\C_0(\T/R)$, is exact. 
\end{pro}
\begin{dem}
	Clearly, Im$(\varphi_3)\subseteq \ker(\varphi_4)$. Let $[\D(\Gm)]\in \C_0(\T/R)$ be such that  $\varphi_4([\D(\Gm)])=[1]$ in $\mathcal{B}(\T/R)$. Then there exists  $[P]\in \Pic_\Z(R)^{(G)}$ such that  $\mathcal{L}([P])=[\D(\Gm)]$. In particular, $\Gm_x\simeq P\ot \T_x\ot P^\m$, for all  $x\in G$. Thus, we have that
	\begin{equation*}
		\T_x\simeq \Gm_x\simeq P\ot \T_x\ot P^\m, \ \mbox{for all} \ x\in G.
	\end{equation*}
	Then $\T_x\ot P\simeq P\ot \T_x$, for all  $x\in G$. Hence, $[P]\in \Pic_\Z(R)\cap\Pics_\Z(R)^{\al^*}$ and  $[\D(\Gm)]\in \mbox{Im}(\varphi_3)$.
\end{dem}

\subsection{The group $\overline{H^1}(G,\al^*,\Pics_0(R))$ and the fourth exact sequence}\label{sec:overline{H}}

Define 
\begin{equation*}
	\Pics_0(R)=\{[P]\in \Pics(R); \ P|R \ \mbox{as bimodules} \}.
\end{equation*}

\begin{lem}\label{propriedadesdePics0}
	\begin{itemize}
		\item[(i)] $\Pics_0(R)$ is commutative.
		\item[(ii)] $\Pics_0(R)\subseteq \Pics_\Z(R)$.
		\item[(iii)] If $[P]\in \Pics_0(R)$, then $[\T_x\ot P\ot \T_{x^\m}]\in \Pics_0(R)$.
		\item[(iv)] $\U(\Pics_0(R))\subseteq \Pic_\Z(R)^{(G)}$.
	\end{itemize}
\end{lem}
\begin{dem}
	Item (i)  directly follows from  Proposition~\ref{isomorfismoT}. 
	
	(ii) Let $[P]\in \Pics_0(R).$ By Remark~\ref{MdivideRimplicaemZbimodulocentral}, we obtain that   $P$ is a central   $\Z$-bimodule.  Consequently,  $\Pics_0(R)\subseteq \Pics_\Z(R)$.
	
	(iii) If $[P]\in \Pics_0(R)$, then $P|R$ and by the    compatibility with the tensor product,   $\T_x\ot P\ot \T_{x^\m}|\T_x\ot \T_{x^\m}\simeq R1_x$. It follows by the transitivity that  $\T_x\ot P\ot \T_{x^\m}|R$, for each $x \in G.$ Therefore, $[\T_x\ot P\ot \T_{x^\m}]\in \Pics_0(R)$.
	
	(iv)  If $[P]\in \U(\Pics_0(R))$, then by  (ii) we have  $[P]\in  \Pic_\Z(R)$. Since $P|R$ and $P^\m|R$, it follows that  $P\ot \T_x\ot P^\m|\T_x$, for each $x \in G$. Hence, $[P]\in \Pic_\Z(R)^{(G)}$.
\end{dem}
%


Items  (ii) and (iii) of Lemma~\ref{propriedadesdePics0} imply that we may  restrict each isomorphism $\al_x^*:\Pics_\Z(R)[R1_{x^\m}]\longrightarrow \Pics_\Z(R)[R1_x]$, defined by  $\al_x^{*}([P][R1_x])=[\T_x\ot P\ot \T_{x^\m}]$, to the ideal   $\Pics_0(R)[R1_{x^\m}]$ and obtain a partial action  of $G$ on $\Pics_0(R)$.

{ Given $[\D(\Gamma)]\in \C(\T/R)$, since $\Gm_x|\T_x$, then $\Gm_x\ot\T_{x^\m}|R$, for all  $x \in G$.    Thus, $[\Gm_x\ot \T_{x^\m}] \in \Pics_0(R)$. Define $f^\Gamma:G\longrightarrow \Pics_0(R)$ by  $f^\Gamma(x)=[\Gamma_x\ot\T_{x^\m}]$, for all  $x \in G$. We have that  $f^\Gm(x)\in \Pics_0(R)[R1_x]$ and 
	\begin{eqnarray*}
		[\Gm_x\ot\T_{x^\m}][\T_x\ot\Gm_{x^\m}] & = & [\Gm_x\ot\T_{x^\m}\ot\T_x\ot\Gm_{x^\m}]\\
		& = & [\Gm_x\ot R1_{x^\m}\ot\Gm_{x^\m}]\\ 
		& = & [\Gm_x\ot\Gm_{x^\m}]=[R1_x].
	\end{eqnarray*}
	Consequently, $f^\Gm(x)=[\Gm_x\ot_R\T_{x^\m}]\in \U(\Pics_0(R)[R1_x])$, for each $x\in G$. 
	Moreover, given $x,y \in G$, we obtain by  (iii) of  Lemma~\ref{isomorfismoscomT} that
	\begin{eqnarray*}
		\al_x^*(f^\Gm(y)[R1_{x^\m}])f^\Gm(xy)^{\m}f^\Gm(x) & = & [\T_x\ot\Gm_y\ot\T_{y^\m}\ot\T_{x^\m}][\T_{xy}\ot \Gm_{(xy)^\m}][\Gm_x\ot\T_{x^\m}]\\
		& = &[\T_x\ot\Gm_y\ot\T_{y^\m}\ot\T_{x^\m}\ot\T_{xy}\ot \Gm_{(xy)^\m}\ot\Gm_x\ot\T_{x^\m}] \\
		& = &[\T_x\ot\Gm_y\ot R1_{y^\m}\ot R1_{(xy)^\m}\ot\Gm_{(xy)^\m}\ot\Gm_{x}\ot\T_{x^\m}]\\ 
		& = &[\T_x\ot\Gm_y\ot\Gm_{(xy)^\m}\ot\Gm_{x}\ot\T_{x^\m}]\\ 
		& = &[\T_x\ot R1_{y}\ot R1_{x^\m}\ot\T_{x^\m}]\\ 
		& = & [R1_{xy}\ot \T_x\ot\T_{x^\m} ]\\
		& = & [R1_{xy}\ot R1_x]=[R1_x][R1_{xy}].
	\end{eqnarray*}
	Hence, $f^\Gm \in Z^1(G,\al^*,\Pics_0(R))$, for all  $[\D(\Gm)]\in \C(\T/R)$.}
\begin{lem}\label{definicaodezeta} The map 
	$$\begin{array}{c c c l}
		\zeta: & \C(\T/R) & \longrightarrow & Z^1(G,\al^*,\Pics_0(R)),\\
		& [\D(\Gamma)]& \longmapsto & 	{f^\Gm,}
	\end{array}$$
	is a group homomorphism whose kernel is  $\C_0(\T/R)$.
\end{lem}
\begin{dem}	
	If $[\D(\Gm)]=[\D(\Om)]$ in $\C(\T/R)$, we have $\Gm_x\simeq \Om_x$, for all  $x \in G$. Then, $\Gm_x\ot\T_{x^\m}\simeq \Om_x\ot\T_{x^\m}$, for all  $x \in G$. Therefore,
	$$f^\Gm(x)=[\Gm_x\ot\T_{x^\m}]=[\Om_x\ot\T_{x^\m}]=f^\Om(x),\ \ \mbox{for all} \ x\in G.$$
	Thus, $f^\Gm=f^\Om$ in $Z^1(G,\al^*,\Pics_0(R))$ and, consequently,   $\zeta$ is well-defined.
	
	Let $[\D(\Gm)],[\D(\Om)]\in \C(\T/R)$. Since $[\D(\Gm)][\D(\Om)]=\left[ \bigoplus_{x\in G}\Gm_x\ot\T_{x^\m}\ot\Om_x\right] $ in $\C(\T/R)$, then
	\begin{eqnarray*}
		f^{\Gm\Om}(x)& = &[\Gm_x\ot\T_{x^\m}\ot\Om_x\ot\T_{x^\m}]= [\Gm_x\ot\T_{x^\m}][\Om_x\ot\T_{x^\m}]=f^\Gm(x)f^\Om(x),
	\end{eqnarray*}
	for each  $x \in G$. Hence, $f^{\Gm\Om}=f^\Gm f^\Om$ in $Z^1(G,\al^*,\Pics_0(R))$.  Therefore, $\zeta$ is a group homomorphism. 
	
	It remains to  verify that  $\ker(\zeta)=\C_0(\T/R)$. If $[\D(\Gm)]\in \C_0(\T/R)$, then $\Gamma_x\simeq \T_x$, for each   $x \in G$. Hence,
	$$f^\Gm(x)=[\Gm_x\ot \T_{x^\m}]=[\T_x\ot\T_{x^\m}]=[R1_x], \ \mbox{for all} \ x\in G.$$
	Therefore, $[\D(\Gm)]\in \ker(\zeta)$. On the other hand,  let $[\D(\Gm)]\in \ker(\zeta)$. Then, $\Gm_x\ot\T_{x^\m}\simeq R1_x$, for each  $x\in G$. Taking the tensor product with  $\T_x$, we obtain $\Gm_x\ot R1_{x^\m}\simeq R1_x\ot\T_x$ {and consequently} $\Gamma_x\simeq \T_x, \ \mbox{for all } x\in G.$ Thus, $[\D(\Gm)]\in \C_0(\T/R)$.
\end{dem}

Define the group  $\overline{H^1}(G,\al^*,\Pics_0(R))$ by the exact sequence 
$$\xymatrix{  \Pic_\Z(R)^{(G)}\ar@{-->}[rr]\ar[rd]_{\mathcal{L}} & & Z^1(G,\al^*,\Pics_0(R))\ar[r] &   \overline{H^1}(G,\al^*,\Pics_0(R))\ar[r] & 1\\
	& \C(\T/R)\ar[ru]_{\zeta}& & &    }$$
that is,
\begin{equation*}
	\overline{H^1}(G,\al^*,\Pics_0(R))=\dfrac{Z^1(G,\al^*,\Pics_0(R))}{(\zeta\circ\mathcal{L})(\Pic_\Z(R)^{(G)})}.
\end{equation*}

\begin{obs}\label{rem:B^1} A natural question is the following: what is the  relation between  $B^1(G,\al^*,\Pics_0(R))$ and the   image of   $\zeta\circ\mathcal{L}$?  The inclusion $B^1(G,\al^*,\Pics_0(R))\subseteq \mbox{Im}(\zeta \circ \mathcal{L})$ always holds.
\end{obs}
Indeed, let  $\rho \in B^1(G,\al^*,\Pics_0(R))$, then there exists  $[P]\in \U(\Pics_0(R))\subseteq \Pic_\Z(R)^{(G)}$ (see Lemma~\ref{propriedadesdePics0}), such that 
\begin{eqnarray*}
	\rho(x) & = & \al_x^*([P][R1_{x^\m}])[P^\m]= [\T_x\ot P\ot \T_{x^\m}][P^\m]\\
	& = & [\T_x\ot P\ot \T_{x^\m}\ot P^\m]=[P^\m\ot \T_x\ot P\ot \T_{x^\m}]\\
	& = & \zeta(\mathcal{L}([P^\m]))(x),
\end{eqnarray*}
for all  $x \in G$. 
Hence, $B^1(G,\al^*,\Pics_0(R))\subseteq \mbox{Im}(\zeta\circ \mathcal{L}),$ as desired. 


We define the homomorphism  $\varphi_5:\mathcal{B}(\T/R)\longrightarrow \overline{H^1}(G,\al^*,\Pics_0(R))$ by means of the following commutative  diagram:

\begin{equation*}
	\xymatrix{ & \C_0(\T/R)\ar[rd]^{\varphi_4}\ar@{_(->}[d] & & \\
		\Pic_\Z(R)^{(G)} \ar[r]^{\mathcal{L}} & \C(\T/R)\ar[r]^{\mathcal{L}^c}\ar[d]_{\zeta} & \mathcal{B}(\T/R) \ar[r]\ar@{-->}[d]^{\varphi_5} & 1\\
		& Z^1(G,\al^*,\Pics_0(R))\ar[r]^{(\zeta \mathcal{L})^c} & \overline{H^1}(G,\al^*,\Pics_0(R))\ar[r] & 1   } \label{definicaodevarphi5}
\end{equation*}

\begin{obs}\label{phi5} The map $\varphi_5$ is well-defined. 
\end{obs}
Indeed, let $[\D(\Gm)]$ and $[\D(\Om)]$ in $\C(\T/R)$ be such that  $\mathcal{L}^c([\D(\Gm)])=\mathcal{L}^c([\D(\Om)])$ in $\mathcal{B}(\T/R)$. Then there exists   $[P]\in \Pic_\Z(R)^{(G)}$ with $[\D(\Gm)]=[\D(\Om)]\mathcal{L}([P])$. By the commutativity of the diagram, we have 
\begin{eqnarray*}
	\varphi_5(\mathcal{L}^c([\D(\Gm)])) & = & ((\zeta\mathcal{L})^c\circ \zeta)([\D(\Gm)])\\
	& = & ((\zeta\mathcal{L})^c\circ \zeta)([\D(\Om)]\mathcal{L}([P]))\\
	& = & (\zeta\mathcal{L})^c(\zeta([\D(\Om)]\zeta\mathcal{L}([P]))\\
	& = & (\zeta\mathcal{L})^c(\zeta([\D(\Om)])\\
	& = & \varphi_5(\mathcal{L}^c([\D(\Om)])).
\end{eqnarray*}
Consequently,  $\varphi_5$ does not depend on the choice of the representative of the class in $\mathcal{B}(\T/R)$.

\begin{pro}\label{quartaseqexata} The sequence 
	\begin{equation*}
		\xymatrix{ \C_0(\T/R)\ar[r]^{\varphi_4} & \mathcal{B}(\T/R)\ar[r]^{\varphi_5\ \ \ \ \ \ \ \ } & \overline{H^1}(G,\al^*,\Pics_0(R))  }
	\end{equation*}
	is exact.
\end{pro}
\begin{dem} Let $[\D(\Gm)]\in \C_0(\T/R)$. By the commutativity of the diagram, we obtain that  $\varphi_5(\varphi_4([\D(\Gm)]))=((\zeta\mathcal{L})^c\circ \zeta)([\D(\Gm)])$. Thanks to  Lemma~\ref{definicaodezeta}, 
	$\zeta([\D(\Gm)])=[1]$. Thus, $\varphi_5(\varphi_4([\D(\Gm)]))=[1]$. Hence, Im$(\varphi_4)\subseteq \ker(\varphi_5)$.
	
	Let $\mathcal{L}^c([\D(\Gm)]) \in \mathcal{B}(\T/R)$, where $[\D(\Gm)]\in \C(\T/R)$, is such that  $\varphi_5(\mathcal{L}^c([\D(\Gm)]))=[1]$. By the  commutativity of the  diagram, we obtain
	\begin{equation*}
		((\zeta\mathcal{L})^c\circ \zeta)([\D(\Gm)])=[1] \ \mbox{em} \ \overline{H^1}(G,\al^*,\Pics_0(R)).
	\end{equation*} 
	Then, $\zeta([\D(\Gm)])\in {\mbox{Im}(\zeta\mathcal{L})}$. Let $[P]\in \Pic_\Z(R)^{(G)}$ be such that  $\zeta([\D(\Gm)])=\zeta(\mathcal{L}([P]))$. Since $\zeta$ is a group homomorphism, it follows that  $[\D(\Gm)]\mathcal{L}[P]^\m\in \ker(\zeta)=\C_0(\T/R)$. Then,
	\begin{eqnarray*}
		\varphi_4([\D(\Gm)]\mathcal{L}([P])^\m) & = & \varphi_4([\D(\Gm)]\mathcal{L}([P^\m]))\\
		& = & \mathcal{L}^c([\D(\Gm)])\mathcal{L}^c(\mathcal{L}([P^\m]))\\
		& = & \mathcal{L}^c([\D(\Gm)]).
	\end{eqnarray*}
	Consequently, $\ker(\varphi_5)\subseteq \mbox{Im}(\varphi_4),$ from which we conclude that the sequence is exact. 
\end{dem}

\subsection{The fifth exact sequence}

Let $g$ be a normalized  $1$-cocycle in $Z^1(G,\al^*,\Pics_0(R))$, that is, 
$$\begin{array}{c c c l}
	g: & G & \longrightarrow & \Pics_0(R),\\
	& x & \longrightarrow & [\nabla_x],
\end{array}$$
where $[\nabla_x]\in \U(\Pics_0(R)[R1_x])$, for all  $x \in G$ and 
\begin{equation}
	\al_x^*(g_y[R1_{x^\m}])g_{xy}^\m g_x=[R1_x][R1_{xy}], \ \ \mbox{for all} \ x,y \in G. \label{gecociclo}
\end{equation}
Observe that 
\begin{equation}
	\al_x^*(g_y[R1_{x^\m}])[\T_x]=[\T_x]g_y, \ \ \mbox{for all} \ x,y \in G. \label{gycomutacomTx}
\end{equation}
Indeed,  
\begin{eqnarray*}
	\al_x^*(g_y[R1_{x^\m}])[\T_x] & = & [\T_x\ot \nabla_y\ot \T_{x^\m}][\T_x]\\
	& = & [\T_x\ot \nabla_y\ot R1_{x^\m}]\\
	& = & [\T_x\ot R1_{x^\m}\ot \nabla_y]\\
	& = & [\T_x\ot \nabla_y]\\
	& = & [\T_x]g_y.
\end{eqnarray*}

\begin{lem}\label{representacaoparcialalteradoporcociclo} Let $g$ be a normalized element in $Z^1(G,\al^*,\Pics_0(R))$. Then 
	\begin{equation*}
		\begin{array}{c c c l}
			U: & G & \longrightarrow & \Pics(R),\\
			& x & \longmapsto & g_x[\T_x],
		\end{array}
	\end{equation*}
	is a unital  partial representation with   $U_x\ot U_{x^\m}\simeq R1_x$ and $U_x|\T_x$, for all  $x\in G$. 
\end{lem}
\begin{dem}
	Since $g$ is normalized, then   $[U_1]=g_1[\T_1]=[R]$. Given $x,y\in G,$ we have
	\begin{eqnarray*}
		[U_x][U_y][U_{y^\m}]& = & g_x[\T_x]g_y[\T_y]g_{y^\m}[\T_{y^{\m}}]\\
		& \stackrel{(\ref{gycomutacomTx})}{=}& g_x\al_x^*(g_y[R1_{x^\m}])[\T_x][\T_y]g_{y^\m}[\T_{y^{\m}}]\\
		& \stackrel{(\ref{gecociclo})}{=}& g_{xy}[R1_x][R1_{xy}][R1_x][\T_{xy}]g_{y^\m}[\T_{y^{\m}}]\\
		& =&g_{xy}[R1_x][\T_{xy}]g_{y^\m}[\T_{y^{\m}}]\\
		&=& g_{xy}[\T_{xy}][R1_{y^\m}]g_{y^\m}[\T_{y^{\m}}]\\
		&=& g_{xy}[\T_{xy}]g_{y^\m}[R1_{y^\m}][\T_{y^{\m}}]\\
		&=& g_{xy}[\T_{xy}]g_{y^\m}[\T_{y^{\m}}]\\
		& = & [U_{xy}][U_{y^\m}].
	\end{eqnarray*}
	Analogously, 
	\begin{eqnarray*}
		[U_{x^\m}][U_x][U_y] & = & g_{x^\m}[\T_{x^{-1}}]g_x[\T_x]g_y[\T_y]\\
		& \stackrel{(\ref{gycomutacomTx})}{=} & g_{x^\m}[\T_{x^{-1}}]g_x\al_x^*(g_y[R1_{x^\m}])[\T_x][\T_y]\\
		&\stackrel{(\ref{gecociclo})}{=} & g_{x^\m}[\T_{x^{-1}}]g_{xy}[R1_x][R1_{xy}][R1_x][\T_{xy}]\\
		& = & g_{x^\m}[\T_{x^{-1}}]g_{xy}[R1_x][\T_{xy}]\\
		& = & g_{x^\m}[\T_{x^{-1}}][R1_x]g_{xy}[\T_{xy}]\\
		& = & g_{x^\m}[\T_{x^{-1}}]g_{xy}[\T_{xy}]\\
		& = & [U_{x^\m}][U_{xy}],
	\end{eqnarray*}
	for all  $x,y \in G$.
	Again, since  $g$ is  normalized, we have: 
	\begin{eqnarray*}
		[U_x][U_{x^\m}]& = & g_x[\T_x]g_{x^\m}[\T_{x^\m}]\\
		& \stackrel{(\ref{gycomutacomTx})}{=}& g_x\al^*_x(g_{x^\m}[R1_{x^\m}])[\T_x][\T_{x^{-1}}]\\
		& = & [g_1][R1_x]=[R][R1_x]=[R1_x],
	\end{eqnarray*}
	for all $x \in G$.
	Moreover, since  $g_x=[\nabla_x] \in \Pics_0(R)$, then $\nabla_x|R$, for all  $x\in G$. Hence, $U_x\simeq \nabla_x\ot\T_x|\T_x,$ for all  $x\in G$.
\end{dem}

Due to the fact that  $U$ is a unital partial representation,  Lemma~\ref{isomorfismoscomT} implies that there exists a family of $R$-bimodule isomorphisms 
\begin{equation}
	f^g=\{f_{x,y}^g:U_x\ot _R U_y \longrightarrow R1_x\ot U_{xy}\simeq 1_xU_{xy}, \ \forall \ x,y \in G\}. \label{familiadeisomorfismosassociadaag}
\end{equation}
As $[U_1]=[R]$, we may choose a family of representatives   $\{U_x\}_{x \in G}$ with  $U_1=R$ and a family of $R$-bimodule isomorphisms  $f^g$ with   $f_{x,1}^g:U_x\ot R\longrightarrow U_x$ and $f_{1,x}^g:R\ot U_x\longrightarrow U_x $ given by the left and right actions of  $R$ on $U_x,$  respectively. Then the commutative diagrams in   (\ref{diagramasdaunidade}) are trivially satisfied. 
By Corollary~\ref{obs3cocicloparasimilares} there exists a  partial  $3$-cocycle   $\widetilde{\beta_{-,-,-}^g}$ in $Z_\T^3(G,\al,\Z)$ associated to  $f^g$.

\begin{lem}\label{definicaodedelta} The map 
	\begin{equation*}
		\begin{array}{c c c l}
			\delta: & Z^1(G,\al^*,\Pics_0(R)) & \longrightarrow & H_\T^3(G,\al, \Z),\\
			& g & \longmapsto& [\widetilde{\beta_{-,-,-}^g}],
		\end{array}
	\end{equation*}
	is a group homomorphism. Moreover, 
	\begin{equation}
		\delta\circ \zeta=[1], \label{deltazetaiguala1}
	\end{equation}
	where $\zeta$ is the homomorphism defined in Lemma~\ref{definicaodezeta}. 
\end{lem}
\begin{dem}
	By Proposition~\ref{proparaboadefinicaodedelta}, $[\widetilde{\beta_{-,-,-}^g}] \in H^3_\T(G,\al,\Z)$ does not depend on the choice of the family of $R$-bimodules  $\{U_x\}_{x \in G},$ nor on the choice of the family of the $R$-bimodule isomorphisms   $f^g$ and, thus,  $\delta$ is well-defined. 
	Let $g,h \in Z^1(G,\al^*,\Pics_0(R))$. Denote $[U_x]=g_x[\T_x]$ and $[V_x]=h_x[\T_x]$, for all $x \in G$. Let 
	\begin{equation*}
		f_{x,y}^g:U_x\ot U_y\longrightarrow1_xU_{xy} \ \ \mbox{and} \ \ f_{x,y}^h:V_x\ot V_y\longrightarrow1_xV_{x,y} \ \ x,y \in G,
	\end{equation*}
	be families of $R$-bimodule  isomorphisms and let  $\widetilde{\beta_{-,-,-}^g}$ and $\widetilde{\beta_{-,-,-}^h}$ be the associated  $3$-cocycles in $Z_\T^3(G,\al,\Z).$ 
	Let $[W_x]=g_xh_x[\T_x]$, for all $x\in G$. Then,
	\begin{eqnarray*}
		[W_x] & = & g_xh_x[\T_x]=g_xh_x[R1_x][\T_x]\\
		& = & g_x[R1_x]h_x[\T_x]=g_x[\T_x][\T_{x^{-1}}]h_x[\T_x]\\
		& = & [U_x\ot \T_{x^{-1}}\ot V_x],
	\end{eqnarray*}
	for each $x \in G$. Let  $f_{x,y}^{gh}$ be the $R$-bimodule isomorphisms defined by
	\begin{equation*}
		\xymatrix{ U_x\ot\T_{x^{-1}}\ot V_x\ot  U_y\ot\T_{y^{-1}}\ot V_y\ar[r]^{T}\ar@/_1.2cm/@{-->}[rdd]_{f_{x,y}^{gh}} &   U_x\ot  U_y\ot\T_{y^{-1}}\ot\T_{x^{-1}}\ot V_x\ot V_y \ar[d]^{f_{x,y}^g\ot f_{x,y}^\T\ot f_{x,y}^h}\\
			& 1_xU_{xy}\ot 1_{y^\m}\T_{(xy)^\m}\ot 1_xV_{xy}\ar@{=}[d]\\
			&  1_xU_{xy}\ot\T_{(xy)^\m}\ot V_{xy}}
	\end{equation*}
	where $T:\T_{x^\m}\ot V_x\ot U_y\ot \T_{y^\m}\longrightarrow U_y\ot \T_{y^\m}\ot\T_{x^\m}\ot V_x$ is the   isomorphism from  Proposition~\ref{isomorfismoT}. By Lemma~\ref{conjuntodefatores}, we obtain that  $\widetilde{\beta^{gh}_{-,-,-}}=\widetilde{\beta_{-,-,-,}^g}\widetilde{\beta_{-,-,-}^h}.$
	Hence,
	\begin{equation*}
		\delta(gh)=[\widetilde{\beta^{gh}_{-,-,-}}]=[\widetilde{\beta_{-,-,-}^g}\widetilde{\beta_{-,-,-}^h}]=[\widetilde{\beta_{-,-,-}^g}][\widetilde{\beta_{-,-,-}^h}]=\delta(g)\delta(h).
	\end{equation*}
	and, consequently,   $\delta$ is group homomorphism. 
	{ As to \eqref{deltazetaiguala1}, } take $[\D(\Gm)]\in \C(\T/R)$ and denote $f^\Gm=\zeta([\D(\Gm)])$, that is,  $f^\Gm_x=[\Gm_x][\T_{x^{-1}}]$, for all  $x \in G$. Thus, 
	\begin{equation*}
		[U_x]:=f_x^\Gm[\T_x]=[\Gm_x][\T_{x^{-1}}][\T_x]=[\Gm_x][R1_{x^\m}]=[\Gm_x], \ \ \mbox{for all} \ x \in G. 
	\end{equation*}
	{ Consider a family of $R$-bimodule isomorphisms}  $f^\Gm=\{f_{x,y}^\Gm:\Gm_x\ot \Gm_y\longrightarrow 1_x\Gm_{xy}, \ x,y\in G\}$ which is a factor set for   $\D(\Gm)$. Then $\widetilde{\beta^{f^\Gm}_{-,-,-}}$ is trivial.  
	Therefore,   $\delta(\zeta([\D(\Gm)]))=[1]$ em $H^3_{\T}(G,\al,\Z)$.
\end{dem}


We  define $\varphi_6:\overline{H^1}(G,\al^*,\Pics_0(R))\longrightarrow H^3_\T(G,\al,\Z)$ via the commutative diagram:
\begin{equation*}
	\xymatrix{ \Pic_\Z(R)^{(G)}\ar[r]^{\zeta\mathcal{L}\ \ \ \ \ \ \ } & Z^1(G,\al^*,\Pics_0(R))\ar[r]^{(\zeta\mathcal{L})^c}\ar[d]_{\delta} & \overline{H^1}(G,\al^*,\Pics_0(R))\ar[r]\ar@{-->}[ld]^{\varphi_6}& 1\\
		& H^3_\T(G,\al,\Z). & &   }
\end{equation*}
Let us see that  $\varphi_6$ is well-defined:  take $h,g \in Z^1(G,\al^*,\Pics_0(R))$ such that 
\begin{center}
 $(\zeta\mathcal{L})^c(h)=(\zeta\mathcal{L})^c(g)$ in $\overline{H^1}(G,\al^*,\Pics_0(R))$. 
\end{center}
Then there exists  $[P]\in \Pic_\Z(R)^{(G)}$ such that  $h=g\zeta(\mathcal{L}([P]))$. By (\ref{deltazetaiguala1}) and by the commutativity of the diagram, we obtain 
\begin{eqnarray*}
	\varphi_6((\zeta\mathcal{L})^c(h)) & = & \delta(h)= \delta(g\zeta(\mathcal{L}([P])))\\
	&= & \delta(g)(\delta\circ\zeta)(\mathcal{L}([P]))\\
	& = & \delta(g) =\varphi_6((\zeta\mathcal{L})^c(g)),
\end{eqnarray*}
which shows that  $\varphi_6$ is well-defined.  

We have the following commutative diagram:  
\begin{equation*}
	\xymatrix{ \Pic_\Z\cap\Pics_\Z(R)^{\al^*}\ar[r]^{\ \ \ \ \ \ \varphi_3}\ar@{_(->}[d] & \C_0(\T/R)\ar[rd]^{\varphi_4}\ar@{_(->}[d] & & \\
		\Pic_\Z(R)^{(G)}\ar@{=}[d] \ar[r]^{\mathcal{L}} & \C(\T/R)\ar[r]^{\mathcal{L}^c}\ar[d]_{\zeta} & \mathcal{B}(\T/R) \ar[r]\ar[d]^{\varphi_5} & 1\\
		\Pic_\Z(R)^{(G)}\ar[r]^{\zeta\mathcal{L}\ \ \ \ \ \ }	& Z^1(G,\al^*,\Pics_0(R))\ar[r]^{(\zeta\mathcal{L})^c}\ar[d]_{\delta} & \overline{H^1}(G,\al^*,\Pics_0(R))\ar[r]\ar[ld]^{\varphi_6} & 1  \\
		& H^3_\T(G,\al,\Z) & &  }
\end{equation*}

\begin{pro}\label{quintaseqexata}
	The sequence  
	\begin{equation*}
		\xymatrix@C=1cm{ \mathcal{B}(\T/R)\ar[r]^{\varphi_5\ \ \ \ \ \ \ \ } & \overline{H^1}(G,\al^*,\Pics_0(R))\ar[r]^{\ \ \ \ \ \ \varphi_6} & H^3_\T(G,\al,\Z)  }
	\end{equation*}
	is exact.
\end{pro}
\begin{dem}
	Let $h \in Z^1(G,\al^*, \Pics_0(R))$ be such that  $(\zeta\mathcal{L})^c(h)\in Im(\varphi_5)$. Then there exists    $[\D(\Gm)]\in \C(\T/R)$, with $\varphi_5(\mathcal{L}^c([\D(\Gm)]))=(\zeta\mathcal{L})^c(h)$. In view of  the commutativity of the   diagram and thanks to (\ref{deltazetaiguala1}),  we have 
	\begin{equation*}
		\varphi_6((\zeta\mathcal{L})^c(h))=\varphi_6(\varphi_5(\mathcal{L}^c([\D(\Gm)])))=(\delta\circ \zeta)([\D(\Gm)])=[1].
	\end{equation*}
	Hence, $Im(\varphi_5)\subseteq \ker(\varphi_6)$.
	
	For the converse inclusion, take    $h \in Z^1(G,\al^*, \Pics_0(R))$ with $\varphi_6((\zeta\mathcal{L})^c(h))=[1]$ in $H^3_\T(G,\al,\Z)$. By the commutativity of the diagram we get   that $\delta(h)=[1]$ in $H^3_\T(G,\al,\Z)$. Denoting $\delta(h)=[\beta_{-,-,-}] \in B^3_\T(G,\al,\Z)$, we have that there exists   $\s_{-,-}:G\times G\longrightarrow \Z$, with $\s_{x,y}\in \U(\Z1_x1_{xy})$, for all  $x,y \in G$, and
	\begin{equation*}
		\beta_{x,y,z}=\al_x(\s_{y,z}1_{x^\m})\s_{xy,z}^\m\s_{x,yz}\s_{x,y}^\m,	\end{equation*}  
	for all  $x,y,z \in G$. Define $[\Om_x]:=h_x[\T_x]$, for all  $x\in G$ and  consider the family of associated isomorphisms   $f^h$  as in  (\ref{familiadeisomorfismosassociadaag}). Then, by definition, we have that  
	\begin{equation*}
		\beta_{x,y,x}\circ f^h_{x,yz}\circ (\Om_x\ot f^h_{y,z})=f_{xy,z}^h\circ(f_{x,y}^h\ot \Om_z),
	\end{equation*}
	for all  $x,y,z \in G$. By  Lemma~\ref{representacaoparcialalteradoporcociclo}, we obtain that  $\Om:G\longrightarrow \Pics(R),$ defined by  $[\Om_x]=h_x[\T_x],$ is a unital  partial representation with  $\Om_x\ot \Om_{x^\m}\simeq R1_x$ 
	and $\Om_x|\T_x$, for all  $x\in G$. By Proposition~\ref{3cocicloquecorrigeaassociatividade}, we have that  
	\begin{equation*}
		\begin{array}{c c c l}
			\bar{f}^h_{x,y}:&  \Om_x\ot \Om_y& \longrightarrow & 1_x\Om_{xy},\\
			& \om_x\ot \om_y & \longmapsto & \s_{x,y}f_{x,y}^h(\om_x\ot \om_y),
		\end{array}
	\end{equation*}
	is a factor set for  $\Om$.  Hence, $\D(\Om)$ is a generalized partial crossed product with   $[\D(\Om)]\in \C(\T/R)$. Moreover,  $\zeta([\D(\Om)])=h$.
	Then, by the commutativity of the diagram, we obtain 
	\begin{equation*}
		\varphi_5(\mathcal{L}^c([\D(\Om)]))=((\zeta\mathcal{L})^c\circ\zeta)([\D(\Om)])=(\zeta\mathcal{L})^c(h).
	\end{equation*}
	Thus, $(\zeta\mathcal{L})^c(h)\in \mbox{Im}(\varphi_5)$. Therefore, $\ker(\varphi_6)\subseteq \mbox{Im}(\varphi_5) $ and, consequently,  the sequence is exact.
\end{dem}

In the following diagram we combine the  results of the section:

\begin{equation*}
	\xymatrix{ H^1_\T(G,\al,\Z)\ar[d]^{\varphi_1}& & \\
		\p_\Z(\D(\T)/R)^{(G)}\ar[d]^{\varphi_2\ \ \ \ } & & \\
		\Pic_\Z(R)\cap\Pics_\Z(R)^{\al^*}\ar[r]^{\ \ \ \ \ \ \varphi_3}\ar@{_(->}[d] & \C_0(\T/R)\ar[rd]^{\varphi_4}\ar@{_(->}[d] &  \\
		\Pic_\Z(R)^{(G)}\ar@{=}[d] \ar[r]^{\mathcal{L}} & \C(\T/R)\ar[r]^{\mathcal{L}^c}\ar[d]_{\zeta} & \mathcal{B}(\T/R) \ar[d]^{\varphi_5} \\
		\Pic_\Z(R)^{(G)}\ar[r]^{\zeta\mathcal{L}\ \ \ \ \ \ }	& Z^1(G,\al^*,\Pics_0(R))\ar[r]^{(\zeta\mathcal{L})^c}\ar[d]_{\delta} & \overline{H^1}(G,\al^*,\Pics_0(R))\ar[ld]^{\varphi_6}\\
		& H^3_\T(G,\al,\Z) &  }
\end{equation*}


Recalling that   $\C_0(\T/R)\simeq H^2_\T(G,\al,\Z),$  thanks to    Proposition~\ref{C0isomorfoH2}, we derive from   Theorems~\ref{primeiraseqexata} and \ref{segundaseqexata} and from  Propositions~\ref{terceiraseqexata}, \ref{quartaseqexata} and \ref{quintaseqexata} our main result:

\begin{teo}\label{theo:main}  The sequence
	\begin{equation*}
		\xymatrix{1\ar[r] &  H_\T^{1}(G,\al,\Z)\ar[r]^-{\varphi_1} &  \p_\Z(\D(\T)/R)^{(G)}\ar[r]^-{\varphi_2} & \Pic_\Z(R)\cap\Pics_\Z(R)^{\al^*}\ar[r]^-{\varphi_3} & H_\T^2(G,\al,\Z)\\ 
			\ar[r]^-{\varphi_4}& \mathcal{B}(\T/R)\ar[r]^-{\varphi_5} & 	 \overline{H}^1(G,\al^*,\Pics_0(R))\ar[r]^-{\varphi_6} & H^3_\T(G,\al,\Z)  }
	\end{equation*}
	is exact. 
\end{teo}

\bibliography{refs}{}
\bibliographystyle{acm}

\end{document}